\documentclass[twoside]{fic-l}
\usepackage{amsthm,amsmath,amssymb,amscd}
\usepackage[greek,english]{babel}
\newtheorem{theorem}{Theorem}[section]
\newtheorem{proposition}[theorem]{Proposition}
\newtheorem{lemma}[theorem]{Lemma}
\newtheorem{corollary}[theorem]{Corollary}
\theoremstyle{definition}
\newtheorem{definition}[theorem]{Definition}

\newtheorem{example}[theorem]{Example}
\theoremstyle{remark}
\newtheorem{remark}[theorem]{Remark}
\newtheorem{Remark}[theorem]{Remarks}

\def\A{{\mathcal A}}
\def\P{{\mathcal P}}
\def\Q{{\mathcal Q}}
\def\Qq{{\mathbf Q}}

\def\B{{\mathcal B}}
\def\Ss{{\mathcal S}}
\def\C{{\mathbf C}}
\def\N{{\mathbf N}}
\def\Z{{\mathbf Z}}
\def\limproj{\mathop{\oalign{lim\cr\hidewidth$\longleftarrow$\hidewidth\cr}}}
\begin{document}
\setcounter{tocdepth}{2}
\title{Valuations, deformations, and toric geometry}
\author{Bernard Teissier}
\maketitle
\markboth{\rm Bernard Teissier}{\rm Valuations, deformations, and toric geometry}
\vskip.3truecm
\centerline{\textbf Summary}
\par\medskip\noindent Given a n\oe therian local integral domain $R$ and a
valuation $\nu$ of its field of fractions which is non negative on $R$, i.e., an inclusion $R\subset R_\nu$ of $R$ in a valuation
ring, I study a geometric specialization of $R$ to the graded ring $\hbox{\rm gr}_\nu R$ determined by the
valuation. If $R_\nu$ dominates $R$  and the residue field extension $k=R/m\to R_\nu/m_\nu$ is trivial, this graded ring corresponds to an essentially
toric variety, of Krull dimension $\leq\hbox{\rm dim}R$ but possibly of infinite embedding dimension; it is of the form:
$$\hbox{\rm gr}_\nu R=k[(U_i)_{i\in I}]/(U^m-\lambda_{mn}U^n)_{(m,n)\in E}\ \ ,\ \ \lambda_{mn}\in k^*,$$ where $I$ and $E$ are
countable sets and
 $U^m=U_{i_1}^{m_{i_1}}\ldots U_{i_r}^{m_{i_r}}$. In order to apply this fact to a characteristic-blind proof of local uniformization by
deformation of a partial resolution of singularities of
$\hbox{\rm Spec}\hbox{\rm gr}_\nu R$, in the case where $R$ is equicharacteristic and excellent, with $k$ algebraically closed, I explore the
following strategy:\par\noindent  1) Extend the valuation $\nu$ to a valuation $\hat\nu$ of a suitable n\oe therian {\it scalewise} $\nu$-adic completion
$\hat R^{(\nu)}$ of
$R$ such that the natural map $$\ \ \ \ \ \\ \ \ \ \ \ \ \ \\ \ \ \ \ \ \ \ \\ \ \ \ \ \ \ \ \ \ \ \ \ \ \ \ \ \ \ \ \hbox{\rm gr}_\nu R\to\hbox{\rm
gr}_{\hat\nu}\hat R^{(\nu)}\ \ \ \ \ \
\hbox{\rm is {\it scalewise} birational}.$$
 2) Obtain a presentation, that is a surjective continuous morphism of $k$-algebras
$$\widehat{k[(w_j)_{j\in J}]}\to\hat R^{(\nu)},$$ where the left hand term is a suitable {\it scalewise} completion of the polynomial ring, and the
kernel is generated up to closure by elements whose initial forms for the term order $t$ deduced from $\hat \nu$ are the binomial equations
defining
$\hbox{\rm gr}_{\hat\nu} \hat R^{(\nu)}$, that is, which are of the form:
$$ w^m-\lambda_{mn}w^n+\sum_{t(s)>t(m)=t(n)}c^{(mn)}_sw^s\ \ \ \ \ \ \ \ \ \ \ \ \ \hbox{\rm with}\ (m,n)\in \hat E,\lambda_{mn}\in k^*,\ \
c^{(mn)}_s\in k.$$
\par\noindent 3) Show that all but finitely many of these equations serve only to express the images of all $w_j$'s in the n\oe therian ring $\hat
R^{(\nu)}$ in terms of finitely many of them. Then show that a toric map in the coordinates $(w_j)_{j\in J}$ which resolves the singularities of the finitely
many binomial equations of $\hbox{\rm gr}_{\hat\nu}\hat R^{(\nu)}$ involving these finitely many variables will uniformize $\hat\nu$ on $\hat R^{(\nu)}$,
and that such resolving toric maps exist.\par\noindent 4) Use the excellence of $R$ and the birationality of the map $\hbox{\rm gr}_\nu
R\to\hbox{\rm gr}_{\hat\nu}\hat R^{(\nu)}$ to lift this to a uniformization of
$\nu$ on
$R$.
\vfill\eject
\tableofcontents
\vfill\eject
\par\bigskip\noindent
\dedicatory{ \textit {This paper is dedicated to David Rees, on the occasion of his eightieth Birthday,\\ and to Pierre Samuel, on the occasion of his
seventy-seventh Birthday}}\par\vskip.3truecm\noindent
\section{\textbf{Introduction}} 
\hfill { \it Lord Calversham: Do you always really understand what you say, sir?}\par\noindent
\hfill{\it Lord Goring (after some hesitation): Yes, father, if I listen attentively.}\par\medskip\noindent
\hfill{Oscar Wilde, {\it An ideal husband}, third Act.}\par\noindent
 \vskip1.2truecm
This text is an attempt to clarify the problems which one meets when trying to
prove local uniformization of a valuation at a singular point of an algebraic or analytic variety along the
lines of what is done for plane branches in [G-T], by re-embedding the singularity in such a
way that the valuation can be uniformized (in an {\it embedded} manner) by a {\it single}
toric modification of the ambient space. It will help the reader to be aware of the basics of
valuation theory contained in [V1] and [Cu], of the theory of toric varieties found in [Cox] and of the
results of [T1] and [G-T].\par Originally intended for the proceedings of ``Mountains
and singularities, \break Working Week on Resolution of singularities,'' and in agreement with the
spirit that Herwig Hauser wished to give to that Working Week, this text is somewhat
programmatic, even speculative in nature since there are several statements which may
not be in their final form, and some proofs are only sketched, or even absent. The correspon\-ding statements are between asterisks. I am grateful to the
Editors of this volume who have kindly agreed to publish it in the same spirit, with some asterisks left. \par 
There is therefore no claim to give here a complete detailed proof of local uniformization for an excellent equicharacteristic local ring
with an algebraically closed residue field, but only to sketch a new and possibly useful way of looking at the problem of local
uniformization and, perhaps more importantly, at a simultaneous uniformization process for all valuations ``close enough''  to a given one. \par Given
an integral domain $R$ with fraction field $K$, a valuation of $R$ is a valuation $\nu\colon K^*\to\Phi$ of $K$, with values in a totally ordered group
$\Phi$, which is non-negative on $R$. It is also the datum of a valuation ring
$R_\nu$ of $K$ containing $R$. To such a datum $(R,\nu)$ are associated the semigroup $\Gamma=\nu(R\setminus\{0\})\subset \Phi_+\cup\{0\}$ and a
filtration by ideals
$\P_\phi(R)=\{x\in R/\nu (x)\geq
\phi\}$, the ``successor'' of $\P_\phi(R)$ being $\P_\phi^+(R)=\{x\in R/\nu (x)> \phi\}$, and therefore a graded algebra $$\hbox{\rm gr}_\nu
R=\bigoplus_{\phi\in \Gamma}\frac{\P_\phi(R)}{\P_\phi^+(R)}.$$ The problem of local uniformization is, given the local ring
$R$ of an algebraic variety or more generally an excellent local ring, and assuming that it is an integral domain, to find for each valuation $\nu$ of $R$ a {\it
regular} local $R$-algebra $R'$, essentially of finite type over $R$ and contained in $R_\nu$. We wish $R'$ to be obtained by localizing an affine chart of a
proper algebraic map $Z\to \hbox{\rm Spec}R$ which is described as precisely as possible, for example as a
composition of blowing ups with non singular centers or, according to what is proposed here, as a proper and birational toric
map with respect to some system of generators of the maximal ideal of
$R$.\par
 The basic idea here is to  view our local ring $R$ as a {\it deformation} of the graded ring $\hbox{\rm gr}_\nu R$
of $R$ with respect to the filtration associated to the valuation, and to obtain the uniformization of the
valuation $\nu$ as a deformation of an {\it initial partial toric uniformization} of the valuation $\nu_{\rm
gr}$ induced by $\nu$ on $\hbox{\rm gr}_\nu R$. In the case where the residue field of $R$ is
algebraically closed, this initial partial uniformization of  $\nu_{\rm gr}$ should contain all the combinatorics
associated to the uniformization of $\nu$. \par
The motivating example is the case of complex plane branches studied as deformations of monomial curves in [T1] and [G-T]. The
analytic algebra of a plane branch has a unique valuation, since its normalization is the valuation ring $\C\{t\}$. It has a semigroup $\Gamma$ and a graded
ring as above,  and we chose to abandon the simplicity of embedding dimension two in exchange for the simplicity of dealing only with the monomial curve
$C^\Gamma$ with algebra
$\C[t^\Gamma]$ in affine
$(g+1)$-space, where $g+1$ is the number of generators of the semigroup of the branch (in characteristic zero,
$g$ is the number of its Puiseux exponents), because the  resolution of singularities of the
monomial curve is a purely combinatorial problem. We proved that after re-embedding a
plane branch $C$ in ${\mathbf C}^{g+1}$ in the way given by elements of its algebra whose valuations generate the semigroup it specializes to the
monomial curve
$C^\Gamma$ with the same semigroup, which corresponds to the graded algebra $\hbox{\rm gr}_\nu R$ of the ring
$R$ of $C$ with respect to its unique valuation.\par The equations of the corresponding family
of curves in ${\mathbf C}^{g+1}$ are of a special form which makes it possible to prove, using the
implicit function theorem, that a toric embedded resolution of $C^\Gamma\subset{\mathbf C}^{g+1}$
is also an embedded resolution for $C\subset{\mathbf C}^{g+1}$.\par\noindent 
In higher dimensions, the first serious difficulty is
that the graded ring $\hbox{\rm gr}_\nu R$ associated to the filtration of $R$ determined by
the valuation $\nu$ is not n\oe therian in general; its structure depends on the semigroup of
the values which the valuation takes on the ring, which is not always finitely generated.\par
 The good news is that for rational valuations, that is for valuations
of a local ring $R$ such that the inclusion $R\subset R_\nu$ in the ring
$R_\nu$ of the valuation induces a trivial residue field
extension\footnote{this meaning is different from that used by Abhyankar,
e.g., in [A3], where it means ``of rational rank one''. Abhyankar writes ``residually rational'' for what I call ``rational''.}, the space $\hbox{\rm
Specgr}_\nu R$ is, up to a very simple deformation, a toric variety, of finite
Krull dimension but possibly of infinite embedding dimension (see Proposition
\ref{binomialideal}); it is this toric variety which generalizes the monomial
curve of [T1] and [G-T].  The uniformization of the valuation $\nu_{\rm gr}$ induced by $\nu$ on $\hbox{\rm gr}_\nu R$ amounts to resolving the
singularities of this toric variety, which is a combinatorial problem (see \S \ref{toricmod}) with a simple solution in the finitely generated case
(see subsection \ref{toricresfin}) and no solution in the usual sense otherwise, but fortunately we need only to partially resolve in that case; see below.
Moreover, it is blind to the characteristic. The combinatorial structure of the binomial equations defining the toric variety reflects the complexity of the
valuation in $R$, including its rational rank and its rank (or height). The fact that however complicated this structure the variety defined
by the binomial equations involving finitely many variables can be resolved by a toric map is crucial. \par In order to
relate the uniformizations of
$\nu$ and
$\nu_{\rm gr}$, the main tool is the {\it  valuation algebra} (see \ref{valalg}) associated to a valuation $\nu$ of a domain $R$; its
spectrum is the total space of  a (faithfully flat) specialization of $\hbox{\rm
Spec}R$ to $\hbox{\rm Spec}\hbox{\rm gr}_\nu R$ (see Proposition \ref{faithflat}). \par 
However, just like in the plane branch case, where we had to re-embed our branch in $\C^{g+1}$ before we could write the
specialization, I need to write explicit equations for this deformation, in suitable variables. \par
This leads to the second serious difficulty, that of finding suitable completions of $R$, inside which one can lift generators of the graded $k$-algebra
$\hbox{\rm gr}_\nu R$ to coordinates, as Cohen's Theorem does for the usual graded algebra $\hbox{\rm gr}_m R$ and the usual $m$-adic completion. Here
the key result is that Laurent monomials
$\overline\eta_j=\overline\xi^{d(j)}$ in the generators
$(\overline\xi_i)_{i\in I}$ of the graded
$k$-algebra
$\hbox{\rm gr}_\nu R$ generate the graded algebra of a suitable completion
$\hat R^{(\nu)}$ of $R$ with respect to an extension $\hat\nu$ of the valuation $\nu$, and can be lifted to elements
$(\eta_j)_{j\in J}$ of
$\hat R^{(\nu)}$, which are in some sense ``coordinates''; the surjection $$k[(W_j)_{j\in J}]\longrightarrow\hbox{\rm gr}_{\hat\nu}\hat R^{(\nu)},\
\
\ W_j\mapsto\overline\eta_j$$  expressing the fact that the $\overline\eta_j$ generate the $k$-algebra
$\hbox{\rm gr}_{\hat\nu}\hat R^{(\nu)}$ lifts to a surjection
$$\widehat{k[(w_j)_{j\in J}]}\longrightarrow\hat R^{(\nu)},\ \ \
w_j\mapsto\eta_j.$$ 
Here the first candidate for $\hat R^{(\nu)}$ is the completion $\hat R^\nu$ of $R$  with respect to
the topology defined by the rational valuation $\nu$. However this is not enough because the quotients of this $\nu$-adic
completion are not complete for the quotient topology in general and therefore $\nu$-adic completeness does not suffice to ensure convergence of series 
with terms of increasing valuation (see subsection \ref{scac}). This difficulty is circumvented by extending the valuation $\nu$ to a valuation $\hat\nu$ of a
suitable quotient $\hat R^{(\nu)}$ of the $m$-adic completion $\hat R^m$ of $R$; it is complete for the $\hat
\nu$-adic topology and its quotients are also complete, while the natural morphism
$\hbox{\rm gr}_\nu R\to\hbox{\rm gr}_{\hat \nu} \hat R^{(\nu)}$ is birational in a \textit{scalewise} sense if the valuation $\nu$
is rational. The ring $\widehat{k[(w_j)_{j\in J}]}$ is a {\it scalewise} completion of the polynomial ring (see \ref{newspec}).\par
The next important fact is that one can describe by equations the specialization of $\hat
R^{(\nu)}$ to $\hbox{\rm gr}_{\hat \nu} \hat R^{(\nu)}$. One can choose generators (up to closure) of the kernel of the surjection
$\widehat{k[(w_j)_{j\in J}]}\to\hat R^{(\nu )}$ such that their initial
forms in a suitable sense are the binomials which generate the kernel in $k[(W_j)_{j\in J}]$ of the surjective homomorphism of
graded $k$-algebras
$k[(W_j)_{j\in J}]\to\hbox{\rm gr}_{\hat\nu}\hat R^{(\nu)}$ sending $W_j$ to the $j$-th generator. \par\noindent These generators of the kernel whose
initial forms generate the kernel of the associated graded map may be interpreted as a generalized Gr\"obner, or standard, basis for this kernel with respect
to the monomial order on the -countably many- variables $w_j$ defined by deciding that
$w^s<w^t$ if the valuation of the image of $w^s$ in $\hat R^{(\nu )}$ is less than the valuation of the image of $w^t$. From the viewpoint taken
here, however, the important thing is that they can be resolved by toric maps.\par This presentation of the scalewise complete n\oe therian
local ring $\hat R^{(\nu)}$ as a quotient of a generalized power series ring may be seen as a valuative
analogue of Cohen's structure theorem, and manifests an important difference with Spivakovsky's and Zariski's approaches. I do not start
from a set of coordinates and equations for the singularity and try to improve them, then follow them through a sequence of blowing-ups,
but I produce directly from the ring of the singularity, given a valuation, a ``natural'' embedding in a (possibly infinite-dimensional) space
equipped with a natural class of coordinates $(w_j)_{j\in J}$, liftings of (Laurent monomials in) the generators of the graded algebra $\hbox{\rm gr}_\nu R$,
and a ``natural'' set of equations, which are as simple as one might wish (deformations of binomials) and in
particular are ``nondegenerate'' in the natural coordinates (see
Propositions \ref{coordinates} and
\ref{completeequations}). A substantial part of the work is then to show that they truly are coordinates and
equations in a useful sense with respect to the given valuation.\par  This viewpoint is illustrated in subsection \ref{chexamples} by a geometric
interpretation of MacLane's ''Key polynomials'' associated to a valuation of $k(x)[y]$ extending the $x$-adic valuation of $k(x)$ and having value group
$\Qq$ as equations of plane approximations of a transcendental  plane curve which is a deformation of an essentially monomial curve of infinite embedding
dimension. 
\par
The specialization of the space defined by all the
equations to the space defined by their binomial initial forms, inside the space with coordinates $(w_j)_{j\in J}$, is
the explicit form of the geometric specialization of the ring $\hat R^{(\nu )}$ to its associated graded ring with respect to the valuation $\hat\nu$.\par
Finally the n\oe therian and henselian properties of the completion $\hat R^{(\nu)}$ are used again to ensure
that the lifting to $\hat R^{(\nu)}$ of finitely many steps of the (possibly infinite) combinatorial
process of resolution of $\hbox{\rm gr}_{\hat\nu} \hat R^{(\nu)}$ suffices to uniformize the extension
$\hat\nu$ of $\nu$ to $\hat R^{(\nu)}$.\par
Here the main point is that all but finitely many of the
equations defining $\hbox{\rm Spec}\hat R^{(\nu)}$ in the possibly infinite-dimensional affine
space with coordinates $(w_j)_{j\in J}$ serve only to express in
$\hat R^{(\nu)}$ the coordinates $(w_j)_{j\notin F}$  in terms of finitely many coordinates $(w_j)_{j\in F}$ which suffice
to generate its maximal ideal and by their valuations, to rationally generate the value group, and have some other properties. If
there were also infinitely many equations not serving this purpose, they would indefinitely add
singularities, and the method could not work. This finiteness implies that it suffices that a toric modification
{\it partially} resolves $\hbox{\rm gr}_{\hat\nu} \hat R^{(\nu)}$, in a precise sense, to uniformize $\hat\nu$. This uniformization occurs in spite of the fact
that
$\hbox{\rm gr}_{\hat\nu}\hat R^{(\nu)}$ may be of a smaller Krull dimension than $\hat R^{(\nu)}$, 
because of the ``abyssal'' phenomenon generated by the form of the equations defining $\hbox{\rm Spec}\hat R^{(\nu)}$; see
example \ref{exzar} in subsection \ref{chexamples}, subsections \ref{abyssalph}, \ref{dimthree} and
\ref{trans}, and Proposition \ref{lineapp}.\par In general, for a rational valuation $\nu$ the
Krull dimension of
$\hbox{\rm gr}_\nu R$ is equal, by a theorem of Piltant (see \ref{piltant}) to the rational rank $\hbox{\rm r}(\nu)$ of $\nu$. By
Abhyankar's inequality
$\hbox{\rm r}(\nu)\leq\hbox{\rm dim}R$. The cases of strict inequality may seem to contradict the fact that there is a specialization of $R$ to $\hbox{\rm
gr}_\nu R$, which is even faithfully flat. But the semicontinuity of fiber dimensions is proved only for morphisms which are essentially of finite
type. The abyssal phenomenon provides explicit examples of failure of this semicontinuity in the general case.
\par It is also necessary to make precise what {\it initial partial resolution} of $\hbox{\rm Spec}\hbox{\rm gr}_\nu R$
means, that is, to approximate in a useful way the (possibly infinite-dimensional) toric map which resolves our toric variety by finite-dimensional ones. This is
possible thanks to the fact that for any valuation ring $R_\nu$ with residue field $k_\nu$ and finite rational
rank $\hbox{\rm r}(\nu)$, the graded algebra $\hbox{\rm
gr}_\nu R_\nu$ is the direct limit of a system of essentially toric injective maps
between polynomial $k_\nu$-subalgebras in $\hbox{\rm r}(\nu)$ variables, and this
system can be described rather explicitly; for archimedian valuations, it is essentially
an avatar of Perron's algorithm (see section \ref{graval}).\par This result may also be seen
as a graded version of (non-embedded) local uniformization,
since local uniformization in the rational case over a field $k$ means that a valuation ring over $k$ with residue field $k$ is
a direct limit of regular local $k$-subalgebra with the same field of fractions which are essentially of finite type.\par\noindent
 Perron's algorithm also appears in Zariski's proof, in an apparently different role. 
\par
Then there remains to show that this process can be
extended to a toric modification of $R$ provided that it is excellent. This last step uses the scalewise birationality of  the map $\hbox{\rm gr}_\nu
R\to\hbox{\rm gr}_{\hat\nu}\hat R^{(\nu)}$ and an observation,
also used by Spivakovsky, which implies in the language used here that if $R$ is excellent, $\hbox{\rm Spec}\hat R^{(\nu)}$ is not contained in the
singular locus of
$\hbox{\rm Spec}\hat R^m$ (see Proposition
\ref{genreg}). 
\par
Another difference with Zariski's approach is that I must go directly to the ``worst''
valuations from the viewpoint of the value group, the rational ones, to take advantage of
the toric nature of the corresponding graded algebra, while  Zariski and Spivakovsky both
begin by the case of height one (archimedian) valuations (see [Z1], top of page 858). That
the center of gravity of the proof, as Zariski puts it, is now the uniformization of rational
valuations has the great advantage that it makes sense to assume that the residue field is
algebraically closed, whereas in Zariski's approach, in order to deal with valuations of
height $>1$, one localizes, and so one cannot make such an assumption. We may therefore, in
this first phase, confine our attention to the case of an algebraically closed residue
field. The rational valuations are also the only ones which appear naturally in the
complex-analytic case.\par The cases of valuations which are not rational are essentially deformations of the rational case,
as I show using the simple behavior of the valuation algebra under composition of
valuations. The reason is that, if we assume that the residue field of $R$ is algebraically
closed, for any valuation
$\nu_0$ there are rational valuations $\nu$ composed with $\nu_0$, of which $\nu_0$ may
be deemed to be a deformation.\par
The specialization to the graded ring
changes the field of fractions, which neither Zariski's proof in [Z1] nor Spivakovsky in [S2] did, but it brings us to a situation where two of the principal
complications have disappeared: they are the behaviour of transformed equations under translation after a birational map (see [Z1], pp. 884-885), and the
proof that one reaches non singularity after finitely many steps, using the fact that some partial derivative is not zero (see [Z1], pp. 886-887). Both are
replaced by the results of section \ref{toricmod}, and the implicit function theorem then allows us, via the abyssal phenomenon, to pass from the graded
ring to a suitable completion of the original ring. From there one finally gets back to the ring itself using excellence.
\par I believe that this strategy clarifies the nature of the
problem by separating as much as possible combinatorial problems from choices of coordinates
and equations. Once the toric setup is understood, finding centers of blowing-up in a given
embedding of the singularity should be considerably easier. \par\noindent  
\par  This work was initially partly motivated by the desire to understand
from the viewpoint of [T1] the inspiring attempt of Spivakovsky in [S2] to prove local uniformization and resolution of singularities;
because of this change in approach, substantial differences soon appeared, which were outlined in [G-T] for the case
of plane branches. The viewpoint developed here is
therefore different from that of [S2], although I am indebted to its author for some specific statements, as the text evidences.  Spivakovsky generalizes to
valuations of height one the construction which for plane branches describes resolution as a composition of {\it plane} toric maps (see sections 7 and 8 of
[G-T]), each of these being in fact described as a composition of blowing ups of points in the plane rather than as a toric map. The graded algebra $\hbox{\rm
gr}_\nu R$ also plays an important, although quite different, role (see the introduction and \S 7 of [G-T], as well as [GP]). Many
difficulties in this approach are due to the fact whose avatar for curves is that after each
plane toric map one must choose again an origin and local coordinates in which the next toric map will be
described. I take a systematically toric approach and allow from the start a change of embedding, even to infinite dimensions, such
that in this new embedding a {\it single} toric map in natural coordinates will uniformize, and I do not seek, in the first phase, to construct a
sequence of blowing-ups with non singular centers. The existence of such a sequence will follow from a theorem of De Concini-Procesi ([DC-P]) on the
domination of toric birational maps by sequences of blowing ups.\par In the redaction, some facts must be applied to valuation rings $R_\nu$, which
are not n\oe therian in general, and so parts are written without that assumption. Basically, 
I try to do here what is needed for local uniformization of an arbitrary valuation on a
local equicharacteristic excellent integral domain with an
algebraically closed residue field. The case where the residue field is not algebraically closed
requires other techniques in this approach. In the first sections I have to assume the existence of
a base field, since it is only after completion that one is sure to have a field of
representatives.\par
The point of view of this paper also avoids the use of the ramification and extension theory of valuations in the study of the local
uniformization problem, replacing it by the relations between the ring and its associated graded ring. The ramification over non singular spaces is in general
much worse for the graded ring, but that is of no importance since it corresponds to a toric variety and as such can be resolved by toric maps unaffected by
the ramification. There have been remarkable advances in the ramification theory (see [C-P], [K1], [K-K]) and the extension theory (see [V2]), some of which
indicate that the graded rings may play a useful role there too, especially when the graded ring extension associated to a separable extension of valued rings
is inseparable.
\par
It was during
Heisuke Hironaka's wonderful Summer meeting at Ohnuma-Hokkaido in 1993 that this train of thought
started. I wish to thank Antonio Campillo, Dale Cutkosky, David Eisenbud, Pedro Gonz\'alez P\'erez, Heisuke Hironaka, Franz-Viktor Kuhlmann, Monique
Lejeune-Jalabert, Olivier Piltant, Mark Spivakovsky, Bernd Sturmfels and Michel Vaqui\'e for stimulating comments, some
important corrections, and suggestions. I am especially grateful to Mark Spivakovsky,  whose
work rekindled the ever living embers of my interest for the role of the monomial curve and
the Riemann-Zariski manifold, and to Franz-Viktor Kuhlmann, Michel Vaqui\'e and the
referees who made very useful comments on a very preliminary version and later on a more elaborated one. I am also grateful to Nicolas
Bourbaki for his incomparably careful and lucid exposition of completion problems. \par
 During the preparation of this work, I enjoyed the
hospitality of Brandeis University, the Fields Institute in Toronto, the University of
Valladolid, Tokyo Metropolitan University, the MSRI at Berkeley and the University of Saskatchewan at Saskatoon.\par
Finally,
I wish to thank warmly Herwig Hauser for conceiving the idea of the Tirol meeting and for his
determination in realizing it, and Murray Marshall, Franz-Viktor and Salma Kuhlmann for
doing the same for the Saskatoon Conference and Workshop. \goodbreak
 \section{\textbf{Valuative sorites}}
This section\footnote{Its title is that of an embryo version of this
paper. The following
 was added at the request of a referee: the origin of the word {\it sorite} is an abbreviation
(found e.g. in Galien) of the Greek {\protect\selectlanguage{greek} swre'ithc sullogism'oc} which
refers to a syllogistic type of reasoning based on the accumulation of premises;
{\protect\selectlanguage{greek} swre'uw} {\selectlanguage{english} means }``I accumulate''. It is
used by some mathematicians to designate such an accumulation, adding a slightly depreciatory
nuance. According to the Concise Oxford Dictionary, it also evokes ``a form of sophism leading by
gradual steps from truth to absurdity''.} introduces the specialization of a ring $R$ to the graded ring
associated to a valuation. It generalizes to filtrations indexed by a totally ordered group a
construction found, for filtrations indexed by
$\Z$, in the work of Gerstenhaber [G1], [G2], Rees [R], and in [T1]. One difference
is that I do not use the valuation algebra as an algebraic construction describing the
interaction of the filtration with the multiplication law, but rather as a way to encode a geometric object. 
\subsection{The valuation algebra}\label{valalg}
 Recall that a valuation of a commutative field $K$ is a map $\nu\colon K^*\rightarrow \Phi$ to a totally
ordered commutative group $\Phi$ (the value group) such that $\nu (xy)=\nu (x)+\nu (y),\
\nu (x+y)\geq \hbox{\rm inf}\bigl(\nu (x),\ \nu (y)\bigr)$, with equality if $\nu (x)\neq \nu
(y)$. One may extend the map to $K$ by giving to the element $0\in K$ a valuation $\infty$ greater than any element of $\Phi$. The set
$\{x\in K\mid \nu (x)\geq 0\}$ is the ring $R_\nu$ of the valuation; it is local with maximal
ideal $m_\nu =\{x\in K\mid \nu (x)> 0\}$. If $\nu\colon R\to \Phi_+\cup \{0\}\cup\{\infty\}$ is defined on a commutative ring $R$ and the first
condition is relaxed to
$\nu (xy)\geq
\nu (x)+\nu (y)$, one has a {\it loose valuation} of $R$; the constructions explained below in this section are valid for $R/\nu^{-1}(\infty)$ if one replaces
``valuation'' by ``loose valuation'', the only difference being that the graded algebras are not integral domains in that case. Loose valuations are useful
because under certain conditions a valuation of a ring $R$ induces a loose valuation on a quotient ring. The trivial valuation has value group
$\{ 0\}$.\par Let $R$ be an integral domain with field of fractions $K$ and let $\nu$ be a valuation of $K$
such that its valuation ring $R_\nu $ contains $R$. Except when there is explicit mention of
the contrary, all valuations are assumed to be non trivial. Let us denote by $\Phi$ the totally
ordered value group of the valuation $\nu$ and by $k_\nu$ its residue field $R_\nu/m_\nu$.
Denote by $\Phi_+$ (resp. $\Phi_-$) the semigroup of positive (resp. negative) elements of
$\Phi$ and set $\Gamma=\nu (R\setminus \{0\})\subset \Phi_+\cup \{0\}$; it is the {\it
semigroup} of $(R,\nu )$; since $\Gamma$ generates the group $\Phi$, it is cofinal in the
ordered set $\Phi_+$.\par\noindent  For $\phi\in \Phi$, set $$\begin{array}{ccl}                        
 \P_\phi(R) &= & \{x\in R\mid\nu (x)\geq \phi \}  \\  &&\\                  
 \P^+_\phi(R) &= & \{x\in R\mid \nu (x)>\phi , \}
\end{array}$$
 where we agree that $0\in \P_\phi$ for all $\phi$, since its value is larger than any $\phi$,  so
that by the properties of valuations the $\P_\phi$ are ideals of $R$. Note that the intersection
$\bigcap_{\phi\in
\Phi_+}\P_\phi=(0)$ and that if
$\phi$ is in the negative part $\Phi_-$ of
$\Phi$, then ${\mathcal P}_\phi (R)={\mathcal P}^+_\phi (R)=R$. Consider the graded algebra {\it \`a la}
Gerstenhaber-Rees (see the papers quoted above), which we may call the {\it valuation
algebra} of $(R,\nu)$: $${\mathcal A}_\nu (R)=\bigoplus_{\phi \in \Phi} {\mathcal P}_\phi (R)
v^{-\phi}\subset R[v^\Phi]$$ where $R[v^\Phi]$ is the group algebra of $\Phi$ with coefficients
in $R$, the element $x_\phi [\phi]$ being written $x_\phi v^{\phi}$.\par{\it When
there is no ambiguity we shall denote this algebra simply by
${\mathcal A}$}; it is a domain, not n\oe therian in general even if $R$ is, as we shall see below. Note
that when $\nu$ is the trivial valuation, with $\Phi =\{0\}$, we have ${\mathcal A}_\nu (R)=R$.
\subsection{The valuation ${\mathbf \mu_\A}$}\label{valmu}
We now define a valuation $\mu_{\mathcal A}$ on ${\mathcal A}$ with value group $\Phi$ by $$\mu_{\mathcal
A}(\sum x_\phi v^{-\phi})=\hbox{\rm min}\bigl(\nu (x_\phi)-\phi \bigr) .$$ By
construction of ${\mathcal A}$ this valuation is non-negative on ${\mathcal A}$. \par\noindent
There are several rings connected with ${\mathcal A}_\nu (R)$: its quotient field is the
quotient field $K(v^\Phi)$ of $K[v^\Phi]$, and it is not difficult to describe the valuation
ring ${\mathcal R}_{\mu_{\mathcal A}}$ of $\mu_{\mathcal A}$; it is the set of quotients $(\sum x_\phi
v^{-\phi}) (\sum y_\psi v^{-\psi})^{-1}$ such that $\hbox{\rm min}\bigl(\nu (x_\phi
\bigr)-\phi )\geq \hbox{\rm min}\bigl(\nu (y_\psi )-\psi \bigr)$. \par
We set ${\mathcal P}_\phi (R_\nu)=\{x\in R_\nu\mid\nu(x)\geq
\phi\}$ as above, and similarly for ${\mathcal P}^+_\phi (R_\nu)$ and ${\mathcal P}_\phi
(K)$, and remark that ${\mathcal P}_\phi (K)$, ${\mathcal P}^+_\phi (K)$ are $R_\nu$-submodules of
$K$, and that for $\phi\in \Phi_+$ we have ${\mathcal P}_\phi (K)={\mathcal P}_\phi (R_\nu)$. With this
notation, we have:
 $${\mathcal R}_{\mu_{\mathcal A}}\cap K[v^\Phi ]={\mathcal A}_\nu (K)=\bigoplus_{\phi\in \Phi_-}{\mathcal
P}_\phi (K) v^{-\phi}\bigoplus R_\nu\bigoplus\bigoplus_{\phi\in \Phi_+}{\mathcal P}_\phi (R_\nu)
v^{-\phi}$$ and 
$${\mathcal R}_{\mu_{\mathcal A}}\cap R_\nu[v^\Phi ]={\mathcal A}_\nu (R_\nu)=\bigoplus_{\phi \in \Phi}
{\mathcal P}_\phi (R_\nu) v^{-\phi}\subset R_\nu[v^\Phi].$$
\par\medskip\noindent  The maximal homogeneous ideal of ${\mathcal A}_\nu (R_\nu)$ containing
$m_\nu$ is:
$${\mathcal M}_\nu=\bigoplus_{\phi\in \Phi_-}R_\nu v^{-\phi}\bigoplus
m_\nu\bigoplus\bigoplus_{\phi\in \Phi_+}{\mathcal P}_\phi (R_\nu)v^{-\phi},$$ so that the
corresponding residue field is the residue field $k_\nu$ of the valuation, while the ideal of
elements of positive value is
$$\tilde{\mathcal Q}_0^+=\bigoplus_{\phi\in \Phi_-}R_\nu v^{-\phi}\bigoplus
m_\nu\bigoplus\bigoplus_{\phi\in \Phi_+}{\mathcal P}^+_\phi (R_\nu) v^{-\phi}.$$ The maximal
homogeneous ideal of ${\mathcal A}_\nu (K)$ containing $m_\nu$ is
$$\bigoplus_{\phi\in \Phi}{\mathcal P}^+_\phi (K) v^{-\phi}$$ which induces $\tilde{\mathcal Q}^+_0$ on
${\mathcal A}_\nu (R_\nu)$.
\medskip\noindent In any case, the ideal of elements of positive $\mu_{\mathcal A}$-value of ${\mathcal
A}_\nu(R)$ is
$${\mathcal Q}_0^+=\bigoplus_{\phi\in \Phi_-}R v^{-\phi}\bigoplus (m_\nu \cap
R)\bigoplus\bigoplus_{\phi\in \Phi_+}{\mathcal P}^+_\phi (R) v^{-\phi},$$ and we remark that it is
exactly the ideal of
${\mathcal A}_\nu(R)$ generated by the family of elements
$(v^{\Phi_+})$. Indeed, it is generated by elements $x_\phi v^{-\phi}$ with $\nu (x_\phi)>\phi$;
but then, 
$$x_\phi v^{-\phi}=v^{\nu (x_\phi )-\phi }x_\phi v^{-\nu (x_\phi)}\in (v^{\Phi_+}){\mathcal
A}_\nu(R).$$
The other inclusion is clear.\par\medskip\noindent
For any commutative ring $A$, let us denote by
$A[v^{\Phi_+}]$ the algebra of the semigroup $\Phi_+\cup \{0\}$ with coefficients in $A$.
\subsection{The specialization to the graded algebra}\label{speciatogr}
 The graded algebra associated with the valuation $\nu$ was introduced in ([L-T], [T1]) for the
very special case of a plane branch (see [G-T]), and in [S1], [S2] in full generality. It is
$$\hbox{\rm gr}_{\nu} R =\bigoplus_{\phi\in \Gamma}{\mathcal P}_\phi (R)/ {\mathcal P}_\phi^+(R).$$\noindent For $\phi\notin \Gamma$,
$\P_\phi(R)=\P_\phi^+(R)$; it is sometimes convenient to view it as graded by $\Phi_+\cup \{0\}$, with a zero homogeneous component for
$\phi\notin
\Gamma$.
\par\noindent For each non zero element $x\in R$, there is a unique $\phi\in\Gamma$ such that
$x\in
\P_\phi\setminus \P_\phi^+$; the image of $x$ in the quotient ${(\hbox{\rm gr}_{\nu} R)}_\phi
={\mathcal P}_\phi (R)/ {\mathcal P}_\phi^+(R) $ is the {\it initial form} $\hbox{\rm in}_\nu x$ of $x$. When there is no
ambiguity, I will often denote it by $\overline x$.\par\noindent
The graded ring $\hbox{\rm gr}_{\nu} R$ inherits naturally from $(R,\nu)$ a valuation $\nu_{\rm gr}$, defined by
$$\nu_{\rm gr} (\sum
\overline{x_\phi})=\hbox{\rm min}(\nu (x_\phi)),$$ where $x_\phi\in R$ is a representative 
of $\overline{x_\phi}$, and the valuation does not depend on its choice. 
\begin{Remark} 1) Since $R$ is an integral domain and $\nu$ is a valuation, the ring $\hbox{\rm
gr}_\nu R$ is an integral domain.
\par\noindent 2) If the value group
$\Phi$ has finite rational rank,  it is countable, since it is torsion free and hence the map
$\Phi\rightarrow {\Qq}\otimes_{\Z}\Phi ={\Qq}^r$ is injective. Therefore the $R/(m_\nu
\cap R)$-algebra $\hbox{\rm gr}_\nu R$ has a countable system of generators, and so does the
$k_\nu$-algebra $\hbox{\rm gr}_\nu R_\nu$. \par\noindent 3) Each ${\mathcal P}_\phi (R)/{\mathcal P}_\phi^+(R)$ is a torsion free $R/(m_\nu\cap
R)$-module, which is of finite type if the ring $R$ is n\oe therian.\par\noindent 4) If $R$ is a local ring with maximal ideal $m$, by definition of
domination, the valuation ring $R_\nu$ {\it dominates} $R$ if
$m_\nu \cap R =m$, and we will sometimes be in that situation. Then $\hbox{\rm gr}_\nu R$ is a $k_R$-algebra, where $k_R$ is
the residue field of $R$. If we want
to emphasize that $R$ and $R_\nu$ have the same field of fractions, we say that $R_\nu$ {\it
birationally dominates} $R$. \end{Remark}
\par\noindent Let us denote by $(v^{\Phi_+}){\mathcal A}_\nu(R)$ the ideal
generated in ${\mathcal A}_\nu(R)$ by the elements of $(v^{\Phi_+})$.    
\begin{proposition}\label{deformation} a) The natural map
$${\mathcal A}_\nu(R)\longrightarrow \hbox{\rm gr}_\nu R$$ defined by
$$x_\phi v^{-\phi}\mapsto x_\phi \ \ \hbox{\rm mod.}{\mathcal P}^+_\phi (R)$$ induces an
isomorphism of graded rings 
$${\mathcal A}_\nu(R)/{\mathcal Q}_0^+= {\mathcal A}_\nu(R)/(v^{\Phi_+}){\mathcal A}_\nu(R)
\stackrel{\simeq}\longrightarrow\hbox{\rm gr}_\nu R.$$\par\noindent  b) The natural inclusion
$R[v^{\Phi_+}]\longrightarrow {\mathcal A}_\nu(R)$ obtained by considering only the part of negative degree
(i.e., $\phi \in \Phi_-$) of ${\mathcal A}$ induces, after taking rings of fractions by the multiplicative
subset
$(v^{\Phi_+})$, an isomorphism
$$(v^{\Phi_+})^{-1}R[v^{\Phi_+}]\stackrel{\simeq}\longrightarrow (v^{\Phi_+})^{-1}{\mathcal A}_\nu(R);$$\noindent the
quotient ${\mathcal A}_\nu(R)/R[v^{\Phi_+}]$ is a torsion $\Z[v^{\Phi_+}]$-module.\par\noindent
c) Given a homomorphism $\phi\mapsto e(\phi )$ from
$\Phi$ to the multiplicative group of units of $R$, it induces a surjection ${\mathcal A}_\nu(R)\to  R$
defined by $\sum x_\phi v^{-\phi}\mapsto \sum x_\phi e(-\phi )$. The kernel of this
surjection is the ideal $\bigl((v^\phi -e(\phi ))_{\phi\in \Phi_+}\bigr){\mathcal A}_\nu (R)$. If
$R$ contains a field
$k$, this is true in particular for homomorphisms $\Phi \to k^*$.
\end{proposition}
\begin{proof}  the map in a) is clearly surjective, and its kernel is homogeneous and generated
by elements $x_\phi v^{-\phi}$ with $\nu (x_\phi)>\phi$; but we saw in \ref{valmu} that it is\break
$(v^{\Phi_+}) {\mathcal A}_\nu(R)=\Q_0^+$. \par\noindent To prove b), consider an element
$v^{-\psi}\sum_{\phi
\in\Phi_-}  x_\phi v^{-\phi}$ of the ring of fractions
\break $(v^{\Phi_+})^{-1}R[v^{\Phi_+}]$, and rewrite it as $v^{-\psi}\sum_{\phi
\in\Phi_-}  v^{-\nu (x_\phi )} x_\phi v^{-\phi+\nu (x_\phi)}$. Since $-\phi \geq 0$ and
$\nu(x_\phi)\geq 0$, this is clearly an element of
$(v^{\Phi_+})^{-1}{\mathcal A}_\nu (R)$, and every element can be written in this way; if $\sum
x_\phi v^{-\phi}\in {\mathcal A}_\nu (R)$, let $\phi_+$ be the largest degree appearing in this sum. If
it is negative, the element is in $R[v^{\Phi_+}]$. Otherwise, rewrite the sum as
$v^{-\phi_+}\sum x_\phi v^{-(\phi-\phi_+)}\in (v^{\Phi_+})^{-1}R[v^{\Phi_+}]$. This proves both
assertions of b). Assertion c) follows from the remarks that any $x\in R$ of valuation $\phi$ is the image of $e(\phi)xv^{-\phi}$ and that if $\sum x_\phi
e(-\phi )=0$, then $\sum x_\phi v^{-\phi}=\sum x_\phi (v^{-\phi}-e(-\phi))$.
\end{proof}
\par\medskip  To interpret this geometrically, let us assume in this section that
$R$ contains\footnote{In this paper, whenever I introduce a base field $k$, it will be tacitly
assumed that the valuations studied are trivial on $k$.} a field $k$ ; we have a
composed map of algebras $k[v^{\Phi_+}]\rightarrow R[v^{\Phi_+}]\rightarrow {\mathcal A}_\nu (R)$.
Then: 
\begin{proposition}\label{faithflat} The  $k[v^{\Phi_+}]$-algebra ${\mathcal A}_\nu(R)$ is faithfully flat. 
\end{proposition}
\begin{proof}
The criterion for faithful flatness of [B3], Chap. I, \S 3, No. 7, Proposition 13 is
that every solution
$(y_k\in \A_\nu (R), \ 1\leq k\leq n)$ of a system
$$\sum _{k=1}^nc_{ki}y_k=d_i\ \ (1\leq i\leq m)\leqno{(*)}$$\noindent of linear equations with
coefficients and right-hand sides in $k[v^{\Phi_+}]$ can be written in the form
$$y_k=x_k+\sum_{j=1}^pz_{jk}e_j\ \ \ (1\leq k\leq n)$$\noindent where $(x_k)$ is a solution of
$(*)$ in $k[v^{\Phi_+}]$, the $e_j$ belong to $\A_\nu (R)$ and the $z_{jk}$ for each fixed $j$ are
solutions in $k[v^{\Phi_+}]$ of the homogeneous system associated to $(*)$.\par\noindent Writing
$y_k=\sum y_{k\psi}v^{-\psi}$, we notice that the elements $y_{k\psi}$ with $\psi>0$ are such that $\nu (y_{k\psi})>0$ and therefore if we write
$$y_k=\sum_{\psi\leq 0} y_{k\psi}v^{-\psi}+\sum_{\psi> 0} y_{k\psi}v^{-\psi}=y_k^{\leq 0}+y_k^{> 0}$$
we must have, since the $d_i$ are in $k[v^{\Phi_+}]$ and so have all their coefficients of valuation zero,
$$\hbox{\rm a)}\ \ \sum _{k}c_{ki}y^{\leq 0}_k=d_i\ \ ,\ \ \ \ \hbox{\rm b)}\ \ \ \sum
_{k}c_{ki}y^{>0}_k=0\ \ (1\leq i\leq m).$$\noindent We may view the system a) as a system
of linear equations with coefficients and right-hand side in
$k[v^{\Phi_+}]$, and for which we have a solution $y_k^{\leq 0}\in R[v^{\Phi_+}]$. By [B3], Chap. 1, \S 3, No.3, Prop. 5, the
$k[v^{\Phi_+}]$-module
$R[v^{\Phi_+}]$ is faithfully flat 
 since it is obtained from the faithfully flat map $k\to R$ by the base extension $k\to k[v^{\Phi_+}]$.  By the criterion quoted above, we may write
$$y^{\leq 0}_k=x^{\leq 0}_k+\sum_{j=1}^pz^{\leq 0}_{jk}e^{\leq 0}_j\ \ \ (1\leq k\leq n)$$\noindent with the same conditions as above except
that $\A_\nu (R)$ is replaced by $R[v^{\Phi_+}]$.\par Let us now consider the second system; here we
have to deal with restrictions on the valuations of the elements; a solution of this system defines
a finite dimensional $k$-vector subspace $V$ of $R$, generated by the coefficients
$y^{>0}_{k\psi}$ of the solutions $y^{>0}_k=\sum_\psi y^{>0}_{k\psi}v^{-\psi}$ of system b). Let us choose a basis $(\tilde e^{>0}_q)_{1\leq q\leq r}$ of this
vector space which is compatible with the finite filtration of $V$ by the vector subspaces
$\P_\phi\cap V$. We can now write $$y^{>0}_{k\psi}=\sum_{q=1}^rw^{>0}_{qk\psi}\tilde e^{>0}_q\
\hbox{\rm with }\ w^{>0}_{qk\psi}=0\ \hbox{\rm if }\nu (\tilde e^{>0}_q)<\psi ,$$ and the coefficients
$w^{>0}_{qk\psi}$ satisfy the system of linear equations $$\sum_k c_{ki}(\sum_\psi
w_{qk\psi}^{>0}v^{-\psi})=0 \ \ \  \forall q,\ \forall i.$$ This means that we may write
$$y^{>0}_k=\sum_q\bigl(\sum_\psi w^{>0}_{qk\psi}v^{\nu (\tilde e^{>0}_q)-\psi}\bigr)\tilde e^{>0}_qv^{-\nu
(e^{>0}_q)},$$ and this is indeed a combination of elements of $\A_\nu (R)$ with coefficients in $
k[v^{\Phi_+}]$ satisfying b). Let us denote by $e^{>0}_q$ the product $\tilde e^{>0}_qv^{-\nu
(\tilde e^{>0}_q)}$ and by $z^{>0}_{qk}$ the sum $\sum_\psi w^{>0}_{qk\psi}v^{\nu (\tilde e^{>0}_q)-\psi}\in k[v^{\Phi_+}]$. Adding up the
decompositions:
$$y_k=y_k^{\leq 0}+y_k^{>0}=x^{\leq 0}_k+\sum_{j=1}^pz^{\leq 0}_{jk}e^{\leq 0}_j+\sum_{q=1}^rz^{>0}_{qk}e^{>0}_q,$$ we see that
we have satisfied the criterion for faithful flatness.
\end{proof}
\begin{remark} Let
$(\xi_i)_{i\in I}$ be elements of $R$ such that their images (initial forms) in $\hbox{\rm
gr}_\nu R$ generate it as a $R/(m_\nu\cap R)$-algebra. Let us denote by $S$ the
$k$-subalgebra of $R$ generated by the elements $(\xi_i)_{i\in I}$ and let $\B_\nu (S)$ the sub
$k[v^{\Phi_+}]$-algebra of $\A_\nu (S)$ generated by the elements $\xi_iv^{-\nu (\xi_i)}$. The
proof we have just seen can be adapted to show that the $k[v^{\Phi_+}]$-algebra $\B_\nu (S)$ is
faithfully flat.
\end{remark}
\par\bigskip The inclusion
$k[v^{\Phi_+}]\rightarrow {\mathcal A}$ defines a map of schemes $$\hbox{\rm Spec}{\mathcal
A}\longrightarrow \hbox{\rm Spec} k[v^{\Phi_+}].$$ Part a) of Proposition
\ref{deformation} is equivalent to the fact that the special fiber of this map is $\hbox{\rm
Spec}\hbox{\rm gr}_\nu R$ and parts b) and c) are ways of stating that its generic and
general fibers are isomorphic to $\hbox{\rm Spec}R$. Since this map is faithfully flat by
Proposition \ref{faithflat}, we can say that we have made $R$ appear as a deformation of its
associated graded ring, in a family parametrized by $\hbox{\rm
Spec}k[v^{\Phi_+}]$.
\begin{remark} Without assuming that $R$ contains a field of representatives, we have an isomorphism
$$\A_\nu(R)/m\A_\nu (R)\simeq k[v^{\Phi_+}]$$ corresponding to an inclusion
$\hbox{\rm Spec}k[v^{\Phi_+}]\subset \hbox{\rm Spec}\A_\nu(R)$.\par\noindent The localization of $\hbox{\rm Spec}\A_\nu(R)$ at the general
point corresponding to localization by $(v^{\Phi_+})$ induces a general point of this subscheme and  produces the ``general fiber'' seen above,
which is essentially 
$\hbox{\rm Spec}R$, while the subscheme of 
$\hbox{\rm Spec}\A_\nu(R)$ defined by the inverse image in $\A_\nu (R)$ of the ideal $(v^{\Phi_+})k[v^{\Phi_+}]$ is $\hbox{\rm
Spec}\hbox{\rm gr}_\nu R$.
\end{remark}
\par\medskip The center of the valuation $\mu_{\mathcal A}$ on $\hbox{\rm Spec}{\mathcal A}_\nu (R)$
is the subscheme $\hbox{\rm Spec}\hbox{\rm gr}_\nu R$; as we saw,  this subscheme is integral since $\nu$ is a valuation and not a
loose valuation.\par\medskip
 Let us now compute the associated graded ring of
${\mathcal A}=\A_\nu (R)$ with respect to the valuation $\mu_{\mathcal A}$. As usual I implicitely add
$\{ 0\}$ to the semigroups.\par\noindent Set $\P_{\phi}=\P_{\phi}(R)$ and
$${\mathcal Q}_\delta=\{a\in {\mathcal A}\mid\mu_{\mathcal A}(a)\geq
\delta\};$$ it is a graded ideal generated by the $\{x_\phi v^{-\phi}\mid\nu (x_\phi)-\phi \geq
\delta\}$, and similarly for
$\Q^+_\delta$. So we have
$$\Q_\delta=\A\cap \bigoplus_{\phi \in \Phi}\P_{\phi +\delta}v^{-\phi},$$
$$\Q^+_\delta=\A\cap\bigoplus_{\phi\in \Phi}\P^+_{\phi +\delta}v^{-\phi}.$$ So $\Q_\delta =\A$
for $\delta\leq 0$ and we can write:
$$\hbox{\rm gr}_{\mu_\A}\A=\bigoplus_{\delta \in \Phi_+\cup\{ 0\}} (\Q_\delta /\Q^+_\delta )
=\bigoplus_{\delta\in
\Phi_+}\bigl(\bigoplus_{\{\phi \mid\phi +\delta \in \Phi_+\cup\{0\}\}}\P_{\phi
+\delta}/\P^+_{\phi+\delta}\bigr) .$$ Thus $\hbox{\rm gr}_{\mu_\A}\A$ has a
$\Phi_+\oplus\Phi_+$-graduation by $(\delta, \phi+\delta)$  if we give degree $\phi +\delta$ to the elements of
$\P_{\phi+\delta}/\P^+_{\phi+\delta}$ as suggested by the graduation of ${\mathcal A}$.\par Let us denote
by $\overline R$ the quotient $R/(m_\nu \cap R)$, and remark that $\hbox{\rm gr}_\nu R$ is an
$\overline R$-algebra.\par We have a map of $\overline R$-algebras $$\hbox{\rm
gr}_{\mu_\A}\A\rightarrow (\hbox{\rm gr}_\nu R)\otimes_{\overline R}\overline R
[v^{\Phi_+}]$$ induced by the map $\Q_\delta \rightarrow (\hbox{\rm gr}_\nu
R)\otimes_{\overline R}\overline R [v^{\Phi_+}]$ defined by:
$$\sum_{\nu (x_\phi)-\phi \geq\delta} x_\phi v^{-\phi}\mapsto
\sum_{\nu (x_\phi)-\phi =\delta}\overline{x_\phi}\otimes v^{\delta},$$ where
$\overline{x_\phi}$ is the initial form in $\P_{\phi +\delta}/\P^+_{\phi+\delta}$ of
$x_\phi$. This map is homogeneous with respect to the $\Phi_+\oplus\Phi_+$ 
graduation of $\hbox{\rm gr}_{\mu_\A}\A$ described above if we give
\break $(\hbox{\rm gr}_\nu R)\otimes_{\overline R}\overline R[v^{\Phi_+}]$ the tensor product
graduation where elements of degree $\delta\oplus
\psi\in \Phi_+\oplus \Phi_+$ are tensor products with $v^{\delta}$ of elements of degree
$\psi$ of $\hbox{\rm gr}_\nu R$.
\begin{proposition} With the gradings just
described, this map is a {\rm ($\Phi_+\oplus \Phi_+$)}- graded isomorphism of $\overline
R$-algebras $$\hbox{\rm gr}_{\mu_\A}\A\stackrel{\simeq}\rightarrow (\hbox{\rm gr}_\nu
R)\otimes_{\overline R}\overline R[v^{\Phi_+}].$$
\end{proposition}
\begin{proof}
The map is clearly surjective, and its injectivity is easily
verified degree by degree. Note that upon taking $\phi +\delta =0$, this map
induces the identity of
$\overline R[v^{\Phi_+}]$. Note also that the valuation $\mu_\A$ takes
its positive values in $\Phi_+$, and not in
$\Gamma$.\end{proof}
\begin{remark} The meaning of this proposition is that the deformation
of graded $\overline R$-algebras corresponding as in Proposition \ref{deformation} to the valuation algebra associated to $\A$ and to
$\mu_\A$ is ``trivial'' in the sense that the total space $\hbox{\rm Spec}\hbox{\rm
gr}_{\mu_\A}\A$ is the product by the base space $\hbox{\rm Spec} \overline R[v^{\Phi_+}]$ of the special fiber $\hbox{\rm
Spec}\hbox{\rm gr}_\nu R$.
\end{remark}
\subsection{The valuation ${\mathbf \nu_\A}$}\label{nuA}
\par\medskip\noindent There  is another natural valuation on the ring $\A_\nu (R)$: it is the
valuation $\nu_\A$ defined by  $$\nu_\A(\sum x_\phi v^{-\phi})=\hbox{\rm min}(\nu(x_\phi
)).$$ Its center is the subscheme of $\hbox{\rm Spec} \A_\nu (R)$ defined by the ideal of
$\A_\nu (R)$ generated by $(m_\nu\cap R,\bigoplus_{\phi\in \Phi_+}\P_\phi v^{-\phi })$. This
subscheme is isomorphic to $\hbox{\rm Spec}\overline R[v^{\Phi_+}]$. In the special case where
$R$ is a local k-algebra dominated by $R_\nu$ and the extension $k\rightarrow R/m$ is trivial,
this means that the center of $\nu_\A$ is the image of a section $$\sigma\colon
\hbox{\rm Spec} k[v^{\Phi_+}]\rightarrow \hbox{\rm Spec} \A_\nu (R).$$
\par\medskip\noindent Let us examine the valuation algebra of $\A_\nu(R)$ with respect to
$\nu_\A$. Set $$\Ss_\delta=\{\sum x_\phi v^{-\phi}\in \A_\nu(R)\mid\hbox{\rm min}(\nu
(x_\phi))\geq \delta \} =\bigoplus_{\phi \in \Phi} \P_{\hbox{\rm max}(\delta ,\phi )}
(R)v^{-\phi} ,$$
$$\Ss^+_\delta=\{\sum x_\phi v^{-\phi}\in \A_\nu(R)\mid\hbox{\rm min}(\nu (x_\phi))> \delta
\}=\bigoplus_{\phi \in \Phi} (\P^+_\delta (R)\cap \P_\phi (R))v^{-\phi} ;$$ we have the
equalities 
$$\hbox{\rm gr}_{\nu_\A}\A_\nu (R)=\bigoplus_{\delta\in \Phi}\Ss_\delta
/\Ss^+_\delta=\bigoplus_{\phi\in \Phi}\bigl(\bigoplus_{\delta\geq
\phi}(\P_\delta (R)/\P^+_\delta (R)) \bigr).$$ Let us define the ``degree algebra'' of a graded
algebra $H=\bigoplus _{\delta \in \Phi_+\cup
\{0\}} H_\delta$ in a way similar to the definition of the valuation algebra, replacing valuation
by degree, and setting $H_\delta=0$ for $\delta\in \Phi_{-}$:
$$\A_{\rm deg}(H)=\bigoplus_{\phi \in \Phi}\bigl(\bigoplus
_{\delta\geq\phi}H_\delta\bigr)v^{-\phi}\subset H[v^\Phi].$$ This algebra is ``trivial'' in view of
the:
\begin{lemma}\label{twist} The map
$$\A_{\rm deg}(H)\rightarrow H\otimes_{H_0}H_0[v^{\Phi_+}]$$ determined by
$$\sum_\phi\bigl(\sum_{\delta\geq \phi}x_\delta\bigr)v^{-\phi}\mapsto \sum_{\phi ,\delta \in
\Phi ,\ \delta \geq \phi} x_\delta\otimes_{H_0}v^{\delta -\phi}$$ is a graded isomorphism of
$H_0$-algebras when one gives the right hand side the twisted $\Phi\oplus
\Phi_+$ graduation for which a homogeneous element of degree $\delta\oplus \psi$ is of the
form $x_\delta\otimes v^{\delta -\psi}$ where $x_\delta$ is of degree
$\delta$.
\end{lemma}
\begin{proof} Indeed, the inverse map is given by $x_\delta\otimes v^\psi \mapsto
x_\delta v^{-(\delta -\psi )}$.
\end{proof}
\par\medskip In order to show that the valuation algebra
$\hbox{\rm gr}_{\nu_\A}\A_\nu (R)$ is ``trivial'', as in the case of the valuation $\mu_\A$, it
only remains to check the:
\par\medskip\noindent
 \begin{lemma} If we view $\hbox{\rm gr}_{\nu}(R)$ as a graded
algebra indexed by  $\Phi_+\cup \{0\}$ with a $0$ component for $\phi\notin \Gamma$, we have
the equality:
$$\hbox{\rm gr}_{\nu_\A}\A_\nu (R)=\A_{\rm deg}(\hbox{\rm
gr}_{\nu}(R)).$$\end{lemma}\begin{proof} It follows directly from the definitions.\end{proof}\par\medskip\noindent Finally, we have proved
\begin{proposition} \label{triviality}With the twisted grading on $(\hbox{\rm gr}_\nu
R)\otimes_{\overline R}\overline R[v^{\Phi_+}]$ described above, the map of Lemma \ref{twist}
is a graded isomorphism of graded $\overline R$-algebras
$$\hbox{\rm gr}_{\nu_\A}\A_\nu (R)\stackrel{\simeq}\rightarrow (\hbox{\rm gr}_\nu R)\otimes_{\overline
R}\overline R[v^{\Phi_+}],$$ determined by the fact that the $\nu_\A$-initial form of $\xi v^{-\nu
(\xi )}$ is mapped to $\overline \xi\otimes 1$, where $\overline \xi$ is the
$\nu$-initial form of $\xi$.\par\noindent Note that $$\A_\nu (R)/\Ss _0=\overline
R[v^{\Phi_+}]=\overline R\otimes_{\overline R} \overline R[v^{\Phi_+}].$$
\end{proposition}
\begin{Remark} 1) The isomorphisms of
the two graded algebras of $\A_\nu (R)$ with respect to $\mu_\A$ and $\nu_\A$ with $
(\hbox{\rm gr}_\nu R)\otimes_{\overline R}\overline R[v^{\Phi_+}]$ have quite different
geometric significations, as one can see by considering the analogous constructions for the Rees
algebra ${\mathcal R}(I)=\bigoplus_{n\in \Z}I^nv^{-n}\subset R[v,v^{-1}]$ of an ideal $I$ of a ring
$R$, where
$I^n=R$ for $n\leq 0$. The valuation $\mu$ corresponds to the $v$-adic filtration; it is the
(loose) valuation associated to this filtration, while the valuation $\nu$ corresponds to the
(loose) valuation induced by the $I$-adic filtration of $R[v,v^{-1}]$. The triviality of the
$\hbox{\rm gr}_{\mu_\A}\A_\nu (R)$ generalizes the fact that since the special fiber
$\hbox{\rm Specgr}_IR$ in this case is a divisor, its normal cone in $\hbox{\rm Spec}{\mathcal R}(I)$
is $\hbox{\rm Spec}(\hbox{\rm gr}_IR)[v]$. The triviality for $\nu_\A$ generalizes the fact, more
specific to the Rees algebra, that the graded ring of ${\mathcal R}(I)$ with respect to the $I$-adic filtration is
trivial in the sense that it is isomorphic to  $\hbox{\rm gr}_IR\otimes_{R/I}R/I[v]$. This reflects the fact that in the specialization of $\hbox{\rm Spec}R$ to
$\hbox{\rm Spec} \hbox{\rm gr}_IR$ provided by $\hbox{\rm Spec}{\mathcal R}(I)$, the graded ring is constant, or more precisely that the induced
specialization on the graded ring is trivial.
\par\noindent 2) We shall see below in subsection
\ref{sec-special} that there is a valuation $\tilde\nu_\A$ on $\A_\nu (R)$ which is composed with the natural
valuation $\nu_{\rm gr}$ on the graded ring $\hbox{\rm gr}_\nu R$. Moreover, this
valuation $\tilde\nu_\A$ induces in a natural way the valuation $\nu$ on $R$. This will
gives a precise meaning to the intuition that $\nu_{\rm gr}$ is the specialization of $\nu$
to $\hbox{\rm gr}_\nu R$.
\par\noindent 3) The graded ring
$\hbox{\rm gr}_{\nu_{\rm gr}}\hbox{\rm gr}_\nu R$ is isomorphic to
$\hbox{\rm gr}_\nu R$.\end{Remark}\noindent
{\textbf {An example:}} Let $R$ be a discrete valuation ring, and $u$ a generator of its maximal ideal. 
The value group of the valuation of $R$ is $\Z$ and we have $\P_n(R)=u^n R$ for $n\geq 0$,
$\P_n(R)=R$ for $n\leq 0$, so that  $$\A_\nu (R)=R[v,uv^{-1}]\subset R[v,v^{-1}].$$
In this case, $\hbox{\rm gr}_\nu R=(R/uR)[U]$, where $U$ is 
the initial form of $u$. I leave it as an exercise to determine the maps
$\A_\nu (R)\to R$ and $\A_\nu (R)\to \hbox{\rm gr}_\nu R$ in this case
(hint: map $u v^{-1}$ to $U$). This applies in particular when $R={\mathbf C}\{
u\}$, as we shall see in subsection \ref{planebranches}, or to $k[[u]]$ where
$k$ is any field. It also applies to the non-equicharacteristic case, for example
to $R=\hat\Z_p$, the ring of $p$-adic integers, and more generally when
$R=W(k)$ is the ring of  Witt vectors associated to a perfect field $k$ of
characteristic $p$, with $u=p$ (see [B3], Chap. IX, \S\S 1,2, [Ei], \S 7). We
shall see below in subsection \ref{Bcomp} that the $\hbox{\rm Spectrum}$ of the
$p$-adic completion of $\A_\nu (W(k))$ contains both
$\hbox{\rm Spec}k[[t]]$ and $\hbox{\rm Spec}W(k)$. This is the very first step
in the extension of the method presented here to the non-equicharacteristic
case. 
\section{ \textbf{Valuation algebras and composition of valuations}}\label{valco} 
In this section I examine the behavior of valuation algebras when one
changes a valuation for an other one with which it is composed, or for the
corresponding residual valuation. This section also contains two results which are very useful in the sequel: Zariski's description of valuation
ideals with respect to two consecutive valuations, and Piltant's theorem on the Krull dimension of $\hbox{\rm gr}_\nu R$.
\subsection{The composition of valuations}\label{compofval}
In this paragraph,
we study the relation of composition of valuations; if $\nu$ is composed with $\nu_1$, both
non-negative on our local ring $R$, it means that we have inclusions of local rings $$R\subset
R_\nu\subset R_{\nu_1},$$ and of ideals $$m_{\nu_1}\subset m_\nu\subset  R_\nu ,$$
$$m_{\nu_1}\cap R\subseteq m_\nu\cap R \subseteq m.$$ I assume known the theory of
composition of valuations, referring to ([V1], section 4, [Z-S], Vol 2, Ch. VI, [A1], [B2]), and fix
notations.\par Let us recall that the rank $\hbox{\rm  r}(\Phi)$ of the
$\Qq$-vector space
$\Qq\otimes_{\Z}\Phi$ is called the {\it rational rank} of $\Phi$; if $\Phi$ is the group
of a valuation $\nu$, it is also denoted by
$\hbox{\rm  r}(\nu )$. Recall also that in the case where $R_\nu$ dominates $R$, the
transcendence degree
$\hbox{\rm t}_{k_R} (\nu )$ of the extension
$(k_\nu :k_R)$ is called the residual transcendence degree, and that the {\it height}
$ \hbox{\rm h} (\nu )$ (often also called the {\it rank}) of the valuation $\nu$ is the number (in
general the ordinal type) of elements of a maximal chain 
$$R_\nu\subset R_{\nu_1}\subset \cdots \subset R_{\nu_{ h -1}}$$  \noindent of valuation rings
of $K$ containing $R_\nu$; it is also the Krull dimension of the ring $R_\nu$. Each inclusion
$R_\nu\subset R_{\nu_i}$ corresponds to a monotone (non-decreasing) map of totally ordered
groups
$\lambda_i\colon
\Phi\to \Phi_i$ such that $\nu_i=\lambda_i\circ \nu$. A monotone non-decreasing map means
here that if
$\phi\leq \psi$, then $\lambda (\phi )\leq \lambda (\psi )$. The kernel
$\Psi_i$ of such a map is an {\it isolated} (or {\it convex}) subgroup of
$\Phi$, so that $\hbox{\rm h} (\nu )$ is also the ordinal type of the chain of isolated subgoups of
$\Phi$; it depends only on $\Phi$ and we will speak of the height of $\Phi$. \par A totally
ordered group of height one is archimedian and can be embedded as an ordered subgroup of
${\mathbf R}$.  Following [A1], we will say that a group of height one is {\it discrete} if it is well
ordered (i.e., isomorphic to $\Z$, see [A1]), and we will say that a group $\Phi$ of 
finite height is discrete if the quotient
of two successive elements of the sequence  $$(0)\subset \Psi_{ h-1}\subset \cdots
\subset \Psi_1\subset \Phi$$ of the isolated subgroups of $\Phi$ is discrete.\par    The additivity of the rational rank in exact sequences
 implies the inequality  $$ \hbox{\rm h} (\nu )\leq  \hbox{\rm  r} (\nu)$$\noindent and for
a valuation ring $R_\nu$ birationally dominating a n\oe therian local ring $R$, we have \it
Abhyankar's inequality \rm (see [V1], section 9, Th. 9.2) $$\hbox{\rm  r}(\nu)+\hbox{\rm
t}_{k_R} (\nu) \leq\hbox{\rm dim}R.$$\par\noindent It is important to note
that given a valuation of height $\nu$ of a ring $R$, the prime ideals of $R_\nu$ form a
chain
$$ m_{\hbox{\rm h} (\nu )-1}\subset m_{\hbox{\rm h} (\nu )-2}\subset \cdots\subset
m_{\nu_1}\subset m_\nu$$ of length $\hbox{\rm h} (\nu )$ but their intersections with
$R$ are not necessarily distinct, and the length of the chain of  distinct prime ideals among the 
$$m_{\hbox{\rm h} (\nu )-1}\cap R\subseteq\cdots \subseteq m_{\nu_1}\cap R\subseteq m_\nu\cap R,$$
which may be called the {\it height of $\nu$ in $R$} and denoted by $\hbox{\rm h}_R(\nu)$, is at most $\hbox{\rm
h}(\nu)$.\par\noindent Abhyankar's inequality has a nice interpretation in terms of the
graded algebra:
\begin{proposition}\label{piltant} {\rm (Piltant, [P])} If $R$ is a n\oe therian local domain
with uncountable residue field $k_R$, birationally dominated by the valuation ring
$R_\nu$, the Krull dimension of $\hbox{\rm gr}_\nu R$ is equal to $\hbox{\rm 
r}(\nu)+\hbox{\rm t}_{k_R} (\nu)$. If $\hbox{\rm t}_{k_R}(\nu)=0$, the assumption on
the cardinality is unnecessary. \end{proposition}
\noindent
In the cases where $k_R$ is uncountable,  Abhyankar's inequality then becomes $$\hbox{\rm dim}\hbox{\rm gr}_\nu R\leq
\hbox{\rm dim}R.$$ Note that if we have $R\subset R_\nu$, the local ring $R_\nu$
dominates $R_{m_\nu \cap R}$.\par\noindent  Spivakovsky had remarked that the
transcendence degree over $R/m$ of $\hbox{\rm gr}_\nu R$  is equal to $\hbox{\rm 
r}(\nu)+\hbox{\rm t}_{k_R} (\nu)$ (see [S3], Remark 3.7).\par\medskip
  I will use Piltant's result only in the case where $\hbox{\rm t}_{k_R} (\nu)=0$.
For the convenience of the reader, I sketch Piltant's proof of his result in this case. In the
beginning, I suppose only that $m_\nu \cap R=m$, i.e. that the valuation ring $R_\nu$
dominates $R$.  Let $S$ be the multiplicative set of homogeneous elements of positive
degree in $\hbox{\rm gr}_\nu R$, and let $[F]_0$ denote the part of degree zero of an
homogenous ring of fractions $F$ of a graded algebra by a multiplicative set of
homogeneous elements.  
\begin{lemma} The map $$[S^{-1}\hbox{\rm gr}_\nu R]_0\to k_\nu,\ \ \ \frac{\overline
x_\phi}{\overline y_\phi}\mapsto  \frac{x_\phi}{y_\phi}\hbox{\rm
mod.}m_\nu$$ is well defined and is an isomorphism.
\end{lemma}
\begin{proof} This follows from a direct verification.
\end{proof}
\par\noindent
Let now $(u_1,\ldots ,u_r)$ be elements of $R$ whose valuations
$(\delta_1,\ldots ,\delta_r)$ rationally generate the group $\Phi$, where $r$ is the rational
rank $\hbox{\rm r}(\nu )$ of $\Phi$; this means that $(\delta_1,\ldots ,\delta_r)$ form a basis
of the $\Qq$-vector space $\Phi\otimes_\Z\Qq$. Let $\Phi^{(u)}=\Z\delta_1 \oplus\cdots
\oplus \Z\delta_r\subseteq \Phi$ be the subgroup generated by the $\delta_i$; set
$\Phi_+^{(u)}=\Phi^{(u)}\cap \Phi_+$ and define the subalgebra 
$$\hbox{\rm gr}^{(u)}_\nu R=\bigoplus_{\phi\in \Phi_+^{(u)}}\P_\phi (R)/\P_\phi^+
(R).$$ Let us denote by $S^{(u)}$ the intersection $S\cap\hbox{\rm gr}^{(u)}_\nu R$ and by
$k^{(u)}$ the induced field $[(S^{(u)})^{-1}\hbox{\rm gr}^{(u)}_\nu R]_0$.
\begin{lemma}\label{integext} The extensions $\hbox{\rm gr}^{(u)}_\nu R\to \hbox{\rm gr}_\nu R$ and
$k^{(u)}\to k_\nu$ are integral.
\end{lemma}
\begin{proof} By [B3], Chap.V, No.1, Corollary to Prop.4, it suffices to check that any
homogeneous element
$\overline x_\phi \in\hbox{\rm gr}_\nu R$ has a power which is in $\hbox{\rm gr}^{(u)}_\nu
R$. Since
$\Phi^{(u)}$ rationally generates $\Phi$, there is an integer $n$ such that $n\phi
\in\Phi^{(u)}$, so that $\overline x^n\in \hbox{\rm gr}^{(u)}_\nu R$. 
\end{proof}

Given an element $\phi\in\Phi^{(u)}$, it can be written uniquely as 
$$\phi=a_1\delta_1+\cdots +a_r\delta_r\ \  \ \hbox{\rm with } a_i\in \Z.$$
\begin{lemma} With the notation just introduced, the map 
$$\iota_u\colon \hbox{\rm gr}^{(u)}_\nu R\to k^{(u)}[t^{\Phi^{(u)}_+}],\ \ \overline
x_\phi\mapsto \frac{\overline x_\phi}{\overline u_1^{a_1}\ldots \overline u_r^{a_r}}t^\phi$$ induces an
isomorphism
$$ (S^{(u)})^{-1}\hbox{\rm gr}^{(u)}_\nu R\stackrel{\simeq}\longrightarrow k^{(u)}[x_1^{\pm 1} ,\ldots, x_r^{\pm 1}].$$
If $k_\nu$ is algebraic over $k_R$, the map $\iota_u$ is integral.\end{lemma}
\begin{proof} For the first part, since $\Phi^{(u)}$ is a free abelian group of rank $r$, it
suffices to show that $ (S^{(u)})^{-1}\hbox{\rm gr}^{(u)}_\nu R\to k^{(u)}[t^{\Phi^{(u)}}]$ is an
isomorphism. Remarking that $\iota_u (S^{(u)})\subseteq \{\lambda t^\gamma ,\ \gamma
\in \Phi^{(u)}_+,\ \lambda\in {k^{(u)}}^*\}$, we see that 
$ (S^{(u)})^{-1}\hbox{\rm gr}^{(u)}_\nu R\subset k^{(u)}[t^{\Phi^{(u)}}]$. The surjectivity follows from the fact that $\Phi_+^{(u)}$ rationally
generates $\Phi^{(u)}=\Z^r$.
\par\noindent If the extension $k_R\to k_\nu$ is algebraic, so is $k_R\to k^{(u)}$. To prove that $\iota_u$
is integral, it suffices to show that any element $\lambda t^\gamma$ with $\lambda\in
k^{(u)},\ \gamma\in\Phi^{(u)}$, is integral. This follows from the fact that $\lambda$ is
algebraic over $k_R$.
\end{proof} 
\par\noindent From these lemmas we deduce
\begin{proposition}{\rm (Spivakovsky, Piltant)}\label{tplusr} The transcendence degree
over $k_R$ of $ \hbox{\rm gr}_\nu R$ is equal to $\hbox{\rm r}(\nu )+\hbox{\rm
t}_{k_R}(\nu)$.
\end{proposition} 
\begin{proof} Because the extension $\hbox{\rm gr}^{(u)}_\nu R\to \hbox{\rm gr}_\nu R$ is integral by Lemma \ref{integext}, this
transcendence degree is equal to that of $\hbox{\rm gr}^{(u)}_\nu R$, hence to that of $
(S^{(u)})^{-1}\hbox{\rm gr}^{(u)}_\nu R$, which is that of
$k^{(u)}[t^{\Phi^{(u)}}]$, and also of $k_\nu [t^{\Phi^{(u)}}]$ since the
extension $k^{(u)}\to k_\nu$ is integral by Lemma \ref{integext}.
\end{proof}
\par\noindent
 Now Piltant uses the following corollary to Cohen's dimension inequality:
\begin{proposition} Let $B$
be an integral domain containing a n\oe therian ring $A$. Then $$\hbox{\rm dim}B\leq
\hbox{\rm dim}A+\hbox{\rm t}_A(B),$$ where $\hbox{\rm t}_A(B)$ is the transcendence
degree of $B$ over $A$.
\end{proposition}
 \begin{proof} The proof is like the classical one
(see [Ei], Chap. 13) except that we do not suppose that $B$ is finitely generated over $A$.
But if $(0)={\mathbf p}_0\subset {\mathbf p}_1\subset {\mathbf p}_2\subset \cdots \subset {\mathbf p}_d$ is
a chain of distinct prime ideals of $B$, choosing for each $i,\  1\leq i\leq d$ an element
$y_i\in {\mathbf p}_i\setminus {\mathbf p}_{i-1}$, we see that the ${\mathbf p}_i\cap A[y_1,\ldots ,y_d]$
form a chain of distinct prime ideals, so that $\hbox{\rm dim}A[y_1,\ldots ,y_d]\geq d$,
and we are reduced to the case of a finitely generated algebra.
\end{proof}
\par\noindent
Applying this result with $A=k_R$ and $B=\hbox{\rm gr}_\nu R$, remarking that we have 
used the n\oe therianity of $R$ only to ensure that $\hbox{\rm r}(\nu )+\hbox{\rm
t}_{k_R}(\nu)<\infty$ and using Proposition \ref{tplusr} gives the first part of the: 
\begin{proposition}\label{pil}{\rm (Piltant, [P])} Whenever the valuation ring $R_\nu$
biratio\-nally dominates the local ring $R$ and $\hbox{\rm r}(\nu )+\hbox{\rm
t}_{k_R}(\nu)<\infty$, the following inequality holds: $$\hbox{\rm dim}\hbox{\rm gr}_\nu
R\leq \hbox{\rm r}(\nu )+\hbox{\rm t}_{k_R}(\nu).$$ If $\hbox{\rm t}_{k_R}(\nu)=0$,
$$\hbox{\rm dim}\hbox{\rm gr}_\nu R= \hbox{\rm r}(\nu ).$$
\end{proposition}
\begin{proof} To prove the second part, we use the lemmas above: the map 
$$ (S^{(u)})^{-1}\hbox{\rm gr}^{(u)}_\nu R\to k^{(u)}[x_1^{\pm 1} ,\ldots, x_r^{\pm 1}]$$ 
is an isomorphism, so these two rings have the same dimension $\hbox{\rm r}(\nu )$. Since
dimension cannot increase by localization, this gives
$\hbox{\rm dim}\hbox{\rm gr}^{(u)}_\nu R\geq \hbox{\rm r}(\nu )$ and finally, because the
extension $\hbox{\rm gr}^{(u)}_\nu R\to \hbox{\rm gr}_\nu R$ is integral, using
([B3], Chap. VIII, No.3. Th.1), $\hbox{\rm dim}\hbox{\rm gr}_\nu R\geq \hbox{\rm r}(\nu )$. But
we have the opposite inequality by the first part of the proposition.
\end{proof} 
\begin{remark} Piltant's result implies that if $\hbox{\rm r}(\nu)<\infty$ we have:
 $$\hbox{\rm dim}\hbox{\rm gr}_\nu R_\nu=\hbox{\rm r}(\nu)$$ 
and if $\hbox{\rm t}_{k_R}(\nu )=0$,
$$\hbox{\rm dim}\hbox{\rm gr}_\nu R=\hbox{\rm t}_{k_R}(\hbox{\rm gr}_\nu R),$$ 
the transcendence degree of $\hbox{\rm gr}_\nu R$ over $k_R$. I will use this last equality
in several examples. The first one should be contrasted
with the fact that $\hbox{\rm dim}R_\nu=\hbox{\rm h}(\nu )\leq \hbox{\rm r}(\nu)$.
\end{remark}\par\medskip\noindent   
\textbf {Question.} Does the proof in the case where $\hbox{\rm t}_{k_R}(\nu)=0$ also show that
$\hbox{\rm gr}_\nu R$ is catenary? \par\noindent  The proof of the fact that equality in the
inequality of Proposition \ref{pil} continues to hold without the assumption that 
$\hbox{\rm t}_{k_R}(\nu )=0$ provided that $k_R$ is uncountable is more delicate. Piltant has
given in [P] an example of a valuation $\nu$ on the ring $R=\overline {\mathbf
Q}[s,t,x,y]_{(s,t,x,y)}$ with value group $\Z^2_{lex}$ and residue field $\overline {\mathbf
Q}(\frac{s}{t})$ such that $\hbox{\rm gr}_\nu R$ has Krull dimension  $3=\hbox{\rm 
r}(\nu)+\hbox{\rm t}_{k_R} (\nu)$ but is not catenary; it possesses a saturated chain of prime
ideals of length two. It would be interesting to decide whether the completion
$\widehat{\hbox{\rm gr}}^{(\nu)}_\nu R$ with respect to $\nu_{\rm gr}$ which we shall meet in
subsection \ref{Bcomp} is catenary.\par\noindent
 Strict inequality in Abhyankar's inequality may happen, but this does not contradict the
semi-continuity theorem for the dimensions of the fibers of the map
$\hbox{\rm Spec}A_\nu (R)\rightarrow  \hbox{\rm Spec}k[v^{\Phi_+}]$ since semicontinuity is
proved only for morphisms which are essentially of finite type. Under the assumptions of
Piltant's theorem, for $\A_\nu (R)$ to be a
$k[v^{\Phi_+}]$-algebra of finite type it is necessary that equality holds in Abhyankar's inequality.
\par\medskip Since $R$ is assumed to be n\oe therian we will in this text remain in the
situation where $\hbox{\rm r}(\Phi), \hbox{\rm t}_{k_R} (\nu )$, etc.. are finite. \par
Recall (see [B3], Chap. VI, \S 10, No.
2) that if $\Phi$ is finitely generated and if $\hbox{\rm r}(\nu)=\hbox{\rm h}(\nu)$, then $\Phi$
is order isomorphic to $\Z^{\hbox{\rm h}(\nu)}$ with the lexicographic
order.\par\medskip\noindent 
I will often make use of the following fact:
\begin{proposition}\label{wellord} {\rm  (Krull, [Kr], Zariski; see ([Z-S], Appendix 3, [S2], remark 2.7, [V1],
Prop. 9.1.)} If $R$ is n\oe therian, the value semigroup $\Gamma =\nu (R\setminus \{
0\})\subset \Phi_+\cup \{ 0\}$ is well ordered and its ordinal type {\rm (see [B5])} is
$\omega^{\hbox{\rm h}_R(\nu)}$, where  $\hbox{\rm h}_R(\nu)$ is the height in $R$ of
the valuation $\nu$.
\end{proposition}
\begin{proof} Zariski's proof is in terms of valuation ideals of $R$, which are the distinct
$\P_\phi (R)$ and therefore in monotonous bijection with the elements of $\Gamma$ by the map $\gamma\mapsto \P_\gamma$ for
the reverse inclusion ordering on ideals; a decreasing sequence of elements of
$\Gamma$ gives an increasing sequence of ideals $\P_\phi (R)$, which has to be stationary since $R$ is n\oe therian. The second part
of the statement follows from [Z-S], Appendix 3, Proposition 2.
\end{proof}
\begin{corollary}\label{gammai} For a valuation $\nu$ of a n\oe therian ring $R$, the value semigroup admits a unique
minimal system of generators
$$\Gamma=\langle \gamma_1,\gamma_2,\ldots ,\gamma_i,\ldots \rangle$$
satisfying $$\gamma_i<\gamma_{i+1}.$$
It is indexed by ordinals $<\omega^{\hbox{\rm h}_R(\nu)}$. If the height of $\nu$ in $R$ is one, in
particular if the height of $\nu$ is one, given $\gamma\in \Gamma$, there are finitely many
generators $\gamma_i\leq \gamma$.\end{corollary}
\begin{proof}Define $\gamma_1$ to be the smallest non zero element of $\Gamma$, and then
inductively $\gamma_{i+1}$ to be the smallest non zero element of $\Gamma$ which is not in
the semigroup  $$\Gamma_i=\langle \gamma_1,\gamma_2,\ldots ,\gamma_i\rangle$$
generated by the previous ones. If after finitely many steps, $\Gamma_i=\Gamma$, our semigroup is finitely generated; otherwise, the elements obtained
after iterating the successor construction countably many times generate a semigroup
$\Gamma_{<\omega}\subseteq
\Gamma$. If it is not equal to $\Gamma$, let $\gamma_\omega$ be the smallest non zero element of $\Gamma$ not contained in $\Gamma_{<\omega}$, set
$\Gamma_\omega=\langle\Gamma_{<\omega},\gamma_\omega\rangle$ and continue in this manner. This is clearly a minimal system of generators of
$\Gamma$. If for some
$i$ we had
$\gamma_{i+1}<\gamma_i$ it would contradict the definition of $\gamma_i$ since
$\gamma_{i+1}\notin \Gamma_{<i}$. The second assertion
also follows [Z-S], Appendix 3, Lemma 4.\end{proof}
\begin{remark} The element $\gamma_\omega$ has no predecessor among the $\gamma_i$. We shall
see below that if $R$ is complete there are only finitely many elements without predecessor. \end{remark}
 Viewing a valuation as a map $\nu\colon K^*\to \Phi$ as in
subsection
\ref{valalg}, we see that a surjective monotone non-decreasing map of ordered groups $\lambda\colon \Phi \rightarrow\Phi_1$ corresponds to a composite
valuation $\nu_1=\lambda\circ \nu$ and another valuation ring $R_{\nu_1}\supset R_\nu$.
However the classical terminology, referring to a geometric notion of composition, is that
$\nu$ is {\it composed} with $\nu_1$ and we will respect it. The map $\lambda$ extends to
a map $$\tilde\lambda\colon k[v^{\Phi_+}]\rightarrow k[w^{\Phi_{1+}}].$$
\noindent We have a map ${\mathcal A}_\nu(R) \rightarrow {\mathcal A}_{\nu_1}(R)$ and in fact a
commutative diagram $$\begin{array}{ccl}\tilde\lambda_{\mathcal A} (R)\colon &{\mathcal
A}_\nu(R)&\longrightarrow {\mathcal A}_{\nu_1}(R) \\ &\ \ \ \uparrow\ \ \ \ &\ \
\ \ \ \ \ \ \uparrow\\ \tilde\lambda \colon &k[v^{\Phi_+}]&\longrightarrow
k[w^{\Phi_{1+}}]\end{array}$$ where the top arrow is induced by the inclusions ${\mathcal
P}_\phi(R)\subset {\mathcal P}_{\lambda(\phi)}(R)$, i.e., $x_\phi v^{-\phi}\mapsto x_\phi
w^{-\lambda (\phi)} $, and the bottom arrow is $\tilde\lambda$. The map
$\tilde\lambda$ is surjective since $\lambda\colon \Phi\rightarrow\Phi_1$ is order preserving,
and its kernel is the ideal generated by the $(v^\psi-1)_{\psi\in \Psi_+}$ where $\Psi$ is the
kernel of the map $\Phi\rightarrow \Phi_1$. Let us show that
$\tilde\lambda_{\mathcal A}(R)$ is also surjective: given $xw^{-\mu}$, where
$\mu=\lambda (\phi)$ and $\nu_1(x)\geq \mu$, either $\nu (x)\geq \phi$ and
$xw^{-\mu}$ is the image of $xv^{-\phi}$, or we have not chosen $\phi$ well, but since
$\lambda$ is order preserving, if $\nu (x)<\phi$, we have $\nu_1(x)=\lambda (\nu (x))\leq
\lambda(\phi)= \mu$ so that $\lambda(\nu(x))=\mu$ and $\phi -\nu(x)\in
\Psi$. Finally, we see that $xw^{-\mu}$ is the image of $xv^{-\nu(x)}$ and we have shown the
surjectivity. The kernel is the ideal of ${\mathcal A}_{\nu}$ generated by
$(v^\psi-1)_{\psi\in\Psi_+}$.\par Finally, we have shown that the formation of
$k[v^{\Phi_+}]\rightarrow {\mathcal A}_\nu(R)$ commutes with the composition of valuations, as follows; 
\begin{proposition}\label{base change} The natural morphism
 $$\begin{array}{ccl}{\mathcal A}_{\nu}(R)\otimes_{k[v^{\Phi_+}]} k[w^{\Phi_{1+}}]&\rightarrow {\mathcal
A}_{\nu_1}(R)\\ x_\phi v^{-\phi}\otimes w^{\phi_1}&\mapsto x_\phi
w^{\phi_1-\lambda(\phi)}\end{array}$$ is an isomorphism of graded
$k[v^{\Phi_+}]$-algebras.
\end{proposition} 
\begin{remark} The fact that
$$\A_{\nu_1}(R)=\A_{\nu}(R)/\bigl((v^\psi -1)_{\psi\in \Psi_+}\bigr)$$ is true without assuming
that $R$ contains a field.\par\noindent
\end{remark}
\begin{corollary}\label{speciagr} The rings $\hbox{\rm gr}_{\nu}R$ and
$\hbox{\rm gr}_{\nu_1}R$ are both specializations of the same ring over the polynomial ring
$k[(t_\psi)_{\psi\in\Psi_+}]$.
\end{corollary}
\begin{proof} The ring $\hbox{\rm gr}_{\nu_1}R$ is isomorphic
to the quotient ring of $\A_{\nu}(R)$ by \break $\bigl((v^\psi -1)_{\psi \in
\Psi_+},(v^\phi)_{\phi\in
\Phi_+\setminus
\Psi_+})\bigr)$, and $\hbox{\rm gr}_{\nu}R=\A_\nu(R)/\bigl((v^\phi)_{\phi\in \Phi_+}\bigr)$.
These rings are specializations of $$\A_{\nu}(R)\otimes_kk[(t_\psi)_{\psi\in\Psi_+}]/\big((v^\psi
-t_\psi)_{\psi\in\Psi_+},(v^\phi)_{\phi\in
\Phi_+\setminus \Psi_+})\big).$$
\end{proof}
\begin{Remark}  1) If we have a sequence of valuation rings as above, with
$R\subset R_\nu$,
$$R_\nu\subset R_{\nu_1}\subset \cdots \subset R_{\nu_{h -1}},$$
\noindent where $h$ is the height of the valuation $\nu$, all the graded algebras $\hbox{\rm
gr}_{\nu_s}R$ are specializations of the same algebra of the form $${\mathcal
A}_{\nu}(R)\otimes_kk[t_\omega]/\big((v^\omega -t_\omega)_{\omega\in \Omega_1},
(v^\omega)_{\omega\in \Omega_2}\big),$$ and we can in a number of cases use a
system of generators of
$\hbox{\rm gr}_{\nu}R$ to produce a system of ``formal'' generators for each $\hbox{\rm
gr}_{\nu_s}R$ (see below).\par\noindent 
2) Denoting by $\lambda\colon \Phi\to \Phi_1$ the morphism of value groups, notice that the graph $(\phi , \lambda
(\phi))\subset \Phi\times \Phi_1$ is isomorphic to $\Phi$ and that we have an isomorphism of graded algebras:
$$\A_\nu (R)\stackrel{\simeq}\rightarrow \bigoplus_{\phi_1\in \Phi_1}\bigl(\bigoplus_{\lambda (\phi)=\phi_1}\P_\phi
(R)w^{-\phi}\bigr)w_1^{-\phi_1}.$$\end{Remark}
\subsection{The valuation $\tilde{\mathbf {\nu}}_\A$}\label{sec-special}
\par\noindent Let us denote by $\tilde{\mathbf {\nu}}_\A$ the valuation on $\A_\nu (R)$ with
values in the group $\Phi\oplus
\Phi$ ordered lexicographically, and defined by
$$\tilde{\mathbf {\nu}}_\A (\sum x_\phi v^{-\phi})=\big(\mu_\A(\sum x_\phi v^{-\phi}), 
\nu_\A(\sum x_\phi v^{-\phi})\big).$$ According to what we saw in \S 1, the center of this
valuation is defined by the ideal $\big( (v^{\Phi_+})\A_\nu (R),
\oplus_{\phi\in\Phi_+}\P_\phi v^{-\phi}\big)$ of $\A_\nu (R)$. Geometrically, this center is
$\hbox{\rm Spec}\overline R$.\par This valuation is composed with $\mu_\A$ in the manner
corresponding to the first projection $\Phi\oplus \Phi
\rightarrow \Phi$, and the residual valuation of $\tilde \nu_\A$ on the center $\hbox{\rm Spec}\hbox{\rm gr}_\nu
R$ of $\mu_\A$ is $\nu_{\rm gr}$.\par This valuation is of interest because it is a
``specialization'' of the valuation
$\nu$ on $R$ to the valuation
$\nu_{\rm gr}$. Indeed, the valuation $\tilde\nu_\A$ extends to a valuation of $(v^{\Phi_+})^{-1}A_\nu
(R)=(v^{\Phi_+})^{-1}R[v^{\Phi_+}]$ (Proposition \ref{deformation}, b)) which is the natural extension of $\nu$ to the algebra of the
``generic fiber'' of the faithfully flat map defined in \S \ref{speciatogr}, and as we saw it induces the
natural valuation on the special fiber. Heuristically, we may think of the inclusion 
$\A_\nu (R)\rightarrow\A_\nu (R_\nu )$ as a uniformization of the valuation $\tilde \nu_\A$ (see
below, \S \ref{graval}), forgetting for a moment the fact that it is not an extension of finite type.
The center of $\tilde \nu_\A$ in $\A_\nu (R_\nu )$ is the origin, and, still heuristically,
$\A_\nu (R_\nu )$ is ``regular''. Then the corresponding map of schemes is a
simultaneous uniformization for $\nu$ and $\nu_{\rm gr}$ since the induced map
$\hbox{\rm gr}_\nu R\rightarrow \hbox{\rm gr}_\nu R_\nu$ is also (in the same heuristic
sense) a uniformization of $\nu_{\rm gr}$.
\subsection{The valuation algebra of the residual valuation}\label{residual} 
Let us now come back to the inclusion
$R_\nu\subset R_{\nu_1}$ associated to an monotone non-decreasing surjection
$\lambda\colon \Phi\rightarrow \Phi_1$ of ordered groups, with kernel $\Psi$; it is a classical 
fact of valuation theory (see [A1], pp. 56-57 and [V1]) that $m_{\nu_1}\subset R_\nu$ and the
image $\overline R_\nu=R_\nu/m_{\nu_1} \subset R_{\nu_1}/m_{\nu_1}$ is a valuation ring
of the field
$R_{\nu_1}/m_{\nu_1}$, corresponding to a valuation
$\overline\nu$ with value group $\Psi$ and center $m_\nu \cap R/m_{\nu_1}\cap R$ in
$R/m_{\nu_1}\cap R$. \par If we endeavor to find the valuation
algebra of $\overline R_1=R/(m_{\nu_1}\cap R)$ with respect to $\overline \nu$, we see that if we denote by $\A_\Psi (R)$ the
subalgebra $\bigoplus_{\phi\in
\Psi}{\mathcal P}_\phi (R)v^{-\phi}$ of $\A_\nu (R)$, then $\A_{\overline
\nu}(\overline R_1)$ is nothing else than the image of
$\A_\Psi (R)$ in $\overline R_1[v^{\Phi}]$ by the canonical map $R[v^{\Phi}] \rightarrow\overline
R_1[v^{\Phi}]$. It is important to note that, setting
$\Gamma =\nu (R\setminus\{ 0\})$, the semigroup $\Gamma
\cap \Psi$ may be zero. This happens if and only if
$m_{\nu_1}\cap R=m_\nu \cap R$; in that case the valuation
$\overline \nu$ is trivial on $\overline R_1$.\par
A typical example is the valuation on
$R=k[u_1,u_2]_{(u_1,u_2)}$ with values in $\Z^2$ ordered lexicographically determined by setting $\nu (u_1)=(1,0),\
\nu (u_2)=(1,1)$ and deciding that the valuation of a polynomial is the smallest valuation of the monomials which it contains; it is composed with
the
$(u_1)$-adic valuation,
$\Gamma$ is the semigroup generated by $(1,0),\ (1,1)$ and $\Psi$ is the subgroup $0\bigoplus\Z$. Note that the two
valuations have the same center in $R$ and that for a given $\phi\in \Phi_+$ there are only
finitely many ideals $\P_{\phi'}$ between $\P_{\lambda (\phi )}$ and $\P^+_{\lambda (\phi )}$. 
\subsection{Comparison of graded algebras for composed valuations}\label{compargr} 
The result of Proposition
\ref{base change} suggests that we seek a closer relation between $\hbox{\rm gr}_\nu R$ and $\hbox{\rm
gr}_{\nu_1} R$ when $\nu$ is composed with $\nu_1$.  First, let us remark that with the usual
notations, we have for $\phi\in \Phi_+\cup\{0\}$ the following inclusions\footnote {The fact that
$\lambda (\phi )\in \Phi_1$ is deemed sufficient, here and in the sequel, to indicate that the
first and last ideals correspond to the valuation $\nu_1$.} $$\P_{\lambda (\phi)}^+\subset
\P_\phi^+\subset
\P_\phi \subset \P_{\lambda (\phi)}, $$ so that we get a filtration of each $(R/m_{\nu_1}\cap R)$-module $$(\hbox{\rm
gr}_{\nu_1} R)_{\phi_1}=\frac{\P_{\phi_1}}{\P_{\phi_1}^+}$$ by the submodules $(\P_\phi/\P_{\lambda
(\phi)}^+)_{\lambda (\phi)=\phi_1}$. In the case where $\phi_1=0$, it is nothing else than the
filtration of $\overline R_1=R/(m_{\nu_1}\cap R)$ by the images of the $(\P_\phi)_{\phi\in \Psi} $
where $\Psi =\hbox{\rm Ker}\lambda$, i.e., the filtration corresponding to the ideals
$(\overline\P_\psi )_{\psi \in \Psi}$ associated to the residual valuation $\overline \nu$ on
$\overline R_1$. Let us denote by $\overline\P (\phi_1)$ this filtration, for each $\phi_1\in
\Phi_{1+}\cup \{ 0\}$. We have immediately
 \begin{lemma}\label{bifiltr} The natural map
$$\P_\phi\rightarrow \frac{\P_\phi}{ \P_{\lambda (\phi )}^+}$$ induces an isomorphism
$$\hbox{\rm gr}_\nu R\stackrel{\simeq}\rightarrow \bigoplus_{\phi_1\in \Phi_{1+}\cup\{ 0\}}\hbox{\rm
gr}_{\overline\P (\phi_1)}(\hbox{\rm gr}_{\nu_1} R)_{\phi_1}.$$ It is graded in the sense that
an element of degree $\phi$ is mapped to an element of bidegree $(\phi, \lambda (\phi))$.
Moreover, for $x\in R$, we have $$\hbox{\rm in}_\nu x=\hbox{\rm in}_{\overline
\P(\nu_1(x))}\big(\hbox{\rm in}_{\nu_1} x \big).$$ \end{lemma}  \par 
In particular,
taking $\lambda (\phi)=0$, we find that it induces a graded isomorphism $$\bigoplus_{\psi\in
\Psi_+\cup\{0\}}(\hbox{\rm gr}_\nu R)_\psi\stackrel{\simeq}\rightarrow \hbox{\rm gr}_{\overline\nu}\overline
R_1,$$ which can also be deduced from \ref{residual}. 
Note that when $R$ is n\oe therian, and so $\Gamma$ is well ordered  by Proposition
\ref{wellord}, the topology on $R/(m_{\nu_1}\cap R)$ corresponding to the $\overline
\nu$-adic filtration is the same as the image by the natural surjection $R\rightarrow
R/(m_{\nu_1}\cap R)$ of the $\nu$-adic topology on $R$. This follows from the inclusions
$\P_\psi(R)\supset  m_{\nu_1}\cap R$ for $\psi\in \Psi$. Thus, the image of $\P_\psi (R)$ in $R/(m_{\nu_1}\cap R)$ is the valuation ideal
$\overline\P_\psi(\overline R)$ for $\overline\nu$. Note that it may be the trivial topology; more generally, if $\nu$ and $\nu_1$ have the same
center in
$R$, the filtrations
$\overline\P (\phi_1)$ are finite ([Z-S], Vol. II, Appendix 3, Corollary).\par\noindent 
I will make considerable use of the following:
\begin{proposition} {\rm (Zariski)}\label{almostone} With the same notations, assume $\nu$ is a valuation on the n\oe therian local
domain $R$ and that the valuation
$\nu_1$ has height one less than
$\nu$, so that $\overline\nu$ is of height one. Given
$\phi_1\in \Phi_{1,+}$, we have, whenever
$\lambda (\phi)=\phi_1$ (see Lemma \ref{bifiltr}) 
$$\P_{\phi_1}^+\subseteq \P_\phi\subseteq \P_{\phi_1}.$$ There are only finitely many
ideals $\P_{\phi'}$ between $\P_{\phi_1}$ and $\P_\phi$. The ideals $\P_\phi$ between
$\P_{\phi_1}$ and $\P_{\phi_1}^+$ are finite in number if the centers of $\nu$ and $\nu_1$ are
equal, and form a simple infinite sequence if the centers of $\nu$ and $\nu_1$ are
different; in this last case, denoting by $\phi (\phi_1)$ the least element of $\Gamma\cap \lambda^{-1}(\phi_1)$, the $R/(m_{\nu_1}\cap
R)$-modules 
$(\P_{\phi (\phi_1)+\psi}/\P_{\phi_1}^+)_{\psi\in\Gamma\cap\Psi_+}$ are cofinal in the sequence of the
$(\P_{\phi}/\P_{\phi_1}^+)_{\lambda (\phi )=\phi_1}$, and for each $\phi\in \Gamma\cap \lambda^{-1}(\phi_1)$ there is an
integer $q$ such that $(m_\nu
\cap R)^q.\P_{\phi_1}\subset \P_\phi$. 
\end{proposition}
\begin{proof} Let $\Gamma\subset \Phi_+\cup\{0\}$ be the semigroup of $R$. Given $\phi_1\in \Gamma_1=\lambda (\Gamma)$, for any
$\phi\in
\Gamma$ such that $\lambda (\phi)=\phi_1$ we can write $\phi=\phi (\phi_1)+\psi$ with $\psi\in \Psi_+\cup\{0\}$, and the result follows
from the fact that if $\overline\nu$ is not trivial, $\Gamma\cap\Psi_+$ is cofinal in
$\Psi_+$ (see subsection \ref{valalg}), and the fact that $\overline\nu$ has center $(m_\nu \cap R)/(m_{\nu_1} \cap R)$ and is archimedian.
See [Z-S], Appendix 3, lemma 4 and its corollary.
\end{proof}
\begin{remark} In fact, whenever we have a
specialization $R\subset R_\nu\subset R_{\nu_1}$ of valuations of $R$, the valuation $\nu$
induces, as described above, a filtration on the ring $\hbox{\rm gr}_{\nu_1} R$, to which we
can associate as in subsection \ref{valalg} a ``filtration algebra'' which determines a
specialization of $\hbox{\rm gr}_{\nu_1} R$ to $\hbox{\rm gr}_\nu R$. This specialization is
the same as that coming from the embedding of $\hbox{\rm Spec}\A_{\nu_1}(R)$ in
$\hbox{\rm Spec}\A_\nu (R)$.
\end{remark} \par\medskip\noindent 
\begin{example}\label{resval} Let $R$ be a local integral domain, and
$f\in R\setminus\{0\}$ be such that $R/(f)$ is again integral and $\bigcap_{n\geq 0} f^nR=\{0\}$. Let us choose a valuation
$\overline \nu$  of $R/(f)$ with group
$\Psi$. Since each element $x$ of $R\setminus \{0 \}$ can be written in a unique way $x=f^ny$ with $y\notin (f)$,
we may define a valuation $\nu$ on $R$ with group $\Z\oplus \Psi$ ordered lexicographically by setting
$$\nu (f^ny)=(n, \overline \nu (\overline y)) \ \ \hbox{\rm where}\ 
\overline y=y\ \hbox{\rm mod.}(f).$$ The map $f^ny\mapsto n\in \N$ defines a valuation
$\nu_1$ on $R$ with value group $\Z$ and with which $\nu$ is composed. With the notations
introduced above, we have $\overline R=R/(f)$, and  $$\hbox{\rm gr}_{\nu_1} R=\overline R
[F]\ \ \ \hbox{\rm with}\ F=\hbox{\rm in}_{\nu_1}f,$$ and from the considerations above
we see that because it is a polynomial ring in one variable, we have:
$$\hbox{\rm gr}_\nu R=(\hbox{\rm gr}_{\overline\nu}\overline R)[F].$$
We shall see the general form of this result in subsection \ref{morestruc}.
\end{example}
\subsection{Comparison of residual valuations}\label{resigroup}
 In this subsection I recall some
direct consequences of the basic facts on the composition of valuations found in ([A1], pp. 56-57,
[V1], section 4, [Z-S]).   Keeping the notations of the beginning of this section, let us
consider a sequence
$$R_\nu\subset R_{\nu_1}\subset R_{\nu_2}$$ of valuation rings; it corresponds to a sequence
$$\Phi\rightarrow \Phi_1\rightarrow
\Phi_2$$ of value groups; let $\Psi_2\subset \Psi_1\subset \Phi$ be the kernels of the maps
$\Phi\rightarrow \Phi_1$ and
$\Phi\rightarrow \Phi_2$ respectively. To this situation is associated an inclusion
$$R_\nu/(m_{\nu_2}\cap R_\nu )\subset R_{\nu_1}/(m_{\nu_2}\cap R_{\nu_1})$$ of valuation
rings of the residue field $k_{\nu_2}$ of $R_{\nu_2}$. Note that in fact
$m_{\nu_2}\cap R_\nu =m_{\nu_2}$ since $m_{\nu_2}\subset R_\nu$ by the general properties
of composed valuations (see [V1]), and similarly for $R_{\nu_1}$; recall also that $R_{\nu_i}$ is the
localization of $R_\nu$ at the prime $m_{\nu_i}$. The value group of the first valuation ring is
$\Psi_1$, and the value group of the second is $\Psi_1/\Psi_2$. If $R_{\nu_1}$ and $R_{\nu_2}$ are
consecutive valuation rings, in the sense that there is no other valuation ring between them, the
group $\Psi_1/\Psi_2$ is of height one.
\par Assuming that this is the case, given a subring
$R\subset R_\nu$, we see that $\nu_1$ induces a valuation of height at most one, with valuation
ring $R_{\nu_1}/(m_{\nu_2}\cap R_{\nu_1})$, on the quotient $R/(m_{\nu_2}\cap R)$; the
center of this valuation is the ideal $(m_{\nu_1}\cap R)/(m_{\nu_2}\cap R)$, which is also the
kernel of the natural morphism
$$R/(m_{\nu_2}\cap R)\rightarrow R/(m_{\nu_1}\cap R).$$ More precisely, it may happen that
$m_{\nu_1}\cap R=m_{\nu_2}\cap R$, and then the residual valuation is trivial on
$R/(m_{\nu_2}\cap R)$.
\subsection{Specialization of valuations and rational valuations}\label{speciaval} 
Let $k$ be a field, and $(R, m)$ a local
$k$-algebra which is an integral domain. Let $R_{\nu_0}$ be a valuation ring in the fraction field of $R$, containing $R$, with
residue field $k_{\nu_0}$, and assume that the transcendence degree over $k$ of $k_{\nu_0}$ is
greater than zero. Then, if
$R/(m_{\nu_0} \cap R) \neq k_{\nu_0}$, there exists a
$k$-valuation ring $V$ of $k_{\nu_0}$ containing $R/(m_{\nu_0}
\cap R)$ and different from $k_{\nu_0}$. If we denote by
$\pi\colon R_{\nu_0}\to k_{\nu_0}$ the canonical projection, we see that $\pi^{-1}(V)$ is a
valuation ring $R_{\nu_1}$ containing
$R$ and contained in $R_{\nu_0}$. The residue field of $R_{\nu_1}$ is equal to the residue field of
$V$, and therefore its transcendence degree over $k$ is less than that
of $k_{\nu_0}$. We may iterate this construction as long as the residue field of $R_{\nu_1}$
is transcendental over $k$ and
$R/(m_{\nu_1} \cap R) \neq k_{\nu_1}$. When we stop, we have built a sequence $$R\subset
R_{\nu_t}\subset R_{\nu_{t-1}}\subset\cdots \subset R_{\nu_0}$$ such that $t\leq
\hbox{\rm t}_k(\nu_0)$ and either $R_{\nu_t}$ dominates $R$ and the residual extension is
trivial, or the residue field of $R_{\nu_t}$ is algebraic over $k$. In this case, the ring\break 
$R/(m_{\nu_t}\cap R)$ contains $k$ and is contained in $k_{\nu_t}$, so it is a field. This means
that $m_{\nu_t}\cap R=m$, so $R_{\nu_t}$ dominates $R$. So at worse we end up the
specialization process with a valuation ring dominating $R$ and an algebraic residual extension.
\par This is a special case of the construction found in [Z-S], Vol. 2, Chap VI, \S
16.\par In the case where $R_\nu$ dominates $R$ and the residual extension is trivial,
one says that the valuation
$\nu$ on $R$ is {\it rational}. So we have:
\begin{proposition}\label{domi} If the residue field of
$R$ is algebraically closed, we can specialize any valuation $\nu_0$ of $R$ to a rational
valuation.
\end{proposition}
\begin{Remark} 1) Assuming that $R$ is excellent, and so in particular a Nagata ring, we may (see [V1], Proposition 10.1), at each step in the preceding
construction, choose a discrete divisorial valuation; then the transcendence degree of the residual valuation drops exactly by one at each
step, while the height of the value group increases exactly by one since the value group of $V$ is
the kernel of the map $\Phi_{\nu_1}\to\Phi_{\nu_0}$. Finally we have:
$$\hbox{\rm h}(\nu_t)=\hbox{\rm h}(\nu_0)+\hbox{\rm t}_{k_{R}}(\nu_0),$$ and since
$\hbox{\rm h}(\nu_t)\leq \hbox{\rm r}(\nu_t)$, in order to prove Abhyankar's inequality in this
case it suffices to check that $\hbox{\rm r}(\nu)\leq \hbox{\rm dim}R$ when
$R_\nu$ dominates $R$, but this is readily done, at least in the algebraic or formal case, if we
present the fraction field $K$ of $R$ as an algebraic extension of a field of rational functions of
transcendence degree $\hbox{\rm dim}R$ over $k_R$, since extending valuations to algebraic
extensions does not affect the rational rank. This suggests a proof of Abhyankar's inequality, not
very different from that described in  [V1]. 
\end{Remark}
\subsection{Monomial valuations}\label{monval}
This section begins to explain\ in what
sense the liftings $(\xi_i)_{i\in I}$ of generators
$(\overline\xi_i)_{i\in I}$ of the graded $R/(m_\nu\cap R)$-algebra $\hbox{\rm gr}_\nu R$
form a system of coordinates ``adapted'' to the valuation. The idea is that, at least in the case of
trivial residual extension, valuations which are monomial in a certain system of coordinates can
be uniformized by a single toric modification of the ambiant space in which this system of
coordinates embeds our singularity.\par Given a valuation $\nu$ of a local or graded ring
$R$, we generalize and modify the definition given in [S2] for a regular local ring as
follows\par\medskip\noindent
\begin{definition}  1) The valuation
$\nu$ of the local (or graded) integral domain
$R$ is {\it algebraically monomial} with respect to a system of generators $(\xi_i)_{i\in I}$ of the center $m_\nu \cap R$ of
$\nu$ in
$R$ if $(\xi_i)_{i\in I}$ is a generating sequence for $\nu$ in the sense (see [S2]) that any element
$x\in R$ can be written as a finite sum 
$$x=\sum_{\nu (\xi^\gamma)\geq \nu (x)}c_\gamma\xi^\gamma$$ with $\gamma\in I ^b,\hbox{
\rm for some} \ b\in \N,\ \ c_\gamma \in R$.\par\noindent 2) If the valuation $\nu$ is rational and $R$ is complete
with respect to the
$\nu$-adic topology and has a field of representatives, we say that $\nu$ is {\it analytically monomial} with respect to
$(\xi_i)_{i\in I}$ if the same holds with ``finite sum'' replaced by ``$\nu$-adically convergent series'' and
``$c_\gamma\in R$'' by ``$c_\gamma\in k$''. In particular, in both cases,
$\nu (x)=\hbox{\rm min}_{\{\gamma/c_\gamma\neq 0\}} (\nu (\xi^\gamma))$.
\end{definition}
\begin{Remark} 1) To say that $\nu$ is algebraically monomial with respect to the $\xi_i$ is equivalent, when $R$ is
n\oe therian, to asking that for each
$\phi\in\Phi_+$ the ideal $\P_\phi (R)$ is generated by monomials in the $(\xi_i)_{i\in I}$, and implies
that $\A_\nu (R)$ is generated as a $R[v^{\Phi_+}]$-algebra by the $\xi_iv^{-\nu(\xi_i)}$ and the initial forms $(\overline \xi_i )_{i\in I}$ in $\hbox{\rm
gr}_\nu R$ of the elements $(\xi_i)_{i\in I}$ generate it as $R/(m_\nu\cap R)$-algebra. Analytic monomiality implies that they generate the
$k$-algebra $\hbox{\rm gr}_\nu R$. Again if $R$ is n\oe therian, since the ideal generated by the exponents of a series is finitely generated, analytic
monomiality implies algebraic monomiality. We shall later meet monomial valuations in non n\oe therian rings. In particular, if the valuation
$\nu$ is algebraically (resp. analytically) monomial with respect to a {\it finite} set of generators of $m_\nu \cap R$ (resp. the maximal ideal $m$), the
$R/(m_\nu\cap R)$-algebra (resp. $k$-algebra) $\hbox{\rm gr}_\nu R$ is finitely generated. Since the expansion $x=\sum_{\nu (\xi^\gamma)\geq \nu
(x)}c_\gamma\xi^\gamma$ is not unique in general, one must stress the fact that the definition requires the existence of one expansion with the stated
properties.\par\noindent 2) If $\Phi$ is a totally ordered group of finite rational rank, if
$R$ a $\Phi_+$-graded $k$-algebra generated by the elements
$\xi_i$ and we assume that the valuation coincides with the grading and that each homogeneous component of $R$
is a finite-dimensional
 vector space, the valuation $\nu$ is monomial
with respect to the
$\xi_i$. Indeed, it suffices to check the definition for homogeneous elements. Every homogeneous element
$\overline x_\phi$ can be written as a finite sum of monomials of the same degree $\overline
x_\phi=\sum_\gamma c_\gamma\xi^\gamma$ with
$c_\gamma\in R/{\mathbf p} \ \hbox{\rm and}\ \nu_{\rm gr}(c_\gamma\xi^\gamma )=\phi$, as
required.\par\noindent
3) If $R$ is a one-dimensional excellent henselian equicharacteristic local domain, and $\nu$ its unique non-trivial valuation,
there are, as we shall see (subsection \ref{branchrev} and Proposition \ref{TC(*)}), finite systems of generators of the maximal ideal of $R$ with
respect to which $\nu$ is analytically monomial; \textit{in this case, the valuation is also algebraically monomial with respect to any such system of
generators in the strong sense given by the theory of semiroots of {\rm [PP]}, where the coefficients in $R$ of the monomials $\xi^\alpha$ depend only on one
variable}.\par The argument in the case where
$R$ is complete is as follows: Let
$k$ be the residue field of
$R$ and let
$(\xi_1,\ldots ,\xi_{g+1})$ be elements of $R$ whose images
$(\overline
\xi_i)$ constitute a system of generators of the graded
$k$-algebra
$\hbox{\rm gr}_\nu R$ (see subsection \ref{branchrev} below).  Let us fix a field of representatives $k\subset R$. Then, by the usual method of successive
approximations, using the fact that the
$\nu$-adic and $m$-adic topologies coincide in this case,  we have a surjection of
$k$-algebras $$k[[u_1,\ldots ,u_{g+1}]]\to R$$ mapping each
$u_i$ to $\xi_i$ and inducing, by passing to the graded rings with respect to $\nu$ and the corresponding monomial order on $k[[u_1,\ldots ,u_{g+1}]]$, the
surjective map
$$k[U_1,\ldots ,U_{g+1}]\to
\hbox{\rm gr}_\nu R $$ mapping $U_i$ to $\overline\xi_i$. Then each element $x$ of $R$ can be written $$\leqno{(**)}\ \ \ \ \ \ \ \  x=\sum_{\alpha
/\nu(\xi^\alpha\geq\nu (x)} c_\alpha\xi^\alpha\ \ \ \hbox{\rm with}\ c_\alpha\in k,$$ which shows that $(\xi_1,\ldots ,\xi_{g+1})$ constitute a
system of generators of the maximal ideal of $R$ with respect to which the valuation
$\nu$ is analyticaly monomial. For each index
$2\leq i\leq g+1$, in the one-dimensional subring
$k[[\xi_1,\xi_i]]\subset R$ with the valuation induced by $\nu$, the minimality of the valuation of
$\xi_1$ implies that $\xi_i$ is integral over the ideal
$\xi_1k[[\xi_1,\xi_i]]$, so that we have a relation $\xi_i^{s_i}+b^{(i)}_1(\xi_1,\xi_i)\xi_i^{s_i-1}+\cdots  +b^{(i)}_{s_i}(\xi_1,\xi_i)=0$ with
$b^{(i)}_\ell(\xi_1,\xi_i)\in (\xi_1^\ell)k[[\xi_1,\xi_i]]$. Using the Weierstrass preparation theorem, we see that we have relations
$$\xi_i^{s_i}+a^{(i)}_1\xi_i^{s_i-1}+\cdots  +a^{(i)}_{s_i}=0\
\ \
\hbox{\rm with }\ a^{(i)}_\ell\in (\xi_1)^{\ell}k[[\xi_1]].$$ Using this relation to eliminate from the representation $(**)$ all terms of order $\geq s_i$ in
$\xi_i$ for
$i\geq 2$, we see that we can now write
$$x=\sum_{\hbox{\rm finite}}a_\alpha (\xi_1)\xi_2^{\alpha_2}\ldots \xi_{g+1}^{\alpha_{g+1}},$$ where $a_\alpha(\xi_1)\in k[[\xi_1]]$ and the sum
involves only terms such that
$$\nu(a_\alpha)+\sum_{i=2}^{g+1}\alpha_i\nu (\xi_i)\geq\nu (x).$$Writing $a_\alpha (\xi_1)=\xi_1^{\alpha_1}v_\alpha (\xi_1)$ with $v_\alpha(\xi_1)$ a
unit in
$k[[\xi_1]]$, this shows that the
$(\xi_i)_{1\leq i\leq g+1}$ form a system of  generators for the maximal ideal of $R$ with respect to which $\nu$ is algebraically monomial in the strong
sense. The proof in the henselian case is similar.  For a different proof in the case of plane branches, see [PP]. 
\par We shall see in Proposition
\ref{coordinates} that if
$\nu$ is a rational valuation, after extending it to a suitable completion $\hat R^{(\nu)}$ one always has a (possibly infinite) system of generators of the
maximal ideal of $\hat R^{(\nu)}$ with respect to which this extension of $\nu$ is analytically monomial. One may then, by analogy with the curve case,
choose a system of parameters
$(\xi_1,\ldots ,\xi_s)$ with
$s=\hbox{\rm dim}\hat R^{(\nu)}$ such that the integral closure of the ideal $(\xi_1,\ldots ,\xi_s)\hat R^{(\nu)}$ is the maximal ideal, and examine under
which conditions the valuation is algebraically monomial in this stronger sense with respect to this system of
generators and this system of parameters, or becomes so after a birational extension.
\end{Remark}
\begin{example} 1) Let $(S,m)$ be a regular local
ring and $\nu$ be a valuation of $S$ having the property of the above definition with
respect to a system of regular parameters $(\xi_0,\ldots ,\xi_n)$ of $S$. If we denote by ${\mathcal
R}_\phi$ the ideal of $R$ generated by the monomials in $(\xi_0,\ldots ,\xi_n)$ of valuation
$\geq \phi$, we have $x\in  {\mathcal R}_{\nu (x)}$ and therefore $\nu$ is monomial in the sense
of [S2]. In this case $\hbox{\rm gr}_\nu R$ is a polynomial ring in $n+1$ variables, and this is clearly the {\it only case} where $\hbox{\rm gr}_\nu
R$ is regular, that is, a polynomial ring.\par\noindent
 2) If $R$ is the algebra of a plane branch
$C$, and $\nu$ its unique valuation, then $\nu$ is monomial with respect to two well chosen generators of the
maximal ideal of $R$ if and only if the semigroup $\nu (R\setminus \{0\})$ of $C$ has two generators. In characteristic zero, this means that $C$ has only
one characteristic pair (see [G-T], [T1]).
\end{example}
\begin{proposition}\label{grismon} Let
$R$ be an integral domain and $\nu$ a valuation on $R$, i.e., $R\subset R_\nu$. Set ${\mathbf
p}=m_\nu \cap R$. The valuation $\nu_{\rm gr}$ on
$\hbox{\rm gr}_\nu R$ is monomial with respect to any system of homogeneous generators
$(\xi_i)_{i\in I}$ of the graded $R/{\mathbf p}$-algebra $\hbox{\rm gr}_\nu R$.
\end{proposition}
\begin{proof} by the definition of $\nu_{\rm gr}$, this follows from one of the remarks made above.
 \end{proof}
 \section{\textbf{The toric structure of the graded algebras of valuation rings and n\oe therian rings}}\label{graval} 
In this section I show that the graded ring $\hbox{gr}_\nu R$ is, when $\nu$ is a rational valuation of finite rational
rank for the local integral domain $R$, a quotient of a polynomial ring in countably many indeterminates over the
residue field by a binomial ideal, and that
$\hbox{gr}_\nu R_\nu$ is the direct limit of a nested sequence of polynomial subalgebras in $\hbox{\rm r}(\nu)$ variables over the
residue field $k_\nu$ of the valuation. This means in particular that it is ``regular'', and hence can be deemed
to be a resolution of singularities of the algebra $\hbox{\rm gr}_\nu R$, if we forget for a
moment that $\hbox{gr}_\nu R_\nu$ is not an extension of finite type of $\hbox{\rm gr}_\nu R$,
and remember only that they have the same field of fractions. It means also that the binomial equations defining
$\hbox{gr}_\nu R_\nu$ are of a very special form.
\par It turns out that the result for $\hbox{gr}_\nu R_\nu$ is equivalent to its
counterpart for $k[t^{\Phi_+}]$ where $k$ is any field. \par The result for $k[t^{\Phi_+}]$ in
turn follows from the stronger statement that the semigroup $\Phi_+$ is a direct limit of a system
of linear maps between finitely generated free monoids, which will give us toric maps between the
corresponding polynomial algebras. \par We shall see that this is true, even with finitely
generated free {\it submonoids}, whenever
$\Phi$ is of finite height.
\subsection{The graded algebra of a rational valuation is essentially toric}\label{essto}
Here and in the sequel, unless otherwise specified, the notation $V^m$ stands for
a monomial $V_{j_1}^{m_{j_1}}\ldots V_{j_s}^{m_{j_s}}$, with $m_{j_i}\in \N\cup \{ 0\}$.
\begin{proposition}\label{dimone} Assume that a valuation ring $R_\nu$ birationally
dominates a local ring $R$ and that the residual extension $k_R=R/m\rightarrow
R_\nu/m_\nu =k_\nu$ is trivial. Then each homogeneous component of the graded
$k_R$-algebra $\hbox{\rm gr}_\nu R$ is a vector space of dimension $\leq 1$ over $k_R$.
\end{proposition} \begin{proof} Indeed, by definition of the ideals $\P_\phi (R)$, we have
the graded inclusion $$\hbox{\rm gr}_\nu R\subset \hbox{\rm gr}_\nu R_\nu$$
of $k_R$-algebras and the result is
true for the second algebra since the initial forms of two elements of
$R_\nu$ with the same valuation differ by a factor which is the initial form of a unit,  i.e., an
element of $k_\nu^*$; each homogeneous component of $\hbox{\rm gr}_\nu R_\nu$ is a
$k_\nu$-vector space of dimension one. 
\end{proof}
\par\noindent 
 From Proposition \ref{dimone} above follows by a direct generalization of a result of Korkina
([Ko], see also [E-S]) the
\begin{proposition}\label{binomialideal} In the situation of Proposition \ref{dimone},
assuming that $\hbox{\rm r}(\nu )$ is finite, the  graded $k_R$-algebra $\hbox{\rm
gr}_\nu R$  is the quotient of a polynomial ring in countably many  indeterminates
$k_R[(U_j)_{j\in J}]$ by a binomial prime ideal. \end{proposition}
\begin{proof} By assumption, $\Phi$ is of
finite rational rank, and therefore countable. Let
$(\overline\xi_j)_{j\in J}$ be countable system of nonzero homogeneous generators for the
$k_R$-algebra $\hbox{\rm gr}_\nu R$ (note that the number of these generators may be finite);
consider the kernel
$K$ of the surjective homomorphism
$$k_R [(U_j)_{j\in J}]\rightarrow \hbox{\rm gr}_\nu R$$ 
determined by $U_j\mapsto \overline\xi_j$. 
It is a
homogeneous ideal for the $\Phi_+$-grading giving to each $U_j$ the degree of $\overline\xi_j$. Let
$P=\sum c_mU^m$ be a homogeneous generator of $K$, say of degree
$\phi$. All the monomials occurring in $P$ have their images in the one
dimensional
$k$-vector space $(\hbox{\rm gr}_\nu R)_\phi$ and the sum of these images is zero.  By
induction on the number of monomials of $P$, we see that $P$ is a sum of binomial terms
$c_{m,m'}(U^m-\lambda_{m,m'}U^{m'})$ with
$c_{m,m'},\lambda_{m,m'}\in k$.\par\noindent The ideal $K$ is generated by binomials of the
form
$U^m-\lambda_{m,m'}U^{m'}$ and there is one such for each $(m,m')$ such that $U^m$ and
$U^{m'}$ have the same degree. By construction, all the $\lambda_{m,m'}$ are $\neq
0$. Since $\hbox{\rm gr}_\nu R$ is integral, the ideal $K$ is prime.
\end{proof}
\begin{corollary}\label{binomials} For any valuation ring $R_\nu$ with countable value group, the graded 
algebra $\hbox{\rm gr}_\nu R_\nu$ is a quotient of a polynomial algebra $k[(U_i)_{i\in I}]$ in countably many
variables over $k_\nu=R_\nu/m_\nu$ by a prime ideal generated by binomials $U^m-\lambda_{mn}U^n$ with $\lambda_{mn}\in k^*$.
\end{corollary}
\begin{remark} This corollary is used in the proof of Proposition \ref{grapp}, but does not
tell us anything precise about the form of the binomial equations.  \end{remark}
Assume now that $R$ is n\oe therian.\par\medskip  Let $\nu (R\setminus
\{0\})=\Gamma \subset \Phi_+ \cup \{ 0\}$ be the semigroup of $R$. Since $R$ is n\oe therian, the
semigroup $\Gamma$ is well ordered (Proposition \ref{wellord}), and we have a countable system of
generators, the valuations of the generators $\overline\xi_i$ of $\hbox{\rm gr}_\nu R$, i.e., $\gamma_i=\nu_{\hbox{\rm gr}}(\overline\xi_i)$: 
$$\Gamma=\langle\gamma_1,\gamma_2,\ldots ,\gamma_j,\ldots\rangle\ \ \ \ \hbox{\rm with}\ \gamma_j\in \Phi_+, \ j\ \hbox{\rm an ordinal}\
<\omega^{\hbox{\rm h}_R(\nu )},\ \ \gamma_j<\gamma_{j+1}.$$ (Corollary \ref{gammai}).\par To this situation
is associated the map of groups   $$\Z^\N\rightarrow \Phi\ \ \ \ \ \hbox{\rm determined
by}\ \ \  \sum\lambda_je_j\mapsto \sum \lambda_j\gamma_j,$$
where $e_j$ is the $j$-th basis vector, such that the image of
$(\Z_{\geq 0})^\N$ is $\Gamma$. It defines a map of $k_\nu$-algebras
$$k_\nu [(U_j)_{j\in J}]\rightarrow k_\nu [t^{\Phi_+}]\ \ \ \ \hbox{\rm determined by}
\ \ \ U_j\mapsto t^{\gamma_j}.$$ The image of this map is the semigroup algebra
$k_\nu [t^\Gamma]\subset k_\nu [t^{\Phi_+}]$, and its kernel is the binomial ideal obtained 
from that of of Proposition \ref{binomialideal} by setting all
$\lambda_{m,m'}$ equal to $1$.
\begin{remark}\label{truebin}We note that since the
$\gamma_i$ are assumed to be a minimal system of generators, none of the binomial equations is of the form
$$U_i-\lambda_iU^{m(i)}=0.$$
\end{remark}
Let $(R,\nu)$ be a valued n\oe therian ring; assume that the center of $\nu$ in $R$ is the maximal ideal $m$; set $k=R/m$ and let $$k[(U_i)_{i\in
I}]\to
\hbox{\rm gr}_\nu R\
\ \ \hbox{\rm determined by}\ \ U_i\mapsto \overline\xi_i$$ be a presentation of the $k$-algebra
$\hbox{\rm gr}_\nu R$ corresponding to a minimal system of generators of the semigroup $\Gamma$ of $\nu$. Let $(\xi_i)\in R$ be a system of
representatives of the $\overline\xi_i$. Let $R_\nu$ be the valuation ring of $\nu$ and consider all the images $\lambda_{mn}$ in $R_\nu/m_\nu$ of
the ratios of valuation zero $\xi^{m-n}$ of monomials in the $\xi_i$. For each $\lambda_{mn}$ which is algebraic over the field $k$,
consider the minimal polynomial $P_{mn}\in k[T] $ of
$\lambda_{mn}$. Setting $T=\frac{U^m}{U^n}$, after multiplication by a suitable power of $U^n$, $P_{mn}$ may be considered as the image in $k[(U_i)_{i\in
I}]$ of an irreducible polynomial. The polynomials $P_{mn}$ belong to the kernel $I$ of the surjective map $k[(U_i)_{i\in I}]\to\hbox{\rm gr}_\nu R$.
After extension of scalars to $k_\nu=R_\nu/\hbox{\rm m}_\nu$, the ideal $Ik_\nu[(U_i)_{i\in I}]$ is contained, in view of
Proposition \ref{binomialideal} in the ideal of $R_\nu/\hbox{\rm m}_\nu [(U_i)_{i\in I}]$ generated by all the binomials
$U^m-\lambda_{mn}U^n$, with $\lambda_{mn}\in k_\nu$ and $\nu (\xi^{m-n})=0$. By the faithful flatness of the extension
$k[(U_i)_{i\in I}]\to k_\nu[(U_i)_{i\in I}]$, the ideal $I$ is contained in the intersection of that ideal with
$k[(U_i)_{i\in I}]$, which is the ideal generated by all the polynomials $P_{mn}(U)$ for $\lambda_{mn}$ algebraic over $k$; the ideal $I$ is equal to the
ideal generated by the $P_{mn}(U)$. For an arbitrary valuation, say with center ${\mathbf p}$, we may apply this to the extension of $\nu$ to the
localization
$R_{\mathbf p}$ with residue field
$\kappa ( {\mathbf p})=R_{\mathbf p}/{\mathbf p}R_{\mathbf p}$:
\begin{corollary}\label{roots} Let $(R,\nu)$ be a valued n\oe therian ring; let ${\mathbf p}$ be the center of $\nu$ in $R$ and extend $\nu$ to the
localization $R_{\mathbf p}$. The kernel
$I$ of the map of $\kappa({\mathbf p})$-algebras
$\kappa({\mathbf p})[(U_i)_{i\in I}]\to \hbox{\rm gr}_\nu R_{\mathbf p}$ determined by $U_i\mapsto \overline\xi_i$  is generated by the
polynomials
$P_{mn}$ obtained by clearing denominators in the minimal polynomials over $\kappa ({\mathbf p})$ of the images in $R_\nu/m_\nu$ of the ratios
$\xi^{m-n}$ of valuation
$0$, viewed as polynomials in $\frac{U^m}{U^n}$.
\end{corollary}
\begin{proposition} Let $(U^m-\lambda_{mn}U^n)_{(m,n)\in E}$ be the defining binomial ideal in $k[(U_i)_{i\in I}]$ for the graded ring $\hbox{\rm gr}_\nu
R$ of a rational valuation of a n\oe therian ring. Let $(\Lambda_{mn})_{(m,n)\in E}]$ be a set of indeterminates and let $\Lambda$ denote the multiplicative
subset of
$k_\nu [(\Lambda_{mn})_{(m,n)\in E}]$ generated by the $\Lambda_{mn}$. The natural map of
$k_\nu$-algebras $$\Lambda^{-1}k_\nu [(\Lambda_{mn})_{(m,n)\in E}]\rightarrow
\Lambda^{-1}k_\nu [(\Lambda_{mn})_{(m,n)\in E}][(U_j)_{j\in
J}]/(U^m-\Lambda_{mn}U^n)_{(m,n)\in E}$$
\noindent  is faithfully flat and corresponds to a flat
family of schemes having one fiber equal to $\hbox{\rm Spec}k_\nu[t^\Gamma]$ and another 
to $\hbox{\rm Spec}\hbox{\rm gr}_\nu R$.
\end{proposition} 
Indeed, if the second
algebra is graded by giving to $U_j$ the degree $\gamma_j$, then by definition of the binomial
relations, each of its homogeneous components becomes a free module of rank one over
$\Lambda^{-1}k_\nu [(\Lambda_{mn})_{(m,n)\in E}]$. \par *This family
parametrized
by the $\Lambda_{m,n}$ has simultaneous resolution along the subspace defined by the vanishing of the $(U_i)_{i\in I}$ in any reasonable sense; it is locally
analytically trivial along this subspace.* We shall see below that it has simultaneous embedded resolution by a
toric map; it follows from the results of subsection \ref{toricresfin}.
\subsection{More on the structure of $\hbox{gr}_\nu R$}\label{morestruc}
Let $\Phi$ be a group of height $h$ associated to a valuation of a local integral
domain $R$ centered at the maximal ideal $m$ of $R$; set $k=R/m$. Let
$$(0)=\Psi_h\subset \Psi_{h-1}\subset\cdots\subset \Psi_1\subset \Psi_0=\Phi$$
be the sequence of isolated subgroups of $\Phi$ (including the trivial ones). Let $\Gamma$ be
the semigroup of $\nu$ on $R$ (note that it may be equal to $\Phi_+$, if $R=R_\nu$) and set
$\Gamma_i=\Gamma\cap \Psi_i$. Considering the valuation $\nu_1$ of height $h-1$ with
center ${\mathbf p}_1$ with which $\nu$ is composed, we have seen in subsection \ref{compargr}
that we could identify $\hbox{\rm gr}_{\overline \nu}\overline R_1$, where $\overline R_1=R/{\mathbf
p}_1$, with the subalgebra $\bigoplus_{\psi\in \Psi_{h-1+}\cup \{0\}}(\hbox{\rm gr}_\nu
R)_\psi$ of $\hbox{\rm gr}_\nu R$. If we distinguish between the generators $(U_i)_{i\in I_1}$, where $I_1\subset I$,
of the $k$-algebra $\hbox{\rm gr}_\nu R$ whose degree lies in $\Psi_{h-1}$ and those whose
degree lies in $\Phi\setminus \Psi_{h-1}$, we have:
\begin{proposition}\label{strucgr} If $R\subset R_\nu\subset R_{\nu_1}$, where $R_\nu$
dominates the local ring $R$ without residual extension and $\nu_1$ is of height one less than $\nu$, the natural map
$$\hbox{\rm gr}_{\overline \nu}\overline R_1[(U_i)_{i\in I\setminus I_1}]\to \hbox{\rm
gr}_\nu R$$ is surjective and its kernel is generated by the images in $\hbox{\rm gr}_{\overline \nu}\overline R_1[(U_i)_{i\in I\setminus I_1}]$ of those
binomials
$U^n-\lambda_{mn}U^n$ which appear in Proposition \ref{binomialideal} and which involve at least one
variable $U_i$ with
$i\in I\setminus I_1$.\end{proposition}
\begin{proof} The subalgebra $\hbox{\rm gr}_{\overline \nu}\overline R_1$ contains all the generators of the $k$-algebra $\hbox{\rm
gr}_ \nu  R$ whose degree is in $I_1$; we need to add only the generators $U_i$ whose degree lies in $I\setminus I_1$, and the relations
between them are those indicated. \end{proof}
\begin{corollary}\label{strugr} If $R_\nu$ dominates $R$ without residual extension, setting $I_t=\{i\in I/\hbox{\rm deg}U_i\in
\Psi_{h-t}\setminus \Psi_{h-t+1}\}$, we can write \font\eightrm=cmr8 $$\begin{array}{lr}\hbox{\rm gr}_\nu
R=\\ \scriptstyle{\bigl(\cdots\bigl( k[(U_i)_{i\in I_1}]/(U^n-\lambda_{mn}U^n)\bigr)[(U_i)_{i\in
I_2}]/(U^n-\lambda_{mn}U^n)\bigr)\cdots\bigr)[(U_i)_{i\in
I_h}]/(U^n-\lambda_{mn}U^n)},\end{array}$$ where at each step $t$ the binomials which are
written involve only variables $U_i$ with degrees in $\Psi_{h-t}$ and must involve one which
has degree in $I_t$. Each monomial of such binomials must then involve variables $U_i$ with $i\in I_t$. \end{corollary}
\begin{proof}The first part is only the iteration of the preceding Proposition. The second part follows from
the fact that there can be no linear expression of  an element of $I_t$ in terms of  those of $\Psi_{h-t+1}$.\end{proof}

\subsection{$\hbox{gr}_\nu R_\nu$ as a limit of polynomial algebras}
In this subsection, I show that $\hbox{gr}_\nu R_\nu$ is the union of a nested sequence of polynomial subalgebras in
$\hbox{\rm r}(\Phi)$ variables over $k_\nu$, generated by homogeneous elements, the inclusions sending variables to terms, i.e., monomials multiplied
by non zero constants. Technically this is essentially Perron's algorithm. This result, however, may be considered as a graded version of local uniformization;
each finitely generated
$k_\nu$-subalgebra of $\hbox{gr}_\nu R_\nu$ is contained in a  polynomial subalgebra in $\hbox{\rm r}(\Phi)$ variables. However this inclusion is
birational only in special cases, for example if our subalgebra contains a finitely generated $\hbox{gr}_\nu R$ with
$R\subset R_\nu$ birational and residually rational.\par\medskip\noindent    I first prove the corresponding result for semigroup algebras
$k[t^{\Phi_+}]$, and I begin with an
example: the case when
$\Phi =\Z^d$ with the lexicographic order. \par\noindent Indeed, consider for each $(d-1)$-uple $(r_1, r_2, \ldots
,r_{d-1})$ of positive integers the map 
$$e_{r_1, r_2, \ldots ,r_{d-1}}\colon \N^d\rightarrow \Z^d$$  defined by $$ e_{r_1, r_2, \ldots
,r_{d-1}}(i_1,\ldots ,i_d)= (i_1,i_2-r_1i_1,\ldots ,i_d-r_{d-1}i_{d-1});$$   its image is in the
semigroup $\Z^d_+$ of positive elements, and it is injective. It defines an injective map of
semigroup algebras $$e_{r_1, r_2, \ldots ,r_{d-1}}\colon k[z_1,\ldots ,z_d] \rightarrow
k[t^{\Z^d_+}],$$ and clearly this last algebra is  the union of these subalgebras as $r_1, r_2,
\ldots ,r_{d-1}$ vary. Note that whenever $(r_1, r_2, \ldots ,r_{d-1})\leq (s_1, s_2, \ldots
,s_{d-1})$ in the sense that $r_j\leq s_j\ ,1\leq j\leq d-1$, the map $e_{r_1, r_2, \ldots ,r_{d-1}}$ factors through a map of algebras
$e_{s_1-r_1,s_2- r_2, \ldots ,s_{d-1}-r_{d-1}}\colon k[z_1,\ldots ,z_d]\rightarrow k[w_1,\ldots ,w_d]$ which is toric,  i.e., given by monomials, and
that $k[t^{\Z^d_+}]$ is the direct limit of this system of toric maps. \par\medskip When $\Phi$ is
a subgroup of ${\mathbf R}$, consisting of numbers of the form $\sum_{i=1}^m a_i\tau_i$ with given
$\tau_i >0$ and $a_i\in \Z$, it is also true that $k[t ^{\Phi_+}]$ is a direct limit of polynomial
subalgebras over $k$.  This is again because the half space in ${\mathbf R}^m$ defined by
$\sum_{i=1}^m a_i\tau_i\geq 0$ can be approximated by regular rational simplicial cones
(generated by integral vectors which form a basis of the integral lattice).\par More precisely,
let $\nu$ be a valuation of height one, i.e., with archimedian value group $\Phi\subset {\mathbf R}$
(see [A1], p.46, [Z-S], Vol. II). Assume first that $\Phi$ is generated by $m$ rationally independent
real numbers $\tau_1,\ldots ,\tau_m$, which we may assume to be positive. I use the Perron algorithm as expounded in
([Z1], B. I, p. 861), but with a somewhat different interpretation. The algorithm consists in
writing $$\tau_1=\tau_m^{(1)} ,\tau_2=\tau_1^{(1)}+a_2^{(0)}\tau_m^{(1)},\ldots
,\tau_m=\tau_{m-1}^{(1)}+a_m^{(0)}\tau_m^{(1)},$$ where $$a_j^{(0)}=[\tau_j/\tau_1],\ \ \
j=2,\ldots ,m,$$ and repeating this operation after replacing $(\tau_1,\ldots ,\tau_m)$ by
$(\tau_1^{(1)},\ldots ,\tau_m^{(1)})$, and so on. After $h$ steps, one has written
$$\tau_i=A_i^{(h)}\tau_1^{(h)}+\cdots +A_i^{(h+m-1)}\tau_m^{(h)}$$ or,  if we denote by $w$ the
(weight) vector $(\tau_1,\ldots ,\tau_m)\in {\mathbf R}^m$ and by $A^{(h)}$ the vector 
$(A_1^{(h)},\ldots ,A_m^{(h)})$, $$w=\tau_1^{(h)}A^{(h)}+\tau_2^{(h)}A^{(h+1)}+\cdots
+\tau_m^{(h)}A^{(h+m-1)}$$ where the $\tau^{(h)}_j$ are positive, the coefficients $A^{(j)}_i$ are
non negative integers, and the matrix of the vectors $$A^{(h)},A^{(h+1)},\ldots ,A^{(h+m-1)}$$
has determinant $(-1)^{h(m-1)}$. Moreover, as $h$ grows the directions in ${\mathbf P}^{m-1}({\mathbf R})$ of the
vectors $A^{(h)}$ tend to the direction of $w$. So we have a sequence of vectors $A^{(h)}$ with
positive integral coordinates whose directions in ${\mathbf
P}^{m-1}({\mathbf R})$ spiral to the direction of $w$ and such that any consecutive $m$ of them as
above form a basis of the integral lattice such that $w$ is contained in the convex cone
$\sigma^{(h)}=\langle A^{(h)},A^{(h+1)},\ldots ,A^{(h+m-1)}\rangle$ which they generate. The
convex dual $\check\sigma^{(h)}$ of $\sigma^{(h)}$ (see [Cox], \S 2, [E], V, 2, p. 149) is contained
in the half space $\sum_{i=1}^ma_i\tau_i\geq 0$, the integral points of which form the
semigroup $\Phi_+$. The algebra of the semigroup $\check \sigma^{(h)}\cap \Z^m$ is a
polynomial algebra $k[x_1^{(h)}, \ldots,x_m^{(h)}]$ ({\it loc.cit.}, VI,2) contained in
$k[t^{\Phi_+}]$, and since by assumption there are no integral points on the hyperplane 
$\sum_{i=1}^ma_i\tau_i= 0$ except the origin, the semigroup $\Phi_+$ is the union of the
$\check \sigma^{(h)}\cap \Z^m$ as $h\rightarrow \infty$. This proves that $k[t^{\Phi_+}]$ is the
union, or direct limit, of these polynomial subalgebras.\par Note that by construction we have
for each $h\geq 1$ the equality  $$A^{(h)}-A^{(h+m)}+a_2^{(h)}A^{(h+1)}+\cdots
+a_m^{(h)}A^{(h+m-1)}=0,$$ which shows, since the $a_j^{(h)}=[\tau_j^{(h)}/\tau_1^{(h)}]$ are
non-negative, that we have $$A^{(h+m)}\in\sigma^{(h)}=\ \langle A^{(h)},\ldots
,A^{(h+m-1)}\rangle$$ and therefore $\sigma^{(h+1)}\subset\sigma^{(h)}$, that is $\check
\sigma^{(h)}\subset \check\sigma^{(h+1)}$ and $$ k[x_1^{(h)}, \ldots,x_m^{(h)}]\subset
k[x_1^{(h+1)}, \ldots,x_m^{(h+1)}],$$ so that our direct system is in fact a nested sequence of
polynomial subalgebras. The morphisms between these polynomial algebras correspond by duality to the expressions 
of the $A^{(h+m)}$ as linear combinations with non negative integral coefficients of $(A^{(h)},\ldots
,A^{(h+m-1)})$ and are therefore monomial as announced.\par If we now consider
a group  with one more generator
$\tau_{m+1}>0$ which is rationally dependent on $\tau_1,\ldots ,\tau_m$, Zariski shows in ([Z1], B. I, p. 862) that the new weight vector 
$w=(\tau_1,\ldots ,\tau_m ,\tau_{m+1})\in {\mathbf R}^{m+1}$ is contained in a rational simplicial  cone $\sigma \subset {\mathbf R}^{m+1}$
generated by $m$ integral vectors $v_1,\ldots ,v_m$ of the first quadrant  forming part of a basis of the integral lattice.
Indeed $w$ is contained in a unique rational hyperplane. The dual cone $\check
\sigma\subset\check {\mathbf R}^{m+1}$  is the product of an $m$ dimensional strictly convex cone generated by vectors $e_1,\ldots, e_m$ by a
1-dimensional vector space (see [E], V, 2), generated by a primitive integral vector $e_{m+1}$, which is the dual of the rational hyperplane containing $w$.
The vectors $e_1,\ldots, e_{m+1}$ are a basis of the integral lattice, and correspond to the variables $x_1,\ldots ,x_{m+1}$ generating a polynomial ring. Note
that the map
$f:
\Z^{m+1}\rightarrow {\mathbf R}$ defined by $(a_1,\ldots ,a_{m+1})\mapsto \sum_{i=1}^{m+1}a_i\tau_i$ is no longer injective; 
the primitive vector $e_{m+1}$ corresponding to the variable $x_{m+1}$ is in the kernel. Let us set
$\widetilde{\Phi_+}=f^{-1}(\Phi_+\cup\{ 0\})$. By refining as above by the Perron algorithm for
$w$ inside the linear span $<\sigma >$ of $\sigma$, starting with the coordinates of
$w$ in $\sigma$,
 we find a sequence of regular simplicial cones $\sigma^{(h)}\subset \sigma $ whose duals
$\check
\sigma^{(h)}\subset \widetilde{\Phi_+}$ correspond ([E], VI, Th. 2.12) to algebras of the
form $k[x_1^{(h)}, \ldots ,x_m^{(h)}, x_{m+1}^{\pm 1}]\subset k[t^{\widetilde{\Phi_+}}]$.  
The free semigroups $\check \sigma^{(h)}\cap \Z^{m+1}$ fill up $\widetilde{\Phi_+}$ as
$h\rightarrow \infty$ since the only  rational points of the hyperplane
$\sum_{i=1}^{m+1}a_i\tau_i=0$ are on the dual of the hyperplane containing $w$, which is
contained in all the $\check\sigma^{(h)}$. So the direct limit of the images of the maps
$k[t^f]\colon k[x_1^{(h)},
\ldots ,x_m^{(h)}, x_{m+1}^{\pm 1}]\rightarrow k[t^{\Phi_+}]$ is $k[t^{\Phi_+}]$. But these images
are isomorphic to $ k[x_1^{(h)},
\ldots ,x_m^{(h)}, x_{m+1}^{\pm 1}]/(x_{m+1}-1)$ so that they are again polynomial rings
 $k[x_1^{(h)}, \ldots,x_m^{(h)}]$. If we have more generators rationally dependent on
$\tau_1,\ldots ,\tau_m$, we can repeat the argument after taking as new generators the
coordinates of the weight vector with respect to the $m$ primitive vectors of
$\sigma$. In both examples, the fact that $k$ is a field plays no role, so we have
proved:
\begin{lemma} Let $\Phi$ be a totally ordered finitely generated
group of height one (i.e., archimedian), or
$\Z^d$ with the lexicographic order. For any commutative ring $A$ the semigroup algebra
$A[t^{\Phi_+}]$ of $\Phi_+$ with coefficients in $A$ is the direct limit of a direct system of graded
subalgebras which are polynomial algebras $A[x_1,\ldots ,x_m]$ over $A$ with
$m=\hbox{\rm r} (\Phi)$.
\par In addition, the maps between these algebras are toric maps,  i.e., each variable of
one is sent to a monomial in the variables of the other, and there is a cofinal subsystem which is
a chain of nested subalgebras.
\end{lemma}
\begin{Remark} 1) In both cases, the smooth
subalgebras  are produced by an algorithm. \par\noindent 2) The result, especially in view of its
algorithmic nature, is much more useful than its consequence that $\hbox{\rm
Spec}A[t^{\Phi_+}]$ is a pro-object in the category of smooth affine toric schemes over
$\hbox{\rm Spec}A$ with toric maps. \par\noindent 3) The proof also shows that the semigroup
algebra of the semigroup  $$\widetilde{\Phi_+}=
\{(a_1,\ldots ,a_m,a_{m+1},\ldots ,a_{m+r})\in \Z^{m+r}\mid \sum_{i=1}^{m+r}a_i\tau_i\geq
0\},$$ where 
$\tau_1,\ldots ,\tau_m$ are rationally independent and the others are rationally dependent upon
them, is a toric direct limit of toric subalgebras of the form 
$$A[x_1,\ldots ,x_m,x_{m+1}^{\pm 1},\ldots ,x_{m+r}^{\pm 1}].$$\par\noindent  4)  The result for
subgroups of ${\mathbf R}$ holds without the finiteness assumption provided the rational rank is finite, since any abelian group
 is a direct limit of its subgroups of finite type.  If we allow
$m$ to vary, even that last assumption is unnecessary.
\end{Remark}
\par\medskip Let now $\Phi$ be a totally ordered group of finite height $h>1$. We have a
surjective monotone non-decreasing map
$\lambda\colon \Phi\rightarrow \Phi_1$ where $\Phi_1$ is of height $h-1$, and the kernel
$\Psi$ of $\lambda$ is of height $1$. By induction on the height we may assume that
$\Phi_{1+}$ is the union of sub-semigroups isomorphic to $\N^m$, and we know from the lemma
above that the same is true for $\Psi_+$. \par Let us denote the free semigroups that fill $\Phi_{1+}$
by $F_i$, and let $\tilde F_i\subset \Phi_+$ be the subsemigroup generated by elements $e_1,
\ldots ,e_{r_i}$ which lift to $\Phi_+\setminus\Psi$ the generators of $F_i$. Similarly let us denote by
$G_j\subset \Psi_+$ free semigroups which fill
$\Psi_+$, generated say by
$f_1,\ldots ,f_{s_j}$. Note that for $\phi\in \Phi_+\setminus\Psi,\ \psi\in \Psi,\ \phi +\psi\in \Phi_+$,  and consider for $r_is_j$-tuples
$n=(n_{s,t},1\leq s\leq r_i,\ 1\leq t\leq s_j)$ of non negative integers, the free semigroups 
$\tilde F_i(n)\subset \Phi_+\cup \{ 0\}$ generated by $e_1-\sum_tn_{1t}f_t,\ldots
,e_{r_i}-\sum_t n_{r_it}f_t$.  Let us check that the direct sums of free semigroups
$ \tilde F_i(n)\oplus G_j$ fill up $\Phi_+$; the proof generalizes that given in the case of the
lexicographic $\Z^d$. Given $\phi \in \Phi_+$, there exists an index $i$ and a $\phi_1\in \tilde
F_i$ such that $\phi-\phi_1\in \Psi$. If we have $\phi -\phi_1\in \Psi_+$, there exists an index
$j$ such that $\phi -\phi_1\in G_j$ so that indeed $\phi \in \tilde F_i\oplus G_j$. This
happens in particular if $\phi\in \Psi$. If
$\phi -\phi_1\in \Psi_-$, there exists an index $j$ such that we can write:
$$\phi=\sum_{s=1}^{r_i}k_se_s-\sum_{t=1}^{s_j}\ell_tf_t\ \hbox{\rm with non negative
integers}\ k_s,\ \ell_t\ \hbox{\rm and}\ f_t\in G_j.$$\noindent  Since $\phi\notin \Psi$, at least one of the $k_s$ is not zero,
say $k_1$. Choosing positive integers $\tilde \ell_t$ such that $k_1\tilde\ell_t >\ell_t$, we may
rewrite $\phi$ as follows
$$\phi=k_1(e_1-\sum_{t=1}^{s_j}\tilde\ell_tf_t)+\sum_{s=2}^{r_i}k_se_s+\sum_{t=1}^{s_j}(k_1\tilde\ell_t-\ell_t)f_t.$$
This shows that $\phi$ is indeed in $ \tilde F_i(n)\oplus G_j$ with $n=\tilde \ell$.\par
Now if we assume, as we may by induction, that the $F_i$ and $G_j$ are nested sequences, and
choose $n(i)$ given by
$(n_{s,t}=i,1\leq s\leq r_i,\ 1\leq t\leq s_j)$, we see that the corresponding groups $ \tilde
F_i(n(i))\oplus G_i$ form a nested sequence which is cofinal in the direct system. So we
have:
\begin{proposition}\label{smgrpapp} For any totally ordered group $\Phi$ of finite height, and
any commutative ring $A$, the  semigroup algebra $A[t^{\Phi_+}]$ is a direct limit of a direct
system of graded polynomial subalgebras over $A$ with monomial maps, and there are cofinal
nested subsystems of such polynomial subalgebras.
 If $\Phi$ is of finite rational rank
$\hbox{\rm r}(\Phi)$, all the polynomial subalgebras may be chosen isomorphic to
$A[x_1,\ldots ,x_{\hbox{\rm r}(\Phi)}]$.
\end{proposition} 
\begin{proof} There remains only to prove the last
sentence. This is done by induction on the height; we have shown above that the result is true for
valuations of height one. Assume that the result is true for $\Psi$ and $\Phi_1$, and so the rank
of the free monoid
$\tilde F_i(n)\oplus G_j$ is the sum of the rational ranks of $\Psi$ and $\Phi_1$. But since the
rational rank is additive in an exact sequence because $\Qq$ is a flat $\Z$-module, this is the
rational rank of $\Phi$. 
\end{proof}

\begin{corollary}\label{goodeq}If the rational rank of $\Phi$ is finite, the semigroup algebra
$A[t^{\Phi_+}]$ endowed with its natural grading is a quotient of a polynomial ring over $A$ in countably
 many indeterminates $A[(V_j)_{j\in J}]$ graded by $\Phi_+$ by a homogeneous binomial ideal of the form 
$$(V_j- V^{m(j)})_{j\in J'},$$\noindent where $J'$ is a subset of $J$ and $\vert m(j)\vert\geq 2$ for $j\in J'$.
\end{corollary}
\begin{proof} We may choose as a system of homogeneous generators of the
$A$-algebra $A[t^{\Phi_+}]$, the union of the generators of the polynomial subalgebras
of which $A[t^{\Phi_+}]$ is the direct limit. The only relations between these generators are those
corresponding to the toric inclusions
$A[x_1,\ldots , x_r]\to A[y_1,\ldots ,y_r]$, and they are of the announced type after we discard
the trivial relations $x_i=y_j$ by removing some generators. 
\end{proof}
\begin{Remark} 1) There is in general no minimal system of generators for $A[t^{\Phi_+}]$ and we can of course discard an arbitrary number of those
generators which appear as $V_j$ in a relation $V_j- V^{m(j)}=0$. For example consider $\Phi=\Z^2_{lex}$. \par\noindent 2) We can apply this to
the subalgebra
$k[v^{\Phi_+}]$ of
$\A_\nu (R)$, and view the $v^\phi$ as coordinates for $\hbox{\rm
Spec}k[v^{\Phi_+}]$ subjected to the binomial relations described in this Corollary.
\end{Remark}
 \begin{proposition}\label{grapp} Given the ring $R_\nu$ of a
valuation of finite rational rank $r(\nu )$:\par\noindent a) the graded
$k_\nu$-algebra  $\hbox{\rm gr}_\nu R_\nu$ is a quotient of a polynomial ring
$k_\nu[(V_j)_{j\in J}]$ in countably many indeterminates over the residue field
$k_\nu$, graded by $\Phi_+$, by a\break homogeneous binomial ideal of the form  $$(V_j-\lambda_j V^{m(j)})_{j\in J'}\ ,\
\lambda_j\in k_\nu^*$$\noindent where $J'$ is a subset of $J$ and $\vert m(j)\vert\geq 2$ for $j\in J'$.\par\noindent  b) The graded $k_\nu$-algebra 
$\hbox{\rm gr}_\nu R_\nu$ is the union of a nested sequence of graded polynomial subalgebras $k_\nu[x^{(h)}_1,\ldots
,x^{(h)}_{\hbox{\rm r}(\nu )}]$, where the inclusions are given by maps sending each variable
$x^{(h)}_i$ to a constant times a monomial in the $x^{(h+1)}_j,\ 1\leq j\leq \hbox{\rm
r}(\nu )$. \end{proposition}
 \begin{proof}
See Proposition \ref{smgrpapp} and the previous Corollary. We have seen in Proposition
\ref{binomialideal} that $\hbox{\rm gr}_\nu R_\nu$ is isomorphic to a quotient of a
polynomial algebra $k[(V_i)_{i\in I}]$ by a binomial ideal whose generators are of the form
$V^m-\lambda_{mn}V^n$. Setting all the constants $\lambda_{mn}$ equal to one gives the
semigroup algebra $k[t^{\Phi_+}]$, hence assertion a). Assertion b)
follows from this and the correspondence between the direct
system of polynomial subalgebras and the binomial equations
exhibited in the proof of the Corollary to Proposition \ref{smgrpapp}.
\end{proof}
\begin{corollary} Given any finite set of homogeneous elements $\overline x_1,\overline x_2,\cdots ,\overline x_s$ in 
$\hbox{\rm gr}_\nu R_\nu$ and a cofinal nested system of polynomial subalgebras as above, there is an algebra in our nested system such that not
only do we have
$\overline x_i\in A[x^{(h)}_1,\ldots ,x^{(h)}_{\hbox{\rm r}(\Phi)}]$
 for $1\leq i\leq s$, but the element of least degree divides all the others in this subalgebra.\end{corollary}\par\noindent
\begin{remark}\label{directsystem} It is important to note that the direct system built in this
way is rather precisely determined by $\Phi$; we shall see on the next examples that at least
in dimension two the structure of the binomial equations, or equivalently that of the direct
system of polynomial subalgebras converging to $\hbox{\rm gr}_\nu R_\nu$, reflects the
structure of the system of points in birational modifications of $\hbox{\rm Spec}R$
corresponding to the valuation $\nu$, as well as the height and rational rank of
$\nu$.\par\medskip\noindent 
 {\textbf {Problem.}} It is an interesting problem to make these correspondences precise,
especially in dimensions $\geq 3$. A very useful reference for dimension two is [A2].
\end{remark}
\subsection{A selection of examples}\label{chexamples}
\begin{example} 1) I give a slightly different proof of the fact, which we saw at the
beginning of this section, that the algebra of the semigroup $\Z^d_+$ for the lexicographic
order is a toric union of polynomial algebras; the relation between the two proofs is the idea
for the induction on the height shown above. Let $\Phi$ be $\Z^d$ with the lexicographic
order; the semigroup $\Phi_+$ is generated by the following elements\par\noindent
 $e_j^{(1)}=(1,-j,0,\ldots ,0), e_j^{(2)}=(0,1,-j,\ldots ,0),\ldots ,e_j^{(d-1)}=(0,0,\ldots ,1,-j),\break 
e_0^{(d)}=(0,0,\ldots ,0,1)\ \ \hbox{\rm for}\ j=0,1,2,\ldots$. To check this, we proceed by
induction on $d$; the case $d=1$ is obvious, so let $d\geq 2$ and
$\phi =(a_1,\ldots ,a_d)\in \Z^d_+$. If we have $a_1=0$ we are reduced to $\Z^{d-1}$, so we may
assume $a_1>0$. If $a_2\geq 0$, we may write $\phi =a_1e_0^{(1)}+(0, a_2,\ldots ,a_d)$ and we
are again reduced to the case of $\Z^ {d-1}$. If $a_2<0$, let $\ell$ be the smallest positive integer
such that $\ell a_1+a_2\geq 0$; we can write $\phi =a_1e_\ell^{(1)}+(0, \ell a_1+a_2, a_3, \ldots
,a_d)$ and again we may conclude by induction. One can check, again by induction, that the
relations are generated by 
$e^{(i)}_j=e^{(i)}_{j+1}+e^{(i+1)}_0\ \hbox{\rm for}\ 1\leq i\leq d-1,\ j\geq 0$.
 The semigroup algebra $A[t^{\Z^d_+}]$ is isomorphic to $$A[(V_{j_1}^{(1)})_{j_1\in \N_0},\ldots
,(V_{j_{d-1}}^{(d-1)})_{j_{d-1}\in \N_0},
V_0^{(d)}]/\bigl((V^{(i)}_j-V^{(i)}_{j+1}V^{(i+1)}_0)_{j\in
\N_0, 1\leq i\leq d-1}\bigr),$$where $\N_0= \N\cup \{ 0\}$. Note that from these
equations we can read the fact that the valuation is of height $d$; they show that
for each $(i,j),\ 1\leq i \leq d-1,\ j\geq 1$ and each integer $n$ we have
$V_j^{(i)}=V_{j+n}^{(i)}(V_0^{(i+1)})^n$, so that the value in $\Phi_+$ of $V_j^{(i)}$ is
greater than $n$ times the value of $V_0^{(i+1)}$ for all $n\geq 1$. This implies the existence
of $d-1$ non trivial convex subgroups of $\Phi$; since we know that $\hbox{\rm h}(\Phi )\leq
\hbox{\rm r}(\Phi )$, we have $\hbox{\rm h}(\Phi )=d$.\par\noindent
 In particular, for $d=2$, we see that
$$A[t^{\Z^2_+}]=A[W, V_0,V_1,\ldots ,V_j,\ldots ]/\bigl((V_j-WV_{j+1})_{j\in
\N\cup\{ 0\}}\bigr).$$ Upon taking for $A$ a field $k$, we see that it corresponds to the sequence of
infinitely near points in the plane ${\mathbf A}^2(k)$ obtained by endlessly repeating the following
process: blow up the origin of an ordered set of coordinates $(u_1,u_2)$, choose the chart of the
blow up where the ideal $(u_1,u_2)$ is generated by $u_1$, and in this chart, with coordinates
$(u_1,u_2/u_1)$, blow up the origin. We see that the direct limit of the monomial maps
$k[u_1,u_2]\to k[u_1,u_2/u_1]$ obtained by repeating this operation is $k[t^{\Z^2_+}]$ presented
as above, with $V_j=u_2u_1^{-j},\ W=u_1$.\par By the correspondence between sequences of
points in blowing ups and valuations ([V1]), this sequence of points corresponds to a valuation
$\nu$ of the local ring $R=k[u_1,u_2]_{(u_1,u_2)}$, the monomial valuation in the coordinates $(u_1,u_2)$
with value group $(\Z^2,\hbox{\rm lex.})$ determined by $\nu (u_1)=(1,0),\ \nu (u_2)=(0,1)$. The
graded algebra $\hbox{\rm gr}_\nu R_\nu$ is equal to $k[t^{\Z^2_+}]$. The valuation $\nu$ is
monomial with respect to the coordinates $(u_1,u_2)$, so that we have $\hbox{\rm gr}_\nu
R=k[U_1,U_2]$, where $U_1,U_2$ are the initial forms of $u_1,u_2$. Things are quite similar for $d>2$.\end{example}
 \begin{example}\label{irrtau} Given a positive real number $\tau$, which I assume to be irrational, let us
denote by $$\tau=[s_1,s_2, \ldots ,s_i,s_{i+1},\ldots ]=s_1+\frac{1}{s_2+\cdots}$$ its continued
fraction expansion. Let $p_i/q_i=[s_1,\ldots ,s_i]$ be the $i$-th approximant. Using the
inductive relations $$p_{i+1}=p_{i-1}+s_{i+1}p_i,\ \ q_{i+1}=q_{i-1}+s_{i+1}q_i\ \ \hbox{\rm
with } p_0=1,\ q_0=0,\ p_1=s_1,\ q_1=1,$$ the reader will see that Perron's algorithm
applied to $(1,\tau )$ produces the vectors $A^{(h)}=(p_h,q_h)\in {\mathbf R}_+^2$ for $h\geq 1$.
If $\Phi$ denotes the subgroup of ${\mathbf R}$ generated by $(1,\tau)$ endowed with the induced order, the
algorithm gives the following presentation, for any commutative ring $A$: $$A[t^{\Phi_+}]=A[(V_i)_{i\in
\N}]/\big((V_i-V_{i+1}^{s_i}V_{i+2})_{i\in \N}\big).$$ It is equivalent to say
that it presents
$A[t^{\Phi_+}]$ as the direct limit of the maps $$A[V_i,V_{i+1}]\rightarrow A[V_{i+1},V_{i+2}]\
\hbox{\rm given by} \ 
 V_i\mapsto V_{i+1}^{s_i}V_{i+2},\  \ V_{i+1}\mapsto V_{i+1}.$$ The case where $\tau\in {\mathbf Q}_+$ is left as an
 exercise.\par\noindent We shall see below in examples \ref{exzar} and \ref{nextex} the interpretation of these maps in terms of blowing ups, the
corresponding valuation of
$k[u_1,u_2]_{(u_1,u_2)}$ and the graded algebra of its valuation ring.\par
Let $\nu$ be the monomial valuation of $R=k[u_1,u_2]_{(u_1,u_2)}$ determined by $\nu (u_1)=1$, $\nu
(u_2)=\tau$. The group of the valuation is $\Phi =\Z+\Z\tau$ with the order induced by that of $\mathbf{R}$, so that $\hbox{\rm gr}_\nu
R_\nu$ is presented as above, with $A=k$:
$$\hbox{\rm gr}_\nu R_\nu=k[(V_i)_{i\in
\N}]/\big((V_i-V_{i+1}^{s_i}V_{i+2})_{i\in \N}\big).$$ The semigroup of the values of $\nu$ on $R$ is $\Gamma=\{a+b\tau/a,b\in \Z_+\cup\{0\}\}$
and the graded algebra
$\hbox{\rm gr}_\nu R$ is the polynomial algebra $k[U_1,U_2]$, where $U_1$ is of degree one and $U_2$ of degree $\tau$.
\end{example} 
\begin{example}\label{exzar} Recall the examples given in [Z-S]. Vol.2, Chap. VI,
\S 15, p.102: let
$k$ be a field, 
$\{s_1,s_2,\ldots \}$ a sequence of positive integers such that the products $s_1s_2\cdots s_i$
tend to infinity with $i$, and $\{c_1,c_2,\ldots \}$ a sequence of elements of $k^*$. We define an
infinite sequence of elements $v_i$ in $k(u_1,u_2)$ by $$v_1=u_1, v_2=u_2, \ldots ,v_{i+2}=
\frac{v_i-c_iv_{i+1}^{s_i}}{ v_{i+1}^{s_i}},\ldots .$$ Setting $R_i=k[v_{i-1},v_i], \ \
m_i=(v_{i-1},v_i)R_i$, Zariski shows that the ring $$R_\nu=\bigcup_{i=2}^\infty (R_i)_{m_i}$$ is
the valuation ring of a valuation $\nu$ of the field $k(u_1,u_2)$ with value group  $$\Phi=\{\frac{n}{s_1\cdots
s_k}\ ,\ n\in \Z ,\ k=1,2,\ldots\},$$ that is, all rational numbers whose denominator is a product
of $s_i$'s. \par One checks that the elements
$\bigl((v_i)_{i\geq 1}\bigr)$ generate the maximal ideal of $R_\nu$, and their initial forms $(V_i)$ generate the $k$-algebra $\hbox{\rm gr}_\nu
R_\nu$. As $\nu (v_{i+2})>0$, we have
$\nu (v_i-c_iv_{i+1}^{s_i})>\nu (v_{i+1}^{s_i})$, which gives the relations $V_i-c_iV_{i+1}^{s_i}=0$ for $i\geq1$; they
generate all relations between the
$(V_i)$, so that the graded algebra of $R_\nu$ is given by
$$\hbox{\rm gr}_\nu R_\nu=k[V_1,V_2,\ldots ,V_i ,\ldots ]/\big((V_i-c_iV_{i+1}^{s_i})_{i\geq
1}\big).$$ Since the rational rank of
$\Phi$ is one, we expect the $k$-algebra $\hbox{\rm gr}_\nu R_\nu$ to be the direct limit of
polynomial algebras in one variable, and indeed this shows that $\hbox{\rm gr}_\nu R_\nu$ is the direct limit of the system of maps $k[V_i]\to
k[V_{i+1}]$ given by
$V_i\mapsto c_iV_{i+1}^{s_i}$.\par\noindent This corresponds, replacing $k$ by a commutative ring $A$ and taking $c_i=1$ for all $i$s, to a presentation of
the semigroup algebra of the semigroup $\Phi_+$ with coefficients in $A$ as
$$A[t^{\Phi_+}]=A[(V_i)_{i\in\N}]/\big((V_i-V_{i+1}^{s_i})_{i\geq 1}\big).$$ This corresponds also to the fact that $\Phi_+$ is the direct limit of its
subsemigroups generated by $\frac{1}{s_1.\ldots s_i}$, the inclusion between two consecutive semigroups being given by multiplication by $s_{i+1}$.\par 
Note that this valuation is not monomial for $R=k[u_1,u_2]_{(u_1,u_2)}$ with respect to the coordinates
$(u_1,u_2)$.\par\medskip\noindent In fact it is an interesting exercise to begin to compute
$\hbox{\rm gr}_\nu R$ as a subalgebra of $\hbox{\rm gr}_\nu R_\nu$ in this case. \par\noindent Recall that
$R=k[u_1,u_2]_{(u_1,u_2)}$, and assume for simplicity that all $s_i$ are $>1$ for $i\geq 2$. From the equations defining the $v_i$, we see that
if we set $ \nu( u_1)=1,\ \hbox{\rm then}\ \nu (u_2)=\frac{1}{s_1}$ and $u^{s_1}_2-c^{-1}_1u_1=-c^{-1}_1u_2^{s_1}v_3\in R$,
I set
$u_3= -c^{-1}_1u_2^{s_1}v_3$, so we have the following equation defining $u_3$ 
$$u^{s_1}_2-c^{-1}_1u_1=u_3, $$ and we see similarly that
$u_3^{s_2}-c_2^{-1}(-c_1)^{s_2}u_2^{s_1s_2+1}=-c_2^{-1}u_3^{s_2}v_4\in R$. I set
$u_4=-c_2^{-1}u_3^{s_2}v_4$ and so we have the equation
$$u_3^{s_2}-c_2^{-1}(-c_1)^{s_2}u_2^{s_1s_2+1}=u_4$$ which we may use to define $u_4$. If we apply the same method to find $u_5$, we first find the
equation
$$u_4^{s_3}-(-1)^{s_3(s_2-1)}c_1^{s_2s_3}c_2^{-s_3}c_3^{-1}u_3^{s_2s_3}v_3=(-1)^{s_3(s_2-1)}c_1^{s_2s_3}c_2^{-s_3}c_3^{-1}u_3^{s_2s_3}v_4^{s_3}v_5.$$
But the monomial $u_3^{s_2s_3}v_3$ is not in $R$; its valuation, however, is in the semigroup $\Gamma$ of $R$; it is the valuation of
$u_2^{s_3(s_1s_2+1)-s_1}u_3$. Working out the computations gives an equation
$$u_4^{s_3}-d_4u_2^{s_3(s_1s_2+1)-s_1}u_3=*u_4^{s_3}v_5+*u_4v_3\prod_{\omega\in \mu_{s_3}\setminus\{1\}}(u_3^{s_2}-*\omega
u_2^{s_1s_2+1}),$$ where $d_4$ and the asterisks are Laurent monomials in $c_1,c_2,c_3$. We choose the right-hand term as $u_5$, note that its valuation is
that of
$u_4^{s_3}v_5$, and continue in this way. The minimal system of generators of the semigroup $\Gamma$ corresponding to the values of
$\nu$ on $R$ is
$$\Gamma=\langle\frac{1}{s_1}, 1+\frac{1}{s_1s_2},
s_2+\frac{1}{s_1}+\frac{1}{s_1s_2s_3},\ldots,\gamma_i,\ldots\rangle ,$$ where $\gamma_i=\nu (u_{i+1})$ and
$$\gamma_{i+1}=s_i\gamma_i+\frac{1}{s_1s_2\cdots s_{i+1}}.$$
From this definition follows that for all $i\geq 2$ $$s_i\gamma_i\in \langle \gamma_1,\ldots ,\gamma_{i-1}\rangle
.$$ We have seen the result for $i=2$. Proceed now by induction on $i$, assuming the result true until $\gamma_{i-1}$ and writing
$s_{i-1}\gamma_{i-1}=\sum_{k=1}^{i-2}\tilde\ell^{(i-1)}_k\gamma_k$, we have
$$\begin{array}{lr}s_i\gamma_i=(s_is_{i-1}+1)\gamma_{i-1}-\gamma_{i-1}+\frac{1}{s_1\ldots
s_{i-1}}=(s_is_{i-1}+1)\gamma_{i-1}-s_{i-2}\gamma_{i-2}\\
= ((s_i-1)s_{i-1}+1)\gamma_{i-1}+(\tilde\ell^{(i-1)}_{i-2}-s_{i-2})\gamma_{i-2}+\sum_{k=1}^{i-3}\tilde\ell^{(i-1)}_k\gamma_k\end{array},$$ the last
equality being provided by the induction assumption with the extra twist that $\tilde\ell^{(i-1)}_{i-2}>s_{i-2}$.\par\noindent Since we know by
our induction hypothesis that
$(s_i-1)s_{i-1}\gamma_{i-1}\in
\langle
\gamma_1,\ldots ,\gamma_{i-2}\rangle$, we see that in fact
$s_i\gamma_i-\gamma_{i-1}\in \langle \gamma_1,\ldots ,\gamma_{i-2}\rangle$ for $i\geq 3$. So we have relations, which we may write
$$s_i\gamma_i=\gamma_{i-1}+\sum_{k=1}^{i-2}\ell^{(i+1)}_k\gamma_k\ \ \ \hbox{\rm with}\ \ell^{(i+1)}_{k}<s_k\ \hbox{\rm for all}\  k\geq 2.$$ The
$\ell$ are indexed by $i+1$ which corresponds to the variable $U_i$ instead of $i$ in order to give naturality to the equations given below.
\noindent These relations give us
the equations for the monomial curve associated to
$\Gamma$ in the infinite-dimensional space with coordinates
$(U_j)_{j\geq 2}$ (note that we have eliminated $U_1$), and therefore also for its avatar $\hbox{\rm Spec}\hbox{\rm gr}_\nu R$; the latter are
$$U_3^{s_2}-d_3U_2^{s_1s_2+1}=0,\ \hbox{\rm and}
\ \ U_{i+1}^{s_i}-d_{i+1}U_i\prod_{k=1}^{i-2}U_{k+1}^{\ell^{(i+1)}_k}=0\ \ \hbox{\rm for}\  i\geq
3,$$ where the $d_i$ are, up to sign, Laurent monomials in the $c_k$'s.\par\noindent Going back to our ring
$R$, note that in the first few relations defining the $u_i$ inductively another expression appeared naturally for the inclusion
$s_3\gamma_3\in
\langle
\gamma_1,\gamma_2\rangle$. In any case the important fact is the following : the equations defining inductively the
$u_j$ are (with a small change in the indexation) of the form 
$$u_j^{s_{j-1}}-d_ju_1^{\ell^{(j)}_1} \ldots u_{j-1}^{\ell^{(j)}_{j-1}}=u_{j+1}$$ with $d_j\in k^*$, the $d_j$'s
being, up to sign, Laurent monomials in the $c_k$'s, and the term on the right of the equation has greater valuation than the
monomials appearing in the left hand side.\par When we compute
$\hbox{\rm gr}_\nu R$ from these equations, we need to take as coordinates $U_2,U_3, \ldots$ but {\it not} $U_1$, because of the
equation $U^{s_1}_2-c^{-1}_1U_1=0$. All the right-hand sides disappear because they have greater
valuation, and the equations give algebraic relations
between the initial forms $U_j$ of $u_j$ and those of $u_2,\ldots ,u_{j-1}$, which
decreases the transcendence degree over $k$, and therefore the Krull dimension, of $\hbox{\rm gr}_\nu R$ to one in agreement with
 Proposition \ref{pil} (see also the remark which follows it). In the ring $R$ however each equation
is captured, via its right-hand side $u_{j+1}$, in an abyssal\footnote{From the Greek
{\protect\selectlanguage{greek} >'abussoc} \selectlanguage{english} meaning \textit{without
bottom}; in spite of a discrepancy in meaning I have preferred this to \textit{abysmal}, which does not exist in French.}
sequence of relations, and never gets a chance to decrease the dimension of
$R$ by one. If we stopped after a finite number of steps, setting $$u_j^{s_{j-1}}-d_ju_1^{\ell^{(j)}_1} \ldots
u_{j-1}^{\ell^{(j)}_{j-1}}=P(u_1,\ldots ,u_j),$$ where $P$ is any polynomial (for example $0$) which is not a zero divisor modulo the ideal generated by
$u_j^{s_{j-1}}-d_ju_1^{\ell^{(j)}_1}
\ldots u_{j-1}^{\ell^{(j)}_{j-1}}$ and the previous equations, the dimension would
decrease to one.\par\noindent This example displays for us a basic phenomenon: the decrease of
dimension of the special fiber in the specialization of $R$ to $\hbox{\rm gr}_\nu R$
in the case where Abhyankar's inequality is strict. In both rings, all
variables of index $\geq 4$ are algebraically dependent upon $u_2$ and $u_3$, but in 
$\hbox{\rm gr}_\nu R$ these two variables are dependent, while in $R$ they are
independent and all others are polynomials in them, which shows ``from the
equations'' that $R$ is regular of dimension two. In fact, we have given a presentation
$$k[u_2,u_3]=k[u_2,u_3,\ldots,
u_i,\ldots]/(u_3^{s_2}-c_1^{-s_2}c^{-1}_2u_2^{s_1s_2+1}-u_4,u_4^{s_3}-\ldots)$$ of our polynomial ring which may seem
complicated, but is in fact an explicit description of the valuation, allowing us to compute the valuation of every element, as follows:\par\noindent
Given a polynomial $P(u_2,u_3)$, replace every occurrence of $u_3^{s_2}$ by $c_1^{-s_2}c^{-1}_2u_2^{s_1s_2+1}+u_4 $, then every occurrence of
$u_4^{s_3}$ by $d_4u_2^{\ell^{(4)}_2}u_3^{\ell^{(4)}_3}+u_5$, still replacing occurences of $u_3^{s_2}$ as above, and so on. This stops after finitely
many steps because the product
$s_2s_3\cdots s_i$ eventually exceed the degree of the polynomial $P$. We obtain finally a polynomial $\tilde P(u_2,u_3,\ldots ,u_k)$ where no variable
$u_i$ appears with an exponent
$\geq s_{i-1}$ for $i\geq 3$, and now the valuation of
$P$ is the smallest valuation of the monomials of $\tilde P$ simply because relations between their initial forms in the graded ring have become
impossible. We shall see below in subsection \ref{abyssalph} that since
$\nu$ is rational of height one, a similar result is always true at least after a suitable completion of the ring $R$, that is
$$\hat R^{(\nu)}=\widehat{k[u_2,u_3,\ldots,
u_i,\ldots]}/\overline{({\mathbf F},(u_j^{n_j}u^{n(j)}-\lambda_{nm}u^{m(j)}+c_{j+1}u_{j+1}+\cdots)_{j\geq j_0})},$$ where the second term is
a quotient of a suitable completion of the polynomial ring, ${\mathbf F}$ is a finitely generated ideal, the $c_{j+1}$ are in $k^*$ and the bar means a
topological closure. If $\hbox{\rm h}(\nu)>1$, a similar result is true, where the variables $u_i$ are indexed by ordinals.
\par\noindent  Another important phenomenon seen here is that because of the minimality the initial form
$U_1$ of $u_1$ does not appear among the generators of $\hbox{\rm gr}_\nu R$, so that if we wish to find generators of the maximal ideal of $R$
whose initial forms are among a minimal system of generators of $\hbox{\rm gr}_\nu R$, we have to write $R=k[u_2,u_3]_{(u_2,u_3)}$. The
valuation of $u_3$ is maximal among those of elements which, together with $u_2$, generate the center of the valuation. Heuristically,
the process leads us to coordinates which have, with respect to the given valuation, a property analogous to Hironaka's maximal
contact (see [H1] and below).
\par\medskip\noindent  We can
describe the valuation in terms of blowing-ups as follows: starting with  ${\mathbf A}^2(k)$ with
ordered coordinates $(v_1,v_2)$, we perform the following sequence of steps. \par\noindent
Step 1: reverse the order of coordinates.\par\noindent Step 2: blow-up the origin and
localize in the chart where $v_2$, now the first coordinate, generates the ideal
$(v_1,v_2)$, taking as coordinates $(v_2, v_1v_2^{-1})$. \par\noindent Step 3: repeat
this second step $s_1-1$ times.\par\noindent Step 4: substract $c_1$ from the second coordinate so
that the point with coordinate $c_1$ on the axis of the  second coordinate (the
exceptional divisor just created) is the origin of coordinates, and localize at this point.
\par\noindent Repeat from step 1, replacing $(s_1, c_1)$ by $(s_2,c_2)$, and so on.
\par\noindent The valuation ring $R_\nu$ is the direct limit of the direct system of local rings corresponding to the sequence of infinitely near
points obtained in this way. Note that the strict transform of $u_3$ is a coordinate at the point obtained after step $4$; and so on; it is one of the attributes
of maximal contact.\par\noindent According to ([Z-S], Vol. 2, p. 104), any subgroup
$\Phi$ of the additive group $\Qq$ of rational numbers can be obtained in this
manner.\par\noindent This applies in particular to $\Phi=\Qq$, obtained by taking $s_i=i\ 
 \hbox{\rm for all} \ i\geq 1$.\par\medskip Another construction of valuations whose
value groups are given subgroups of $\Qq$ is given by Zariski in [Z3], p.648, which
leads to the distinction between those prime numbers which divide only a finite
number of the $s_i$'s and those which divide an infinite number of them. If we
worked over a field whose characteristic belongs to the second group and tried to
effectively build inductively elements like the $v_i$ by solving equations of the type
$c_iv_{i+1}^{s_i}-v_i+\cdots =0$, we would run into substantial difficulties. However,
this is not at all what I do here; I only use the direct system as a means to
approximate $\hbox{\rm gr}_\nu R_\nu$ by polynomial subalgebras.  Note anyway that in these examples, $\hbox{\rm gr}_\nu R$ is not regular.
\end{example}
\begin{example}\label{nextex} One may ask what happens if, in the sequence of steps just described, one makes no
translation; it amounts to setting all the $c_i$ equal to zero. Since the equations
describing the formation of the $v_i$ are  $$v_i=v_{i+1}^{s_i}(c_i+v_{i+2}),$$ we see
that the graded algebra of our new valuation $\nu_0$ is $$\hbox{\rm
gr}_{\nu_0} R_{\nu_0}=k[V_1,V_2,\ldots ,V_i ,\ldots
]/\big((V_i-V_{i+1}^{s_i}V_{i+2})_{i\in \N}\big).$$ It follows from what we saw in
example \ref{irrtau} that $\nu_0$ is the monomial valuation of height one and rational rank
two defined by $\nu_0 (u)=1, \nu_0 (v)=\tau$, where $\tau$ is the number with
continued fraction expansion $$\tau =[s_1,s_2,\ldots ,s_i,s_{i+1},\ldots ].$$ Its value
group is  $$\Phi_\tau=\Z+\Z \tau\ \subset {\mathbf R}$$ with the order induced by that of
${\mathbf R}$. We have $\hbox{\rm gr}_{\nu_0} R=k[U_1,U_2]$, with $U_1$ of degree 1 and $U_2$ of degree $\tau$; it is regular. \par\noindent In
particular, the valuation with value group
$\Qq$ constructed above becomes, upon setting all the $c_i$ equal to zero, the monomial valuation
corresponding to $\nu_0 (u)=1,\ \nu_0 (v)=\rho$, where $$\rho =[1,2,3,\ldots ,i, i+1,\ldots ].$$\end{example}
\par\medskip\noindent
{\textbf {Problem.}} It would be interesting if, given
any rational valuation $\nu$ of a n\oe therian local integral domain $R$, one could
specialize it (in the sense used just above, which differs from the usual one) to a
valuation of $R$ of rational rank equal to the dimension of $R$. Heuristically this
means we can always ``suppress the translations'' continuously. As we shall
see in Section \ref{complepro} it is reasonable to try first for a valuation $\nu_0$ of
rational rank equal to the dimension of the scalewise $\nu$-adic completion of $R$.
Can one use the results of [A2] to verify this statement in dimension
two?\par\medskip\noindent \begin{example}
In [Z3], Zariski gives, in the case where the characteristic of $k$ is zero, yet another construction, also found in 
[McL-S], of valuations on $k[u_1,u_2]_{(u_1,u_2)}$ with a value group $\Phi\subset \Qq$ given in
advance. Let $k[[u^{\Qq_+}]]$ be the ring of power series with exponents forming a well
ordered subset of $\Qq_+$, and with the same $s_i$ as above, given a sequence of integers
$m_i$ and a sequence of elements $e_i\in k^*$, consider the series 
$$w(u)=e_1u^{\frac{m_1}{s_1}}+e_2u^{\frac{m_2}{s_1s_2}}+\cdots+e_ju^{\frac{m_j}{s_1\ldots s_j}}+\cdots \ .$$ Since $k$ is of characteristic zero, the series
$w(u)$ cannot be algebraic over $k(u)$ and the map
$ u_1\mapsto u,\ u_2\mapsto w(u)$ induces an injection  $$k[u_1,u_2]_{(u_1,u_2)}\subset k[[u^{\Qq_+}]]$$ and the $u$-adic valuation
induces on $k(u_1,u_2)$ a valuation having as group of values the group of rational numbers
whose denominator is a product of $s_i$'s.\par\noindent The correspondence between this construction and the
previous one is as follows. \par\noindent In our set of equations $u_j^{s_{j-1}}-d_ju_1^{\ell^{(j)}_1} \ldots
u_{j-1}^{\ell^{(j)}_{j-1}}=u_{j+1}$ for $k[[u_1,u_2]]$, stop at $j=\ell$ and set $u_{\ell+1}=0$, then using the previous equations, eliminate the
variables $u_i,\ i>2$ to finally get an equation $Q_\ell(u_1,u_2)=0$ in $u_2$ and $u_1$ which has as solution, if the characteristic of $k$ is zero, a Puiseux
series $u_2=w_\ell(u_1)$ whose exponents have denominators
$s_1\ldots s_{\ell-1}$. As we let $\ell$ go to infinity, the series $w_\ell(u_1)$ converge to a series of the form $w(u_1)$ in
$k[[u^{\Qq_+}]]$. This is a generalization of the procedure given in [T1] to compute the equation of a branch having a given
semigroup. \par What I have just done is to describe the transcendental branch defined by $u_1= u,\ u_2= w(u)$ as a plane
deformation of the avatar corresponding to the constants $d_j$ of the monomial curve associated to the semigroup $\Gamma$ described above,
which lives in an infinite-dimensional space. One can see, with a little work, that the $Q_\ell(u_1,u_2)$, viewed as polynomials in $u_2$, are ``key
polynomials'' in the sense of [McL2] and [Ka] for the valuation $\nu$ viewed as extending to $k(u_1)(u_2)$ the $(u_1)$-adic valuation of
$k(u_1)$; the irreducibility of the $Q_\ell$ in the sense of \textit{loc. cit.} follows from their construction and the results of [T1] on deformations of
monomial curves, and the minimality of their degrees follows from the easily verified fact that for each index $i$, the smallest non zero integer 
$d$ such that $d\gamma_i\in\langle\gamma_1,\ldots ,\gamma_{i-1}\rangle$ is $s_i$. This provides a geometric
interpretation for these polynomials and also of the construction of
\S 5 of [McL-S] which builds, in arbitrary characteristic, a rational valuation of $k[u_1,u_2]_{(u_1,u_2)}$ with a given group of values contained in ${\mathbf
Q}$, and for some of the constructions of [F-J] and [V2].
\par\noindent Note that, as $\ell$ goes to infinity, the parametric representation $u_2=w_\ell (u_1)$ of our approximating plane branches tend to a
transcendental series, while their equations $Q_\ell$ tend to zero in the $(u_1,u_2)$-adic topology; this last fact explains why the algebraic description,
which is blind to transcendental series, gives us the whole (formal germ of the) affine plane as a limit.\end{example}\par\noindent
\begin{example}A power series in one variable with rational exponents having unbounded denominators is not necessarily very transcendental in
positive characteristic:\par\noindent
Let $k$ be a perfect field of finite characteristic $p$ and $u$ an indeterminate; set $K=\bigcup_{n\geq 1}k(u^{\frac{1}{p^n}})$, the perfect closure of $k(u)$.
There is a unique extension to $K$ of the $u$-adic valuation of $k(u)$, and its valuation ring is \textit{not} n\oe therian. Consider the series
$$v=\sum_{i=1}^{\infty}u^{1-\frac{1}{p^i}}\in k[[u^{{\mathbf Q}_+}]];$$
it is a solution of the polynomial equation
$$v^p-u^{p-1}(1+v)=0.$$ This equation is an Artin-Schreier equation: it is obtained from the standard Artin-Schreier
$y^p-y=\frac{1}{u}$ by replacing $y$ by $\frac{v}{u}$. If we set $L=K(v)$, it is shown in [K3] that the extension $L/K$ has degree $p$ and
defect $p$. More precisely, the unique extension $\nu$ to $K$ of the $u$-adic valuation of $k(u)$ has a unique extension $\nu'$ to $L$, with the same
group of values, so that the ramification index
$e=[\Phi':\Phi]$ is equal to one, and no residual extension so that the inertia degree $f=[\kappa(\nu'):\kappa (\nu)]$ is also equal to one. The extension is
of degree $p$ so that the Ostrowski ramification formula (see [K3], [Roq]), which is $[L:K]=def$, where $d$ is the defect, gives $d=p$. This defect
complicates the parametrization but does not make it more difficult to create a non-singular model. We remark that our curve is a deformation of the
monomial curve
$v^p-u^{p-1}=0$, and apply to this monomial curve the toric resolution process of [G-T] and subsection \ref{toricmod} below: it gives us a chart
$u=y_1^py_2,v=y_1^{p-1}y_2$. Our equation then becomes $y_1^{p(p-1)}y_2^{p-1}(y_2-1-y_1^{p-1}y_2)$, so that the strict transform
$y_2-1-y_1^{p-1}y_2=0$ is non singular. It can be parametrized in a  neighborhood of the exceptional divisor $y_1y_2=0$ by
$y_2=\frac{1}{1-y_1^{p-1}}$, so that we have the following parametrization of our curve in a neighborhood of the origin:
$$u=\frac{y_1^p}{1-y_1^{p-1}};\ \ v=\frac{y_1^{p-1}}{1-y_1^{p-1}}.$$
The fact that the extension $K\subset K(v)$ has defect seems to be related to the fact that while the extension of fields $k(u)\to k(u)[v]/(v^p-u^{p-1}(1+v))$
is separable, the extension of graded rings associated to the $u$-adic valuation of $k[u]$ and its extension to $k[u,v]/(v^p-u^{p-1}(1+v))$,
which is
$$k[U]\to k[U,V]/(V^p-U^{p-1})$$ is purely inseparable of degree $p$.
\end{example}
\section{\textbf{Completion problems}}\label{complepro} 
In this section I study the $\nu$-adic
completion of a n\oe therian local ring $R$ and show that it is n\oe therian in certain cases.
 Although completions of valued fields are a classical topic, the algebraic structure of the  $\nu$-adic completions of n\oe
therian rings does not seem to have been studied. Perhaps the reason is that the n\oe therian hypothesis is not natural from the viewpoint of
the theory of valued fields, while completions other than those with respect to an ideal are not
commonplace in algebraic geometry. This is a meeting place for the number-theoretic/henselian and the algebro-geometric traditions of
valuation theory.\par The general fact is that the completion of
$R$ for the
$\nu$-adic topology is a quotient of the completion of $R$ for the topology of the symbolic powers of the center of the height one valuation with
which $\nu$ is composed.  If
$R$ is analytically irreducible this topology is finer than the $m$-adic topology.\par Being complete for the $\nu$-adic topology, however, does
not suffice to imply the convergence of the sequences which occur when one tries to lift elements from the graded ring to $R$. In order to
overcome this difficulty I consider, in the case  of a rational valuation, a quotient
$\hat R^{(\nu )}$ of the $m$-adic completion $\hat R^m$, which is complete for a valuation $\hat \nu$
extending $\nu$ and such that the inclusion $R\to\hat R^{(\nu )}$ induces a scalewise birational morphism $\hbox{\rm gr}_\nu R\to\hbox{\rm
gr}_{\hat\nu} \hat R^{(\nu )}$. This ring $\hat R^{(\nu )}$ is not only complete for
the $\nu$-adic topology but also \textit{scalewise} complete, which ensures the convergence of the sequences mentioned above. It is also henselian,
and in the equicharacteristic case admits a field of representatives.
\par I also show how to specialize the ring
$\hat R^{(\nu )}$ to a suitable (scalewise) $\hat\nu_{\rm gr}$-adic completion of $\hbox{\rm gr}_{\hat\nu} \hat R^{(\nu )}$ by a
suitable completion of the valuation algebra, and the existence of a coordinate system and
equations for the ring of this specialization, which is a fundamental fact in this approach. 
\subsection{The $\nu$-adic completion of $R$}\label{nuadiccomp}
 In general, even if the local ring $R$ has a field of representatives, it is not possible to lift generators of the graded algebra
$\hbox{\rm gr}_\nu R$ to generators of $R$; we need a hypothesis of completeness (see
[B3], Chap. III). The arguments using the fact that a deformation preserves smoothness and
transversality also use some form of the implicit function theorem, or Hensel's Lemma, the
validity of which is another attribute of complete rings. This leads us to the study of the
completion of a ring $R$ with respect to the $(\P_\phi)_{\phi\in \Phi_+}$-filtration. We
note that the quotients $R/\P_\phi$ form a countable projective system since the
semigroup $\Phi_+$ is countable and totally ordered. One can then define the $\nu$-adic
completion of $R$ as $$\hat R^\nu=\displaystyle\limproj_{\phi\in \Phi_+} R/\P_\phi .$$
There is a natural map $R\rightarrow \hat R^\nu$, which is injective since the filtration is
separated. The ring $\hat R^\nu$ is an integral domain and the valuation $\nu$ has a canonical extension $\hat \nu$ to a valuation of $\hat
R^\nu$ with values in $\Phi$. To see this it is enough to apply directly the explicit form of
the definition of the $\nu$-adic completion as a subset of the product of the quotients
$R/\P_\phi$. We note that $\hat R^\nu$ is the closure of $R$ in the completion $\hat K^\nu$ of
its field of fractions with respect to the topology ${\mathcal T}_\nu$ defined by $\nu$ (see [B3],
Chap VI, \S 5, No. 3).
 In particular, we have:
$$R_{\hat \nu}=\hat R^\nu_\nu.$$  One sees also immediately that if we set 
$$\hat \P_\phi=\{x\in \hat R^\nu\mid\hat \nu (x)\geq \phi\},\ \ \hat \P^+_\phi=\{x\in \hat
R^\nu\mid\hat \nu (x)>
\phi \},$$ they are respectively the closure of $\P_\phi ,\ \P_\phi^+$ in $\hat R^\nu$ ([B3], Chap.
VI, \S 5 No. 2)  and the natural inclusion $\P_\phi/\P^+_\phi\subset \hat\P_\phi/\hat\P^+_\phi$
is an equality, so that the inclusion $R\subset \hat R^\nu$ induces an equality of graded rings
$$\hbox{\rm gr}_\nu R=\hbox{\rm gr}_{\hat \nu}\hat R^\nu .$$
 We note that the semigroup of values $\hat \nu (\hat R^\nu\setminus \{0\})\subset
\Phi_+\cup \{ 0\}$ is the same as that of $\nu$ on $R$, and so is well ordered when $R$ is
n\oe therian (Proposition \ref{wellord}). \par\medskip\noindent
\par\medskip
\par Let us assume that $R$ is a n\oe therian local ring with maximal ideal $m$ and
$\nu$ is a valuation of $K$ of height one such that its valuation ring $R_\nu$ contains $R$. Set
${\mathbf p}=m_\nu\cap R$; since $R$ is n\oe therian, we have $\nu ({\mathbf p})>0$, say $\nu ({\mathbf
p})=\phi_0$, and since $\nu$ is of height one, its group
$\Phi$ is archimedian, so that for any $\phi\in \Phi$ there exists an integer $N(\phi )$ such that
$N(\phi )\phi_0\geq \phi$, that is ${\mathbf p}^{N(\phi )}\subset \P_\phi$. Then we have a
projective system of surjective maps
$$R/{\mathbf p}^{N(\phi )}\rightarrow R/\P_\phi ,$$ which induces a continuous morphism of
completions
$\hat R^{\mathbf p}\rightarrow \hat R^\nu$, where $\hat R^{\mathbf p}$ is the completion of $R$ for the
${\mathbf p}$-adic topology. Since $\hat R^\nu$ is an integral domain, the kernel of this map is a
prime ideal of $\hat R^{\mathbf p}$.  By the definition of completions, this kernel is
$\bigcap_{\phi\in\Phi_+}\P_\phi\hat R^{\mathbf p}$, so that in the case where $\nu$ is of height one
and $R_\nu$ dominates $R$, i.e., ${\mathbf p}=m$, it coincides with the ideal
 called the implicit ideal of $\hat R$ by Spivakovsky ([S2]). \par\noindent
Recall that the symbolic power ${\mathbf p}^{(n)}={\mathbf p}^nR_{\mathbf p}\cap R$ is the set of elements $x\in R$ such that there
exists
$s\notin {\mathbf p}$ such that $sx\in {\mathbf p}^n$. Because our ring $R$ is a n\oe therian domain, the symbolic powers define a separated topology on
$R$ (see [Z-S], Chap. IV, Th. 23), for which there is a completion
$\hat R^{({\mathbf p})}$. Note that, since when $s\notin {\mathbf p}$ we have $\nu (sx)=\nu (x)$, with the notations just introduced we
have
${\mathbf p}^{(N(\phi ))}\subset
\P_\phi$, so that we have natural continuous maps 
$$\hat R^{\mathbf p}\to \hat R^{({\mathbf p})}\to \hat R^\nu.$$ Let us denote by $H$ the kernel of the map
$\hat R^{({\mathbf p})}\to \hat R^\nu$ and remark that $H\cap R=(0)$.  We have the following
\begin{proposition}\label{complht1} For a valuation $\nu$ of
height one on a n\oe \-therian local ring $R$, with the notations introduced above, the natural continuous injection
\footnote{This
corrects an error in the statement of Proposition 1.3 of [T2].}
$$\hat R^{({\mathbf p})}/H\to\hat R^\nu$$\noindent is an isomorphism of topological rings inducing the identity on $R$. If the center of
$\nu$ in $R$ is the maximal ideal $m$, the natural map $\hat R^m\to \hat R^{(m)}$  is an isomorphism, so that $\hat R^\nu $ is a quotient of  $\hat R^m$
and is n\oe therian.\par\noindent If
$R$ is a regular local ring, the map $\hat R^{\mathbf p}\to \hat R^{({\mathbf p})}$ is an isomorphism for all prime ideals,
$\hat R^\nu$ is a quotient of $\hat R^{\mathbf p}$ and therefore is n\oe therian.
\end{proposition} 
\begin{proof} The injection
$$\hat R^{({\mathbf p})}/H\subset \hat R^\nu $$ implies that the valuation $\nu$ extends to
$\hat R^{({\mathbf p})}/H$ as a valuation $\hat \nu$ of height one.\par\noindent The valuation $\nu$ is
still $\geq 0$ on the localization $R_{\mathbf p}$, and so we have an injection
$\hat R^{{\mathbf p}R_{\mathbf p}}_{\mathbf p}/H_{{\mathbf p}}\subset \hat R_{\mathbf p}^\nu$ corresponding to $R_{\mathbf p}$.
Let us replace
$R$ by $R_{\mathbf p}$ for a moment, so that ${\mathbf p}$ is the maximal ideal $m$. By [Z-S], Vol. 2,
Appendix 3, Lemma 3 p. 343, the distinct valuation ideals (i.e., our
$\P_\phi (\hat R^m/H)$ written without repetition) form a simple infinite descending chain of
ideals, which are primary for the maximal ideal $\hat m$ of
$\hat R^{m}/H$ by the archimedian property, and have intersection zero. I denote them by $\hat\P_j$. Now we can apply
Chevalley's Theorem ([Z-S], Vol. 2, Chap. VIII, \S 5, Th. 13, p. 270), which asserts that there exists
an integer-valued function
$s(n)$ tending to infinity with $n$ and such that for each valuation ideal $\hat\P_j$, we have
$\hat\P_j\subset
\hat m^{s(j)}$. This, added to the fact that the $\hat\P_j$ are primary for
$\hat m$, proves that in the ring
$\hat R^{m}/H$ the $\hat\nu$-adic topology coincides with the
$\hat m$-adic topology, so that it is complete for both, and therefore has to be equal to
$\hat R^\nu$. \par Now coming back to the case where the center of $\nu$ is not necessarily
$m$, we can remark that by definition of the symbolic powers the inclusion $R\to R_{\mathbf p}$ extends to an inclusion $$\hat
R^{({\mathbf p})}\to
\hat R_{\mathbf p}^{{\mathbf p}R_{\mathbf p}}.$$ In fact $\hat R^{({\mathbf p})}$ is the closure of $R$ in $\hat R_{\mathbf p}^{{\mathbf
p}R_{\mathbf p}}$. The ${\mathbf p}^{(n)}$-adic topology  of the first ring is
induced by the ${\mathbf p}^n$-adic topology of the second since $ {\mathbf p}^{(n)}\hat
R^{({\mathbf p})}={\mathbf p}^n\hat R_{\mathbf
p}^{{\mathbf p}R_{\mathbf p}}\cap\hat R^{({\mathbf p})}$ as one can check with Cauchy sequences. The valuation
$\nu$ extends to
$R_{\mathbf p}$, and so by the special case where the center is the maximal ideal, the $\nu$-adic completion of $R_{\mathbf p}$ is
$\hat R_{\mathbf p}^{{\mathbf p} R_{\mathbf p}}/H_{\mathbf p}$ where $H_{\mathbf p}=\bigcap_{\phi\in \Phi_+}
\P_\phi(R_{\mathbf p})\hat R_{\mathbf p}^{{\mathbf p}R_{\mathbf p}}$. Setting
$H=H_{\mathbf p}\cap
\hat R^{({\mathbf p})}$, we have an injection $$\hat R^{({\mathbf p})}/H\to \hat R_{\mathbf p}^{{\mathbf p}R_{\mathbf
p}}/H_{\mathbf p}.$$  The quotient topology on $\hat R^{({\mathbf p})}/H$ is the ${\mathbf p}^{(n)}$-adic topology and is induced by
the quotient topology of
$\hat R_{\mathbf p}^{{\mathbf p}R_{\mathbf p}}$; this last one coincides with the $\hat\nu$-adic topology by the special case, so that in the
end the $\hat\nu$-adic topology and the ${\mathbf p}^{(n)}$-adic topology coincide on $\hat R^{({\mathbf p})}/H$. From this follows that
this quotient is equal to $\hat R^\nu$. Since $\P_\phi (R)$ is finitely generated, its closure in $\hat R^{({\mathbf p})}/H$ is $\P_\phi (R)\hat
R^{({\mathbf p})}/H$ and we have 
$H=\bigcap_{\phi\in\Phi_+}\P_\phi (R)\hat R^{({\mathbf p})}$. This gives the first result in the general case.
\par\noindent We have already seen that
$m^{(n)}=m^n$ and finally, in a regular ring the adic and symbolic topologies coincide for all ideals ([Ve], Th. 3.5). \end{proof}
\begin{Remark} 1) This
shows that the ring $\hat R^m/H$ to which Spivakovsky ([S2]) extends the valuation $\nu$ in the
special case of a valuation of height one such that $R_\nu$ dominates $R$ is nothing but
$\hat R^\nu$. \par\medskip\noindent 2) In particular if $R$ is a regular local ring, since its
completion $\hat R^{\mathbf p}$ for the ${\mathbf p}$-adic topology is a  n\oe therian ring ([B3], Chap.
III, \S 3, No.4, Prop.8), for every valuation of height one 
$\nu$ of the field of fractions of $R$ whose valuation ring contains $R$,  the completion $\hat
R^\nu$ of $R$ for the $\nu$-adic valuation is n\oe therian; as a topological ring it is in fact a
Zariski ring with respect to the center of the valuation $\hat \nu$.\end{Remark}
\begin{proposition}\label{oneequicompletion} Let $R$ be a n\oe therian local ring contained in a valuation ring $R_\nu$ of its
field of fractions. The completion $\hat R^\nu$ is isomorphic to
the completion $\hat R^{\nu_1}$ of $R$ with respect to the height one valuation $\nu_1$ with
which the valuation $\nu$ is composed, and is a quotient of the symbolic ${\mathbf p}$-adic completion of
$R$, where ${\mathbf p}$ is the center of $\nu_1$ in $R$.  \end{proposition} 
\begin{proof} We saw
above that is true for valuations of height one. Because of Abhyankar's inequality, the
valuation $\nu$ is of finite height. Let us assume that $ h (\nu)>1$.  By assumption we have a
surjective monotone map $\lambda\colon \Phi\rightarrow \Phi_1$ with $\Phi_1$ of height
one and $\hbox{\rm Ker}\lambda$ of height $h-1$, a maximal isolated subgroup of
$\Phi$.\par Denote by $\nu_1$ the valuation $\lambda\circ\nu$. Then we have the inclusion
$\P_\phi\subset\P_{\lambda (\phi )}$. Whenever $\phi_1>\lambda (\phi)$, we have
$\P_{\phi_1}\subset\P_\phi\subset \P_{\lambda (\phi)}$, which shows that the $\nu$-adic
topology on $R$ coincides with the $\nu_1$-adic topology, hence the equality of the
completions.                                                                                                                                                                                                                    
We are thus reduced to the case of height one treated above, and we see that the completion
$\hat R^\nu$ is a quotient of the completion of $R$ with  respect to the symbolic $(m_{\nu_1}\cap
R)$-adic filtration.\end{proof} 
\begin{Remark} 1) In a n\oe therian domain $R$ with a non negative valuation $\nu$, the ideal
of elements of $R$ which are topologically nilpotent with respect to the $\nu$-adic topology is
the center in $R$ of the height one valuation with which $\nu$ is
composed.\par\medskip\noindent
 2) The equality of the topologies determined on a ring by two valuations, one of
which is composed with the other, is well known; in the n\oe therian case, a more precise result
can be found in [Z-S], Vol. 2, Lemma 4 of Appendix 3, p. 344.\end{Remark}
 \begin{corollary}\label{equicompletion} Given a n\oe therian local integral domain
$R$  and a valuation $\nu$ of its field of fractions, non negative on $R$, if $\nu$ is composed
with a valuation $\nu_1$, the natural map of completions
$$\hat R^\nu\longrightarrow \hat R^{\nu_1}$$ induced by the inclusions $\P_\phi\subset
\P_{\lambda (\phi )}$ is an isomorphism of topological rings. In particular, the
valuation $\nu_1$ extends to a valuation $\hat\nu_1$ on $\hat R^\nu$.
\end{corollary}
\begin{corollary} If $R$ is a regular excellent local ring, for any valuation $\nu$
of its field of fractions which is non negative on $R$, the ring $\hat R^\nu$ is
excellent.
\end{corollary} 
\begin{proof} By Proposition
\ref{oneequicompletion}, and the second part of  Proposition \ref{complht1}, the $\nu$-adic completion $\hat R^\nu$ is a quotient of the
completion $\hat R^{\mathbf p}$ of the local ring $R$ with respect to an ideal, so it is excellent as a
quotient of an excellent ring (see [Ro1], [Ro2]).
\end{proof}
When trying to understand how to complete some non-n\oe therian rings, essentially
semigroup algebras, one encounters questions such as this:
\par\medskip\noindent
\textbf {Question.} (Exercise; Compare with [Ha], [Ka]) Let $A$ be a
commutative ring and $\Phi_+$ (resp. $\Gamma$) the positive semigroup of a totally ordered archimedian abelian group of finite rational rank (resp. the
value semigroup of a valuation of height one on a n\oe therian ring
$R$). When is it true that the completion of the semigroup algebra $A[t^{\Phi_+}]$ (resp. $A[t^\Gamma]$) with respect to the
topology deduced from the canonical valuation $\nu (\sum a_\phi t^\phi )=\hbox{\rm
min}\{\phi | a_\phi \neq 0\}$ is the ring of power series $A[[t^{\Phi_+}]]$ (resp. $A[[t^\Gamma]]$) in the sense of
[McL-S]? \par\noindent
Remark that
$A[[t^\Gamma]]$, which a priori is only a group, is a ring, since in $\Gamma$ there are only
finitely many ways to write an element as a sum of other elements (because
$\Gamma$ is well ordered by Proposition \ref{wellord}; see [Ka], [B2], [C-G], Theorem 1),
so that we can multiply two series. 
\par\medskip\noindent \begin{example}\label{excomp} 1) Take the usual example of
a rank one valuation not satisfying Abhyankar's equality (see above and [Z-S], Chap
VI, \S 15, Example 2, [S1], [S2]): let $k$ be a field, and take $R=k[u_1,u_2]_{(u_1,u_2)}$ so that
$\hat R=k[[u_1,u_2]]$. Let $\sum_{i\geq 1} a_iu_1^i \in k[[u_1]]$ be transcendental over $k(u_1)$ and set $f(u_1,u_2)=u_2-\sum_{i\geq 1} a_iu_1^i \in
\hat R$. Define a valuation $\nu$ on $R$ by means of the injective map $R\rightarrow k[[t]]$ determined by $u_1\mapsto t\ , u_2\mapsto \sum
a_it^i$. The residual transcendence degree of $\nu$ is $0$ and its rational rank is one. The
completion $\hat R^\nu$ is isomorphic to the quotient $k[[u_1,u_2]]/(f)$; it is
one-dimensional. \par\noindent This shows that the injective map $R\to \hat R^\nu$
is not flat in general. The group of the valuation is $\Z$ and for each $n\in \N\cup\{0\}$ we have $\P_n(R)=(t^n)k[[t]]\cap
R$. We see that for each $n\geq 1$ we have $u_2-\sum_1^{n-1}a_iu_1^i\in \P_n(R)$, which shows that there is no inclusion
$\P_n(R)\subset m^{s(n)}$ with $s(n)>1$. Note that we could have taken any irreducible
element
$f\in (u_1,u_2)k[[u_1,u_2]]$ not having a root algebraic over $k(u_1)$ and the composed map
$k[u_1,u_2]_{(u_1,u_2)}\to k[[u_1,u_2]]/(f)\to k[[t]]$ where the last map is the normalization.
\par\medskip\noindent 2) On the same ring $R$, consider the valuation $\nu$ with
value group $\Z^2$ ordered lexicographically, which attributes value $(1,0)$ to $u_1$,
value $(0,1)$ to $u_2$, and such that the valuation of a polynomial is the infimum
of the valuations of its monomials. There is a natural monotone surjective map
$\lambda=\hbox{\rm pr}_1\colon \Z^2\rightarrow \Z$, hence a valuation $\nu_1$, by
the $u_1$-adic order. The ideal $\P_{(i,j)}$ is the product ideal $u_1^i(u_1,u_2^j)$ of $R$, while
the ideal $\P_i$ for $\nu_1$ is $(u_1^i)$. Since $(u_1^{i+1})\subset u_1^i(u_1,u_2^j)\subset (u_1^i)$,
these two filtrations define the same topology on $R$, so that in this case, $$\hat
R^\nu=\hat R^{\nu_1}.$$\noindent Remark that the residual ideal
$(\P_{(i,j)}:\P_i)=(u_1,u_2^j)$ is a simple $\nu$-ideal in $R$ in the sense of Zariski's theory
(see [Z-S], Vol.2, Appendix 3). Notice however that, if we look at the primary decomposition as in \textit{loc.cit.}, we have $\P_{(i,j)}=(u_1^{i+1},u_2^j)\cap
(u_1^i)$.\par\noindent In this case, the completions of
$k[u_1,u_2]_{(u_1,u_2)}$ with respect to $\nu$ and $\nu_1$ are isomorphic to 
$$\widehat{k[u_1,u_2]_{(u_1,u_2)}}^\nu=\limproj_{i\geq 0}k[u_1,u_2]_{(u_1,u_2)}/(u_1)^i=k[u_2]_{(u_2)}[[u_1]];$$ the ideal ${\mathbf
p}$ is $(u_1)k[u_1,u_2]_{(u_1,u_2)}$ and the ideal $H$ is zero.\par\noindent
Remark that the the quotient $\hat R^\nu/{\mathbf p}\hat R^\nu=k[u_2]_{(u_2)}$ is not
complete for the $\overline \nu$-adic (in this case the $(u_2)$-adic)
topology.\par\noindent 3) (a non-discrete case) Given a subgroup $\Phi\subset {\mathbf R}$ ordered by the order of ${\mathbf R}$,
 consider, following H. Hahn ([Ha]), the ring $k[[t^{\Phi_+}]]$ of formal power series with coefficients in $k$ and exponents in
$\Phi_+$, such that the set of their exponents is a well ordered subset of $\Phi_+$ (see also [Ka], [B2]). This
ring is naturally equipped with the $t$-adic valuation, with values in $\Phi_+$.
 Choose algebraically independent series without constant term $u_i(t)\in k[[t^{\Phi_+}]] ,\ 1\leq i \leq r$ (see [McL-S]) and
let $R$ be the localization of the sub $k$-algebra of $ k[[t^{\Phi_+}]]$ generated by $u_1(t),\ldots
,u_r(t)$ at the maximal ideal generated by these elements. Since $R$ is the image of an injective
morphism $\iota\colon k[u_1, \ldots ,u_r]_{(u_1, \ldots ,u_r)}\rightarrow  k[[t^{\Phi_+}]]$ determined by
$u_i\mapsto u_i(t)$,  the $t$-adic valuation of $k[[t^{\Phi_+}]]$
induces a rational valuation of height one on $k[u_1, \ldots ,u_r]_{(u_1, \ldots ,u_r)}$. The morphism $\iota$ extends, by substitution in formal
power series, to a map\break $k[[u_1, \ldots ,u_r]]\rightarrow  k[[t^{\Phi_+}]]$, with a
certain kernel $H$, which describes the formal relations between the $u_i(t)$. By
Proposition \ref{complht1}, the completion $\hat R^\nu$ is the image of this map.
\end{example} 
\subsection{Scalewise completions of $R$}\label{scac}
The gist of this subsection is the construction of $\nu$-adic completions of $R$ such that their quotients by the centers of the
valuations with which $\hat\nu$ is composed are complete for the image of the $\hat\nu$-adic topology. The natural choice is a quotient
of the $m$-adic completion of $R$.\par\noindent The connexion with the previous subsection is provided by the following:
\begin{proposition}\label{zarsymb} {\rm (Zariski, [Z5], p.29)} If the n\oe therian local ring $R$ is analytically irreducible, for any prime ideal
${\mathbf p}$ high symbolic powers of ${\mathbf p}$ are contained in high powers of the maximal ideal $m$.\qed\end{proposition}
\begin{corollary}\label{finer} If the n\oe therian local ring $R$ is analytically irreducible, for any non-negative valuation $\nu$, the $\nu$-adic topology of
$R$ is finer than its $m$-adic topology.\end{corollary}
what I shall really use is the following
\begin {proposition}\label{nucomplete} If the n\oe therian local integral domain $R$ is complete for the $m$-adic topology, for any non
negative valuation $\nu$ it is complete for the 
$\nu$-adic topology as well. \end{proposition}
\begin{proof} One can use Corollary \ref{finer} but also argue directly as follows: by Proposition \ref{oneequicompletion} we may assume that $\nu$ is
of height one, and then the distinct valuation ideals $\P_j$ form a simple sequence having intersection zero; since $R$ is complete, by Chevalley's theorem we
have
$\hat\P_j\subset m^{s(j)}$ with $s(j)$ tending to infinity with $j$. Let $(x_n)_{n\in\N}$ be a Cauchy sequence in $R$ for the 
$\nu$-adic topology. By the preceding observation, it is also a Cauchy sequence for the
$m$-adic topology, and so converges, since $R$ is complete, to some $x\in R$. Let us show that it converges to $x$ for the
$\nu$-adic topology as well. Say that we have $x-x_n\in m^{s(n)}$, with $s(n)$ tending to infinity with $n$. Since $(x_n)$ is a Cauchy
sequence for the $\nu$-adic topology, for
$j\geq n$ we have
$x_j-x_n\in \P_{t(n)}$, with $t(n)$ tending to infinity with $n$. Since $x-x_j\in m^{s(j)}$, for all $j\geq n$ we
have $x-x_n\in \P_{t(n)}+m^{s(j)}$. But $\P_{t(n)}$ is closed for the $m$-adic topology since $R$ is n\oe therian,
so letting $j$ increase gives $x-x_n\in \P_{t(n)}$, which shows that $(x_n)$ converges
to $x$ in the $\nu$-adic topology.\end{proof}
 Since the property of being $m$-adically complete is preserved under
taking quotients, it is reasonable to try to extend the valuation
$\nu$ to an $m$-adically complete ring; it will then be complete for $\nu$, and is quotients as well. The case of rational valuations of height
one suggests to try to extend $\nu$ to a quotient of $\hat R^m$.
\begin{definition}
 Let $R$ be a n\oe therian local domain. Given a family $(\xi_i)_{i\in I}$ of elements  of
the maximal ideal $m$ of $R$, we say that a birational map $R\to R'$ to a local ring $R'$ is a \textit{birational toric extension} of $R$
associated to the family if there exists a finite subset $F\subset I$ and monomials $(\xi^{b^j})_{j\in F}$ with $b^j\in \Z^F$ such that $R'$
is a localization of
$R[(\xi^{b^j})_{j\in F}]$ at a maximal ideal.
\end{definition}
\begin{Remark}
1) In this text the elements $\xi_i$ will be representatives in $R$ of a system of generators $\overline\xi_i$ of the $k$-algebra $\hbox{\rm gr}_\nu R$ for a
rational valuation $\nu$ of $R$. We will consider only toric extensions of the following form: there exists a regular rational cone $\sigma=\langle a^1,
\ldots ,a^N\rangle$ in ${\mathbf R}_+^N$ with $N=\#F$, such that if we denote by $M$ the matrix with columns $(a^1,\ldots ,a^N)$, the
matrix with columns
$(b^1,\ldots ,b^N)$ is $M^{-1}$. \par\medskip\noindent 2) Each vector $b^j$ can be written as a difference $b^j=b^j_+-b^j_-$ of vectors with non negative
coordinates.  The equality $\nu (\xi^{b^j})=0$ means that there is a constant $c_j\in k^*$ such that $\overline\xi^{b^j_+}-c_j\overline\xi^{b^j_-}=0$ in
$\hbox{\rm gr}_\nu R$. If we have presented, according to Corollary \ref{binomials}, the graded algebra $\hbox{\rm gr}_\nu R$ as
$k[(U_i)_{i\in I}]/\bigl((U^m-\lambda_{mn}U^n)_{(m,n)\in L}\bigr)$, this shows that the binomial $U^{b^j_+}-c_jU^{b^j_-}$
must belong to the binomial ideal, and $c_j$ is a monomial in the
$\lambda_{mn}$. Heuristically, this means that some of the binomial relations in $\hbox{\rm gr}_\nu R$ correspond to translations of the center of the
valuation
$\nu$ with respect to the origin of the new coordinates in toric extensions subordinate to $\nu$, and the constants $\lambda_{mn}$
correspond to the coordinates of these translations.\end{Remark}\par\noindent
If $R$ is endowed with a rational valuation $\nu$, it will be assumed implicitely that $\nu (\xi^{b^j})\geq 0$ for $j\in F$; the maximal
ideal with respect to which we localize
\linebreak $R[\xi^{b^1},\ldots,
\xi^{b^N}]$ will then be the center of $\nu$ in this ring. In this manner we obtain a birational toric extension $$R\subset R'\subset R_\nu .$$ In this
case we say that the toric extension is \textit{subordinate} to $\nu$.
Let us now consider a valuation $\nu_1$ of $R$, with which $\nu$ is composed;  denote by $\Phi_1$ the group of $\nu_1$, by
${\mathbf p}_1$ the center of $\nu_1$ in $R$ and by $\overline \nu$ the valuation on $R/{\mathbf p}_1$ which is the image of $\nu$. 
\begin{lemma}\label{torictensor} Let $\nu$ be a valuation of $R$ and $R\to R'$ be a birational toric extension subordinate to $\nu$ associated to
a finite family
$(\xi_j)_{j\in F}$ of elements such that
$\nu_1 (\xi_j)=0$ for $j\in F$. For each $\phi_1\in \Phi_1\cup\{0\}$, each element of the the kernel and cokernel of the natural map of
$R/{\mathbf p}_1$-modules
$$\P_{\phi_1}(R)/\P^+_{\phi_1}(R)\otimes_{R/{\mathbf p}_1}R'/{\mathbf p}_1R'\to \P_{\phi_1}(R')/\P^+_{\phi_1}(R')$$ is annihilated by a
monomial in the $(\xi_j)_{j\in F}$. In particular they are torsion modules, and the kernel is the torsion submodule of
$\P_{\phi_1}(R)/\P^+_{\phi_1}(R)\otimes_{R/{\mathbf p}_1}R'/{\mathbf p}_1R'$.
\end{lemma}
\begin{proof} Let $$x=(\sum_i c_i\xi^{\alpha_{i+}-\alpha_{i-}})(\sum_k d_k\xi^{\beta_k})^{-1},$$ with $c_i,d_k\in R$, $\alpha_{i+}\geq0$, $\alpha_{i-}\geq
0$,  and
$\nu(\sum_k d_k\xi^{\beta_k})=0$ be a representative of an element
$\overline x$ of
$\P_{\phi_1}(R')/\P_{\phi_1}^+(R')$. Set $\alpha_-=\hbox{\rm max}_i\alpha_{i-}=\alpha_{i_0-},\ \alpha_+=\alpha_{i_0+}$.\par\noindent The element
$\xi^{\alpha_-}(\sum_k d_k\xi^{\beta_k})x$ is in $\P_{\phi_1}(R)$, so that
$$\xi^{\alpha_+}x=\xi^{\alpha_+-\alpha_-}.\xi^{\alpha_-}(\sum_k d_k\xi^{\beta_k})(\sum_k d_k\xi^{\beta_k})^{-1}x$$ is in $\P_{\phi_1}(R)R'$. This shows
that the cokernel of our map is a torsion module. The proof for the kernel is analogous if we notice that
$\P_{\phi_1}(R')/\P_{\phi_1}^+(R')$ is a torsion-free $R/{\mathbf p}_1$-module.\par\noindent This follows from the fact that the $R'/\P_0^+(R')$-module
$\P_{\phi_1}(R')/\P^+_{\phi_1}(R')$ is torsion-free (see subsection \ref{speciatogr}) and the map $R/{\mathbf p}_1\to R'/\P_0^+(R')$ is
injective since by construction
$\P_0^+(R')\cap R=\P_0^+(R)={\mathbf p}_1$. The last sentence also follows from this remark.
\end{proof}
\begin{remark} For
$\phi_1=0$, our map becomes
$R'/{\mathbf p}_1R'\to R'/\P_0^+(R')$. This map is clearly surjective;  its kernel consists of the torsion elements of the $R/{\mathbf
p}_1R$-module $R'/{\mathbf p}_1R'$. The ideals ${\mathbf p}_1R'$ and $\P_0^+(R')$ correspond respectively to the total transform and the
strict transform under our toric extension of the center of $\nu_1$ in $\hbox{\rm Spec}R$.\end{remark}
\begin{definition}
Assertion TLU(d) (Toric local uniformization in dimension $\leq d$) is the following: For every local equicharacteristic excellent integral domain
$R$ of dimension
$\leq d$, given a rational valuation $\nu$ on $R$, and elements $(\xi_i)_{i\in I}$ in $R$ whose initial forms generate the graded $k$-algebra $\hbox{\rm
gr}_\nu R$, there is a birational toric (in the $\xi_i$) extension
$R\subset R'\subset R_\nu$  such that:\par\noindent a) The local ring $R'$ is regular, with a system of coordinates containing
elements whose images in $\hbox{\rm gr}_\nu R'$ are part of a minimal system of generators.\par\noindent
b) The kernel of the scalewise completion map $\hat R{'}^{m'}\to \hat R{'}^{(\nu)}$ (see Proposition \ref{madicomp} below and the remarks which
follow)  is generated by part of a regular system of generators of the maximal ideal of the regular local ring $\hat R{'}^{m'}$; in particular the complete
local ring
$\hat R{'}^{(\nu)}$ is regular.
\end{definition}
\begin{definition}
Assertion TP(d) (Toric principalization in dimension $\leq d$) is the following: For every local equicharacteristic excellent integral domain $R$ of dimension
$\leq d$, given a rational valuation $\nu$ and an ideal $I$ of $R$ there is a birational toric extension $R'\subset R_\nu$ such that $IR'$ is of the
form $(\xi')^\delta u'$, where $u'$ is a unit of $R'$, the $\xi'_j$ are elements of the maximal ideal $m'$ and their images in
$\hbox{\rm gr}_\nu R'$ are part of a system of generators of this $k$-algebra.
\end{definition}
\begin{definition} Let $k$ be a ring and $\Phi$ be a totally ordered group of height $\hbox{\rm h}$ and $$(0)\subset \Psi_{h-1}\subset
\cdots \subset \Psi_1\subset \Psi_0=\Phi$$ the ordered sequence of its convex subgroups. A graded map
$G\to G'$ of $\Phi$-graded $k$-algebras is said to be \textit{scalewise graded} if for each $\Psi_k$ the image of the subalgebra $G_k\subset
G$ generated by homogeneous elements whose degree is in $\Psi_k$ is contained in the subalgebra $G'_k$ generated by homogeneous
elements of $G'$ whose degree is in
$\Psi_k$.
\end{definition}
\begin{lemma}\label{scalbir} If $k$ is a field, if $G$ and $G'$ are $\Phi$-graded integral $k$-algebras, if all homogeneous
components of $G$ and $G'$ are $k$-vector spaces of dimension $\leq 1$ and if the semigroups $\Gamma=\{\phi\in \Phi/G_\phi\neq (0)\}$ and
$\Gamma'=\{\phi\in \Phi/G'_\phi\neq (0)\}$ generate $\Phi$ as a group, then a graded inclusion $G\to G'$ is birational.\end{lemma}
\begin{proof} It suffices to show that any non zero homogeneous element of $G'$ is the quotient of the images of two homogeneous elements of $G$. Such an
element is of the form $cy_\phi$ where $y_\phi$ is a generator of the $k$-vector space $G'_\phi$ and $c\in k^*$. Since $\phi$ is the difference of two
elements of $\Gamma$, we may write $y_\phi=d_\phi x_{\phi_1}x^{-1}_{\phi_2}$ where $x_{\phi_i}$ is a generator of $G_{\phi_i}$ and $d_\phi\in k^*$. The
result follows.\end{proof}
 \par A rather special
case (valuations of height one and center
$m$) of the first part of the following Proposition was proved \textit{en passant} by Zariski (see [Z6], from the bottom of p. 63) and has been proved
again by Spivakovsky ([S2]) with a different perspective, in both cases without the assumption that the residue field is algebraically closed; see the
remark after Proposition \ref{complht1}. Spivakovsky has also communicated to me that he had stated a different extension of his result to arbitrary finite
heights.\par\medskip\noindent
\begin{proposition}\label{madicomp} {\rm *}Let $R$ be a local equicharacteristic excellent integral domain with maximal ideal $m$ and an
algebraically closed residue field, and let
$\nu$ a rational valuation on $R$ with value group $\Phi$.\par\noindent
 a) There exists a prime ideal $H$ of the $m$-adic completion $\hat R^m$ of $R$ such that $\nu$ extends to a valuation $\hat\nu$ of 
$\hat R^m/H$ in such a way that:\par\noindent b) the associated graded map
$$\hbox{\rm gr}_{\nu} R\to \hbox{\rm gr}_{\hat {\nu}}\hat R^m/H$$ is a scalewise birational extension of $\Phi$-graded $
R/m$-algebras, uniquely determined by $(R,\nu)$.\par\noindent If the valuation $\nu$ is of height one, the map $\hbox{\rm gr}_{\nu} R\to \hbox{\rm
gr}_{\hat {\nu}}\hat R^m/H$ is an isomorphism.{\rm *}
\end{proposition}
The proof is not yet complete and should be the subject of a subsequent paper.
\begin{Remark} 1) The fact that $\hat\nu$ extends $\nu$ implies that
$H\cap R=(0)$.\par\noindent
2) The fact that the extension of the graded rings is scalewise birational implies that any toric modification we make in $\hat R^{(\nu)}$ with
respect to representatives $(\eta_j)_{j\in J}$ of the generators of the graded algebra $\hbox{\rm gr}_{\hat\nu}\hat R^{(\nu)}$ can be viewed as
coming from a toric modification in $R$ with respect to the chosen system of representatives
$(\xi_i)_{i\in I}$.
3) Although the ideal $H$ is unique in the case considered by Zariski and
Spivakovsky, it is not uniquely determined by $(R,\nu)$ in general; it is only the graded map associated to the map
$R\to
\hat R^m/H$ which is uniquely determined. Once $H$ is fixed, the extension of $\nu$ is unique.
\par\noindent  The ideal $H$ depends on the choice of the representatives
$\xi_i\in R$ of the generators of  $\hbox{\rm gr}_\nu R$ and of the complementary generators of the maximal ideal of $R$.\par  Let
$R=k[x_1,x_2,x_3]_{(x_1,x_2,x_3)}$, let
$w(x_1)=\sum_{k=1}^\infty c_kx_1^k$ be a series which is transcendental over $k(x_1)$, consider the ring $\overline R=k[x_1,x_2]_{(x_1,x_2)}$, its injection
$\iota\colon\overline R\subset k[[x_1]]$ defined by $x_1\mapsto x_1,\ x_2\mapsto w(x_1)$ and the valuation $\overline\nu$ of $\overline R$ defined
by the restriction via the injection 
$\iota$ of the $x_1$-adic valuation. Let $\nu$ be the rational valuation of $R$ composed of $\overline\nu$ and the $x_3$-adic valuation. As we have
seen in example \ref{excomp} of subsection \ref{nuadiccomp}, the ideal
$\overline H$ is the ideal of $k[[x_1,x_2]]$ generated by $x_2-w(x_1)$.
For given series without constant terms $y(x_3),z(x_3)$ we may choose as representatives of $x_1,x_2$ the elements
$x_1+y(x_3),\  x_2+z(x_3)$. The corresponding ideal $H_{y,z}$ of $k[[x_1,x_2,x_3]]$ is
$(x_2+z(x_3)-w(x_1+y(x_3)))k[[x_1,x_2,x_3]]$. For each pair $y,z\in k[[x_3]]$ of series without constant term the valuation $\hat\nu$ on
$k[[x_1,x_2,x_3]]/H_{y,z}$ which is composed of the valuation
$\hat{\overline\nu}$ (which is the $x_1$-adic valuation on $k[[x_1]]$), and the $\tilde x_3$-adic valuation, where $\tilde x_3$ is the image of $x_3$
modulo $H_{y,z}$, induces the valuation $\nu$ on $R$.\par\medskip\noindent
\textit {In spite of this non uniqueness, I shall denote by $\hat R^{(\nu)}$ a quotient $\hat R^m/H$ as in the Proposition and sometimes even speak of ``the''
scalewise completion of $(R,\nu)$.} In a sense they are all ''equisingular'' with respect to $\hat\nu$, just as general sections by non singular spaces
transverse to an equisingularity stratum through a singular point are all equisingular.\par\medskip\noindent 
4) A completion which does not induce an isomorphism on the associated graded rings may be somewhat surprising; let us take example \ref{resval} and
assume that $\hat R^m$ is an integral domain but $f\hat R^m$ is not a prime ideal; say that we have $f\hat R^m= \hat g\hat h\hat R^m$ with $\hat g\hat
R^m$ and
$\hat h\hat R^m$ prime ideals. By Proposition \ref{genreg} below, an ideal $\overline H$ of $\hat R^m/f\hat R^m$ as in the Proposition contains a unique
minimal prime of
$\hat R^m/f\hat R^m$. Assume that we are in the simple situation where $\overline H=\hat g(\hat R^m/f\hat R^m)$. The valuation $\nu$ extends to a
valuation on
$\hat R^m$ which is composed of the $\hat g$-adic valuation, of height one, of $\hat R^m$, and the valuation $\hat{\overline \nu}$ on $\hat R^m/\hat g\hat
R^m$ which extends the valuation $\overline\nu$ on $R/(f)$. Remark that this is a case where $\overline H\neq (0)$ but $H=(0)$. Now by example
\ref{resval}, the graded ring
$\hbox{\rm gr}_{\hat\nu}\hat R^m$ is equal to $\hbox{\rm gr}_{\hat{\overline\nu}}(R^m/\hat g\hat R^m)[\hat G]$, where $\hat G$ is the initial form of
$\hat g$ with respect to $\hat\nu$. The inclusion $$\hbox{\rm gr}_\nu R=\hbox{\rm gr}_{\overline\nu}(R/fR)[F]\subset\hbox{\rm
gr}_{\hat{\overline\nu}}(R^m/\hat g\hat R^m)[\hat G]$$ is described by saying that $\hbox{\rm gr}_{\overline\nu}(R/fR)\to \hbox{\rm
gr}_{\hat{\overline\nu}}(R^m/\hat g\hat R^m)$ is a scalewise birational map by induction on the dimension, and $F\mapsto (\hbox{\rm
in}_{\hat{\overline\nu}}\overline{\hat h})\hat G$, where
$\overline{\hat h}$ is the image of $\hat h$ in $\hat R^m/\hat g\hat R^m$. This last step introduces a scalewise
birational character in the graded rings extension if it was not already present for $R/fR$.
\end{Remark}

\begin{corollary}\label{Abhyin} Abhyankar's inequality can be made more
precise in the case of a rational valuation $\nu$ of an excellent local ring $R$; with the notations of section \ref{valco}, we have:
$$\hbox{\rm r}(\nu )\leq \hbox{\rm dim}\hat R^{(\nu)}\leq \hbox{\rm
dim}R.$$ In particular, if $\hbox{\rm r}(\nu )=\hbox{\rm
dim}R$ the ideal $H$ must be a minimal prime of $\hat R^m$ and if in addition $R$ is analytically irreducible, we have $\hat R^{(\nu)}=\hat R^m$. In all
cases, in view of Proposition
\ref{genreg}, the ideal
$H$ contains a unique minimal prime of $\hat R^m$.
\end{corollary} 
\begin{Remark}
 1) If $\nu$ is a valuation of height one whose ring birationally dominates $R$,
we have $\hat R^{(\nu)}=\hat R^\nu$. One may ask,  if we only assume that $\nu$ is of height one, what is the relationship of the ring $\hat R^{(\nu)}$ 
with the $m$-adic Hausdorff completion of $\hat R^\nu$. \par\noindent
2) The Proposition can be used to study extensions of a rational valuation $\nu$ of an excellent local ring $R$ to the completion $\hat
R^m$; by composition with $\hat\nu$, every valuation of the regular local ring $\hat R^m_H$ centered at is maximal ideal gives rise to a valuation of
$\hat R^m$ extending
$\nu$. As we saw in Corollary \ref{Abhyin}, if
$\hbox{\rm r}(\nu )=\hbox{\rm dim}R$, the ideal $H$ must be a minimal prime of $\hat R^m$ and if $R$ is analytically irreducible, $\nu$ extends to $\hat
R^m$ with the same group. This extension is then unique since $\hat\P_\phi=\P_\phi\hat R^m$.\par If $\hbox{\rm r}(\nu )=\hbox{\rm dim}R-1$, and we
assume that $H$ is not a minimal prime of $\hat R^m$, then the localization
$\hat R^m_H$ is one-dimensional because of the Proposition, and it is a regular local ring by Proposition \ref{genreg} below. It is a valuation ring with group
$\Z$. Denote by
$\mu$ its valuation. The valuation $\tilde\nu$ on $\hat R^m$ with values in $\Z\oplus \Phi$ ordered lexicographically, defined by $\tilde\nu (x)=(\mu
(x),\hat\nu (x \ \hbox{\rm mod.}H))$ is an extension of $\nu$ to $\hat R^m$.\par In order to check the uniqueness of the extension,
assuming that $R$ is analytically irreducible, remark that $\bigcap_{\phi\in \Phi_+} \P_\phi\hat R^m\subseteq H$ by construction of $H$.
The first ideal is prime because the graded ring associated to the filtration of $\hat R^m$ by the $\P_\phi\hat R^m$ is equal to
$\hbox{\rm gr}_\nu R$, which is a domain. Since $R$ is analytically irreducible, and in view of the dimension condition, we have
$\bigcap_{\phi\in \Phi_+}
\P_\phi\hat R^m= H$. Then, given an extension $\tilde\nu\ '$ of $\nu$ to $\hat R^m$ with group $\tilde\Phi$, the ideal $H$ appears as
the set of elements
$z\in\hat R^m$ such that $\tilde\nu\ '(z)>\phi\ \ \forall\ \phi\in \Phi\subset \tilde\Phi$, and $\Phi\subset \tilde\Phi$ is a convex
subgroup of the group
$\tilde\Phi$. Then $\tilde\nu\ '$ is composed of a valuation on the localization
$\hat R^m_H$, which is regular of dimension one, and a valuation on $\hat R^m/H$, which is uniquely determined since its valuation
ideals must be the $\P_\phi\hat R^m/H$, and coincides with the extension $\hat\nu$ of the Proposition. This shows that if we assume
$R$ analytically irreducible, in the case where $\hbox{\rm r}(\nu )=\hbox{\rm dim}R-1$ there is a unique extension of $\nu$ to $\hat R^m$, described as
above. These two results on extensions of valuations to
$\hat R^m$ in the case where it is a domain are due to Heinzer and Sally ([H-S]) under an assumption weaker than excellence.
\par\noindent  3)
After replacing $R$ by
$\hat R^{(\nu)}$ and $\nu$ by $\hat \nu$, we may assume that the local ring $R$ is complete for the $m$-adic topology, with a valuation
$\nu$, so that it is henselian by Hensel's lemma, admits (in the equicharacteristic case) a field of representatives by Cohen's theorem and is complete for
the $\nu$-adic topology. 
\end{Remark} 
\begin{definition} Given $R$ and a rational valuation as above, we call a map $R\to \hat R^{(\nu)}$ such as we have built in Proposition
\ref{madicomp} a \textit{scalewise}
$\nu$-adic completion of $(R,\nu )$.\end{definition}
\begin{definition} We say that the valuation $\nu$ is \textit{normally flat in degree $\phi$} if either $\P_\phi/\P^+_\phi$ is
zero, or it is a flat $\overline R$-module, where $\overline R=R/{\mathbf p}$ with ${\mathbf p}$ the center of $\nu$ in $R$.\end{definition}
\noindent \textbf {Question.}: Given a rational valuation $\nu$ on $R$,
and a valuation $\mu$ with semigroup $\Delta$, with which $\nu$ is composed, is it true that for any
$\delta\in \Delta$ there exits a birational toric extension $R\to R'$ compatible with $\nu$, determined by elements of $\mu$-valuation zero, and such that
$\mu$ is normally flat in degree $\delta$ on $R'$?
\par \medskip
I now show that with this definition of scalewise completion one has a Cohen-type theorem.
\begin{lemma}\label{height-one finiteness}
Let $R$ be an excellent equicharacteristic local integral domain, and $\nu$ a rational valuation on $R$. Assume that the valuation $\nu$ is of height one.
Lift a system of generators
$(\overline\xi_i)_{i\in I}$ of the $k$-algebra $\hbox{\rm gr}_{\hat\nu}\hat R^{(\nu)}$ to elements
$\xi_i\in R$ and choose a field of representatives $k\subset\hat R^{(\nu )}$. The $k$-algebra $k[(\xi_i)_{i\in I}]$ is dense in $\hat R^{(\nu )}$ for the
$\hat\nu$-adic topology. 
\end{lemma}
\begin{proof} We retain the notations of Proposition \ref{madicomp}, remarking that our assumption implies that $\hat\nu$ is trivial on $k$ and
remembering that in this case $\hbox{\rm gr}_{\hat\nu}\hat R^{(\nu)}=\hbox{\rm gr}_\nu R$. If we pick $x\in \hat R^{(\nu )}$, its $\hat\nu$-adic initial
form is in $\hbox{\rm gr}_\nu R$, and therefore $\hbox{\rm in}_{\hat\nu} x=P_1(\overline \xi_i)\ ,\  \hbox{\rm with}\ P_1\in k[(X_i)_{i\in I}]$. The
polynomial $P_1$ is in fact a term $cX^a$ with $c\in k^*$ and $X^a$ a monomial, but I will not use this here. Abbreviating $(\xi_i)_{i\in I}$ to $(\xi )$, we set
$x^{(0)}=x$,
$x^{(1)}=x-P_1(\xi)$, and iterate this process. We have for all $j\geq 1$ the inequality $\hat\nu (x^{(j+1)})>\hat\nu (x^{(j)})$
and thus $\hat\nu (x^{(j+1)}-x^{(j)})=\hat\nu (x^{(j)})$. Now remark that since $\hat R^{(\nu )}$ 
is n\oe therian by Proposition \ref{madicomp}, the valuation $\hat\nu (\hat m)$ is $>0$, and since
$\hat \nu$ is of height one, for any $\phi\in \Phi_+$ there is an integer $N(\phi)$ such that
$\hat m^{N(\phi)}\subset \hat\P_\phi$. Then there are only finitely many distinct
$\hat\P_{\phi '}$ containing $\hat\P_\phi$, so that in a finite number of steps we have $\hat\nu
(x^{(j+1)})>\phi$, which proves that the sequence $x^{(j)}$ converges to $x$, and the Lemma.
\end{proof}
\begin{remark}\label{oneapprox} Without the height one assumption, it could happen that the
valuations of the differences $x-\sum_{j=1}^rP_j(\xi )$ increase at each step of the
construction but remain bounded in the semigroup $\Phi_+$.
\end{remark}
\par\noindent \textit{I assume until the next Corollary that the local
ring $R$ is complete for its $m$-adic topology; see remark 2) at the end of subsection
\ref{nuadiccomp}.} 
\begin{example} Let $k$ be a field, let $R$  be a complete n\oe therian local $k$-algebra with
residue field $k$, and let $f\in R$ generate a prime ideal. We go back to example \ref{resval} of subsection \ref{compargr} 
and its notations. Let $x\in R$ and let us write it in the form $x=af^\ell$ with $a\notin (f)$. We have $\nu_1(x)=\ell$ and its $\nu$-initial form  is $(\hbox{\rm
in}_{\overline \nu}\ \overline{a})F^\ell$. Let us fix a field of representatives $k\subset R$; in general,  the element $a$ will not be a polynomial
with coefficients in $k$ in the liftings $(\xi_i)$ to $R$ of the generators of
$\hbox{\rm gr}_{\overline\nu}\overline R$, so that we will have, denoting by $\tilde{a}$ a
lifting in $R$ of 
$\hbox{\rm in}_{\overline \nu}\overline a$ which is a polynomial in the
$(\xi _{i})$, and choosing $f$ as representative of $\xi$, a lifting $\widetilde{\hbox{\rm in}}_\nu x$ of ${\hbox{\rm in}}_\nu x$ which is polynomial in
the $\xi_i$ and $f$, and such that
$$x-\widetilde{\hbox{\rm in}}_\nu x=(a-\tilde {a})f^\ell $$
and the $\nu_1$-value
will not have changed. \par However, if we approximate $\overline a$ in $\overline R$ by a
sequence of polynomials $ P_s((\overline\xi_i)_{i\in I_1})$ in
$k[(\overline\xi_i)_{i\in I_1}]$, we will have after replacing the $\overline\xi_i$ by their representatives $\xi_i$ in $R$ $\hbox{\rm
lim}_{s\to\infty}\bigl(a-P_s((\xi_i)_{i\in I_1})\bigr)\in fR$ so that we find as limit of a sequence of polynomials $ P_s((\overline\xi_i)_{i\in
I_1})f^\ell$,  with the notations of the first remark of this subsection, an element
$z^{(1)}\in R$ congruent to $a f^\ell$ modulo the ideal $f^{\ell +1}R$. We may write $z^{(1)}=a_\ell f^\ell$. Then we have to repeat the process with
$x_1=x-z^{(1)}$ and the next term. At each stage we build by the convergence of a series an element $z^{(t)}\in R$ and $x$ is the limit of the series
$\sum_{t\geq1} z^{(t)}$. In fact, the element
$x$ then appears as the sum of a series
$\sum_{k=\ell}^\infty a_k f^k$. I call this situation ``scalewise density''\footnote{F.-V. Kuhlmann has suggested to me that this idea could be related
to Ostrowski's pseudo-Cauchy sequences (see [Roq], [Ka]) and Ribenboim's notion of ``complet par \'etages'' ([Ri]) in the theory of completion of valued
fields. The problem is similar, but I did not compare the two definitions; there is a major difference in that I work with n\oe therian rings, not valuation
rings.} for the subalgebra $k[(\xi_i)_{i\in I};\ f] $ in $R$.
\end{example}\par\noindent
Now for the general case:
\par\medskip\noindent Let $R$ be a \textit{complete} integral n\oe therian local ring with a field of
representatives $k$, and $\nu$ a $k$-valuation of its field of fractions, positive
on $R$ and whose center is the maximal ideal $m$ of $R$. Let $(\xi_i)_{i\in I}$ be a set of
elements of $R$ whose images in $\hbox{\rm gr}_\nu R$ generate it as a $k$-algebra. Consider
the set of isolated subgroups of the value group $\Phi$ of $\nu$ $$(0)=\Psi_h\subset
\Psi_{h-1}\subset \cdots \subset \Psi_t\subset \cdots\subset\Psi_1\subset \Psi_0=\Phi ,$$
and the corresponding valuations $\nu_t$ with value group $\Phi/\Psi_t$ indexed by their
height. Let $${\mathbf p}_r \subset {\mathbf
p}_{r-1}\subset\cdots\subset {\mathbf p}_1\subset m$$
be the sequence of the distinct centers of the
valuations $\nu_t$ in $R$. Let us partition the set $I$ as follows: $I=\bigcup_{k=1}^{h}I_t$
where
$$I_t=\{i\in I\mid\nu (\xi_i)\in\Psi_{h-t}\setminus\Psi_{h-t+1}\},$$
remarking that some of the sets $I_t$ may be empty. Let us also define the following
sequence of subrings of $R$: \par\noindent
$\hat R_1$ is the closure of 
$k[(\xi_i)_{i\in I_1}]$ in $ R$ for the topology defined by the filtration by the
ideals $(\P_\psi)_{\psi\in\Psi_{h-1}}$. We can
continue and define inductively $\hat R_t$ to be the closure of $\hat R_{t-1}[(\xi_i)_{i\in I_t}]$
in $ R $ for the topology defined by the filtration by the ideals $(
\P_\psi)_{\psi\in\Psi_{h-t}}$.
\begin{proposition}\label{allocal} For each $t,\ 1\leq t\leq h$, the ring $\hat R_t$ is local and its maximal ideal is generated by the
$(\xi_i)_{i\in\bigcup_{\ell =1}^tI_\ell}$. \end{proposition}
\begin{proof} The ring
$\hat R_1$ is local because if
$y=\sum a_\alpha\xi^\alpha$ is not in the closure of the ideal
$(\xi_i)_{i\in I_1}$, we can write $y=a_0(1+\sum_{\alpha>0} a_0^{-1}a_\alpha\xi^\alpha)$ with $a_0\in k^*$, and to show that $y^{-1}$
is in $\hat R_1$ we have to show that the series $$1-\sum_{\alpha>0} a_0^{-1}a_\alpha\xi^\alpha+(\sum_{\alpha>0}
a_0^{-1}a_\alpha\xi^\alpha)^2+\cdots +(\sum_{\alpha>0}
a_0^{-1}a_\alpha\xi^\alpha)^n+\cdots $$ converges, which is clear because $\Psi_{h-1}$ is a group of height one. The passage from
$\hat R_{t-1}$ to
$\hat R_t$ is similar, if we notice the following \begin{lemma}\label{archi+} If $\Psi$ is an isolated subgroup of $\Phi$ such that $\Phi/\Psi$ is of height
one, the semigroup $\Phi_+\setminus\Psi_+$ is archimedian.\end{lemma}\begin{proof} It follows directly from the fact that $\Phi/\Psi$ is
archimedian.\end{proof} So, if $\nu (\xi)\in I_t$, then as $n\in \N$ tends to infinity, $n\nu (\xi)$ eventually exceeds any element of $\Psi_{h-t}$.
\end{proof}
\par\noindent \begin{definition} We say that $k[(\xi_i)_{i\in I}]$ is
scalewise dense in $  R $ for the $ \nu$-adic topology if $\hat R_h$ is equal to $ R
$.\end{definition}
\begin{proposition}\label{denseinscales} Let $R_{\nu}$ be the ring of a rational valuation of a
complete equicharacteristic n\oe therian local integral domain $R$. If we fix a field of
representatives $k$ of $ R$ and choose representatives $(\xi_i)_{i\in I}$ in
$R$ of a system of generators $(\overline\xi_i)_{i\in I}$ of the $k$-algebra
$\hbox{\rm gr}_\nu R$, where $k=R/m=R/(m_\nu\cap R)$, the
algebra $k[(\xi_i)_{i\in I}]$ is scalewise dense in $ R $ for the $\nu$-adic
topology. \end{proposition} 
\begin{proof} It is useful to keep in mind Corollary \ref{strugr}, which gives us the structure of the initial forms
 which appear here. We have to prove that, given
$x\in  R$, there is a sequence
$(r_j)_{j\in \N}$ of elements of $\hat R_{h-1}[(\xi_i)_{i\in I_h}]$ converging to $x$ in the
$\tilde{\nu_h}$-adic topology. Define the sequence of $x^{(j)}$ and $P_j(\xi_i)$ as in the proof
of Lemma \ref{height-one finiteness}, and let us build from it a sequence $(r_j)_{j\in \N}$ of elements
of $\hat R_{h-1}[(\xi_i)_{i\in I_h}]$ such that the valuations $\nu (x-r_j)$ exceed, for large
$j$, any $\phi\in \Phi_+$ given in advance. Let $s$ be the least integer such that for
infinitely many values of the integer $j$ we have $\nu_s(x^{(j+1)})>\nu_s(x^{(j)})$. If $s=1$,
since the valuation $\nu_1$ is of height one, for any $\phi_1\in \Phi_{1+}$, after finitely
many steps we have $\nu_1 (x^{(j)})>\phi_1$, hence for any $\phi\in \Phi_+$, after finitely
many steps we have $\nu(x^{(j)})>\phi$, so we may take $r_j=x^{(j)}$ and get $\nu
(x-\sum_{r=0}^{j-1}P_r(\xi_i))>\phi$ and this shows that $x\in \hat R_h$. If $s>1$, by
definition, there exists a $j_0$ such that $\nu_{s-1}(x^{(j+1)})= \nu_{s-1}(x^{(j)})$ for $j\geq
j_0$. By substracting from $x$ a polynomial in the $\xi_i$, we may assume that
$j_0=1$.\par\noindent Note that the center of $\nu_{s-1}$ is necessarily distinct from the
center of $\nu_s$, since if the two centers are equal, there are only finitely many distinct
$\nu_s$-ideals between two consecutive $\nu_{s-1}$-ideals of $R$, according to [Z-S], Vol.
II, Appendix 3, Corollary p.345. \par\noindent We consider the initial forms
$\hbox{\rm in}_{\nu_{s-1}} x^{(j)}\in \hbox{\rm gr}_{\nu_{s-1}} R $; for $j\geq j_0$, they all
have the same degree, say $\phi_{s-1}$. Since the residual valuation $\overline \nu_s$ on
$R/(m_{\nu_{s-1}}\cap R)$ is of height one, with group $\Psi_{s-1}/\Psi_s$, and each
homogeneous component of $\hbox{\rm gr}_{\nu_{ s-1}} R $ is a $R/(m_{\nu_{s-1}}\cap
R)$-module of finite type, and therefore complete for the $m/(m_{\nu_{s-1}}\cap R)$-topology,
so that it is also complete for the $\overline\nu_s$ topology by Proposition \ref{nucomplete}, the sequence of the initial forms
$\hbox{\rm in}_{\nu_{s-1}}(x-x^{(i)})$ converges in $(\hbox{\rm gr}_{\nu_{s-1}}
R)_{\nu_{s-1}(x)}$ to a unique limit $\overline x_1^{(1)}$, which we may lift to an element
$x_1^{(1)}\in
\hat R_h$. This element has the property that for all $j\geq 1$, we have 
$\nu_{s-1}(x_1^{(1)})=\nu_{s-1}(x^{(j)})$ and
$\nu_{s-1}(x-x_1^{(1)})>\nu_{s-1}(x-x^{(j)})$.
 Replacing now $x$ by $x-x_1^{(1)}$, and so on, we build a sequence of elements
$x_1^{(q)}\in R_h$ such that $\nu_{s-1} (x-x_1^{(q)})$ increases at each step. 
We may now repeat
the whole process,  replacing the sequence $x^{(j)}$
by $x_1^{(q)}$; $s$ is replaced by $s-1$ 
and ultimately we are reduced to the case where $s=1$,
which is already settled. This proves the Proposition.
\end{proof}
\textit{I stop assuming that the local ring $R$ is complete, but I assume until the end of this section that its residue field is algebraically
closed.} 
\begin{corollary}\label{alldense} Let $\nu$ be a rational valuation on the local equicharacteristic excellent integral domain $R$; let $k$ be the residue field
of
$R$ and let
$(\xi_i)_{i\in I}$ be representatives in $R$ of a system of generators $(\overline
\xi_i)_{i\in I}$ of the graded $k$-algebra $\hbox{\rm gr}_\nu R$. Let $k\subset \hat R^{(\nu)}$ be a field of representatives in the scalewise completion of
$R$ built in Proposition \ref{madicomp} and let $(\eta_j)_{j\in J}$ be representatives in $\hat
R^{(\nu)}$ of a system of generators
$(\overline
\eta_j)_{j\in J}$ of the graded $k$-algebra $\hbox{\rm gr}_{\hat\nu} \hat R^{(\nu)}$, which are Laurent monomials in the $(\overline\xi_i)_{i\in I}$ (see
\ref{madicomp}). The subalgebra
$k[(\eta_j)_{j\in J}]$ is scalewise dense in 
$\hat R^{(\nu)}$. The local subrings $\hat R_t\subset \hat R^{(\nu)}$ of Proposition \ref{allocal} are complete and
n\oe therian, and the valuation $\nu$ restricted to $\hat R_t$ is analytically monomial with respect to the
$(\eta_j)_{j\in\bigcup_{\ell=1}^tI_\ell}$.\end{corollary}
\begin{proof} The first part follows from the Proposition in view of the fact that $\hat R^{(\nu)}$ is complete for its $\hat m$-adic topology. To prove that
the rings $R_t$ are n\oe therian  the last statement, use descending induction on
$k$: the ring
$\hat R_h=\hat R^{(\nu)}$ is n\oe therian, and each $\hat R_{t-1}$ is isomorphic to the quotient of $\hat R_t$ by the closure of the ideal generated by
$(\xi_i)_{i\in I_t}$. It also reproves that the $\hat R_t$ are local.\par To
prove the last two assertions, we can check that $m_{\nu_{h-1}}\cap \hat R_1=(0)$ so that an element of $R/(m_{
\nu_{h-1}}\cap   R)$, which is in the closure of $k[(\overline \xi_i)_{i\in I_1}]$, lifts uniquely to an
element of this closure in $R$, and we have seen in Proposition \ref{grismon} that the valuation
$\overline \nu_1$ is monomial with respect to the $(\overline \eta_j)_{j\in J_1}$. One passes from $R_{t-1}$ to $R_t$ in the same way.\end{proof}
\subsection{The $\nu_\A$-adic completion of $\A_\nu (R)$}\label{Bcomp}
\par\medskip\noindent  The main result of this subsection is, given a system $(\overline \xi_i)_{i\in I}$ of generators of the $k$-algebra
$\hbox{\rm gr}_\nu R$, giving rise to a surjective map $k[(U_i)_{i\in I}]\to \hbox{\rm gr}_\nu R$, a field of representatives $k\subset \hat
R^{(\nu)}$ and representatives
$(\xi_i)_{i\in I}$ in $R$ of the $\overline\xi_i$, the extension of the map $k[(u_i)_{i\in I}]\to \hat R^{(\nu)}$ mapping $u_i$ to $\xi_i$ to a continuous
surjective map
$$\widehat{k[(u_i)_{i\in I}]}\to \hat R^{(\nu)}$$ from the scalewise completion of the polynomial ring to the scalewise $\nu$-adic
completion of $R$. When $\widehat{k[(u_i)_{i\in I}]}$ is endowed with the term order obtained by giving to $u_i$ the valuation of its image
$\xi_i\in \hat R^{(\nu)}$, the associated graded map is the map 
$$k[(U_i)_{i\in I}]\to \hbox{\rm gr}_\nu R$$
mapping $U_i$ to $ \overline\xi_i$. \par\noindent This is a
$\nu$-adic analogue, called here $TC(*)$, of Cohen's structure theorem ([B3], Chap. IX, \S 3, No.3, Th. 2, a)). Thus, $TC(*)$ is the
generalization of the embedding of a plane complex branch in $\C^{g+1}$, while the similar result for the valuation algebra corresponds to
the construction of the family degenerating the re-embedded branch to the monomial curve.
\par\noindent For brevity I construct the map directly in the case of the valuation algebra (see Proposition \ref{coordinates}), which is only slightly
more complicated than the
$\nu$-adic Cohen theorem.
\begin{proposition}\label{completion} In the situation of \ref{nuadiccomp}, the completion
$\widehat {\A_\nu (R)}^{\nu_\A}$ is the closure of $\A_{\hat\nu } (\hat R^\nu)$ in the
algebra of restricted power series $\hat R^\nu\{v^\Phi\}$ defined by \rm $$\hat
R^\nu\{v^\Phi\}=\displaystyle\limproj_{\phi\in \Phi_+} (R/\P_\phi )[v^\Phi].$$
 \end{proposition}
 \begin{proof} This follows from the definitions, the fact that $\Ss_\delta=\P_\delta
R[v^\Phi]\cap \A_\nu (R)$ and  [B4], Chap II, \S 3, No. 9, Cor.1.
 \end{proof}
\par Note that the
closure in $\hat R^\nu$  of $\P_\phi$ is $\hat \P_\phi$ and that, with the notations of subsection \ref{nuA}, we
have:$$\widehat {\A_\nu (R)}^{\nu_\A}=\displaystyle\limproj_{\delta\in \Phi_+} \A_\nu
(R)/{\mathcal S}_\delta=
\displaystyle\limproj_{\delta\in \Phi_+}\bigoplus_{\phi \in \Phi} \big({\mathcal P}_\phi (R)/{\mathcal
P}_\phi (R)\cap {\mathcal P}_\delta (R)\big) v^{-\phi}.$$ \par So we may view $\widehat {\A_\nu
(R)}^{\nu_\A}$ as the subring of the $\hat R^\nu$-module $\hat R^\nu[[v^\Phi]]$ of formal power
series $\sum x_\phi v^{-\phi }$, consisting of power series $\sum x_\phi v^{-\phi}$ such that
$\hat \nu (x_\phi )
\geq \phi$ and which are restricted, i.e., such that for any $\psi \in \Phi_+$, all but a finite
number of the coefficients $x_\phi$ are in $\P_\psi$. This last condition implies that when we
formally multiply two series 
$$\big(\sum x_\phi v^{-\phi}\big)\big(\sum x_\psi v^{-\psi}\big)=\sum_{\eta
\in \Phi}\big(\sum_{\phi +\psi =\eta} x_\phi y_\psi \big) v^{-\eta},$$ the coefficient of
$v^{-\eta}$ is a convergent series in $\hat R^\nu$.
\par\medskip\noindent
\begin{lemma}\label{closurefiber} The closure $$\overline{(v^{\Phi_+})\widehat{{\mathcal
A}_\nu(R)}^{\nu_\A}}\ \ \ \hbox{\rm in}\ \ \ \widehat{\A_\nu (R)}^{\nu_\A}\ \ \hbox{\rm of the
ideal }(v^{\Phi_+})\widehat{\A_\nu (R)}^{\nu_\A}$$ is the set of series
$$\sum x_\phi v^{-\phi}\in\widehat{\A_\nu (R)}^{\nu_\A}\ \ \ \hbox{\rm such that}\ \
\hat\nu(x_\phi )>\phi .$$
\end{lemma}
\begin{proof} The proof follows immediately from \ref{completion}.\end{proof}
\begin{definition}\label{completloc} With the notations of subsection \ref{speciatogr}, let us define on the model of
{\rm ([B4], exerc. 27 for
\S 2)} the algebra
$$\widehat{\A_\nu (R)}^{\nu_\A}\{(v^{\Phi_+})^{-1}\}=\displaystyle\limproj_{\delta \in
\Phi_+}(v^{\Phi_+})^{-1}\A_\nu (R)/\Ss_\delta=\widehat{\A_{\hat\nu} (\hat
R^\nu)}^{\hat\nu_{\hat\A}}\{(v^{\Phi_+})^{-1}\}.$$
\end{definition}
 This ring is also the Hausdorff
completion of  $(v^{\Phi_+})^{-1}\A_\nu (R)$ for the topology having as a fundamental system of
neighborhoods the $(v^{\Phi_+})^{-1}\P_\delta$. I follow Bourbaki's convention according to
which writing $A\{S^{-1}\}$ entails the completion of $A$, i.e., $A\{S^{-1}\}=\hat A\{{S'}^{-1}\}$
where $S'$ is the image of $S$ in $\hat A$. The inclusion
$\A_\nu (R)\longrightarrow 
\A_\nu (R)\{(v^{\Phi_+})^{-1}\}$ is universal with respect to continuous maps $\A_\nu (R)\to B$
for linearly topologized Hausdorff and complete rings $B$ such that the images in $B$ of the
elements of the multiplicative set $(v^{\Phi_+})$ are invertible.\par\noindent I also consider
the Hausdorff completion $\widehat{\hbox{\rm gr}_\nu}R$ of the graded algebra $\hbox{\rm
gr}_\nu R$ with respect to the valuation $\nu_{\rm gr}$. 
\par\noindent
 According to the general properties of completions of linearly topologized modules ([Ma], Prop.
in 23.I, p.167), the completion of $\widehat{{\mathcal A}_\nu(R)}^{\nu_\A}/(v^{\Phi_+})\widehat{{\mathcal
A}_\nu(R)}^{\nu_\A}$ for the quotient topology is $\widehat{{\mathcal
A}_\nu(R)}^{\nu_\A}/\overline{(v^{\Phi_+})\widehat{{\mathcal A}_\nu(R)}^{\nu_\A}}$; it is also the
completion, for the quotient topology, of ${\mathcal A}_\nu(R)/(v^{\Phi_+}){\mathcal A}_\nu(R)$. I can
now state
\begin{proposition}\label{complete deformation} a) The natural map
$${\mathcal A}_\nu(R)\longrightarrow \hbox{\rm gr}_\nu R$$ defined by
$$x_\phi v^{-\phi}\mapsto x_\phi \ \ \hbox{\rm mod.}{\mathcal P}^+_\phi (R)$$ induces an
isomorphism of filtered rings 
$$\widehat{{\mathcal A}_\nu(R)}^{\nu_\A}/\overline{(v^{\Phi_+})\widehat{{\mathcal
A}_\nu(R)}^{\nu_\A}}\stackrel{\simeq}\longrightarrow\widehat{\hbox{\rm gr}_\nu} R.$$\par\noindent  b) The
natural inclusion $R[v^{\Phi_+}]\longrightarrow {\mathcal A}$ obtained by considering only the part of
negative degree ( i.e., $\phi \in \Phi_-$) of ${\mathcal A}$ induces an isomorphism
$$R[v^{\Phi_+}]\{(v^{\Phi_+})^{-1}\}= \hat R^\nu
\{v^{\Phi_+}\}\{(v^{\Phi_+})^{-1}\}\stackrel{\simeq}\longrightarrow
\widehat{{\mathcal A}_\nu(R)}^{\nu_\A}\{(v^{\Phi_+})^{-1}\}.$$\noindent Here, as above, $\hat R^\nu
\{v^{\Phi_+}\}$ denotes the ring of restricted power series in $v^\phi ,\ \phi\in \Phi_+$ and
coefficients in $\hat R^\nu$.\par\noindent  c) Given a homomorphism $\phi\mapsto c(\phi )$
from $\Phi$ to the multiplicative group of units of $R$, it induces a surjection $\widehat{{\mathcal
A}_\nu(R)}^{\nu_\A}\to \hat R^\nu$ defined by $\sum x_\phi v^{-\phi}\mapsto \sum x_\phi
c(-\phi )$. The kernel of this surjection is the closure of the ideal $$\bigl((v^\phi -c(\phi
))_{\phi\in\Phi_+}\bigr)\widehat{{\mathcal A}_\nu(R)}^{\nu_\A}.$$
\end{proposition}
 \begin{proof}  Part a) follows from \ref{completion} and the general facts on
completion, and  Part b) follows from \ref{deformation} and the remark following definition
\ref{completloc}. Part c) follows from \ref{deformation} which contains the fact that the
corresponding result for $\A_\nu (R)$ is true. 
\end{proof}
\begin{corollary}\label{completecompo} If $\nu$ is composed with a
valuation $\nu_1$, according to a surjective monotone non-decreasing group  homomorphism
$\lambda\colon \Phi\rightarrow \Phi_1$, setting $\Psi=\hbox{\rm Ker}\lambda$, we have: 
$$\widehat{{\mathcal A}_{\nu_1}(R)}^{\nu_{1,\A}}=\widehat{{\mathcal A}_{\nu}(R)}^{\nu_\A}/
\overline{\big((v^\psi -1)_{\psi\in \Psi_+}\big)},$$ where the bar denotes the closure.
\end{corollary}
\begin{proof} This follows from Proposition \ref{completion} and the fact that in view
of \ref{equicompletion} the continuous map $\hat R^\nu
\{v^\Phi \}\to \hat R^{\nu_1} \{v^{\Phi _1}\}$ is surjective and its kernel is the ideal
$\overline{((v^\psi -1)_{\psi\in \Psi_+})}$.
\end{proof}
\par\noindent In fact this corollary is
the generalization to valuation algebras of the first corollary of Proposition
\ref{oneequicompletion}. \par\noindent  In spite of the fact that none of this applies when
$\nu_1$ is the trivial valuation, we have:
$$\widehat{{\mathcal A}_{\nu}(R)}^{\nu_\A}/\overline{\big((v^\phi -1)_{\phi\in \Phi_+}\big)}=\hat
R^\nu ,$$  because the natural map ${\mathcal A}_{\nu}(R)\longrightarrow R$ with kernel $(v^\phi
-1)_{\phi
\in \Phi_+}$ is continuous not only when $R$ has the trivial topology given by the trivial
valuation, but also when $R$ has the $\nu$-adic topology.\par\noindent
Assuming now, as we did in subsection \ref{speciatogr}, that $R$ is a $k$-algebra for some field
$k$, we have:
\begin{proposition}\label{complete flatness} {\rm *}The $k[v^{\Phi_+}]$-algebra
$\widehat {\A_\nu (R)}^{\nu_\A}$ is faithfully flat.{\rm *}
\end{proposition}
\begin{proof} This is
proved exactly like Proposition \ref{faithflat} except that we have to accept infinite sums $\sum
y_{k\psi}v^{-\psi}$ and an infinite-dimensional vector space $V$ in $\hat R^\nu$, but
everything is convergent in $\widehat {\A_\nu (R)}^{\nu_\A}$ because we deal with restricted power series. 
\end{proof} 
\par\noindent Remark that since $\nu$ is trivial on $k$, we have
$k[v^{\Phi_+}]=k\{v^{\Phi_+}\}$.\par\medskip\noindent
In order to control the specialization we need in fact to consider the scalewise completion of both $\A_\nu (R)$ and
$\hbox{\rm gr}_\nu R$, which  correspond to the scalewise completion $\hat R^{(\nu)}$ of $R$.\par\noindent
I only give the definitions and state the  properties we need.\par\noindent
Let us denote by
$$\Phi\stackrel{\lambda_1}\rightarrow\Phi_1\stackrel{\lambda_2}\rightarrow\ldots\stackrel{\lambda_{h-1}}\rightarrow\Phi_{h-1}$$
the sequence of groups of the valuations with which $\nu$ is composed. Notice that the subset
$$\Lambda=\{(\phi,\phi_1,\ldots,
\phi_{h-1})\in \Phi\times\Phi_1\times\cdots \times \Phi_{h-1}/\phi_i=\lambda_i(\phi_{i-1})\hbox{\rm for}\ 1\leq i\leq
h-1\}$$ is isomorphic to $\Phi$. Iterating the last remark before subsection
\ref{sec-special}, we can write:
$$\A_\nu (R)=
\bigoplus_{\phi_{h-1}\in
\Phi_{h-1}}\bigl(\bigoplus_{\lambda_{h-1} (\phi_{h-2})=\phi_{h-1}}\bigl(\ldots  \P_\phi
(R)w^{-\phi}\ldots\bigr)w_{h-2}^{-\phi_{h-2}}\bigr)w_{h-1}^{-\phi_{h-1}}.$$
\begin{definition}
The scalewise completion of $\A_\nu (R)$ consists in completing first with respect to the $\nu_{h-1, \A}$-adic topology, then
completing the space of coefficients of each $w^{-\phi_{h-1}}$ with respect to the $\nu_{h-2,\A}$-adic topology, and so on.
Note that each behaves like the filtration associated to a valuation of height one.  I will denote the result by 
$$\widehat{\A_\nu(R)}^{(\nu_\A)}.$$\end{definition}
I define similarly the scalewise completion of $\hbox{\rm gr}_\nu R$ after iterating the construction of Lemma \ref{bifiltr}
to obtain an isomorphism
$$\hbox{\rm gr}_\nu R\stackrel{\simeq}\rightarrow\bigoplus_{(\phi,\phi_1,\ldots ,\phi_{h-1})\in \Lambda}\bigl(\hbox{\rm gr}_{\overline\P}\hbox{\rm
gr}_{\overline\P_1(\phi_1)}\ldots (\hbox{\rm gr}_{\nu_{h-1}}R)_{\phi_{h-1}}\bigr)_\phi,$$ where $\hbox{\rm gr}_{\overline\P_i}$ is the 
graded ring corresponding to the filtration  $\overline\nu_i$ on each homogeneous component of $\hbox{\rm gr}_{\nu_{i+1}}$ induced by
$\nu_i$. Again each behaves like in the case of height one.
\begin{definition} The scalewise completion of $\hbox{\rm gr}_\nu R$ is obtained by successively completing with respect to
the filtrations $\nu_{h-1,\hbox{\rm  gr}}$, $\overline \nu_{h-2,\hbox{\rm  gr}}$, ..., $\overline\nu_{\hbox{\rm  gr}}$. It will
be denoted by $\widehat{\hbox{\rm gr}}^{(\nu)}_\nu R$.
\end{definition}
\begin{proposition} {\rm *}The $k[v^{\Phi_+}]$-algebra $\widehat{\A_{\hat\nu}(\hat R^{(\nu)})}^{(\nu_\A)} $ is faithfully flat; the
special fiber of the map
$$\hbox{\rm Spec} \widehat{\A_{\hat\nu}(\hat R^{(\nu)})}^{(\nu_\A)} \to \hbox{\rm Spec}k[v^{\Phi_+}]$$ is $\hbox{\rm
Spec}\widehat{\hbox{\rm gr}}_\nu ^{(\nu)}R$ and its ``general'' fiber, over $\hbox{\rm Spec}(v^{\Phi_+})^{-1}k[v^{\Phi_+}]$, is isomorphic to
$\hbox{\rm Spec}(v^{\Phi_+})^{-1}\hat R^{(\nu)}[v^{\Phi_+}]$.{\rm *}\end{proposition}
\par\medskip\noindent In any case, we shall see below why it is for the family $$\hbox{\rm
Spec}\widehat{\A_{\hat\nu}(\hat R^{(\nu)})}^{(\nu_\A)}\rightarrow\hbox{\rm Spec}k[v^{\Phi_+}]$$
that we may indeed hope to have simultaneous \textit{embedded} uniformization; the
$k[v^{\Phi_+}]$-algebra $\widehat{\A_{\hat\nu}(\hat R^{(\nu)})}^{(\nu_\A)}$ is a quotient of the scalewise completion  of
$k[v^{\Phi_+}][(w_j)_{j\in J}]$, which corresponds to a simultaneous embedding of its fibers. This
will put us in a position to apply the implicit function theorem to extend a partial uniformization of
$\nu_{\rm gr}$ on $\hbox{\rm gr}_\nu R$ to a uniformization of $\hat\nu$ on $\hat R^{(\nu )}$. 
\subsection{The specialization of $\hat R^{(\nu )}$ to $\widehat{\hbox{\rm gr}}_\nu ^{(\nu)}R$}\label{newspec}
\par
By Proposition \ref{triviality} we have a graded isomorphism 
$$\hbox{\rm gr}_{\nu_\A}\A_\nu (R)\stackrel{\simeq}\rightarrow (\hbox{\rm gr}_\nu R)\otimes_{\overline
R}\overline R[v^{\Phi_+}].$$
If we take a system of generators $\overline \xi_i$ of the 
$\overline
R$-algebra $\hbox{\rm gr}_\nu R$,
we obtain generators $\overline\upsilon_i=\overline \xi_i
\otimes 1$ of the $\overline R[v^{\Phi_+}]$-algebra
$\hbox{\rm gr}_{\nu_\A}\A_\nu (R)= \hbox{\rm gr}_{\hat\nu_\A}\widehat{{\mathcal
A}_{\nu}(R)}^{(\nu_\A )}$
and we can lift them to homogeneous elements
$\upsilon_i\in {\mathcal A}_{\nu}(R)\subset\widehat{{\mathcal A}_{\nu}(R)}^{(\nu_\A )}$.\par\noindent
Note that according to Proposition \ref{triviality}, each element $\upsilon_i$ is of the form
$\xi_iv^{-\nu (\xi_i)}$ where $\xi_i\in R$ has initial form $\overline \xi_i$.\par
\noindent 
Let us begin by assuming that $R$ is local and has a field of representatives $k$.
Then we may consider the morphism of $k[v^{\Phi_+}]$-algebras 
$$k[v^{\Phi_+}]{[(u_i)_{i\in I}]\to {\mathcal
A}_{\nu}(R)}$$ determined by $u_i\mapsto \upsilon_i=\xi_iv^{-\nu (\xi_i)}$. We can consider
the \textit{scalewise completion} of $k[v^{\Phi_+}][(u_i)_{i\in I}]$ defined as follows: 
define, with the notations
introduced in Subsection \ref{morestruc}, the sequence of rings 
$$\hat S_1=k[v^{\Phi_+}][(u_i)_{i\in I_1}]^{\widehat {}\ \tilde\nu_1},\ldots ,
\hat S_t=\hat
S_{t-1}[(u_i)_{i\in I_t}]^{\widehat {}
\ \tilde\nu_t},\ldots
,\hat S_h=\hat S_{h-1}[(u_i)_{i\in I_h}]^{\widehat {}\ \tilde\nu_h}$$ where $\tilde\nu_t$ is the monomial valuation on $\hat S_{t-1}[(u_i)_{i\in I_t}]$
whose values lie in \linebreak $\Psi_{h-t}/\Psi_{h-t+1}$ defined by 
$\tilde\nu_t (\sum a_\alpha u^\alpha )=\hbox{\rm
min}(\nu_t (\xi^\alpha ))$, and where the hat means the completion. 
\begin{definition} We can inductively define a valuation $\tilde\nu$ on $\hat S_h$ as follows: it
is trivial on $k[v^{\Phi_+}]$ and, assuming that it has been defined on $\hat S_{t-1}$, for an
element $\sum  a_\alpha u^\alpha\in \hat S_t$ we set
$$\tilde\nu ( \sum  a_\alpha u^\alpha)=\hbox{\rm min}_\alpha\{ \tilde\nu
(a_\alpha)+\nu (\xi^\alpha )\}\in \Psi_{h-t},$$ 
which is well defined because the semigroup $\Gamma=\nu (R\setminus\{ 0\})$ is well
ordered. It extends the monomial valuation on the polynomial ring with the $(u_i)_{i\in I}$ as
system of generators and $\tilde\nu (u_i)=\nu(\xi_i)=\gamma_i$, the $i$-th generator of the semigroup $\Gamma$.\end{definition} 
Note that $\tilde\nu$ is {\it not} a rational valuation in  general.\par\medskip\noindent We can apply the same procedure to $k[(u_i)_{i\in
I}]$ instead of
$k[v^{\Phi_+}][(u_i)_{i\in I}]$:
\begin{definition} The ring $\widehat{k[(u_i)_{i\in I}]}$ obtained from the polynomial ring by the process just described is called the {\it
scalewise completion} of the polynomial ring $k[(u_i)_{i\in I}]$ .\end{definition}
\begin{proposition}\label{analmonom}The scalewise completion $\widehat{k[(u_i)_{i\in I}]}$ is a local ring and its maximal ideal is the
closure of the ideal generated by the $(u_i)_{i\in I}$. The valuation $\tilde\nu$ is rational in this case and analytically monomial with respect to
$(u_i)_{i\in I}$.\end{proposition}
\begin{proof} The proof of the first two statements is the same as that of Proposition \ref{allocal}, and the last one follows from
the very definition of $\tilde\nu$. The graded algebra with respect to $\tilde\nu$ is given by
$$\hbox{\rm gr}_{\tilde\nu}\widehat{k[(u_i)_{i\in I}]}=k[(U_i)_{i\in I}].$$\end{proof}
\begin{Remark}1) If the $k$-algebra $\hbox{\rm gr}_\nu R$ is finitely generated, say by $(\overline \xi_i)_{i\in F}$, the ring
$\widehat{k[(u_i)_{i\in I}]}$ is the usual power series ring $k[[(u_i)_{i\in F}]]$, and $\tilde\nu$ is the usual monomial valuation
determined by $\tilde\nu(u_i)=\nu (\xi_i)=\gamma_i$, the $i$-th generator of the semigroup $\Gamma$.\par\noindent
2) In the general case, but with $\Phi$ of height one say, $\widehat{k[(u_i)_{i\in I}]}$ contains elements such as $u_1+u_2+\cdots
+u_i+\cdots$.\end{Remark}
\begin{proposition}\label{scacomp} The ring $\hat S_h$ is complete for the topology defined by $\tilde\nu$.\end{proposition}
\begin{proof} Again this is proved by induction on $t$. For $t=1$ is follows from the definitions. Assume the proposition true up to $t-1$ and
let $\sum_\alpha a_\alpha^{(i)}u^\alpha\in \hat S_t$ be a Cauchy sequence. The minimum over $\alpha$ of the $\tilde\nu$ valuations of the
differences
$(a_\alpha^{(i)}-a_\alpha^{(j)})u^\alpha$ must become arbitrarily large in $\Psi_{h-t}$ as $\hbox{\rm min}(i,j)$ grows. For each $\alpha$ the
$a_\alpha^{(i)}$ must form a Cauchy sequence for the $\tilde\nu$ topology, which converges by the induction hypothesis.\end{proof}
\begin{definition} We will denote by $\widehat{k[v^{\Phi_+}][(u_i)_{i\in I}]}$ the ring $\hat S_h$ endowed
with the topology determined by the $\tilde\nu$-adic filtration, and call it the
{\it scalewise completion} of the ring $k[v^{\Phi_+}][(u_i)_{i\in I}]$. \end{definition}
\noindent Assume that the local ring $R$ is complete and let $k\subset R$ be a field of representatives. The map
$k[(u_i)_{i\in I}]\rightarrow R$ defined by $u_i\mapsto \xi_iv^{-\nu (\xi_i)}$
induces a map $$k[v^{\Phi_+}][(u_i)_{i\in I}]\rightarrow \A_\nu (R);\ \ \ 
v^{-\phi}\mapsto v^{-\phi},\ \ u_i
\mapsto \xi_iv^{-\nu (\xi_i)},$$
which in turn induces a surjective map of
$k[v^{\Phi_+}]$-algebras
$$\widehat{k[v^{\Phi_+}][(u_i)_{i\in I}]}\longrightarrow 
\widehat{\A_{\nu}(R)}^{(\nu_\A)},$$ 
where the first hat designates the scalewise completion
defined above. The surjectivity follows from Proposition \ref{denseinscales} or rather
from repeating its proof in this situation where $k$ is replaced by $k[v^{\Phi_+}]$. This map
induces, modulo the {\it closures} of the ideals  generated by $(v^{\Phi_+})$, the surjective
map  $$\widehat{k[(U_i)_{i\in I}]}\to\widehat{\hbox{\rm gr}}^{(\nu)}_\nu R $$
which extends to the completions the original map
$U_i\mapsto
\overline\xi_i$.\par\noindent We can summarize this as follows
\begin{proposition}\label{coordinates}{\rm *} Let $R$ be a complete equicharacteristic
n\oe therian local integral domain with residue field $k$, and $\nu$ 
 a rational valuation of $R$ with value group $\Phi$. Let $k\subset R$ be a field of
representatives for $R$ and $(\xi_i)_{i\in I}$ be representatives in
$R$ of a system of homogeneous generators
$(\overline\xi_i)_{i\in I}$ of the $k$-algebra $\hbox{\rm gr}_\nu R$.
Then the map of $k[v^{\Phi_+}]$-algebras
$k[v^{\Phi_+}][(u_i)_{i\in
I}]\rightarrow \A_{ \nu}(R)$ determined by the injection
$k[v^{\Phi_+}]\subset
\A_{\nu}(R)$ deduced from $k\subset R$ and the applications
$u_i\mapsto\xi_iv^{-\nu (\xi_i)}$ extends to a continuous surjective map of topological
$k[v^{\Phi_+}]$-algebras
$$\widehat{k[v^{\Phi_+}][(u_i)_{i\in I}]}\rightarrow\widehat{
{\mathcal A}_{\nu}(R)}^{(\nu_\A)},$$
which induces, modulo the closures in each ring  of the ideal $(v^{\Phi_+})$, the map
$$\widehat{k[(U_i)_{i\in I}]}\rightarrow\widehat{\hbox{\rm gr}}^{(\nu )}_\nu
R.\hbox{\rm *}$$
\end{proposition}\hfill $\square$\par
 (I write capital $U$'s for the initial forms of the $u$'s).\par
 The kernel of this last map is the closure in $\widehat{k[(U_i)_{i\in I}]}$ of
the binomial ideal
$(U^m-\lambda_{mn}U^n)k[(U_i)_{i\in I}]$ defining $\hbox{\rm gr}_\nu R$; it is the 
ideal
consisting of possibly infinite sums $\sum_{m,n} A_{m,n}(U)(U^m-\lambda_{mn}U^n)\in
\widehat{k[(U_i)_{i\in I}]}$, with
$A_{m,n}(U)\in\widehat{k[(U_i)_{i\in I}]}$.\par\noindent Using the fact that ${\mathcal A}_{\nu}(R)$ is a
flat
$k[v^{\Phi_+}]$-algebra, one can show 
\begin{proposition}\label{completeequations}{\rm *} a) In the situation of the preceding Proposition, there exist
elements of the form 
$$F_{mn}=u^m-\lambda_{mn}u^n+\sum_sc^{(mn)}_s(v^\phi)u^s\in
\widehat{k[v^{\Phi_+}][(u_i)_{i\in I}]},$$
 with, for all $s$ appearing in $F_{mn}$, $\tilde\nu
(u^s)>\tilde\nu (u^n)=\tilde\nu (u^m)$ and $c^{(mn)}_s$ a term (=constant times a monomial) in
the $v^\phi$, and such that the closure in $\widehat{k[v^{\Phi_+}][(u_i)_{i\in I}] }$ of the ideal
generated by the $F_{mn}$ is the kernel of the surjection $$\widehat{k[v^{\Phi_+}][(u_i)_{i\in
I}]}\to \widehat{\A_{\nu}(R)}^{(\nu_\A)}.$$
b) The $F_{mn}$ may be chosen so that, with the notations of  Proposition \ref{allocal}, if\linebreak $U^m-\lambda_{mn}U^n\in k[(U_i)_{i\in
\Psi_{h-t}}]$, then $F_{mn}\in \widehat{k[v^{\Phi_+}][(u_i)_{i\in\Psi_{h-t}}]}$.{\rm *}
\end{proposition}
 \begin{proof} By the remark following Proposition
\ref{faithflat}, and with its notations, after replacing $R$ by $\hat R^{(\nu)}$, the
$k[v^{\Phi_+}]$-algebra ${\mathcal B}(S)$, which is the image of the map of
$k[v^{\Phi_+}]$-algebras
$$k[v^{\Phi_+}][(u_i)_{i\in I}]\to \A_\nu (R),\ \ \ u_i\mapsto \xi_iv^{-\nu (\xi_i)}$$ 
is  faithfully
flat. It follows that, denoting by $J$ the kernel of this map, we can lift the generators
$U^m-\lambda_{mn}U^n$ of the kernel of the map
$k[(U_i)_{i\in I}]\to \hbox{\rm gr}_\nu R$ to
elements $F_{mn}\in k[v^{\Phi_+}][(u_i)_{i\in I}] $ generating an ideal $J' \subseteq J$ 
such that
$J\subseteq J'+(v^{\Phi_+}).J$. The equality $(v^{\Phi_+})\cap J =(v^{\Phi_+}).J$ holds by 
flatness
since it is equivalent to $$\hbox{\rm Tor}_1^{k[v^{\Phi_+}][(u_i)_{i\in I}]}(k[v^{\Phi_+}]
[(u_i)_{i\in
I}]/(v^{\Phi_+}),k[v^{\Phi_+}][(u_i)_{i\in I}]/J)=0.$$ 
It means that we can write in $\A_\nu (R)$
generators of $J$ which begin with
$$u^m-\lambda_{mn}u^n+\sum_sc_s(v^\phi )u^s+\sum_i v^{\phi_i}A_i
(u^{m^i}-\lambda_{m^in^i}u^{n^i}+\cdots ).$$  Now we
observe that since this sum must give zero in $\widehat{\A_{\nu}(R)}^{(\nu_\A)}$, the
$w$-order of $\sum_i v^{\phi_i}A_i (u^{m^i}-\lambda_{m^in^i}u^{n^i}+\cdots )$ 
must be
greater than that of $u^m$.  Following through the process of scalewise completion,
one sees that
the closure in $\widehat{k[v^{\Phi_+}][(u_i)_{i\in I}] }$ of the ideal $J'$ is the kernel of the
surjection $$\widehat{k[v^{\Phi_+}][(u_i)_{i\in I}]}\to \widehat{\A_{\nu}(R)}^{(\nu_\A)}.$$
To prove part b), one needs only apply this argument to each $\hat R_t$.\end{proof}
Before the next corollary, let us record the
\begin{definition} Assertion TC(d) (Toric coordinatization in dimension $\leq d$) is the following: For every excellent equicharacteristic local domain
$R$ of dimension $\leq d$, given a rational valuation $\nu$, there exists a quotient $\hat R^{(\nu)}$ of $\hat R^m$ by an ideal $H$ such that
$H\cap R=(0)$ and an extension $\hat\nu$ of $\nu$ to $\hat R^{(\nu)}$ inducing a scalewise birational map $\hbox{\rm gr}_\nu R\to \hbox{\rm
gr}_{\hat\nu}\hat R^{(\nu)}$ of graded algebras such that, given a field of representatives
$k\subset \hat R^{(\nu)}$ of the residue field of
$R$ and representatives $(\eta_j)_{j\in J}$ in $\hat R^{(\nu)}$ of a system of generators of the $k$-algebra
$\hbox{\rm gr}_{\hat\nu} \hat R^{(\nu)}$, the map of $k$-algebras $k[(w_j)_{j\in J}]\rightarrow \hat R^{(\nu)}$ determined by $w_j\mapsto\eta_j$
extends to a continuous surjection for the $\tilde \nu$ and $\hat\nu$-adic topologies 
$$\widehat{k[(w_j)_{j\in J}]}\to \hat R^{(\nu)}$$ with a surjective associated graded map $$ k[(W_j)_{j\in J}]\to \hbox{\rm gr}_{\hat\nu}
\hat R^{(\nu)}, $$ where the ring $\widehat{k[(w_j)_{j\in J}]}$ is the scalewise completion of the polynomial ring (see the definition before Proposition
\ref{analmonom}). We denote by $(W^m-\lambda_{mn}W^n)_{(m,n)\in\hat E}$ a system of generators for the kernel of the map $k[(W_j)_{j\in J}]\to
\hbox{\rm gr}_{\hat\nu}
\hat R^{(\nu)}$.
\end{definition}\par\noindent
\begin{corollary}\label{TC(*)}{\rm Toric coordinatization:} {\rm *}{\rm(assuming that TLU(d-1) holds)}\par\noindent
Given a rational valuation $\nu$ of the excellent local domain $R$ of dimension $\leq d$, let $c$ be a character
of
$\Phi$ i.e., an homomorphism $c\colon\Phi\to k^*$;  replacing $R$ by  $ \hat R^{(\nu )}$ and 
$v^\phi$ by $c(\phi )$ in \ref{completeequations} gives, in view of Propositions \ref{complete deformation}, c) and \ref{nucomplete} a continuous
surjection as in TC(d)
$$\hat c\colon \widehat{k[(w_j)_{j\in J}]}
\to \hat R^{(\nu )} ,\ \ \ w_j\mapsto c(-\gamma_i)\eta_j$$ 
where the hat means scalewise completion.
The kernel of
this map is the closure of the ideal generated by the elements
$$c(G_{mn})=w^m-\lambda_{mn}w^n+\sum_sc^{(mn)}_s(c(\phi ))w^s\ \ (m,n)\in\hat E.$$
In
particular, taking the trivial character $c_1$ determined by $c_1(\phi )=1\ \ \forall\phi\in\Phi$, we have a surjective map
$$\hat c_1\colon \widehat{k[(w_j)_{j\in J}]}
\to \hat R^{(\nu )}\ \ \ \ w_j\mapsto\eta_j$$ whose kernel is generated up to closure by elements
$$c_1(G_{mn})=w^m-\lambda_{mn}w^n+\sum_sc^{(mn)}_s(1)w^s.\hbox{\rm *}$$
\end{corollary}
As I have already noted, this may be thought of as a valuative analogue of Cohen's theorem stating that
every complete
equicharacteristic n\oe therian local ring is a quotient of a  power series ring 
in finitely many variables over its residue field, remembering that the associated graded map must then be surjective. 
The advantage here is that, when the residue field is algebraically closed, we know that the
kernel of the map is generated, up to closure,
 by deformations of binomials.\par\noindent One may summarize a part of my approach by
 saying that the elements $\eta_j\in \hat R^{(\nu )}$ representing generators $\overline \eta_j$ of the 
$k$-algebra $\hbox{\rm gr}_{\hat\nu}\hat R^{(\nu )}$ form a system of generators of the maximal
ideal of  $\hat R^{(\nu)}$ in which the valuation $ \hat\nu$ is analytically monomial in the sense of section
\ref{monval}. The elements $c(G_{mn})$ may be thought of as an infinite standard, or
Gr\"obner, basis for the kernel of the map
$\widehat{k[(w_j)_{j\in J}]}\to \hat R^{(\nu)}$ with respect to the the monomial order on the -countably many- variables $w_j$ defined by deciding that 
$w^s<w^t$ if the valuation of the
image of $w^s$ in $\hat R^{(\nu )}$ is less than the valuation of the image of $w^t$.
\begin{corollary}\label{present}  a) With the notations of subsection \ref{nuA}, each ideal\linebreak ${\mathcal S}_\delta(\A_{\hat\nu} (\hat
R^{(\nu)}))$ is generated as a
$k[v^{\Phi_+}]$-module, up to closure, by monomials in the
$\eta_jv^{-\hat\nu (\eta_j)}$.
\par\noindent b) Each ideal $\P_\phi (\hat R^{(\nu)})$ is generated as a $k$-vector space, up to closure, by monomials in the $\eta_j$.
\par\noindent c) If the
$k$-algebra $\hbox{\rm gr}_{\hat\nu}\hat R^{(\nu )}$ is finitely generated,  say by $(\overline \eta_j)_{j\in
F}$, the $k[v^{\Phi_+}]$-algebra 
$\widehat{\A_{\hat\nu}(\hat R^{(\nu)})}^{(\hat\nu_\A)}$ is a quotient of $k[v^{\Phi_+}][[(w_j)_{j\in F}]]$.\end{corollary}
\begin{proof}: a) is a direct consequence of \ref{coordinates} and Proposition
\ref{analmonom}. Statement b) is a consequence of a) by evaluation via a character of $\Phi$ in $k^*$. For c), one checks by
induction on the height that in the finitely generated case, the scalewise completion of
$k[v^{\Phi_+}][(w_j)_{j\in F}]$ coincides with
$k[v^{\Phi_+}][[(w_j)_{j\in F}]]$.\end{proof}
\begin{corollary}\label{relatfin} Let $R$ be a complete n\oe therian local ring endowed with a rational valuation $\nu$ and let $\mathbf p$ be the center
of  the valuation $\nu_1$ of height one with which $\nu$ is composed. For each $\phi_1\in \Phi_{1+}$ there are finitely many elements $\gamma_i$ of our
minimal system of generators for the semigroup $\Gamma$ of $\nu$ which are the valuations of elements of $\P_{\phi_1}\setminus \P^+_{\phi_1}$. 
\end{corollary}
\begin{proof} The $R/{\mathbf p}$-module $\P_{\phi_1}/\P^+_{\phi_1}$ is finitely generated. Each of the generator can be written as a series in the
$\xi^\alpha$ with coefficients in $k$. Each of the monomials in the $\xi_i$ which appear in this way can be written ${\xi'}^{\alpha'}{\xi''}^{\alpha''}$ with
$\nu_1({\xi'}^{\alpha'})=0$ and $\nu_1({\xi''}^{\alpha''})\geq\phi_1$. All the monomials ${\xi''}^{\alpha''}$ which are of valuation $\phi_1$ generate a
$R/{\mathbf p}$-submodule of $\P_{\phi_1}/\P^+_{\phi_1}$. Since $R/{\mathbf p}$ is n\oe therian, this submodule is generated by finitely many of these
monomials, and the $\nu$-valuation of every element of $\P_{\phi_1}\setminus
\P^+_{\phi_1}$ is a linear combination with coefficients in $\N\cup\{0\}$ of elements of the semigroup $\overline \Gamma$ of the residual valuation
$\overline\nu$ and the valuations of these monomials. \end{proof}
\begin{corollary}\label{finpredes} Let $R$ be a complete n\oe therian local ring endowed with a rational valuation $\nu$. In the minimal system of
generators of the semigroup $\Gamma$ of $(R,\nu)$ there are only finitely many elements without predecessor.\end{corollary}
\begin{proof} The proof is by induction on the height of $\nu$. The result is true if $\nu$ is of height one because in this case every element has
a predecessor by the second part of Corollary \ref{gammai}. Keeping the same notations as in Corollary \ref{relatfin} and assuming the result to be true for
the valuation $\overline\nu$ of $R/{\mathbf p}$, we see that if $\gamma_i$ has no predecessor and $\nu_1(\gamma_i)=\phi_1>0$, the only possibility is
that there is an infinite sequence $\gamma_k<\gamma_{k+1}<\ldots$ of our generators of $\Gamma$ whose $\nu_1$-valuation is the predecessor of
$\phi_1$ in the image $\Gamma_1$ of $\Gamma$ in $\Phi_1$, which exists since $\nu_1$ is of height one. By Corollary \ref{relatfin} this is impossible
unless the predecessor of $\phi_1$ is zero and therefore $\gamma_i$ has to be the smallest element of $\Gamma$ which is not in $\Gamma\cap\Psi$ where
$\Psi$ is the convex subgroup corresponding to $\nu_1$. The number of elements of the well ordered set $(\gamma_1,\ldots, \gamma_i,\ldots)$ without
predecessor in the set is therefore at most equal to the height of $\nu$. 
\end{proof}
\begin{Remark} 1) Part b) of the Corollary means that the valuation $\hat\nu$ on $\hat R^{(\nu)}$ is analytically monomial in the
variables $(w_j)_{j\in J}$.\par\noindent 2) Assume that the ring $R$ contains a field of representatives. In spite of the flatness of the
$k[v^{\Phi_+}]$-algebra
$\A_\nu (R)$, it may happen that the
$k$-algebra
$\hbox{\rm gr}_\nu R$ is finitely generated while the $k[v^{\Phi_+}]$-algebra $\A_\nu (R)$ is not; if we go back to example \ref{excomp} of
subsection
\ref{nuadiccomp}, we can see that for each integer $n$, the ideal $\P_n(R)$ is generated by $(u_2^n, u_1-\sum_1^{n-1}a_iu_2^i)$, while
$\P_n(\hat R^{(\nu)})=u_2^nk[[u_2]]$ and
$\hbox{\rm gr}_\nu R=\hbox{\rm gr}_{\hat\nu}\hat R^{(\nu)} =k[U_2]$.\par\noindent 3) It would be interesting to relate the fact that $\hbox{\rm
gr}_{\hat\nu}\hat R^{(\nu )}$ is finitely generated with the same condition for $\hbox{\rm gr}_\nu R$.\end{Remark}
\subsection{The abyssal phenomenon}\label{abyssalph}
Let us keep the notations of Corollary \ref{TC(*)}. The key fact now
is that $\hat R^{(\nu)}$ is
n\oe therian and complete for the $m$-adic topology. The  images by $\hat c$ of the $w_j$ generate the maximal ideal of $\hat
R^{(\nu)}$, and therefore there exists a finite subset $L\subset J$, which in the case of a valuation of height one we may assume to
consist of
$\{1,2,\ldots ,q\}$, such that the images of
$(w_\ell)_{\ell\in L}$ generate the maximal ideal of $\hat R^{(\nu)}$. This means that for
$i\notin L$ and any character $c\colon \Phi_+\to k^*$, we have a relation
$$w_j-\sum_{\ell\in L} R_\ell w_\ell\in
\overline{\big(\{w^m-\lambda_{mn}w^n+\sum_s c_s^{(mn)}(c(\phi ))w^s\}_{(m,n)}\big)}\ \ \hbox{\rm in}\
\widehat {k[(w_j)_{j\in J}]}, $$  and this implies that one of the equations
$w^m-\lambda_{mn}w^n+\sum_s c_s^{(mn)}(c(\phi ))w^s$ must contain $w_j$ linearly (compare
with example \ref{planebranches}). According to Remark \ref{truebin}, this $w_j$ does not
appear in the binomial term of $G_{mn}$. \par
The abyssal phenomenon is that, whenever we have among the equations defining $\hat R^{(\nu)}$ as a quotient of $\widehat {k[(w_j)_{j\in J}]}$ 
infinitely many consecutive indices $j$ for which such an equation exists, these equations cannot create singularities, and at the same time they produce no
decrease of the Krull dimension of $\hat R^{(\nu)}$ (compare with example \ref{exzar}). The creation of singularities,
as well as the decrease in dimension, are pushed away to infinity, since there is ''no equation'', but an endless sequence of substitutions. The proof of the
abyssal phenomenon relies on an infinite dimensional implicit function theorem which is \textit{not} proved here. It states that, using the fact that
$\hat R^{(\nu)}$ is henselian, we may solve this equation for $w_j$ in
$\hat R^{(\nu)}$. 
\par\noindent More precisely, what we need has the following form, where each variable $w_j$ has weight $\nu(\eta_j)\in \Phi$:
\begin{theorem}[Implicit functions Theorem ]* Let $(G_h)_{h\in H}$ be elements in $\widehat{k[((w_j)_{j\in J}]}$ of the form
$$G_h=w^{n(h)}-\lambda_hw^{m(h)}+c_hw_h+\sum_sc_sw^s,$$  indexed by a subset $H\subset J$, with $c_h\in k^*,\ c_s\in k$ and the weight of
$w_h$ and all the $w^s$ is greater than the weight of the terms of the initial isobaric binomial .\par\noindent
Then, these series generate, up to closure, the same ideal as the following:
$$c_hw_h-(w^{n(h)}-\lambda_hw^{m(h)}+\sum_td_tw^t)$$ where the weight of all the monomials $w^t$ is greater than the weight
of the initial binomial and no variable $w_k$ for $k\in H$ appears in the monomials $w^t$.*\end{theorem} 
\begin{proof} This requires an infinite-dimensional implicit function theorem in\break $\widehat{k[((w_j)_{j\in J}]}$, somewhat similar to that of [VdH],
see also ([B1], Chap. IV), to use the fact that all the $\partial_{u_{h+1}}c(G_h)$ are invertible in 
$\widehat{k[((w_j)_{j\in J}]}$. This implicit function theorem should also allow us to recognize from the partial derivatives of the $G_h$ when a n\oe
therian quotient of $\widehat{k[((w_j)_{j\in J}]}$ by the {\it closure} of the ideal generated  by elements $(G_h)_{h\in H}$ is regular. \par
Because of part b) of Proposition \ref{completeequations} and Lemma \ref{archi+}, the implicit function theorem can be proved by induction on the height,
and at each step we have only a situation of height one type.
\end{proof}
Let us now consider a finite set of representatives $(\eta_j)_{j\in F}\in \hat R^{(\nu)}$ of generators $\overline\eta_j$ of the $k$-algebra $\hbox{\rm
gr}_{\hat\nu}\hat R^{(\nu)}$ which has the following properties (by definition, as in Corollary \ref{gammai} of subsection \ref{compofval}, the
ordinal $i+1$ is the index of the generator
$\hat\gamma_{i+1}$ of the semigroup $\hat\Gamma$ following $\hat\gamma_i$):\par\noindent a) the representatives $(\eta_j)_{j\in F}$ form a system
of generators of the maximal ideal of
$\hat R^{(\nu)}$,\par\noindent b) their valuations rationally generate the valuation group $\Phi$ of $\nu$,\par\noindent
c)With the usual notation for convex subgroups of $\Phi$, whenever the set $J_t=\{j\in J\vert \nu(\eta_j)\in \Psi_{h-t}\setminus\Psi_{h-t+1}\}$ is finite, it
is contained in
$F$.\par\noindent 
d) For any $k\notin F$, the valuation $\hat\gamma_k$ of $\eta_k$ has an immediate predecessor among the $\hat\gamma_i$ which is rationally
dependent upon the $(\hat\gamma_j)_{j\in F}$.
\par\noindent  e) If $i\in J$ satisfies $i\notin F$, then $i+1\notin F$.
\begin{definition}\label{defsuff} Such sets of generators of the maximal ideal of $\hat R^{(\nu)}$ will be called \textit {sufficient sets of
generators}.\end{definition}
\begin{lemma}\label{Peano} Let $\nu$ be a valuation of a n\oe therian local ring $R$, and $\Gamma=\langle
\gamma_1,\ldots ,\gamma_i,\ldots\rangle$ its semigroup of values. Let $\Delta=(\delta_1<\ldots <\delta_j<\delta_{j+1}<\ldots)$ be an
increasing sequence of elements of $\Gamma$. Given a finite subset
$F$ of $\Delta$, there is a finite subset $\tilde F$ of $\Delta$ containing $F$ and such that $\delta_j\notin \tilde F$ implies $\delta_{j+1}\notin \tilde F$.
\end{lemma}
\begin{proof}  If $\nu$ is of height one in $R$, it
suffices to take for $\tilde F$ the set of elements of
$\Delta$ smaller than or equal to the largest element of $F$; it is a finite set. Assume that the result is true for valuations
of height $\leq h-1$, and let $\nu$ be a valuation of height $h$. Let $\lambda\colon\Phi\to\Phi_1$ be the map of groups
corresponding to the valuation $\nu_1$ of height $h-1$ with which $\nu$ is composed.  Set $\Delta_1=\lambda (\Delta)$ (as a set) and
$F_1=\lambda (F)$; by the induction hypothesis, we have a finite set $\tilde F_1$ containing $F_1$ and with the property of the Proposition with
respect to $\Delta_1$. Define $\tilde F$ as follows: it is the union of finite subsets $\tilde F_{\phi_1}$ of the $\lambda^{-1}(\phi_1)$ for $\phi_1\in \tilde
F_1$. If
$\Delta\cap\lambda^{-1}(\phi_1)$ is finite (which is always the case if the center of
$\nu_1$ is equal to the center of
$\nu$, by [Z-S], Appendix 3, Corollary to Lemma 4, or Proposition \ref{almostone}), set $\tilde F_{\phi_1}=\Delta\cap\lambda^{-1}(\phi_1)$. If
$\lambda^{-1}(\phi_1)\cap \Delta$ is infinite, define $\tilde F_{\phi_1}$ to be the set of elements of $\lambda^{-1}(\phi_1)\cap
\Delta$ which are smaller than or equal to the largest element of $F$ contained in $\lambda^{-1}(\phi_1)\cap \Delta$. By {\it loc. cit.},
Lemma 4, or Proposition \ref{almostone}, each of these sets if finite. If $\delta_j\notin \tilde F$, either $\lambda (\delta_j)\notin \tilde F_1$, and then
$\lambda (\delta_{j+1})\notin\tilde F_1$ so that $\delta_{j+1}\notin \tilde F$, or $\lambda (\delta_j)=\phi_1\in \tilde F_1$ and $\lambda^{-1}(\phi_1)\cap
\Delta$ is infinite, so that $\lambda (\delta_{j+1})=\phi_1$ and $\delta_{j+1}\notin \tilde F$ by definition of $\tilde F_{\phi_1}$.
\end{proof}
\begin{proposition}\label{suffig} Sufficient sets of generators for the maximal ideal of $\hat R^{(\nu)}$ exist.\end{proposition}
\begin{proof} That condition a) can be satisfied follows from the n\oe therianity of $\hat R^{(\nu)}$. For condition b) it is the finiteness of the rational rank
$\hbox{\rm r}(\nu)$, and for condition c) the finiteness of the height $\hbox{\rm h}(\nu)$, which shows that there are finitely many sets $J_t$. That
condition d) can be satisfied follows from the finiteness of the number of generators of $\Gamma$ without a predecessor; see Corollary \ref{finpredes}. The
only thing left to check is that e) can be satisfied; it is a consequence of the preceding Lemma.
\end{proof}
\par\medskip\noindent  Let
$$\hat\Gamma=\langle \hat\gamma_1,\hat\gamma_2,\ldots ,\hat\gamma_i,\ldots\rangle\subset
\Phi_+$$ be the semigroup of the values taken on $\hat R^{(\nu )}\setminus\{0\}$ by the valuation $\hat\nu$. For each $i\geq 1$, set 
$\hat\Gamma_i=\langle
\hat\gamma_1,\hat\gamma_2,\ldots ,\hat\gamma_i\rangle$, and similarly, denote by $\hat\Gamma_{<i}$ the semigroup
generated by elements $\hat\gamma_k$ with $\hat\gamma_k<\hat\gamma_i$; I will abbreviate this last inequality to $k<i$.\par\noindent  If the height
of
$\Phi$ is one, by Corollary \ref{gammai}, there is a smallest integer
$p_0$ such that
the group generated by
$\hat\Gamma_{p_0}$ is equal to $\Phi$. This means that for
each $i>p_0$ we have an expression $t_i\hat\gamma_i=\sum_{1\leq j\leq p_0}s^{(i)}_j\hat\gamma_j$
with $s^{(i)}_j\in \Z$. In the general case, we still have the following:
\begin{proposition}\label{ratdep} a) With the notations just introduced, there is a finite set of generators
$(\hat\gamma_j)_{j\in F}$ which rationally generates $\Phi$. For every $i\notin F$ there is a finite  relation with non
negative integral coefficients
$$n_i\hat\gamma_i +\sum_{j\in F}n^{(i)}_j\hat\gamma_j=\sum_{j\in F}s^{(i)}_j\hat\gamma_j,$$ such that $n_i > 0$ is
minimal among all such relations. In addition, each relation between the
$\hat\gamma_j$ is equivalent, modulo these relations, to a relation between the $\hat\gamma_j, \ j\in F$. If the height of $\nu$
in $R$ is one, in particular if the height of $\nu$ is one, there is an integer $p_0$ such that one may choose $F=\{1,\ldots
,p_0\}$. \par\noindent
b) Given any relation between the generators $\hat\gamma_s$ of $\hat\Gamma$, and denoting by $\hat\gamma_i$ the generator of
highest degree which appears in it, it is equivalent, modulo the relations which exist between the $(\hat\gamma_k)_{k<i}$, to
a relation of the form:
$$n_i\hat\gamma_i +\sum_{k<i}n^{(i)}_k\hat\gamma_k=\sum_{k<i}\ell^{(i)}_k\hat\gamma_k,$$ such that $n_i >0$ is
minimal among all such relations.\par\noindent
c) We may assume that no relation beween the $(\hat\gamma_k)_{k<i}$ appears in the difference
$\sum_{k<i}n^{(i)}_k\hat\gamma_k-\sum_{k<i}\ell^{(i)}_k\hat\gamma_k$, in the sense that no restriction of the sums gives a difference equal to zero.
\end{proposition} 
\begin{proof} The first two statements follow from the fact that the rational rank of $\Phi$ is finite; the cardinality of
$F$ may be taken to be $\hbox{\rm r}(\Phi)$. More precisely, we may take $\hat\gamma_1$ as the first element of $F$, then
$\hat\gamma_{i_2}$ to be the smallest element of $\hat\Gamma$ rationally independent of $\hat\gamma_1$, $\hat\gamma_{i_3}$ to be the smallest
element
$>\hat\gamma_{i_2}$ which is rationally independent of $\hat\gamma_1,\hat\gamma_{i_2}$, and so on. If $$(0)=\Psi_h\subset
\Psi_{h-1}\subset\cdots\subset \Psi_1\subset \Psi_0=\Phi$$ are the convex subgroups of $\Phi$, the first $\hbox{\rm
r}(\Psi_{h-1})$ elements of this sequence will be in $\Psi_{h-1}\cap \hat\Gamma$, the next $\hbox{\rm
r}(\Psi_{h-2})-\hbox{\rm r}(\Psi_{h-1})$ will be in $(\Psi_{h-2}\setminus \Psi_{h-1})\cap \hat\Gamma$, and so
on. In this way we build a {\it finite} sequence $ \hat\gamma_1<\hat\gamma_{i_2}<\cdots <\hat\gamma_{i_r}$ which rationally generates $\hat\Gamma$.
The notation is chosen to stress the fact that it is a subset of a generating sequence of the semigroup.\par\noindent For
$i\notin F$, if $\hat\gamma_{i_j}$ is the largest element in our sequence such that $\hat\gamma_{i_j}<\hat\gamma_i$, then $\hat\gamma_i$ is rationally
dependant on $(\hat\gamma_1,\hat\gamma_{i_2},\ldots ,\hat\gamma_{i_j})$ and we may write 
$$n_i\hat\gamma_i+\sum_{k<i}n^{(i)}_k\hat\gamma_k=\sum_{k<i}\ell^{(i)}_k\hat\gamma_k$$ a relation expressing the rational dependance of
$\hat\gamma_i$ on the $(\hat\gamma_k)_{k<i}$; we may choose such a relation with minimal $n_i$. Let
$$m_i\hat\gamma_i+\sum_{k<i}m^{(i)}_k\hat\gamma_k=\sum_{k<i}p^{(i)}_k\hat\gamma_k$$ be another relation with
non negative coefficients. By the minimality of $n_i$, we have $m_i\geq n_i$, and we may write
$m_i=qn_i+r$ with $q\geq 1$ and $0\leq r<n_i$. Adding $q\sum_{k<i}n^{(i)}_k\hat\gamma_k$ to both sides
of the second relation, and substracting to each side of the equality $q$ times the corresponding  side of the
first gives
$$r\hat\gamma_i+\sum_{k<i}m^{(i)}_k\hat\gamma_k+q\sum_{k<i}\ell^{(i)}_k\hat\gamma_k=
\sum_{k<i}p^{(i)}_k\hat\gamma_k+q\sum_{k<i}n^{(i)}_k\hat\gamma_k.$$
By definition of $n_i$, we must have $r=0$, so that we have proved that modulo the first relation, the second
one becomes a relation between $(\hat\gamma_k)_{ k<i}$. The last statement of a) follows from the fact
that in the height one case, there are only finitely many $\hat\gamma_k$ smaller than any given
element of $\hat\Gamma$.\par\noindent
Statement  c) follows from the fact that among all relations with minimal $n_i$, we may choose one such that the number of $\ell^{(i)}_k$ which
it contains is minimal.
\end{proof}
\par\medskip\noindent  It is important in the sequel to remember only the fact that for
$i\notin F$ some multiple of $\hat\gamma_i$ to which we add  a sum of the form
$\sum_{k<i}t^{(i)}_k\hat\gamma_k$ is in $\hat\Gamma_{<i}$; choosing the smallest such
multiple, we will then write: 
$$(R_i)\ \ \ \ \ \ \ \ \ \ \ \ \ \ \ \ \ \ \ \ \ \ \ \ \ 
n_i\hat\gamma_i+\sum_{k<i}n^{(i)}_k\hat\gamma_k=\sum_{k<i}\ell^{(i)}_k\hat\gamma_k .\ \ \ \ \ \ \ \ \ \ \ \ \ \ \ \ \ \ $$ The coefficients of each
such relation have no common divisor. The reason why it is important is that all relations between the $\hat\gamma_i$'s appear in this way, and the
coefficient
$n_i$ is really the smallest which can occur in a relation involving $\hat\gamma_i$ with smaller terms. If we require
that these smaller terms are actually in $F$, the exponent $n_i$ is in general larger.
\par\medskip
 With the notations of subsection \ref{approxpro}, given a finite sufficient set of generators $(w_j)_{j\in F}$ for the maximal ideal of $\hat R^{(\nu)}$, let us
choose $h$ large enough for\break
$k[x_1^{(h)},\ldots ,x_r^{(h)}]$ to contain the images in $\hbox{\rm gr}_{\hat\nu}\hat R^{(\nu)}$ of the $(W_j)_{j\in
F}$. \par Remark that the relations
$R_i$ just written provide us with the following binomial relations among those defining $\hbox{\rm gr}_{\hat\nu}\hat R^{(\nu)}$; we
have, for $i\notin F$:
$$(B_i)\ \ \ \ \ \ \ \ \ W_i^{n_i}\prod_{k<i}W_k^{n^{(i)}_k}-\lambda_i\prod_{k<i}W_k^{\ell^{(i)}_k}=0\ \ \
\hbox{\rm in}\ 
\hbox{\rm gr}_{\hat\nu}\hat R^{(\nu)}\ \ \hbox{\rm with}\ \lambda_i\in k^*.$$ \noindent
Recall that the products appearing here and below are finite.\par\noindent
Let us denote by ${\mathbf B}(F)$ the binomial equations between $(W_j)_{j\in F}$ which
define the graded algebra $\hbox{\rm gr}^{(F)}_{\hat\nu}\hat R^{(\nu)}\subset \hbox{\rm
gr}_{\hat\nu}\hat R^{(\nu)}$. The ideal ${\mathcal B}$ in $k[(W_j)_{j\in J}]$ which defines
$\hbox{\rm gr}_{\hat\nu}\hat R^{(\nu)}$ contains the ideal ${\mathcal B}'$ generated by $\bigl({\mathbf B}(F),
(B_i)_{i\notin F}\bigr)$. \par\medskip\noindent 
From now on, we take a finite subset
$(\delta_j)_{j\in F}$ of the generators of the semigroup $\hat\Gamma$ of the values which $\hat\nu$ takes on $\hat R^{(\nu)}$, such that
$(\eta_j)_{j\in F}$ is a sufficient set of generators for the maximal ideal of $\hat R^{(\nu)}$. Such a set exists by Proposition \ref{suffig}.
\begin{proposition}\label{predeceq} Let $j\notin F$ be the index of one of the variables $w_j$. By the definition of $F$ the element $j$ has a predecessor
in $J$ and we write $j=i+1$. The generator $\hat\gamma_i$ is rationally dependent upon the $(\hat\gamma_k)_{k\in F}$; there exists an equation among
the
$G_{mn}$ of  Corollary \ref{TC(*)}, and in which
$w_{i+1}$ appears linearly. The equation can be chosen to be such that its initial binomial is of the form
$W_i^{n_i}W^{n(i)}-\lambda_iW^{m(i)}$.
\end{proposition}
\begin{proof} The first observation is that in $\hat R^{(\nu)}$, by Cohen's theorem, each of the elements $\eta_\ell; \ \ell\notin F$, is a power series in
the
$(\eta_j)_{j\in F}$, say $\eta_\ell=\sum_pd^{(\ell)}_p\eta^p$. Since it involves only finitely many variables, this series makes sense in
$\widehat{k[(w_j)_{j\in J}]}$ and so we get elements $w_\ell-\sum_pd^{(\ell)}_pw^p$ which belong to the kernel $T$ of the map $\widehat{k[(w_j)_{j\in
J}]}\to\hat R^{(\nu)} $. We may, in the series $G_{mn}$ which generate this kernel up to closure, substitute freely the series $\sum_pd^{(\ell)}_pw^p$ for
every occurrence of $w_\ell$ outside of the initial binomial, without changing this kernel. This shows that we may choose the $G_{mn}$ corresponding
to binomials containing only variables with indices in $F$ in such a way that all the
variables $w_j$ which they contain have indices $j\in F$ as well. Of course the initial forms of these generators of the kernel no longer necessarily generate
the initial ideal $(W^m-\lambda_{mn}W^n)_{(m,n)\in \hat E}$. Another remark is that any system of generators, up to closure, of the ideal $T$ must contain
for each $\ell\notin F$ a generator in which the variable $w_\ell$ appears linearly; otherwise there could not be a relation such as
$\eta_\ell=\sum_pd^{(\ell)}_p\eta^p$ in $\hat R^{(\nu)}$.\par\noindent
\par\noindent Now we use again the fact that each element
$\eta_j$ for $j\notin F$ belongs in
$\hat R^{(\nu )}$ to the maximal ideal, generated by $(\eta_j)_{j\in F}$. Let now $i$ denote the smallest index of a generator of the semigroup $\hat
\Gamma$ such that $i\notin F$. By the Proposition and our choice of generators for the ideal defining
$R^{(\nu )}$, we may choose a relation among those which are such that
$w_i$ appears for the first time in their initial (binomial) part, with minimal $n_i$; it is of the form: 
$$G_i=w_i^{n_i}\prod_{k<i}w_k^{n^{(i)}_k}-\lambda_i\prod_{k<i}w_k^{\ell^{(i)}_k}+\sum_sc_s^{(i+1)}(v^\phi )w^s=0,$$
with the condition that the weight of each $w^s$ must be greater than the weight of
$w_i^{n_i}\prod_{k<i}w_k^{n^{(i)}_k}$, which is equal to that of
$\prod_{k<i}w_k^{\ell^{(i)}_k}$, so that the variable $w_{i+1}$ cannot appear in an equation $G_{mn}$ whose initial binomial contains a variable of
index
$>i$. By what we have just seen, we may assume that $w_{i+1}$ does not appear in the equations $G_{mn}$ whose initial form contains only variables
of index in $F$. We may assume that
$w_{i+1}$ appears linearly in the equation
$G_i$, since it must appear in one of the equations whose initial binomial contains
$w_i$ and other variables of index $<i$. Now for any index $i\notin F$ we may assume that
$w_{i+1}$ has not appeared in any of the equations $G_{mn}$ whose initial binomial contains variables of index $<i$ and consider similarly equations
whose initial binomial comes from the expression of $\hat\gamma_i$ in terms of $\hat\gamma_j$ with $j<i$. 
Then $w_{i+1}$ must appear linearly in at least one of these equations; we choose one where this occurs, which we call $G_i$.
\end{proof} 
Turning now to $\hat R^{(\nu)}$, we note that the result just proved provides us with a totally ordered subset $(G_i)_{i\notin F}$ of the set of equations
$G_{mn}$.\par\noindent 
\begin{proposition}\label{lineapp}{\rm(The abyssal phenomenon)} {\rm *}With the same notations, there exists a finite set $F'$ containing $F$ such that
the kernel $T$ of the map
$$\widehat{k[(w_j)_{j\in J}]}\to\hat R^{(\nu)} $$ can be generated up to closure by :\par\noindent
- a finite set of equations involving the variables $(w_j)_{j\in F'}$ and whose initial
forms are binomials generating the kernel of the natural map $k[(W_j)_{j\in F'}]\to \hbox{\rm gr}_{\hat\nu}\hat
R^{(\nu)}$, and\par\noindent - equations
$G_i=w_i^{n_i}w^{n(i)}-\lambda_iw^{m(i)}+c_{i+1}w_{i+1}+\sum_pc_pw^p\ ,\ i+1\notin F',\ c_{i+1}\in
k^*$.{\rm *}\par\noindent \end{proposition}
\begin{proof} We keep the same notations, but please remark that we passed from $i\notin F'$ to $i+1\notin F'$. To achieve this is easy:  first add to our
sufficient set $F$ the smallest element which is not in $F$. We may also
transform, without modifying the ideal, all the other equations
$G_{mn}$ whose initial binomial contains only variables of index $\leq i$ into equations in the $(w_j)_{j\in F}$. Then we add the ordinal $j$ to our set $F$
and continue with its successor. By transfinite induction on
$j\notin F$, using Proposition \ref{predeceq}, we build a set of equations 
$$G_i=w_i^{n_i}w^{n(i)}-\lambda_iw^{m(i)}+c_{i+1}w_{i+1}+\sum_pc_pw^p$$ in bijection with the variables $w_i$ such that $i+1\notin F$. Let us now consider
the ideal ${\mathbf F}$ of $\widehat{k[(w_j)_{j\in J}]}$ generated by the equations $G_{mn}$ whose initial binomial contains only variables with
index in $F$. By the implicit function theorem, the quotient of $\widehat{k[(w_j)_{j\in J}]}$ by the closure of the ideal $({\mathbf F},(G_i)_{i\notin F})$
is a quotient of
$k[[(w_j)_{j\in F}]]$ which maps onto $\hat R^{(\nu)}$. The kernel of this map is generated by the images of the equations $G_{mn}$ whose initial form
involves some
$W_j$ with $j\notin F$ and which we have not used in the construction of the $G_i$. Since the ring $k[[(w_j)_{j\in F}]]$ is n\oe therian, the images of
finitely many of the equations
$G_{mn}$ suffice to generate the kernel. We choose such a finite set and add to our set $F$ the variables involved in the initial forms of these equations, and
consider the finite set
$F'$ containing this new set and having the property that $i\notin F'$ implies $i+1\notin F'$. We consider the new ideal 
${\mathbf F}'$ obtained in this manner. The closure of the ideal generated by ${\mathbf F}'$ and the equations $(G_i)_{i\notin F}$ is now the
kernel $T$. By construction, this set of generators of $T$  has the
property that the initial forms of those among its elements which involve only variables of index
$i\in F'$ generate the prime ideal which is the kernel of the map $k[(W_j)_{j\in F'}]\to \hbox{\rm
gr}_{\hat\nu}\hat R^{(\nu)}$.
\end{proof}
This is a
fundamental fact: all but finitely many of the equations defining $\hat R^{(\nu)}$ serve only to express $w_{i+1}$, for $i+1\notin F$, in terms
of $(w_j)_{j\in F}$ and the initial binomials of the finitely many equations generate a prime ideal in a polynomial ring in finitely many variables. More
precisely, keeping the notations introduced at the beginning of this subsection, we have:
\begin{corollary}\label{genid} Let us
denote by ${\mathbf F'}$ the finite set consisting of those equations among the $G_{mn}$ whose
initial forms depend only on the $(W_j)_{j\in F'}$. Given
an homomorphism $c\colon \Phi\to k^*$, the kernel of the surjective map $$\hat c\colon\widehat
{k[(w_j)_{j\in J}]}\to \hat R^{(\nu)},\ \ \ \  w_j\mapsto \eta_jc(-\hat\gamma_j)$$ is generated, up to closure, by the ideal $\bigl(c( {\mathbf F'}),
c(G_i)_{i\notin F}\bigr)$. The initial binomials of these generators which depend only on the $(W_j)_{j\in F'}$ generate the kernel of the
map of $k$-algebras
$k[(W_j)_{j\in F'}]\to
\hbox{\rm gr}_{\hat\nu}\hat R^{(\nu)}$ determined by $W_j\mapsto \overline\eta_jc(-\hat\gamma_j)$.
\end{corollary} 
  \begin{definition} Given an excellent local domain $R$ with a rational valuation $\nu$, a finite finite subset of the index set of the minimal system of
generators of the semigroup $\hat \Gamma=\hat\nu (\hat R^{(\nu)}\setminus\{0\})$  having the property described in Proposition
\ref{lineapp} will be called a \textit {quite sufficient} set. \end{definition}
\begin{corollary}\label{semir}* With the same notations, let $F$ be a quite sufficient set of generators for the maximal ideal of $\hat R^{(\nu)}$. For each $j$
such that
$j+1\notin F$, the image $\eta_{j+1}$ in $\hat R^{(\nu)}$ of the element $w_{j+1}$ can be written as the image of a series in the
images of $(w_j)_{j\in F'}$; such a series, viewed in $\widehat{k[(w_j)_{j\in J}]}$, may
be called an $j$-th {\it semiroot} for $\hat R^{(\nu)}$; semiroots are now indexed by ordinals $<\omega^{\hbox{\rm h}_R(\nu)}$. Its expression in $\hat
R^{(\nu)}$ as 
the image of a series in $(w_k)_{k\leq j}$ is: $$\eta_{j+1} =d_{j+1}\big( \eta_j^{n_j}\prod_{k<j}\eta_k^{n^{(j)}_k}-\lambda_j
\prod_{k<j}\eta_k^{\ell^{(j)}_k}+\sum_sd_s^{(j)}\eta^s\big),$$ with
$d_{j+1}\in k^*$ and $\hat\nu (\eta^s)>\hat\nu (\prod_{k<j}\eta_k^{\ell^{(j)}_k})=\hat\nu (\eta_j^{n_j}\prod_{k<j}\eta_k^{n^{(j)}_k})$, obtained by
solving for the $(w_{j+1})_{i\notin F}$ the equations
$$\bigl(c(G_j)=0\bigr)_{j\notin F}.\hbox{\rm *}$$
\end{corollary}
This suggests a way to generalizes the theory of approximate roots, or rather of semiroots in the sense of [PP]. \par\noindent From a
valuation theoretic viewpoint the idea of approximate roots appears already in the ``key polynomials'' of MacLane (see [McL1], [McL2] and [V2]),  but in a
different guise it was developed by Abhyankar-Moh, (see [A-M1], [A-M 2]), Cossart-Moreno (in preparation), Lejeune-Jalabert, (see [L]) (for curves, the
idea being  that of branches having ``critical'' contact with a given plane curve singularity),
and related notions are studied by Spivakovsky (in [S1] for surfaces, and in [S2] for valuations of height one). I refer to [PP], especially to Corollary 1.5.4, to be
compared with Proposition \ref{coordinates}, and Corollary 1.5.5. The connection between the various viewpoints is explained in a special case at the end of
subsection \ref{chexamples}.\par\medskip\noindent\par\noindent
\begin{corollary} For $i+1\notin F$,
we have the inequality $$\gamma_{i+1}>n_i\gamma_i .$$ 
\end{corollary} 
\begin{remark}The choice of the finite set $F$ made above is simple but far from economical; one could begin by taking elements whose valuations
rationally generate $\Gamma$, in minimal number $\hbox{\rm r}(\nu)$. Then Lemma \ref{integext} implies that the elements of $\hbox{\rm gr}_\nu R$
are connected to them by binomial equations, and we can add enough such elements for their representatives to generate the maximal ideal. Whether
one of the corresponding equations $F_{mn}$ will contain a linear term then depends on the structures of $\Gamma$ and $R$. The simplest exemple is the
case of a plane branch.
\end{remark}
\subsection{An example: complex plane branches}\label{planebranches} 
\begin{example} Let $R$ be the analytic algebra of a complex analytic plane
branch 
 (i.e., germ of analytically irreducible curve)
$(X,0)$; then it is dominated by only one valuation, the $t$-adic valuation
$\nu$ induced by the injection $R\subset {\mathbf C}\{t\}$ of $R$ into its normalization, 
and the
valuation ring is ${\mathbf C}\{t\}$. Let $\Gamma=\nu (R\setminus \{0\})\subset {\mathbf N}\cup\{ 0\}$ 
be the
semigroup of values of $\nu$ on $R$; it has $g+1$ generators\footnote{ In the complex plane branch case, the generators of the semigroup are often
denoted by
$\overline\beta_i$ to underline their kinship with the Puiseux  exponents $\beta_i$, according to a tradition initiated by Zariski;
we followed that tradition in [G-T].} $(\overline{\beta_0},\ldots
,\overline{\beta_g})$ and if we choose for each $0\leq i\leq g$ an element $\xi_i(t)\in R$ 
with valuation $\overline{\beta_i}$, then the completion $\hat {\mathcal A}(R)$ of ${\mathcal
A}(R)$ with  respect to the $(t)$-adic valuation (our $\nu_\A$) (resp. the analyticization
${\mathcal A}(R)^{an}$ of ${\mathcal A}(R)$) is isomorphic to the $(t)$-adic completion $\hat{\mathcal
A}$ (resp. the  analyticization) of the sub
$\C[v]$-algebra of $R[v,v^{-1}]$ generated by the elements $\xi_i(t)v^{-\overline
{\beta_i}}$.
This follows from ([T1]) and what we have seen above.  Anyway, this is the algebra of 
a family $M$, parametrized by $v$, of branches, each fiber $M(v)$ for $v\neq 0$ being
isomorphic to the branch $(X,0)$, and the special fiber being isomorphic to the monomial
curve described  parametrically by
$u_i=t^{\overline{\beta_i}};\ \ 0\leq i\leq g$. There is a natural map which resolves its
singularities; it is the map ${\mathcal A}^{an}\rightarrow {\mathbf C}\{ v,t'\} $ given by
$$\xi_i(t)v^{-\overline{\beta_i}}\mapsto \xi_i(vt')v^{-\overline{\beta_i}}.$$ Indeed, 
the natural
resolution of singularities for $\hbox{\rm Spec}{\mathcal A}(R)$ is  $\hbox{\rm Spec}
\A ({\mathbf
C}\{ t\})$. This last ring is a genuine Rees algebra for the ideal $(t)$ in ${\mathbf C}\{t\}$
$$ \A({\mathbf C}\{ t\})= {\mathbf C}\{ t\}[v]\bigoplus\bigoplus_{n\in \N}(t)^nv^{-n},$$
and its
analytization maps isomorphically to ${\mathbf C}\{v,t'\}$ by $v\mapsto v,\ t\mapsto vt'$.
\par\noindent This shows that the parametric specialization  described in 
([T1], 1.10, p. 167) and
([G-T], \S 3) coincides with the specialization given by the valuation algebra.
 \par\noindent 
It was shown in [T1] that the binomial equations
$$u_i^{n_i}-u_0^{\ell^{(i)}_0}\ldots u_{i-1}^{\ell^{(i)}_{i-1}}=0\ \ ,\ 1\leq i\leq g,$$ 
defining the
monomial curve in $ {\mathbf C}^{g+1}$ extend to equations defining the family $M$ in 
$ {\mathbf C}\times
{\mathbf C}^{g+1}$, which are of the form 
$$
\begin{array}{lr}u_1^{n_1}-u_0^{\ell^{(1)}_0}+b_2vu_2+\sum_{w(s)>n_1\overline
\beta_1}c^{(1)}_s(v)u^s &=0\\
u_2^{n_2}-u_0^{\ell^{(2)}_0}u_1^{\ell^{(2)}_1}+b_3vu_3+\sum_{w(s)>n_2\overline
\beta_2}c^{(2)}_s(v)u^s&=0\\ \vdots &\vdots\\  
u_{g-1}^{n_{g-1}}-u_0^{\ell^{(g-1)}_0}\cdots
u_{g-2}^{\ell^{(g-1)}_{g-2}}+b_gvu_g+\sum_{w(s)>n_{g-1}\overline
\beta_{g-1}}c^{(g-1)}_s(v)u^s&=0\\ u_g^{n_g}-u_0^{\ell^{(g)}_0}\cdots
u_{g-1}^{\ell^{(g)}_{g-1}}+\sum_{w(s)>n_g\overline
\beta_g}c^{(g)}_s(v)u^s&=0 \end{array}
$$
 where the $c^{(j)}_s(v)$ are in $(v){\mathbf C}\{ v\}$,
$w(s)=\sum_0^g\overline \beta_js_j$ is the weight of the monomial $u^s$ with respect 
to the weight vector $w=(\overline\beta_0, \ldots ,\overline\beta_g)$, i.e., $w(s)=\langle w,s\rangle$, and with all $b_j\neq 0\
\hbox{\rm in}\ {\mathbf C}$, which corresponds to the fact that our branch is plane. The other 
$u^s$
appearing in the $j$-th equation are different from
$u_{j+1}$. Note that since $n_i\overline \beta_i<\overline\beta_{i+1}$ (see [T1]), we add to
each
binomial only terms of higher weight, one of which is linear except in the bottom equation, which ``creates'' the singularity.\par\noindent We see that 
we have
a surjective map $$\pi\colon \C\{v \}[[u_o,\ldots ,u_g]]\rightarrow \hat{\mathcal A}(R),$$ and its
kernel is generated by the equations just written, which illustrates Proposition
\ref{completeequations}. \par\noindent It was shown in [G-T] that a toric map $\pi\colon
Z\to{\mathbf C}^{g+1}$ resolving the monomial curve extends to a toric map $
\hbox{\rm Id}_{\mathbf
C}\times\pi \colon {\mathbf C}\times Z\to{\mathbf C}\times{\mathbf C}^{g+1}$ inducing a resolution of 
$M$
which is a strong simultaneous embedded resolution along ${\mathbf C}\times \{0\}$. In particular the map $\pi$ gives an embedded resolution for our
germ of curve $C\subset \C^{g+1}$, the embedding being given by the generators $(\xi_0,\ldots ,\xi_g)$ of the maximal ideal of $R$, or alternatively by
the equations written just above, with $v=1$.
\par\noindent
In this example, the valuation $\tilde \nu$ of subsection~\ref{sec-special} has values in the
group ${\mathbf Z}\oplus {\mathbf Z}$ ordered lexicographically, with the values
$\tilde \nu (v)=(1,0),\ \ \tilde \nu (t)=(1,1)$, this last value corresponding to $\tilde \nu
(t')=(0,1)$ on ${\mathbf C}\{v,t'\}$, and so to the natural height two valuation on this last ring
adapted to the coordinates $(v,t')$. \par
By successive elimination in the equations shown above, and thanks to the usual implicit function theorem, we can express the image of $u_2$  in the
ring $R$ of our branch in the form
$$u_2=b_2^{-1}(u_1^{n_1}-u_0^{\ell_0^{(1)}}+\sum d^{(2)}_su_0^{s_0}u_1^{s_1}),$$ and then
$$u_3=b_3^{-1}((u_1^{n_1}-u_0^{\ell_0^{(1)}}+\sum
d^{(2)}_su_0^{s_0}u_1^{s_1})^{n_2}-u_0^{\ell^{(2)}_0}u_1^{\ell^{(2)}_1}+\sum_sd^{(3)}_su_0^{s_0}u_1^{s_1}),$$ and so on. The plane branches with
equations $u_2=0,\ u_3=0,\ldots ,u_k=0\ \hbox{\rm for}\ 2\leq k \leq g$ correspond to semiroots in the sense of [PP] of the Weierstrass
polynomial defining our plane branch. Compare with Corollary \ref{semir}; the fact that no abyssal phenomenon is involved here does not
prevent the fact that all variables but two appear linearly in the equations from having its consequences.
\end{example}
\section{\textbf{Toric modifications}}\label{toricmod}  
This section begins the study of partial uniformization, by a combinatorial process, of the valuation
$\nu_{\rm gr}$ on $\hbox{\rm gr}_\nu R$ when $\nu$ is a rational valuation of the n\oe 
therian local ring $R$. I show how to resolve singularities of an irreducible variety
defined by a binomial ideal in finite-dimensional space, which is the case for $\hbox{\rm Specgr}_\nu R$ in the situation
just described if in addition $\hbox{\rm gr}_\nu R$ happens to be a finitely generated $k$-algebra, and is the case for a finitely generated approximation of
$\hbox{\rm Specgr}_\nu R$ in general. As far as I know, existing literature deals only with toric varieties associated to fans, which are in particular normal. I
begin with an example.
\subsection{Valuative motivation: an example} 
Assume now that
$k=k_R$ is algebraically closed and  that $\nu$ is a valuation of rational rank and height equal to 
$d=\hbox{\rm dim} R$ dominating $R$. By
Abhyankar's inequality, the extension $k\rightarrow k_\nu$ is trivial, so we are in the situation
of Proposition \ref{binomialideal}, and  the ordered group $\Phi$ is isomorphic
to ${\mathbf Z}^d$ with the lexicographic order (see the text before Proposition
\ref{wellord}). \par Then we can write: $$\Gamma=\langle\gamma_1,\gamma_2,\ldots
,\gamma_j,\ldots\rangle \ \hbox{\rm with}\ \gamma_j\in {\mathbf Z}^d_{\geq 0}.$$
The semigroup $\Gamma$ is the image of $
{\mathbf N}^{\mathbf N}$ by
a surjective map
$$b\colon {\mathbf Z}^{\mathbf N}\rightarrow {\mathbf Z}^d$$
defined by sending the $i$-th basis vector to $\gamma_i$. Passing to semigroup algebras over $k$ givess a map 
$$k[(U_j)_{j\in J}]\rightarrow k[t_1^{\pm 1},\ldots ,
t_d^{\pm 1}]\ \hbox{\rm determined by}\ \  U_j\mapsto t^{\gamma_j}$$
having $k[t^\Gamma]$ for
image. This describes $\hbox{\rm Spec}k[t^\Gamma]$ as the closure of an orbit of the torus 
$(k^*)^d$ in 
the ${\mathbf N}$-dimensional affine space ${\mathbf A}^{{\mathbf N}}(k)$.\par\noindent  Note that since 
$R$ and
$R_\nu$ have the same field of fractions, the semigroup $\Gamma$ generates $\Phi$ 
as a group. It will not be the case in general for a valuation with value group $\Phi$ that there is a finite index $i_0$ such that the finitely
generated semigroup
$\langle\gamma_1,\ldots ,\gamma_{i_0}\rangle$ generates $\Z^d$ as a group. Note also that even in this special case
where $\Phi=\Z^d$ we have when $d$ is $>1$ a difficulty which does not appear when
$\Phi =\Z$; the
$\gamma_j$ do not have  all their
coordinates positive, so that our orbit closure is not in general the image of a map ${\mathbf A}^d(k)
\rightarrow
{\mathbf A}^{\mathbf N}(k)$. \par By dualizing the map $b$, we obtain a map 
$$w\colon \check \Z^d\rightarrow \check \Z^{\mathbf N}.$$ If the semigroup $\Gamma$ was finitely generated and therefore the right-hand group was finite
dimensional, say $\check\Z^N$, we would know how to find an embedded resolution of the singularities of
$\hbox{\rm Spec}k[t^\Gamma]\subset \hbox{\rm Spec}k[U_1,\ldots ,U_N]$ by a single toric map, as is explained in the next section. 
\subsection{Embedded toric resolution of orbit closures and binomial ideals in finite
dimensions}\label{toricresfin}
 I describe here the embedded resolution of  singularities of torus orbit
closures (in finite dimensions) as a generalization of what is done in [G-T] for monomial
curves. I assume that the reader is aware of the basics of toric
geometry, and especially of the fact that any fan  can be refined to a regular simplicial fan, by iterated stellar
subdivisions. I refer to  David Cox's notes [Cox] and to G\"unter Ewald's book ([E], VI, No.8, p. 253), and I use
the description of orbit closures given by Sturmfels in [St1], [St2]; see [Cox], \S 4. I keep the
notations introduced above, so that our  orbit closure is described by a map of semigroups
$b\colon\N^N\rightarrow \Z^d$, where $N$ is finite and the image of $b$ generates the group $\Z^d$. Choosing generators 
$(m^\ell-n^\ell)_{\ell\in \{1,\ldots , L\}}$ for $\hbox{\rm Ker}b$  gives us an exact sequence
$$\Z^L\rightarrow \Z^N\rightarrow \Z^d\rightarrow  0\ .$$\noindent Note 
that the
generators $m^\ell-n^\ell$ are necessarily primitive vectors in
$\Z^N$, and we may choose them in such a way that $m^\ell$ and $n^\ell$ both have
non-negative coordinates (see [St]). I assume from now on that they are chosen in this
way.\par\noindent
 By dualizing, we obtain an exact sequence $$0\rightarrow\check \Z^d\rightarrow
\check\Z^N\rightarrow \check \Z^L,$$ where the image of
 the $i-th$ basis vector of $\check\Z^d$ is the vector
$$w^i=(\gamma_{1i},\gamma_{2i},\ldots ,\gamma_{Ni}).$$\par\noindent
 Choosing one of the $m^\ell-n^\ell$ gives us an injective map $\Z\rightarrow \Z^N$, and a
surjective $\check\Z^N\rightarrow \check\Z$ having kernel
$H_\ell=\{a\in \check \Z^N\mid \langle a,m^\ell-n^\ell\rangle=0\}$.\par\noindent The
intersection
$\bigcap_{\ell\in \{1,\ldots , L\}}(H_\ell\otimes_{\Z}{\mathbf R})$ is the ${\mathbf R}$-linear span $W$ of the
$w^i$. From now on I will
 write $H_\ell$ indifferently for $H_\ell\otimes_{\Z}{\mathbf R}$ or its integral points. The interpretation will be clear from the context.
\par\noindent  By construction, we have 
$\langle w^i,m^\ell-n^\ell\rangle =0$ for $i=1,2,\ldots ,d,\  \hbox{\rm and}\ \ell\in \{1,\ldots , L\}$, hence also
$\langle w,m^\ell-n^\ell\rangle =0$ for $\ell\in \{1,\ldots , L\}$ and all vectors
 $w\in W$. Conversely, $W$ is the dual of the vector space generated by the vectors $m^\ell-n^\ell$ for $
\ell\in\{1,\ldots , L\}$. Remark that we may assume that $W$ contains no basis vector; if it contained the $i$-th
base vector, the variable $U_i$ would appear in none of the equations, so we might as well remove it. The
complement in
$\check {\mathbf R}_+^N$ of the union of the hyperplanes
$H_\ell$ is a  union of rational $N$-dimensional convex cones, and the union of hyperplanes itself is the 
reunion of rational convex cones of smaller dimensions. Note that $H_\ell$ meets the interior
of $ \check {\mathbf R}_+^N$ since the coordinates of the vectors $m^\ell -n^\ell$ are not all of the same
sign.\par\noindent  By the fundamental result on resolution of toric varieties (see [O], [Cox]), we can  find a regular fan
$\Sigma$ with  support  $\check {\mathbf R}_+^N$ which is not only compatible with the $w^i$ 
but also with the union of the hyperplanes $H_\ell$ in the sense that any $\sigma\in\Sigma$ 
meets any one of the $H_\ell$ along one of its faces. Then, $\Sigma$ is also compatible with the
intersection  $W$ of the $H_\ell$. In the sequel I tacitly assume $\Sigma$ to have these properties.\par I am going to show that for a suitable choice of
such a regular fan $\Sigma$, the toric map
$$\pi (\Sigma)\colon Z(\Sigma)\to {\mathbf A}^N(k)$$ induces an embedded resolution of the orbit closure $X\subset{\mathbf
A}^N(k)$ corresponding to $b$.\par Let 
$k$ be an algebraically closed field and let
$(U^{m^\ell}-U^{n^\ell})_{\ell\in \{1,\ldots , L\}}$ be generators of the kernel of  the map 
$$k[b ]\colon k[U_1,\ldots ,U_N]\rightarrow k[t_1^{\pm 1},\ldots ,t_d^{\pm 1}]\ \ \hbox{\rm  given by }\ \
U_i\mapsto t^{\gamma_i},$$ where $U^m=U_1^{m_1}\cdots U_N^{m_N}$. \par\noindent
Since the ideal $\hbox{\rm ker}k[b]$ defines an orbit closure, in view of ([E-S], Corollary 2.5) and the identity
$$U^{m+m'}-\lambda_{mn}\lambda_{m'n'}U^{n+n'}=U^{m'}(U^m-\lambda_{mn}U^n)+\lambda_{mn}U^n(U^{m'}-\lambda_{m'n'}U^{n'}),$$
the ideal $\hbox{\rm ker}k[b]$ is generated by binomials $U^m-U^n$ such that the exponents $m-n$ belong to the lattice ${\mathcal L}\subset
{\Z}^N$ generated by the $m^\ell-n^\ell$ for $\ell\in \{1,\ldots , L\}$, and this lattice is a direct factor in
${\Z}^N$, which corresponds to the fact that the ideal $\hbox{\rm ker}k[b]$ is prime. If we take $N-d$ of the
$m^\ell -n^\ell$, they generate a sublattice ${\mathcal L}_1\subset {\mathcal L}$. We can choose these elements in such a
way that they rationally generate the lattice ${\mathcal L}$, which means that ${\mathcal L}/{\mathcal L}_1$ is a torsion
$\Z$-module. Since ${\mathcal L}$ is a direct factor, it is the saturation of ${\mathcal L}_1$ in the sense of [E-S]; it is the
lattice of elements of $\Z^N$ which possess a multiple in ${\mathcal L}_1$. \par\medskip\noindent If we take a regular
simplicial cone of maximum dimension
$\sigma =\langle a^1,\ldots ,a^N\rangle\in \Sigma$, in the affine chart $U_\sigma\subset Z(\Sigma )$ with coordinates
 $Y_1,\ldots ,Y_N$ associated to $\sigma$ by definition of the toric map $\pi (\Sigma )\colon Z(\Sigma )\rightarrow {\mathbf A}^N(k)$,
we have the following expression for $\pi (\Sigma )\vert U_\sigma\colon U_\sigma \rightarrow  {\mathbf A}^N(k)$
$$U_i\mapsto Y_1^{a^1_i}\cdots Y_N^{a^N_i}\ \ \ \hbox{\rm and so}\ \ U^m\mapsto
Y_1^{\langle a^1,m\rangle}\cdots Y_N^{\langle a^N,m\rangle},$$ where $\langle a^i,m\rangle =\sum_{j=1}^Na^i_jm_j$. \par\noindent Let us
compute the transform by this map of one of our binomial generators, denoted by $U^m-U^n$. We may assume that
$a^1,\ldots ,a^t$ are those among the $a^j$ which lie on the hyperplane $H$ dual to
$m-n$,  i.e., such that $ \langle
a^j,m-n\rangle=0\ ,\ 1\leq j\leq t$. Because our fan is compatible with $H$,  all the other $\langle a^j,m-n\rangle$ are of the same sign, say
$\langle a^j,m-n\rangle>0$. We have then
$$U^m-U^n\mapsto Y_1^{\langle a^1,n\rangle} \cdots Y_N^{\langle a^N,n\rangle}(Y_{t+1}^{
\langle
a^{t+1},m-n\rangle}\cdots Y_N^{\langle a^N,m-n\rangle}-1).$$ Compare with [Z1], Th. 2, 
p.863 and
[S2]. This essentially means that in toric geometry, there is no need for Hironaka's game (see [H1]); 
it is replaced by the existence of a fan compatible with the
$H_\ell$ as above. Hironaka's game comes into play when one wants to dominate our toric 
map by the composition of a sequence of blowing-ups with non singular
centers.\par\noindent   The exceptional divisor of the map
$\pi (\sigma )$ is the union of the $y_j=0$ for those $j$ such that
$a^j$ is not a basis vector of $\check\Z^N$ (see [G-T]).
 For simplicity I will assume that none of the $a^j$ is a basis vector, so that the exceptional
divisor  is $Y_1\ldots Y_N=0$. If $t=0$, the strict transform $Y_{t+1}^{\langle
a^{t+1},m-n\rangle}\cdots Y_N^{\langle a^N,m-n\rangle}-1=0$ does not meet
 the exceptional divisor; therefore  the strict transform of $U^m-U^n=0$ meets the
exceptional divisor only in those charts 
$U_\sigma$ for which at least one of the primitive vectors 
$a^i$ of $\sigma$ is in $H$. More generally, we have.
\begin{proposition}\label{strict} If the strict transform by $\pi (\Sigma )$ of the subspace 
$X\subset
{\mathbf A}^N(k)$ defined by the ideal $(U^{m^\ell}-U^{n^\ell})_{\ell\in \{1,\ldots , L\}}$ is to meet the
exceptional divisor in the chart $\pi (\sigma )\colon U_\sigma \rightarrow {\mathbf A}^N(k)$, where
$\sigma=\langle a^1,\ldots ,a^N\rangle$, there must be a vector
$a^j$ such that 
$\langle a^j,m^\ell-n^\ell\rangle =0\ \hbox{\rm for}\ \ell\in\{1,\ldots , L\}$,
 i.e., $a^j\in W$.\end{proposition}\begin{proof} Look at the equations of the strict transform.\end{proof} This generalizes the result of [G-T],
5.2.\par\medskip\noindent The
${\mathbf R}$-vector space
$\check W$ generated by the $m^\ell -n^\ell$ is of dimension $N-d$; let us choose $N-d$ of the vectors $m^\ell -n^\ell$
 which generate $\check W$, say $m^1-n^1,\ldots ,m^{N-d}-n^{N-d}$. \par\noindent I now
have to check that  whenever the strict transform of the  subspace of ${\mathbf A}^N(k)$ defined by
$$U^{m^1}-U^{n^1}=\cdots = U^{m^{N-d}}-U^{n^{N-d}}=0$$ is not empty it is non singular and
transverse to the exceptional divisor, and that the strict transforms of all the other $
U^{m^\ell}-
U^{n^\ell}=0$ vanish on one of its irreducible components. \par The strict transforms of the equations have the form
$$\begin{array}{ccl}Y_1^{\langle a^1,m^1-n^1\rangle}\cdots Y_N^{\langle
a^N,m^1-n^1\rangle}-1&=&0\\ Y_1^{\langle a^1,m^2-n^2\rangle}\cdots Y_N^{\langle
a^N,m^2-n^2\rangle}-1&=&0\\ \vdots&\vdots&\vdots\\ Y_1^{\langle
a^1,m^{N-d}-n^{N-d}\rangle}\cdots Y_N^{\langle a^N,m^{N-d}-n^{N-d}\rangle}-1&=&0
\end{array}$$ Note that this strict transform meets the component $Y_j=0$ of the exceptional
divisor if and only if $a^j\in W$;  let us renumber the $a^j$ so that $a^1,\ldots ,a^t$ are those which
lie in $W$, i.e., $\langle a^i,m^s-n^s\rangle=0\ \hbox{\rm for} \ 1\leq i\leq t\ \hbox{\rm and all}
\ s$; if none of the $a^i$ is in $W$, set $t=0$; in this case, the strict transform does not meet any of the coordinate hyperplanes in the  chart
$U_\sigma$. As in ([G-T], \S 4),
 we can compute the jacobian matrix $J$ of these equations by logarithmic differentiation, and
find the equality of $N\times (N-d)$ matrices
$$Y_{t+1}\ldots .Y_N J=Y_{t+1}^{\sum_s\langle a^{t+1},m^s-n^s\rangle}\ldots Y_N^{\sum_s
\langle
a^N,m^s-n^s\rangle}(\langle a^j,m^s-n^s\rangle).$$ 
\begin{proposition}\label{binor} Given an irreducible binomial variety $X\subset {\mathbf A}^N(k)$, with the notations just introduced, the
rank of the image in
${\rm Mat}_{N\times  L}(k)$ of the matrix
$(\langle a^j,m^s-n^s\rangle)\in {\rm Mat}_{N\times  L}(\Z)$ is $N-d$. The strict transform by $\pi (\Sigma
)$ of the subspace 
$X\subset
{\mathbf A}^N(k)$ defined by the ideal $(U^{m^\ell}-U^{n^\ell})_{\ell\in \{1,\ldots , L\}}$ is regular and
transversal to the exceptional divisor; it is also irreducible in each chart.
\end{proposition}
\begin{proof} Let $\sigma$ be a cone of of maximal dimension in the fan $\Sigma$. Let us denote as above by $t$ the dimension of
$\sigma\cap W$. Since
$\hbox{\rm dim}W=d$, we must have
$t\leq d$; since the vectors $a^j$ form a basis of $\Qq^N$, and the space $\check W$  generated
by the $m^s-n^s$ is of dimension $N-d$, the rank of the matrix
$(\langle a^j,m^s-n^s\rangle)$ is $N-d$, which proves the lemma if $k$ is of characteristic zero.
\par\noindent To prove the general case, view our matrix $(\langle a^j,m^s-n^s\rangle)$ as  describing the composed map
$$\Z^L\to \Z^N \stackrel{M(\sigma)}{\longrightarrow }\Z^N;$$ we have to check that the $(N-d)\times
(N-d)$-minors of this matrix have no common divisor. This follows from the:
\begin{lemma}\label{mincopr} In a sequence of maps of free $\Z$-modules as above, assume that the image of the first map is a direct factor of rank
$N-d$ and that the determinant of $M(\sigma)$ is $\pm 1$. Then the $(N-d)\times (N-d)$-minors of the matrix describing the composed map
have no common divisor.\end{lemma}
 \begin{proof}Let us consider the sequence of maps obtained by taking $(N-d)^{\hbox{\rm th}}$ exterior powers:
$$\stackrel{N-d}\Lambda \Z^L\to \stackrel{N-d}\Lambda\Z^N \stackrel{\stackrel{N-d}\Lambda
M(\sigma)}{\longrightarrow }\stackrel{N-d}\Lambda\Z^N.$$ 
Since the image of the first map is a direct factor, so is the image of its exterior power, which is
of rank one and generated by the $(N-d)\times (N-d)$-minors of the matrix defining the first map. This image is therefore a primitive vector, which implies
that the minors are coprime. We have to prove that the image of this primitive vector by $\stackrel{N-d}\Lambda M(\sigma)$ is again a primitive vector.
This will be true if the determinant of that matrix is $\pm 1$. Remark that its entries are the $(N-d)\times (N-d)$-minors extracted from
$M(\sigma)$ and suitably ordered. By a result of Sylvester (see [Ai], p.92), nowadays an exercise in linear algebra, this implies that its
determinant is the ${N-1}\choose {d}$-th power of the determinant of $M(\sigma)$, and the result.
\end{proof}
 So we see that the rank of $J$ is $N-d$ everywhere on the
strict transform, and by Zariski's jacobian criterion ([Z4]) this strict transform is smooth and transversal to the exceptional divisor.
 Note however that it is not necessarily irreducible; let us show that the strict transform
of our orbit closure is one of its irreducible components. Since the differences of the exponents in the total transform and the strict transform of a
binomial are the same, the lattice of exponents generated by the exponents of all the strict transforms of the binomials
$(U^{m^\ell}-U^{n^\ell})_{\ell\in \{1,\ldots , L\}}$ is the image
$M(\sigma){\mathcal L}$ of the lattice
${\mathcal L}$ by the linear map
$\Z^N\to \Z^N$ corresponding to the matrix
$M(\sigma)$ with rows $(a^1,\ldots ,a^N)$. Similarly the exponents of the strict transforms of $U^{m^1}-U^{n^1}, \ldots
,U^{m^{N-d}}-U^{n^{N-d}}$ generate the lattice $M(\sigma){\mathcal L}_1$. The lattice $M(\sigma){\mathcal L}$ is the saturation of $M(\sigma){\mathcal
L}_1$, and so according to [E-S], since we assume that $k$ is algebraically closed, the strict transform of our orbit closure is one of the irreducible
components of the binomial variety defined by the $N-d$ equations displayed above.
\par\noindent
The charts corresponding to regular cones $\sigma\in \Sigma$ of dimension $<N$
are open subsets of those which we have just studied, so they contribute nothing new.\end{proof}
\par\medskip I now show that our
toric map can be chosen so that it induces an isomorphism outside of the singular locus of our irreducible binomial variety, so that it
is an embedded resolution and not only a pseudo-resolution in the sense of [G-T].  A part of this is the fruit of common reflexions with Pedro
Gonz\'alez-P\'erez (see [GP-T]). The same argument may also be used to replace the jacobian rank
argument just given, since it reduces the general toric case to the case of normal toric varieties, where the jacobian
rank argument is unnecessary. See [GP-T] for a proof of the resolution of singularities of binomial varieties along these lines.\par\medskip\noindent
With the notations introduced at the beginning of this subsection, remembering that the image of the map $b$ generates
$\Z^d$, let us denote by $\gamma_1,\ldots ,\gamma_N$ the images by the map $b$ of the basis vectors $e_1,\ldots ,e_N$ of $\Z^N$, and
by
$\check\sigma$ the cone which they generate in ${\mathbf R}^d$. It depends only on $b$ and we may also
denote it by $\check \sigma(b)$. Its convex dual
$\sigma$ (or $\sigma(b)$) in $\check{\mathbf R}^d$ is a strictly convex cone of dimension $d$ whose image in
$\check{\mathbf R}^N$ by the map $\check b$ is contained in the intersection of
$W$ with the first quadrant of $\check {\mathbf R}^N$. Indeed, we know by basic properties of duality that the image of $\check {\mathbf
R}^d$ in $\check {\mathbf R}^N$ consists of those linear forms on ${\mathbf R}^N$ which vanish on the kernel of $b$, and is
therefore equal to $W$. The image of $\sigma$ consists of those linear forms $\tilde f$ on ${\mathbf R}^N$ which are in $W$ and such that
$\tilde f (e_i)=f(\gamma_i)\geq 0,\ i=1,\ldots ,N$, where $f\in\sigma$. Therefore $\tilde f$ is in the dual of the first
quadrant of ${\mathbf R}^N$, which is the first quadrant of $\check {\mathbf R}^N$.
\begin{definition} Given a map $b\colon \N^N\to\Z^d$ of semigroups, the strictly convex cone $\sigma(b)\subset W\cap\check {\mathbf R}_+^N$ is
called the \textit {weight cone} associated to $b$.\end{definition}
\begin{definition} (Condition RES($b$)) Let $X$ be an orbit closure corresponding to a map $b\colon \N^N\to
\Z^d$ as above. Let
$\sigma (b)$ be the weight cone associated to $b$ just above. A fan
$\Sigma$ with support ${\mathbf R}_+^N$ satisfies condition RES($b$) if for some system of generators $(m^\ell-n^\ell)_{\ell\in
\{1,\ldots , L\}}$ of the kernel of the extended map
$b\colon\Z^N\to \Z^d$, it is compatible with the hyperplanes
$H_\ell$ dual to the vectors $(m^\ell-n^\ell)_{\ell\in \{1,\ldots , L\}}$, and in addition all the regular faces of the
cone $\sigma (b)$ belong to $\Sigma$.\end{definition}
Remark that by the classical results on resolution of normal toric varieties (see [Cox]), fans satisfying (RES) do exist.
\begin{proposition}\label{truereso} Assuming that the field $k$ is algebraically closed, let $X\subset
{\mathbf A}^N(k)$ be an orbit closure corresponding to a morphism $b\colon\Z^N\rightarrow \Z^d$. The singular locus of $X$ is a union of intersections
of $X$ with linear coordinate spaces. The toric map $\pi (\Sigma )\colon Z(\Sigma )\rightarrow {\mathbf A}^N(k)$
associated to a regular simplicial fan $\Sigma$
satisfying the condition {\rm RES($b$)} just defined is an embedded resolution of the singularities of $X$.
\end{proposition}
\begin{proof} Let $d$ be the dimension of our binomial
variety $X$ and $(m^\ell-n^\ell)_{\ell\in \{1,\ldots , L\}}$ and $\Sigma$ be a system of generators and a fan as given by condition RES($b$). A
straightforward computation using logarithmic differentials shows that the jacobian determinant
$J_{I,L'}$ of rank
$c=N-d$ of the generators $(U^{m^\ell}-\lambda_{m^\ell
n^\ell}U^{n^\ell})_{\ell\in\{1,\ldots , L\}}$ of our prime binomial ideal $ P\subset k[U_1,\ldots
,U_N]$ associated to $I=(i_1,\ldots , i_c)$ and a subset $L'\subseteq \{1,\ldots , L\}$ of
cardinality $c$ satisfies the congruence 
$$U_{i_1}\ldots U_{i_c}.J_{I,L'}\equiv
\big(\prod_{\ell\in L'}U^{m^\ell}\big)\hbox{\rm Det}_{I,L'}\big( (\langle m-n\rangle )\big) \ \ 
\hbox{\rm mod.} P,$$ where $\big(\langle m-n\rangle\big)$ is the matrix of the
vectors $(m^\ell-n^\ell)_{\ell\in \{1,\ldots , L\}}$,  and $\hbox{\rm Det}_{I,L'}$ indicates the minor
in question. By Proposition \ref{binor}, the rank of the image in $k^{N\times L}$ of the jacobian matrix is
equal to $c$. This proves the first part of the Proposition. For example, the only possibility for the intersection of
our variety with
$U_i=0$ not to be in the singular locus is that there exists such a minor such
that all the $m^\ell_i$ (resp. $n^\ell_i$) for $\ell\in L'$ are zero except one, which is equal to one. In
that case one sees that one of the equations of our binomial variety is of the form
$$U_i{U'}^r-{U'}^s=0,$$ and all the others are independent of $U_i$. \par\noindent
  By ([St3], Corollary 13.6), the
normalisation of our binomial variety
$X$ is
$X_{\sigma (b)}=\hbox{\rm Spec}k[\check \sigma (b)\cap \Z^d]$. The proof is as follows: let $P$ be the prime binomial ideal which is the kernel of
the map of semigroup algebras $$k[b]:k[U_1,\ldots ,U_N]\to k[t_1^{\pm 1},\ldots ,t_d^{\pm 1}]$$ corrresponding to $b$. By construction of
$\check\sigma (b)$ the image of this map is contained in
$k[\check\sigma (b)\cap\Z^d]$. Finally we have an injection
$$ k[U_1,\ldots ,U_N]/P\to  k[\check\sigma (b)\cap\Z^d]$$ which is birational since the image of $b$ generates $\Z^d$, and the algebra on the right is
normal since the semigroup $\check\sigma (b)\cap\Z^d$ is saturated. A toric resolution of singularities of
$X$ is obtained by a regular refinement of the fan in $W\cap \check {\mathbf R}_+^N$ consisting of $\sigma (b)$ and its faces which does
not affect the regular faces of
$\sigma (b)$, since that toric map resolves the singularities of $X_{\sigma (b)}$ and is an isomorphism outside of the singular
locus. Now we can choose a regular fan supported in $\check {\mathbf R}_+^N$, compatible with the hyperplanes $H_{mn}$
and which contains the cones of our regular subdivision of $\sigma (b)$. Any such fan corresponds to a toric map
$Z(\Sigma)\to {\mathbf A}^N(k)$ which is a toric embedded resolution of $X$, since the strict transform of $X$ is
transversal to the exceptional divisor by the argument given above, and the restriction to this strict transform
is an isomorphism outside of the singular locus of $X$. This proves the result.
\end{proof}
\begin{Remark}1) At least in characteristic zero, the equivariant resolution theorem implies that an irreducible binomial variety $X$ as above can
be resolved by a sequence of equivariant permissible blowing ups. The theorem of De Concini-Procesi ([DC-P]) implies that our toric resolution
can be dominated by a sequence of equivariant blowing-ups. It would be interesting to prove directly that these blowing ups can be chosen to be
permissible, in the sense that at each step the space blown up is normally flat along the center.\par\noindent 2) It has been remarked by P.
Gonz\'alez P\'erez (see [GP]) that any fan $\Sigma$ which is compatible with $W$ and contains the regular cones of $\sigma (b)$ will also provide a
resolution of singularities of the orbit closure correponding to $b$; this description depends only on $b$ and not on the choice of equations. This is
illustrated by the fact that in [G-T] a resolution of the monomial curve is obtained precisely with such fans.
\end{Remark}
\begin{proposition} \label{TP} Let $(u^{m_k})_{k\in K}$ be a finite collection of monomials in a polynomial ring $R=k[(u_i)_{i\in I}]$
endowed with a rational valuation $\mu$. There is a birational toric extension $R\to R'=R[(u^{\alpha_j})_{j\in J}]_{m'}$ such that  in $R'$ the ideal
$(u^{m_k})_{k\in K}R'$ is principal and generated by the monomial with the least valuation.\end{proposition}
\begin{proof} Let $J\subset I$ be the finite set of the variables appearing in the monomials considered. Consider the toric modification associated to a fan of
the first quadrant of ${\mathbf R}_+^{\vert J\vert}$ which is compatible with the hyperplanes dual to the vectors
$(m_k-m_{k'})_{k\neq k'}$; by the valuative criterion for properness the valuation $\mu$ picks a point in a chart of that modification, and the
corresponding toric extension satisfies the condition of the proposition.\end{proof}
 Let us now suppose that we
have an  irreducible subspace of
${\mathbf A}^N(k)$ defined by binomials 
$U^{m^\ell}-\lambda_{m^\ell n^\ell }U^{n^\ell}$ for ${\ell \in \{1,\ldots , L\}}$, where as usual
$U^m=U_1^{m_1}\cdots U_N^{m_N}$ and $\lambda_{mn}\in k^*$. By ([E-S],  Cor.
2.3), at least if our field $k$ is algebraically closed, the subspace defined by the binomials 
obtained by changing each $\lambda_{mn}$ to $1$ is a torus orbit closure, and by Proposition \ref{truereso} has an
embedded  resolution  by a toric map $\pi (\Sigma)\colon Z(\Sigma)\rightarrow {\mathbf A}^N(k)$
as above.  In one of the charts of $Z(\Sigma)$, the ideal $(U^{m^\ell}-U^{n^\ell})_{\ell \in \{1,\ldots , L\}}$ becomes the ideal
generated (up to  permutation of
the coordinates) by  
$$Y^{e(n^k)}(Y_{t+1}^{\langle a^{t+1},m^k-n^k\rangle}
\cdots Y_N^{\langle
a^N,m^k-n^k\rangle}-1);\ \ \ \ \ \ 1\leq k\leq N-d,$$
where $e(n)=(\langle a^1,n\rangle ,\ldots ,\langle a^N,n\rangle )$, and we are assuming as above that $$\langle a^i,m-n\rangle =0,\ 1\leq
i\leq t\ \ \hbox{\rm and}\ \ \langle a^i,m-n\rangle >0,\ t+1\leq i\leq N.$$ The nature of the computation shows that we
have:
$$U^m-\lambda_{mn}U^n\mapsto Y^{e(n)}(Y_{t+1}^{\langle a^{t+1},m-n\rangle}\cdots
Y_N^{\langle a^N,m-n\rangle}-\lambda_{mn}).$$
The relations between the $\lambda_{mn}$
associated to linear relations between the $m-n$, of the form (with the notations used above)
$$\lambda_{m^\ell n^\ell}=\prod_{k=1}^{N-d}\lambda_{m^k n^k} ^{b^\ell_k},$$
are exactly what is
needed to ensure that once we have chosen a system of generators $m^k-n^k$ as above, the
transforms of the remaining binomials $U^m-\lambda_{mn}U^n$ will define an irreducible component of the
 non singular variety defined by
$$Y_{t+1}^{\langle a^{t+1},m^k-n^k\rangle}\cdots
Y_N^{\langle a^N,m^k-n^k\rangle}-\lambda^{\epsilon (k)}_{m^kn^k}=0 \ ;\ \ \ \ \ 1\leq k\leq N-d,$$
where $\epsilon (k)=\pm 1$ is the common sign which the vectors $a^i,\ t+1\leq i\leq N$ take on $m^k-n^k$, so that the
subspace defined by our binomial ideal is also resolved by $\pi(\Sigma)$. Keeping the
notations just introduced, we see that we have proved:
\begin{proposition}\label{resbin}
Assuming that the field $k$ is algebraically closed, let $X\subset {\mathbf A}^N(k)$ be an
irreducible binomial variety, defined by the ideal of $k[U_1,\ldots ,U_N]$ generated by $(U^{m^\ell}-\lambda_{m^\ell
n^\ell }U^{n^\ell})_{\ell \in \{1,\ldots , L\}},\
\lambda_{m^\ell n^\ell}\in k^*$. For any regular simplicial fan $\Sigma$ satisfying the condition {\rm RES($b$)} with respect to the orbit closure defined
by the binomials $(U^{m^\ell}-U^{n^\ell})_{\ell \in
\{1,\ldots , L\}}$, the corresponding toric map
$\pi (\Sigma )\colon Z(\Sigma )\rightarrow {\mathbf A}^N(k)$ is an embedded resolution of the singularities of
$X\subset {\mathbf A}^N(k)$.
\end{proposition} 
\begin{remark} After the description of the singular locus of a binomial variety given in the proof of Proposition \ref{truereso}, Zariski's
jacobian criterion implies that the
$(N-d)\times (N-d)$-minors of the matrix 
 $\big( m^\ell-n^\ell\big)$ are coprime integers if the binomial variety determined by the vectors $(m^\ell-n^\ell)$ is  reduced and irreducible
over any algebraically closed field. The geometric fact that this is true if it is reduced and irreducible over $\C$ also follows, in view of
[E-S], from the combinatorial fact expressed in Lemma \ref{mincopr}. \par\noindent  In the case where
$d=1$ and our (equations-wise) binomial variety is the (parametri\-zation-wise) monomial curve corresponding to the semigroup
$\Gamma=\langle
\overline{\beta_0},\ldots ,\overline{\beta_g}\rangle$ and is a complete intersection, as in the example of subsection
\ref{planebranches}, this Lemma \ref{mincopr} solves the minor difficulty mentioned at the
end of \S 5 of [G-T] and shows that the minors of the matrix coincide, up to sign, with the generators of the
semigroup.  \end{remark}

\subsection{An example: branches revisited}\label{branchrev}
 Let $R$ be a one-dimensional excellent equicharacteristic analytically irreducible local ring,
i.e., the local ring of a branch. Let us assume that the residue field $k$ of $R$ is algebraically
closed. The ring $R$ has only one non trivial valuation, whose valuation
ring is the integral closure $R_\nu$ of $R$ in its field of fractions. Since $R$ is excellent,
$R_\nu$ is a finite $R$-module ([EGA], Scholie 7.8.3, vii)). Therefore there exists an element $d\in R,\ d\neq 0$, such that $dR_\nu\subset R$ and
$\N\setminus \Gamma$ is finite, so that the semigroup $\Gamma$ of values of $\nu$ on
$R\setminus\{ 0\}$ is finitely generated. Let us write
$\Gamma=\langle
\gamma_1,\gamma_2,\ldots , \gamma_N\rangle$. Since $k$ is algebraically closed, the valuation $\nu$ is rational, and the
completion $\hat R^m$ of $R$ is equal to $\hat R^\nu$ and is a subring of $k[[t]]$, the valuation
$\hat \nu$ being induced by the $t$-adic valuation. Let us choose elements
$(\xi_i )_{1\leq
i\leq N}$ in $R$ such that $\nu (\xi_i)=\gamma_i$. Their initial forms for the $\nu$-adic
valuation generate the graded ring 
$$\hbox{\rm gr}_\nu R=k[t^{\gamma_1},\ldots ,t^{\gamma_N}]\subset k[t]=\hbox{\rm gr}_\nu
R_\nu , $$ which is the affine algebra corresponding to the monomial curve $C^\Gamma$ associated to $\Gamma$. Since the $(\xi_i )_{1\leq
i\leq N}$ in $R$ are finite in number, they form a minimal quite sufficient set of generators for the maximal ideal of $\hat R^{(\nu)}=\hat
R^m$.\par\noindent By Proposition \ref{coordinates}, the map $u_i\mapsto
\xi_iv^{-\nu (\xi_i)}$ determines a surjective map $$k[v][[u_1,\ldots ,u_N]]\rightarrow \widehat{\A_\nu
(R)}\subset k[v][[t]].$$ The kernel of this map is generated by series
$$F_{mn}=u^m-u^n+\sum_s c^{(mn)}_s(v)u^s$$ which are deformations of  binomial generators
$u^m-u^n$ of the kernel of the map $$k[u_1,\ldots ,u_N]\to k[t^{\gamma_1},\ldots
,t^{\gamma_N}],\ \ \ u_i\mapsto t^{\gamma_i}.$$
Denoting by $w=(\gamma_1,\gamma_2,\ldots
,\gamma_N)\in \N^N$ the weight vector, for each of these equations, we add to the binomial
$u^m-u^n$ only terms $c^{(mn)}_p(v)u^s$ such that
$w(s)>w(m)=w(n)$. If now we choose a regular fan $\Sigma$ in ${\mathbf R}_+^N$ compatible
with the dual hyperplanes of the vectors $m-n$, and so in particular with the vector $w$,
we obtain a toric map
$$\pi (\Sigma)\colon Z(\Sigma )\to {\mathbf A}^N (k)$$
which is a toric embedded resolution of the monomial curve $C^\Gamma \subset {\mathbf A}^N
(k)$ according to what we saw in subsection \ref{toricresfin}. The map $$\hbox{\rm Id}_{{\mathbf
A}^1(k)}\times \pi (\Sigma )\colon {\mathbf A}^1(k)\times Z(\Sigma )\to  {\mathbf A}^1(k)\times {\mathbf
A}^N (k)$$ then induces an embedded resolution of singularities for the formal subspace
defined in $k[v][[u_1,\ldots ,u_N]]$ by the ideal generated by the $F_{mn}$, which is a
simultaneous resolution for all the fibers of the projection to ${\mathbf A}^1(k)$.
This is checked exactly as in
[G-T]. It implies that the map $\pi (\Sigma )$ gives an embedded resolution of the
singularities of the image of the formal completion of our curve embedded in
${\mathbf A}^N (k)$ by using the elements $\xi_i\in R\subset \hat R^\nu$ as coordinates.
Now since $\pi (\Sigma )$ is a monomial map, it corresponds to adding to $\hat R^\nu$
certain monomials in the $\xi_i$, with some negative exponents, and then localizing the
$\hat R^\nu$-algebra which they generate. This operation has a meaning
in $R$ itself, and builds a local ring $R'$ which is essentially of finite type over $R$
and birationally equivalent to $R$. The local ring $R'$ is excellent since $R$ is, and its
completion is regular since it is the completion of the transform of $\hat R^\nu$, so $R'$ itself is regular. \par
To illustrate the difference between this approach and Abhyankar's, consider, over an algebraically closed field 
$k$ of characteristic $p$, the curve given parametrically by:
$$u_0=t^p+t^{p+1},\ \ \ \ \ \ \ u_1=t^{p^2+1}+t^{p^2+p+1}.$$
It has no Puiseux presentation of the form $u_0'=t^p,\ u_1'=u_1'(t)$. However, we may consider it as a deformation of a plane monomial curve, with an
equation of the form:
$$u_1^p-u_0^{p^2+1}-u_0^{p^2-p+1}u_1+\cdots$$
Here we can remark that the irreducible initial binomial describes a curve which is purely inseparable over $k((u_0))$, while this is not the case
for our original curve. However, if we resolve this binomial, say by a toric modification of ${\mathbf A}^2(k)$ having as one of its charts
$$u_0=y_0^py_1,\ \ \ \ \ u_1=y_0^{p^2+1}y_1^p,$$ we find that the transform of the equation is
$$y_0^{p(p^2+1)}y_1^{p^2}(1-y_1-y_0y_1+\cdots )$$
and therefore we have an embedded resolution of our curve.\par
The usual relations between Puiseux exponents and semigroup may fail in characteristic $p$. Consider, over an algebraically closed field $k$ of
characteristic
$p$, the curve (see [Ca], p.114):
$$u_0=t^{p^3},\ \ \ \ \ \ u_1=t^{p^3+p^2}+t^{p^3+p^2+p+1}.$$
One can check that although only three exponents are visible on the expansions, its semigroup has four generators:
$$\Gamma=\langle p^3,\ p^3+p^2,\ p^4+p^3+p^2+p,\ p^5+p^4+p^3+p^2+p+1\rangle,$$ and the corresponding monomial curve has
equations 
$$U_1^p-U_0^{p+1}=0,\ \ \ U_2^p-U_0^{p(p+1)}U_1=0,\ \ \ U_3^p-U_1^{p^3}U_2=0,$$ while the curve itself has, in its natural
coordinates, equations
$$u_1^p-u_0^{p+1}-u_2=0,\ u_2^p-u_0^{p(p+1)}u_1+u_3=0,\ u_3^p-u_1^{p^3}u_2-u_0^{p+1}u_3^p=0.$$\par\medskip
The following question and example are motivated by two problems: \par\noindent 1) If there is a natural toric specialization of a
singularity, where are the regular points of the strict transform of the singularity by a toric resolution of the toric variety?\par\noindent
2) Of what nature is the complexity of the binomial relations defining $\hbox{\rm gr}_\nu R$? this is of course essential for the description of its
initial partial toric resolution.
\subsection{A question on surfaces}\label{exersurf}  Let $k$ be an algebraically closed field of characteristic $p$. Consider in ${\mathbf
A}^3(k)$ the surface
$$ F=z^{p+1}+y^{p-1}z+x^{p+1}=0$$ studied by Abhyankar in ([A4], p.589); it is quasi-ordinary with respect to the projection onto the $(x,y)$-plane; the
$z$-discriminant is $x^{p+1}=0$. For $p\geq 5$, however, its Galois group with respect to this projection is ``very large and quite
complicated'' ({\it loc. cit.}) and not a cyclic group as one would expect in characteristic zero. Indeed, Abhyankar later showed (see [A5]) that the Galois
group over $\hat K=k((y,x^{p+1}))$ of the polynomial $ F\in \hat K[z]$  is $\hbox{\rm PGL}_2(\hbox{\rm\textbf{F}}_p)$.\par\noindent Study the rational
valuations, if any, for which the associated graded ring of this surface singularity is
$$k[X,Y,Z,U]/(Z^p+Y^{p-1}, ZU+X^{p+1}),$$ so that $k[x,y,z]/(z^{p+1}+y^{p-1}z+x^{p+1})$ appears as a deformation of its associated graded ring; its
equations in natural coordinates are
$$ z^p+y^{p-1}-u=0, \ \ \ zu+x^{p+1}=0.$$
Note that the associated graded ring is {\it not} quasi-ordinary for the projection onto the $(X,Y)$-plane and that the field extension associated to this
projection displays inseparability.
\subsection{Another example, in dimension three}\label{dimthree}
\begin{example} (partly inspired by an example of Spivakovsky) Let us give
$\Z^2$ the lexicographic order and denote by $k((t^{\Z^2_{\rm lex}}))$ the field of
Puiseux
series associated to the ordered group $\Z^2_{\rm lex}$ and the field $k$ as in [Ka], [B2]; it is the
field of formal series with exponents forming a well ordered subset of
$\Z^2_{\rm lex}$. It is naturally endowed with the $t$-adic valuation with values in $\Z^2_{\rm lex}$. Let us denote by
$k[[t^{\Z^2_+}]]$ the corresponding valuation ring. Choose a 
sequence
of pairs of positive integers $(a_i,b_i)_{i\geq 3}$ and a sequence of elements $(\lambda_i\in k^*)_{i\geq 3}$ such that $b_{i+1}>b_i$, the series
$\sum_{i\geq 3}\lambda_i u_2^{b_i}$ is not algebraic over $k[u_2]$, and the ratios
$\frac{a_{i+1}-a_i}{b_{i+1}}$ are positive and increases strictly with $i$. Let $R_0$ be the
$k$-subalgebra of $k[[t^{\Z^2_+}]]$ generated by
$$u_1=t^{(0,1)}, u_2=t^{(1,0)},u_3=\sum_{i\geq 3}\lambda_i u_1^{-a_i}u_2^{b_i}.$$
There
cannot be an algebraic relation between $u_1,u_2,$ and $u_3$, so the ring
$R_0=k[u_1,u_2,u_3]$ is the polynomial ring in three variables. It inherits the $t$-adic 
valuation of
$k[[t^{\Z^2_+}]]$. One checks that this valuation extends to the localization
$R=k[u_1,u_2,u_3]_{(u_1,u_2,u_3)}$; it is a rational valuation of height two and rational rank
two.  Let us try to
compute the semigroup $\Gamma$ of the values that it takes on $R$. We have
$\gamma_1=(0,1),\ \gamma_2=(1,0),\gamma_3=(b_3,-a_3)\in \Gamma$. Set
$\Gamma_3=\langle \gamma_1,\gamma_2,\gamma_3\rangle$. Then we have
$u_1^{a_3}u_3-\lambda_3 u_2^{b_3}=\sum_{i\geq 4}\lambda_i u_1^{a_3-a_i}u_2^{b_i}\in R$,  so that
$\gamma_4=(b_4,a_3-a_4)$ is in $\Gamma$. It is easy to deduce from our assumptions
that  {\it no
multiple of $\gamma_4$ is in $\Gamma_3$}, and that it is the smallest element of $\Gamma$
which is not in $\Gamma_3$. We set $u_4=u_1^{a_3}u_3-\lambda_3 u_2^{b_3}$, and continue in
the same manner: $u_1^{a_4-a_3}u_4-\lambda_4 u_2^{b_4}=u_5$, ...,
$u_1^{a_i-a_{i-1}}u_i-\lambda_i u_2^{b_i}=u_{i+1}$,... with
the generators $\gamma_i=\nu (u_i)=(b_i,a_{i-1}-a_i)$ for $i\geq 4$. Finally we have:
$$\Gamma=\langle \gamma_1,\gamma_2,\ldots ,\gamma_i ,\ldots\rangle ,$$ the initial forms of the $u_i$ constitute a minimal system of generators of
the graded $k$-algebra
$\hbox{\rm gr}_\nu R$, and the equations (setting $a_2=0$) $$u_1^{a_i-a_{i-1}}u_i-\lambda_i u_2^{b_i}=u_{i+1},\ \ i\geq 3$$
above  describe $\hat R^{(\nu )}$. In fact they even describe $R$; it is clear that from
them we can reconstruct the value of $u_3$ as a function of $u_1,u_2$ by (infinite)
elimination. The binomial equations for $\hbox{\rm gr}_\nu R$ are the 
$$U_1^{a_i-a_{i-1}}U_i-\lambda_i U_2^{b_i}=0,\ \ i\geq 3,$$ 
showing that all
the $U_i$ for $i\geq 3$ are  rationally dependent on $U_1,U_2$. Using the remark following
proposition \ref{pil}, we see that the  Krull dimension of $\hbox{\rm gr}_\nu R$ is two. The
fact that $R$ is regular of dimension three is due to the abyssal phenomenon that we met in
example \ref{exzar}. Remark that in this example $\hbox{\rm
gr}_\nu R$  is not regular.\par\noindent
From our assumption on the growth of the ratios we see moreover that {\it no multiple of}
$\gamma_i$ is in $\Gamma_{i-1}=\langle \gamma_1,\ldots ,\gamma_{i-1}\rangle$. In fact
$\gamma_i$ is outside of the cone with vertex $0$ generated by
$\Gamma_{i-1}$ in ${\mathbf R}^2$.\par\noindent
It is interesting to study in this case the construction of a system of regular fans which will
resolve the singularities of $\hbox{\rm gr}_\nu R$.\end{example}\par\noindent 
\begin{remark} Olivier Piltant has communicated to me an example of a valuation of height one on the ring
$R=k[u_1,u_2,u_3]_{(u_1,u_2,u_3)}$ where it is also the case that no multiple of $\gamma_i$ is in $\Gamma_{i-1}$. The valuation
is induced on
$R$ from a rational height one valuation, of rational rank three, on the power series ring
$k[[u_1,u_2/u_1,u_3/u_1]]$ which is monomial in the variables $u_1,u_2/u_1$ and $u_3/u_1-(u_2/u_1)^2-(u_2/u_1)^3$. It is
related to the example of [C-G-P].\end{remark}
\noindent {\textbf {Problem.}} If $\Gamma=\langle \gamma_1,\ldots, \gamma_i,\gamma_{i+1},\ldots\rangle$ is the semigroup of a
valuation of a n\oe therian local integral domain $R$, and $\Gamma_i=\langle \gamma_1,\ldots,
\gamma_i\rangle$, can the fact that even for large $i$ no multiple of $\gamma_i$ is in $\Gamma_{i-1}$ occur with
$R$ of dimension two? Hint: use [A2]. Compare with subsection \ref{exersurf}.\par 
\section{\textbf{Local uniformization} (speculative sketch)}\label{proof}
{\par\vskip.2truecm \hfill {\it ...and he will be too dazzled by the light to look at the
objects whose shadows he was seeing a moment ago...}\par\hfill Plato, {\it  The Republic}, VII,
515\par\medskip \par\medskip\noindent  Assume in this section that the local equicharacteristic excellent integral domain $R$ has an
algebraically closed residue field $k$.\par  Given a valuation $\nu_0$ on the ring $R$, we may specialize $\nu_0$ to a valuation $\nu$
whose ring $R_\nu$ dominates $R$ with trivial residue field extension (see subsection
\ref{speciaval}, Proposition \ref{domi}).\par\noindent We first uniformize the valuation  $\hat
\nu$ on $\hat R^{(\nu)}$. In the last subsection I will sketch the descent from
$\hat R^{(\nu)}$ to $R$.\par\noindent After this descent, by the first Corollary to Proposition
\ref{complete deformation}, which ensures the nice behaviour of the specialization above
with respect to composition of valuations, the same toric map will also simultaneously 
uniformize the valuations
in the fibers of the total space of the specialization of $ R$ to 
$\hbox{\rm
gr}_{\nu_0}R$, which is the subspace of the preceding specialization defined by the ideal $(v^\psi -1)_{\psi \in \Psi_+}$, and hence also
uniformize 
$ \nu_0$ on $R$. The proof proceeds by induction in the dimension, since in order to prove the existence of a scalewise completion
for a ring of dimension $d$ I have used local uniformization in dimension $<d$. The proof actually follows
the scheme $TLU(d-1)\implies TC(d)\implies TLU(d)$. On the way we use that $TLU(d-1)\implies TP(d-1)$ by Corollary
\ref{monomialization} below.\par\medskip\noindent Since
$\nu$ is rational, by Proposition
\ref{binomialideal} the graded algebra $\hbox{\rm gr}_\nu R$ is a quotient of a polynomial
ring by a binomial ideal. Using Proposition
\ref{madicomp}, we may assume that the completion
$\hat R^{(\nu)}$ has a field of representatives $k\subset \hat R^{(\nu)}$ which we take  as our base field. According to
Propositions \ref{coordinates} and \ref{completeequations}, we have a surjective map
$$\widehat{k[v^{\Phi_+}][(w_j)_{j\in J}]}\rightarrow\widehat{{\mathcal A}_{\hat\nu}(\hat R^{(\nu
)})}^{(\nu_\A)},$$ the kernel of which is the closure in  $\widehat{k[v^{\Phi_+}][(w_j)_{j\in J}]}$ of the ideal generated by the elements 
$$\tilde G_{mn}=w^m-\lambda_{mn}w^n+\sum_sc^{(mn)}_s(v^\phi )w^s,$$ 
where the tilda means that we have taken the varying (with coefficients in $k[v^{\Phi_+}]$) version of the equation.
The idea for the proof is that
if the $k$-algebra $\hbox{\rm gr}_\nu R$ is finitely generated, say with $N$ 
generators, it is also finitely presented and one can
proceed essentially as in [G-T]: we have in this case by Lemma \ref{fingen} below $\hbox{\rm dim}\hat R^{(\nu)}=\hbox{\rm
dim}\hbox{\rm gr}_\nu R=\hbox{\rm r}(\nu )$, and a toric map $Z\to {\mathbf A}^N(k)$ which
resolves $\hbox{\rm Spec}\hbox{\rm gr}_\nu R$ will also uniformize $\hat \nu$ on $\hat
R^{(\nu )}$. If $\hbox{\rm gr}_\nu R$ is not finitely generated, we have to show that it suffices to
resolve a finite approximation to $\hbox{\rm gr}_\nu R$, the subalgebra $\hbox{\rm
gr}^{(h)}_\nu R$ of $k[x_1^{(h)},\ldots ,x_r^{(h)}]$ generated by the images of $(U_j)_{j\in F}$ for sufficiently large $h$,
where $h$ is a level of approximation for $\hbox{\rm gr}_\nu R_\nu$  in the sense of
Proposition \ref{grapp} and $F$ is chosen according to Proposition \ref{Peano} and the text which follows it as a quite sufficient initial generating set in
$\Gamma$. Here the key point is that  the abyssal phenomenon of subsection \ref{abyssalph}) tells us that if we consider the family defined in
$\widehat{k[v^{\Phi_+}][(w_j)_{j\in J}]}$ by the equations  $(\tilde {\mathbf F},(\tilde G_i)_{i\notin F})$ it has ''generic fiber'' $\hbox{\rm Spec}\hat
R^{(\nu)}$ and a special fiber which may be larger that $\hbox{\rm Spec}\hat{\hbox{\rm gr}}_{\hat\nu}R^{(\nu)}$, but coincides with it if one
considers only the variables with index in $F'$.
Then, because the equations $\tilde G_i$ create no new singularity it is
sufficient to resolve the space defined by $\tilde {\mathbf F}$, and for that, it is sufficient to resolve the ideal generated by its initial binomials, which
reduces us to the case of the finitely generated approximating algebra. 
\subsection{The approximation process}\label{approxpro} 
{\it A priori}, we must examine whether it makes sense to say that we resolve $\nu_{\rm gr}$  by a single infinite dimensional
toric modification of $k[(U_i)_{i\in I}]$, hence uniformize the completion $\hat\nu_{\rm gr}$ on
$\widehat{\hbox{\rm gr}}^{(\nu)}_\nu R$. This toric modification should then uniformize
$\hat\nu$ on
the completion $\hat R^{(\nu)}$. Fortunately, as we shall see below, we do not have to overcome this difficulty, because a
finite initial partial toric resolution of $\hbox{\rm gr}_{\hat\nu}\hat R^{(\nu)}$ suffices to uniformize $\hat\nu$ on $\hat R^{(\nu )}$. 
\par\noindent Let us keep the notations of the preceding section for a valued ring $(R,\nu)$; we will apply the results to $(\hat R^{(\nu)},\hat\nu)$. 
The variables
$U_i$ are ordered by the valuations
$\gamma_i\in
\Gamma$ of their images in $\hbox{\rm gr}_\nu R$.  Let us consider the
nested sequence of polynomial algebras $k[x^{(h)}]=k[x_1^{(h)},\ldots ,x_r^{(h)}]$ approximating
$\hbox{\rm gr}_\nu R_\nu$ according to Proposition \ref{grapp}, b).
Remember that $r$ is the
rational rank of the value group $\Phi$ of $\nu$. For each finite set $F$ of generators of the $k$-algebra $\hbox{\rm gr}_\nu R$, there is an integer $h$
such that $F$ is contained in
$\hbox{\rm gr}_\nu R\cap k[x^{(h)}]$. The degrees of the elements of this subalgebra are in a free subsemigroup
$\N^r\subseteq\Phi_+\cup\{0\}$ and
 the subalgebra $\hbox{\rm gr}^{(F)}_\nu R$ of $\hbox{\rm gr}_\nu R$ generated by $F$ is the image of the map $$k[(U_j)_{j\in F}]\to k[x^{(h)}]$$
sending $U_j$ to its image in $k[x^{(h)}]$. The kernel of the map
$k[(U_j)_{j\in F}]\to \hbox{\rm gr}^{(F)}_\nu R$ is a prime binomial ideal, and corresponds to a torus orbit closure in ${\mathbf A}^F(k)$ characterized
by the map $\omega (F)\colon\Z^F\to \Z^r$ sending the $j$-th basis element to $\nu (\xi_j)\in
\N^r\subseteq\Phi_+\cup\{0\}$.  I apply the results of subsection \ref{toricresfin} to this binomial ideal 
in finitely many
variables, obtaining a regular fan in $\check {\mathbf R}_+^F$. Adding an element to $F$ to get
$F'=F\cup \{ g \}$, we have
a commutative diagram 

\begin{equation*}\begin{CD}\Z^F@ >>>\Z^r\\ @VVV @VVV\\
\Z^{F'}@>>>\Z^r \end{CD} \end{equation*} 

and after dualizing  

\begin{equation*} \begin{CD}
\check\Z^r@ >>>\check\Z^{F'}\\ @VVV @VV{p_{F',F}}V\\ \check\Z^r@>>>\check\Z^F 
\end{CD}
\end{equation*} 
Note that in the case of curves, $r=1$ and, taking for $F$ the finite set of generators of the
algebra $\hbox{\rm gr}_\nu R$, the image of $\check \Z$ in $\check \Z^F=\check \Z^N$ is the
weight vector. Although no approximation is needed, we can illustrate the procedure by setting
$e_i=\hbox{\rm gcd}(\gamma_1,\ldots ,\gamma_i)$ and considering the nested sequence of
subalgebras  $$k[t^{e_1}]\subset  k[t^{e_2}] \subset
\cdots\subset k[t^{e_{N-1}}]\subset k[t].$$ 
The construction shown above then builds a sequence
of vector spaces \break $({\mathbf  R}^i\to {\mathbf  R}^{i-1})_{2\leq i\leq N}$ and a weight vector in
each. In the
case of plane branches, where the relations between the $\gamma_i$ are particularly simple, 
this suggests a procedure to build inductively a fan compatible with the hyperplanes 
corresponding to these relations; see subsection \ref{planebranches} and [G-T].\par\noindent
  For each finite set $F$ of variables $U_j$ as above, we take a fan
$\Sigma_F$ in $\check {\mathbf  R}_+^F$ compatible with the  hyperplanes dual
to generators of the kernel of
$\omega (F)$. For each inclusion $F\subset F'$, we consider similarly fans $\Sigma_{F'}$
compatible with the generators of the kernel of $\omega (F')$  and with the inverse image by
$p_{F',F}$ of $\Sigma_F$. In this manner, letting $F$ grow, we obtain a projective system of
vector spaces indexed by the filtering system of finite subsets of $I$, and a subordinate system
of regular fans, which is what we would like to call a regular fan for the projective system of
vector spaces. \par\noindent We are
led to the following definition of a possibly infinite dimensional fan:
\par\medskip\noindent
\begin{definition}  Let $P$ be a projective system of ${\mathbf R}$-linear maps between
finite-dimen\-sio\-nal real vector spaces, indexed by a filtering ordered set $I$. 
$$P\colon
\big((p_{i,j}\colon {\mathbf R}^{n(i)}\rightarrow {\mathbf R}^{n(j)})_{i,j\in I, i>j}\big),$$
where $n(i)$ is an
integer valued function. A fan $\Sigma$ subordinate to $P$ is the datum for each $i\in I$ of a
rational fan $\Sigma _i$ in ${\mathbf R}^{n(i)}$ in such a way that for each $i>j$, the fan $\Sigma_i$ is
a refinement of $p_{i,j}^{-1}(\Sigma_j)$.\par\noindent We say that the fan is of finite type if
there exists a finite subset $F\subset I$ such that for $i,j \notin F\ ,\ i>j$, the image of
$\Sigma_i$ is equal to $\Sigma_j$, and if we denote by $K_i\subset {\mathbf R}^{n(i)}$ the kernel of
$p_{i,j}$, the map $\sigma\mapsto (p_{i,j}(\sigma ), \sigma\cap K_i)$ is a bijection from
the set
of cones in $\Sigma_{n(i)}$ to the set of pairs consisting of a cone in $\Sigma_j$ and one of the
cones in $\Sigma_{n(i)}$ whose image by $p_{i,j}$ is $0$. We say that $\Sigma$ is simplicial
(resp. regular) if  all $\Sigma_i$ are. We say that $\Sigma$ is {\it rational} if there exists a
finite subset $F\subset I$ such that for $i,j \notin F\ ,\ i>j$, the map $p_{i,j}$ sends a basis
of
the integral lattice of ${\mathbf R}^{n(i)}$ to part of a basis of the integral lattice of
${\mathbf R}^{n(j)}$.
\end{definition}
 To such a system of fans is associated a system of toric maps
$\pi (\Sigma_i)\colon Z(\Sigma_i )\to {\mathbf A}^{n(i)}(k)$; if we have a simple system of
surjective maps
$p_i\colon {\mathbf R}^{n(i+1)}\to {\mathbf R}^{n(i)}$, we get inclusions ${\mathbf A}^{n(i)}(k)\subset {\mathbf
A}^{n(i+1)}(k)$, so that we get a compatible system of toric maps of an increasing sequence of
affine spaces, which we may call an infinite-dimensional toric map. I leave it to the reader to
interpret geometrically the various conditions defined for the system of fans. 
\subsection{Toric modifications and transforms of $\hat R^{(\nu)}$}\label{trans}
Finally, in order to uniformize the valuation $\hat \nu$ on $\hat R^{(\nu )}$ it
should suffice to resolve the singularities of a finite number of the equations $F_{mn}$, so that
we may consider only finitely many elements of the projective system of toric blowing-ups which ``resolves'' $\hbox{\rm gr}_\nu R$.\par\noindent
To check this, we must now compute the effect on $\hat R^{(\nu)}$ of the toric modifications associated to our regular fan. \par
 Using the specialization of $\hat R^{(\nu)}$ to $\widehat{\hbox{\rm gr}}^{(\nu)}_\nu R$, which we
can write in a very explicit form (see
\ref{complete deformation}) as the space associated to
$$\widehat {k[v^{\Phi_+}][(w_j)_{j\in IJ}]}/\overline{\big(w^m-\lambda_{mn}w^n+\sum_s
c_s^{(mn)}(v^\phi)w^s\big)},$$ with $\hbox{\rm weight} (w^s)>\hbox{\rm weight}(w^n)=\hbox{\rm weight}(w^m)$, the weight of a monomial being the
valuation of its image in $\hat R^{(\nu)}$. For any homomorphism 
$c\colon \Phi_+\to k^*$, say
$\phi\mapsto c(\phi )$, we have an isomorphism $$ \hat c\ \colon \ \widehat
{k[v^{\Phi_+}][(w_j)_{j\in J}]}/\overline{\big(w^m-\lambda_{mn}w^n+\sum_s
c_s^{(mn)}(v^\phi)w^s,\ (v^\phi-c(\phi ))_{\phi\in \Phi_+}\big)}\buildrel
\simeq\over\longrightarrow \hat R^{(\nu)} .$$ 
We will denote by $c(G_{mn})$ the series $$c(G_{mn})=w^m-\lambda_{mn}w^n+\sum_s
c_s^{(mn)}(c(v^\phi ))w^s\in \widehat {k[(w_j)_{j\in J}]}.$$\par
Let us first assume that the $k$-algebra $\hbox{\rm gr}_{\hat\nu}\hat R^{(\nu)}$ is of finite type, generated say
by $(\overline \eta_1,\ldots ,\overline \eta_N)$. We have
\begin{lemma}\label{fingen} If the $k$-algebra $\hbox{\rm gr}_{\hat\nu}\hat R^{(\nu)}$ is finitely generated, the equality\break  $ \hbox{\rm
dim}\hat R^{(\nu)}=\hbox{\rm r}(\nu)$ holds.\end{lemma} 
\begin{proof} Using Propositions \ref{coordinates}, \ref{completeequations}, and their corollaries (see Corollary \ref{present}, c)), we are
reduced to the case of a family of algebras which are ``formally of finite type''. 
This shows that we can apply the semi-continuity of the
dimensions of fibers to the map 
$$\hbox{\rm Spec}\widehat{\A_{\hat\nu} (\hat R^{(\nu)})}^{(\nu_\A)}\to\hbox{\rm Spec} k[v^{\Phi_+}]$$
so that by Piltant's Theorem (see \ref{piltant}),
$\hbox{\rm r}(\nu)=\hbox{\rm dim}\hbox{\rm gr}_{\hat\nu}\hat R^{(\nu)}\geq \hbox{\rm dim}\hat R^{(\nu)}$.
The reverse inequality is Abhyankar's inequality applied to $\hat R^{(\nu)}$.\end{proof}
\begin{Remark} 1) By classical results on curves (see [T1]) and the results of Spivakovsky in [S1], when $\hbox{\rm dim}\hat R^{(\nu)}\leq 2$, for a
rational valuation $\nu$ the equality
$\hbox{\rm r}(\nu)= \hbox{\rm dim}\hat R^{(\nu)}$ implies that the $k$-algebra $\hbox{\rm gr}_\nu R$ is finitely generated and the
equality $\hbox{\rm r}(\nu)= \hbox{\rm dim} R$ implies that the $k[v^{\Phi_+}]$-algebra $\A_\nu (R)$ is finitely generated.\par\noindent 2)
There remains the problem of characterizing among rational valuations of an excellent equicharacteristic local domain the Abhyankar
valuations (characterized by the equality $\hbox{\rm dim} R=\hbox{\rm r}(\nu)$) and the ''weakly Abhyankar''  valuations
(characterized by $\hbox{\rm dim}\hat R^{(\nu)}=\hbox{\rm r}(\nu)$) in terms of the abyssal presentation of $\hat R^{(\nu)}$ and
the construction of the ideal $H$ of $\hat R^m$. One may also ask whether the $k$-algebra $\hbox{\rm gr}_{\hat\nu}\hat R^{(\nu)}$ is finitely
generated if $\hbox{\rm gr}_\nu R$ is.
\end{Remark}
 Let us keep the notations of the preceding subsection, with a quite sufficient initial set of generators $(\hat\gamma_j)_{j\in F}$ (see
subsection \ref{abyssalph}) and the corresponding variables $W_j$ of degree $\hat\gamma_j$, while $h$ is large enough for $k[x_1^{(h)},\ldots ,
x_r^{(h)}]$  to contain the images in $\hbox{\rm gr}_\nu R_\nu$ of the $(W_j)_{j\in F}$, and keep the notations used in the proof of Proposition
\ref{resbin} as well. I number the variables
$(W_j)_{j\in F}$ as $(W_{j_1},\ldots ,W_{j_N})$, with $N=\#F$. We must compute the strict transforms of the equations
$G_{mn}$ under the toric map associated to a regular fan $\Sigma$ as in \ref{resbin}.  By assumption $(\gamma_{j_1},\ldots ,\gamma_{j_N})$ are in the
free submonoid
$\N^r\subseteq \Phi_+$ and we may consider the map $$\Z^N\to \Z^r,\ \ (n_1,\ldots ,n_N)\mapsto
\sum_{t=1}^Nn_t\gamma_{j_t}\in\Z^r.$$ The image $E\subset \check\Z^N$ of the dual map is the
$r$-dimensional free submodule generated by the $r$ elements
$e^i=(\gamma_{j_ti})_{1\leq t\leq N}\in\check\Z^N$.\par\noindent We want to show that if
$w^s$ has value greater than $w^n$, after a toric modification corresponding to a regular fan
compatible with the hyperplanes defining $W$, the transform of $w^n$ divides the
transform of
$w^s$.\par\noindent More generally, we have the:
\begin{proposition}\label{principalization} Let $(w^{p_1},\ldots ,w^{p_c})$ be a finite set of monomials in $\widehat {k[(w_j)_{j\in J}]}$, arranged in
increasing term order. There is a toric modification in the variables $w_j$ appearing in these monomials such that, at the point of the strict transform of
$\hbox{\rm Spec}\hat R^{(\nu)}$ picked by $\hat\nu$, the image of $w^{p_1}$ divides the images of all the other $w^{p_j}$.
\end{proposition}
\begin{proof}
 By construction, and in view of Proposition \ref{grapp}, the assumption means
that at least if we have taken $N$ large enough and $1\leq i\leq c-1$,
$\langle w^j,p_{i+1}-p_i\rangle\geq 0$ for $1\leq j\leq r$ (see the Corollary following Proposition \ref{smgrpapp}).  In view of the fact that
$U^{p_i-p_1}\in \hbox{\rm gr}_\nu R_\nu$ is in $k[x^{(h)}]$ for large enough $h$, after refining our fan
$\Sigma$ to a regular fan $\Sigma '$ without changing its intersection with $W$, we may assume that these
inequalities imply for each cone $\sigma =\langle a^1,\ldots ,a^N\rangle\in \Sigma$ such that
$\sigma\cap W\neq\emptyset$ the inequalities
$\langle a^k,p_i-p_1\rangle \geq 0,\ \ 1\leq k\leq
N$. These inequalities precisely mean that in the chart of $Z(\Sigma ')$ corresponding to
$\sigma$, the image in the transform $\hat {R'}^{m,\nu }$ of $\hat R^{(\nu)}$ by  the toric map corresponding to $\Sigma'$ of the transform
of $u^{p_1}$ divides the image of the transform of the
$u^{p_i}$ (compare with [Z1], Th.2, p.
863).
\end{proof}
\begin{corollary} Let $h_1,\ldots ,h_p$ be elements of $\hat R^{(\nu)}$. There is a toric modification in finitely many of the variables $w_j$ such that at
the point of the strict transform of $\hbox{\rm Spec}\hat R^{(\nu)}$ picked by $\hat\nu$, the image of the ideal $(h_1,\ldots ,h_p)$ is generated by the
element of least valuation.\end{corollary}
\begin{proof} The ideal of $\hat R^{(\nu)}$ generated by all the monomials appearing in the $h_k$'s is finitely generated since $\hat R^{(\nu)}$ is n\oe
therian. It suffices to apply the Proposition to a finite set of generators of this ideal.\end{proof}
\begin{corollary}{\rm *}\label{monomialization} {\rm TLU(d) implies TP(d):}   Given finitely many elements $g_k\in R$, there exists a birational
toric extension $R\to  R'$ subordinate to $\nu$ such that in $R'$ the ideal generated by the elements $g_k$
is equal to a monomial in
$\xi'_i$ times a unit. Here the $\xi'_i\in R'$ lift generators of $\hbox{\rm gr}_\nu R'$. In particular, given finitely many
elements of $R$,  we can achieve that in the ring $R'$ obtained from $R$ by a birational toric extension the element with  the
smallest valuation divides the others.{\rm *}
\end{corollary}
\begin{proof} Let us first assume that $R$ is complete: by the preceding Corollary, it suffices to remark that in $R$ the ideal generated by the
monomials in the
$\xi_i$  which appear in the series representing the $g_k$ is generated by finitely many of them.\par\noindent
In the general case, assuming TLU(d) for $d=\hbox{\rm dim}R$, we have a birational toric extension $R'$ of $R$ such that $\hat R'$ is
regular, and $\nu$ extends to a valuation $\hat\nu$ of a regular quotient $\hat R'/H$ of $\hat R'$ with the same graded ring. We may
assume, since this is a matter of refining a fan to make it compatible with hyperplanes, that in this toric extension the images of the
elements $g_k$ in $\hat R'/H$ generate an ideal of the form ${\xi'}^{\alpha}u$ where $u$ is a unit in $\hat R'/H$. We need to show that
$u$ is in fact in $R'$. Since $\hat R'/H$ is 
regular we may view it as a subring of $\hat R'$, and argue as follows: ${\xi'}^{\alpha}u$ is in $R'$, and so is
${\xi'}^{\alpha}$, so $u\in K\cap \hat R'$, which is equal to $R'$ by the flatness of $\hat R'$ as an $R'$-module ([B3], Chap. 3, \S 3, No. 5,
Corollary 4).
\end{proof}
\begin{corollary}\label{finit} If the graded $k$-algebra $\hbox{\rm gr}_\nu \hat R^{(\nu)}$ is finitely generated, say by the initial forms
$(\overline\eta_j)_{1\leq j\leq N}$ of elements $\eta_j\in \hat R^{(\nu)}$, a toric map $Z(\Sigma)\to {\mathbf A}^N(k)$ 
corresponding to a
regular fan compatible with the hyperplanes
$H_{m-n}$ associated to the binomial equations defining
$\hbox{\rm gr}_{\hat\nu}\hat R^{(\nu)} $ will produce, upon taking the strict transform of $\hbox{\rm Spec}\hat R^{(\nu)}$ embedded in
${\mathbf A}^N(k)$ via the elements $\eta_j$, an embedded uniformization of the valuation
$\hat\nu$ on
$\hat R^{(\nu)}$.\end{corollary}
 The $k$-algebra $\hbox{\rm gr}_{\hat\nu}\hat R^{(\nu)} $ is also finitely presented and we have to
deal with finitely many equations $G_{mn}$. We keep the fan
$\Sigma'$ of Proposition
\ref{principalization} and apply this to compute the strict transform of the
$G_{mn}$. By the same argument as in the proof of Corollary \ref{monomialization}, we need concern ourselves with only a finite set of the monomials
appearing in $F_{mn}$.
 We may assume without loss of generality that the intersection of $\sigma$ with $W$ is
of dimension $r$, so that $(a^1,\ldots ,a^r)$, say, are in $W$. Then in our monomial transformation
$$w_i=y_1^{a^1_i}\ldots y_N^{a^N_i},\ \ 1\leq i\leq N,$$  we find that the strict transform of each
$G_{mn}$ in that chart has the following expression, where I write $y^{\langle a,m\rangle}$
for $y_1^{\langle a^1,m\rangle}\ldots y_N^{\langle a^N,m\rangle}$ except when
I want to single out some of the $y_i$'s
$$\begin{array}{lr} G_{mn}\circ
\pi (\sigma)=\\ y_1^{\langle a^1,m\rangle}
\ldots y_r^{\langle a^r,m\rangle}\bigl(y_{r+1}^{\langle
a^{r+1},n-m\rangle}\ldots y_N^{\langle a^N,n-m\rangle}-
\lambda_{mn}+\sum_pc^{(mn)}_p(v^\phi
)y^{\langle a,p-m\rangle}\bigr).\end{array}$$
It shows that the strict transform of $\hat
R^{(\nu )}$ is non singular, as a deformation of the strict transform of $\hbox{\rm gr}_{\hat\nu}\hat R^{(\nu)}$.
Now the refinement of $\Sigma$ to $\Sigma '$ corresponds to a birational map $Z(\Sigma ')\to
Z(\Sigma )$ which induces an isomorphism outside the divisor
$\prod_{j=r+1}^Ny_j=0$, which the strict transform of our space does not meet, so that we may descend the
result of simultaneous resolution from
$Z(\Sigma ')$ to
$Z(\Sigma )$. This is the same process as in [G-T].
\begin{remark}I do not know whether the toric local uniformization map can be chosen in such a way that it induces an isomorphism outside of the singular
locus of $\hbox{\rm Spec}\hat R^{(\nu)}$.\end{remark}
\par\medskip
In the general case the first problem is to show that our valuation $\hat\nu$ on $\hat
R^{(\nu)}$ can be uniformized after a toric modification which is sufficiently far in the
projective system, but occurs after finitely many steps.
\par\medskip
We keep the notations of  subsection \ref{abyssalph}. Let us now consider the binomial ideal in $k[(W_j)_{j\in F'}]$ which is the kernel of the map to
$\hbox{\rm gr}^{(h)}_{\hat\nu}\hat R^{(\nu)}$ sending $W_j$ to $\overline \eta_j$. As we saw, it is a prime ideal, so that it
corresponds to a torus orbit closure since our residue field is algebraically closed ([E-S]). We apply to it
Proposition \ref{resbin}, and obtain a regular fan $\Sigma$ of ${\mathbf R}_+^{F'}$ compatible with the hyperplanes dual to the
$m-n$ where $(W^m-\lambda_{mn}W^n)_{mn}$ is a system of generators of the binomial ideal. We
may choose a regular simplicial cone $\sigma\in \Sigma$ whose intersection with the space of weights
$E\subset {\mathbf R}_+^{F'}$ is of the largest possible dimension, that is, the rational rank $r$ of
$\Phi$. This corresponds to a monomial map  $$w_i\mapsto y_1^{a^1_i}\ldots y_N^{a^N_i},\ \
1\leq i\leq N=\sharp F',$$ which we may extend by setting  $$w_k=y_k,\ \ k\notin F'.$$ Now we compute the
transforms of the equations of $\hat R^{(\nu)}$ under this transformation. By construction, using
the same argument as we used in the case where the graded algebra is finitely generated, the equations 
$G_{mn}$ whose initial forms involve only $(W_i)_{i\in F'}$, become 
$$\begin{array}{cl} y^{\langle a,m^1-n^1\rangle}(y^{\langle
a,n^1\rangle}-\lambda_{m^1n^1}+\cdots) &=0\\ \vdots\ \ \ \ \ \ \ \ \ \  &\vdots\\ y^{\langle
a,m^t-n^t\rangle} (y^{\langle a,n^t\rangle}-\lambda_{m^tn^t}+\cdots) &=0, \end{array}$$ while
the $G_i$ become $$
w_i^{n_i}\prod_{k\notin F',k<i}w_k^{n^{(i)}_k}y^{t(n^{(i)})}-\lambda_{i}\prod_{k\notin
F,k<i}w_k^{\ell^{(i)}_k}y^{t(\ell^{(i)})}+\cdots +c^{(i)}(v^\phi )w_{i+1}+\cdots =0$$
where $y^{t(n^{(i)})}$ (resp. $y^{t(\ell^{(i)})}$) is a monomial  which is the transform of the monomial in the variables
$(w_j)_{j\in F'}$ which affected $w_i^{n_i}$ (resp. appeared in the other term) in the original relations. 
Note that one can compute the valuations of the $y_k$ from those of the $w_j$. Some of the $y_k$, say for
$k\in T$, will have valuation zero, so that for $k\in T$ there is a constant $c_k\in k_R^*$ such that $y_k-c_k$
has positive valuation; the equations $(y_k-c_k=0)_{k\in T}, (y_k=0)_{k\in F\setminus T}, (w_j=0)_{j\notin F}$
define the closed point which is the center of the valuation in the space obtained by the
toric modification. By openness of transversality, if our toric
map has resolved the $\hbox{\rm gr}^{(h)}_{\hat\nu}\hat R^{(\nu)}$, the first set of equations
defines a regular $\nu$-complete local ring at that point. By the nature of the equations, and the implicit function theorem mentioned in subsection
\ref{abyssalph}, it will be the case for the whole set of equations and the strict transform of $\hat R^{(\nu)}$ will be regular.\par\noindent
Note that by the faitfull flatness of $\widehat{\A_{\hat\nu}(\hat R^{(\nu)})}^{(\nu_\A)}$ over $k[v^{\Phi_+}]$, if the initial binomials of the equations
$G_{mn}$ do not constitute a regular sequence, relations between them will lift to relations between the $G_{mn}$ so that the strict transforms of some of the
$G_{mn}$ will vanish on the nonsingular space defined by the strict transforms of the others, as in the proof of the second part of Proposition \ref{binor}.
\begin{Remark} 1) If one compares with the case where $\hbox{\rm
gr}_{\hat\nu}\hat R^{(\nu)}$ is finitely generated, as in the example of plane branches, one see that the effect
of the ``abyssal'' phenomenon is also to send singularities away to infinity: if we stopped at
any finite step, replacing $u_{i+1}$ by $0$ or a polynomial in $(u_k)_{k\leq i}$ of higher
weight than the binomial initial form, the new equation thus obtained from $G_i$ would
be singular.\par\noindent
b) The last equation written underlines the role of the $(u_i)_{i\in I}$ as semiroots in the sense of [PP]; after our toric modification the $(u_k)_{k\notin F}$
are still part of a natural coordinate system on the transformed ring for our valuation.\par\noindent
c) Note that there are at least $\hbox{\rm r}(\nu )$
independent variables, and possibly more if $\hbox{\rm dim}\hat R^{(\nu)} >\hbox{\rm r}(\nu)$, and that after substitution of the expressions for the
$u_j,\ j\notin F$, some of the equations may become identities. 
\end{Remark}
 \par\noindent At this point we have proved local uniformization by
a finite dimensional toric map for the rational valuation $\hat \nu$ on $\hat R^{(\nu)}$.
\par\medskip\noindent To deduce from this local uniformization for $R$, the fact
that $R$ is excellent is again crucial. In [S2] Spivakovsky used in a similar situation the following \par\noindent
\begin{proposition}\label{genreg}{\rm (see [Ma], \S 33, Spivakovsky [S2], [L-M], Chap. 15, Proof of Th. 15.7, p. 202)}  If $R$ is an excellent local
integral domain, for any prime ideal ${\mathbf q}$ of $R$ and any prime ideal
 ${\mathcal H}$ of $\hat R^{\mathbf q}$ such that ${\mathcal H}\cap R=(0)$, the localization $\hat R^{\mathbf
q}_{\mathcal H}$ is regular.\end{proposition} 
\begin{proof} Since we have ${\mathcal H}\cap
R=(0)$, the ring $\hat R^{\mathbf q}_{\mathcal H}$ is a localization of the ring $\hat R^{\mathbf q}\otimes_R
K$, where $K$ is the field of fractions of $R$. It suffices to prove that this tensor product is
regular. By ([EGA], Scholie 7.8.3, (v)), if $R$ is excellent, the map $\hbox{\rm Spec}\hat R^{\mathbf
q}\to \hbox{\rm Spec}R$ is regular. Our tensor product corresponds to the generic fiber of
this map, so it is regular; it is even a geometrically regular $K$-algebra.  \end{proof}
 Since the $\hat\nu$-initial forms of the elements $\eta_j$ which we used to  generate
topologically the $k$-algebra $\hat R^{(\nu)}$ are Laurent monomials in the $\overline\xi_i$, with $\xi_i\in R$, we may perform in $R$ the trace of
the toric modification described above; it consists in adding to $R$ the monomials (with
some negative exponents) in the $\xi_i$ which are the $y_i$ and localizing at the point
specified  by the values of the $y_i$: this gives us a local ring $R'$ with maximal ideal $m'$. If we denote by $\hat {R'}^{(\nu)}$ the
scalewise completion of $R'$, we have maps of $R$-algebras
$$R'\to R'\otimes_R\hat R^m\to \hat{R'}^{m'}\to \hat{R'}^{(\nu)}.$$
Following through constructions, we see that the image in $\hat{R'}^{m'}$ of the strict transform in $R'\otimes_R\hat R^m$ of the ideal $H$ is contained
in the ideal $H'$ corresponding to $\nu$ in $\hat {R'}^{m'}$, so that we have
$$\hbox{\rm dim}\hat{R'}^{(\nu)}\leq \hbox{\rm dim}\hat{R}^{(\nu)}.$$
That strict inequality may occur has been shown by Spivakovsky, who discovered this ``subanalytic'' phenomenon in the case of valuations of height one
(see [S2]). From this inequality follows the
\begin{proposition}\label{transform}{\rm *} After replacing $R$ by a toric transform, me may assume that the map $$R'\otimes_R\hat R^{(\nu)}\to
\hat{R'}^{(\nu)}$$ is the completion with respect to the maximal ideal $m'\otimes 1+1\otimes \hat m$, and $\hat{R'}^{(\nu)}$ is the
completion of the transform of $\hat{R}^{(\nu)}$. In that case, the ideal $H'$ of $\hat {R'}^{m'}$ is the strict transform of $H$.{\rm *}\end{proposition}
\begin{proof} Let $R\to R'$ be a toric modification, and consider the transform ${\hat R^{(\nu)}}{'}$ of $\hat R^{(\nu)}$, which is $R'\otimes_R\hat
R^{(\nu)}$ localized at the prime ideal $n=m'\otimes 1+1\otimes m\hat R^{(\nu)}$. The extension $\hat R^{(\nu)}\to {\hat R^{(\nu)}}{'}$ is birational
since
$R$ is a subring of $\hat R^{(\nu)}$, and the graded ring $\hbox{\rm gr}_{\hat\nu} {\hat R^{(\nu)}}{'}$ is the strict transform of $\hbox{\rm
gr}_{\hat\nu}R^{(\nu)}$ by the toric map defining our transform. We may choose as representatives in ${\hat R^{(\nu)}}{'}$ of generators of
$\hbox{\rm gr}_{\hat\nu} {\hat R^{(\nu)}}{'}$ the transforms of positive valuation by our toric modification, as in subsection \ref{scac}, of elements
$\eta_j\in R^{(\nu)}$ whose initial forms generate
$\hbox{\rm gr}_{\hat\nu}R^{(\nu)}$. Now the scalewise completion (Proposition \ref{madicomp}) of
${\hat R^{(\nu)}}{'}$ may be a nontrivial quotient of its $n$-adic completion, and in that case $\hbox{\rm dim}{R'}^{(\nu)}<\hbox{\rm
dim}\hat R^{(\nu)}$, but the remarks above show that given a toric transform
$R'$ of $R$, and another one $R^*$ such that the dimension of $\hat {R^*}^{(\nu)}$ is minimal among all toric modifications of $R$, there is a toric
modification of
$R'$ which dominates $R^*$ by a toric map and therefore also reaches the minimum dimension of its scalewise completion.\end{proof}
This implies that if the transform ${\hat R^{(\nu)}}{'}$ of $\hat R^{(\nu)}$ by a toric map is regular, so is $\hat{R'}^{(\nu)}$.
 Let $\hat R^{(\nu)}$ be the scalewise $\nu$-adic completion of a local equicharacteristic excellent integral domain. We know by \ref{nuadiccomp} that
$\hat R^{(\nu)}$ is a  quotient of the completion $\hat R^m$ by an ideal $H$. By Proposition
\ref{genreg}, the assumption that $R$ is excellent implies that the localization $\hat
R^m_H$ is regular. If we have equality in Abhyankar's inequality, the kernel of the map $\hat
R^m\to \hat R^{(\nu)}$ is a minimal prime and the same is true for the map $\hat {R'}^{m'}\to\hat
{R'}^{(\nu )}$. We conclude by case 1) of the general case below.\par\noindent   For the
general case, where Abhyankar's inequality is not assumed to be an equality, I use a variant of the method described by Spivakovsky ([S2], 1997,
Lemma 6.4). We are reduced to the case where
$R$ is such that $\hat R^m_H$ and $\hat R^{m}/H$ are both regular, the first
one by excellence and the second one because we have uniformized $\hat \nu$ on
$\hat R^{(\nu )}$. We want to show that we can make $\hat R^{m}$ (i.e., its strict
transform) non singular by a toric modification defined in $R$. Then $R$
itself will be made non singular. As we have seen in
Proposition \ref{transform}, we may assume that the strict transform of $H$ is equal to $H'$. \par\noindent -case 1) If the ideal $H$ contains only zero
divisors of
$\hat R^m$, it is a minimal prime ideal of $\hat R^m$, say $\hat q_0$. The blowing up in $\hbox{\rm Spec}\hat R^m$ of the sum $J$ of the minimal
primes of $\hat R^m$ separates the strict transforms of the irreducible components. Let
$J_0$ be the ideal of
$\hat R^m/H$ which is the image of the intersection of the other minimal primes of $\hat R^m$. According to toric principalization, which we may apply in
$\hat R^m/H$ because whe know that it is a quotient of
$\widehat {k[(w_j)_{j\in J}]}$, there is a toric modification in the $(\xi_i)_{i\in I}$ which principalizes $J_0$. The corresponding toric map dominates the
strict transform of $\hbox{\rm Spec}\hat R^m/H$ by the blowing up of $J$. Therefore the completion of the transform $R'$ of $R$ by this toric map is an
integral domain, the ideal $H'$ associated to $\nu$ in $ \hat {R'}^{m'}$ is zero and the kernel of the map $\hat R^m\to \hat {R'}^{m'}$ is the ideal $H$. Now
$R'$ has to be regular.\par\noindent -case 2)  If the ideal $H$ contains non zero divisors, the regularity of $\hat R^m/H$ and of $\hat R^m_H$ implies
the regularity of
$\hat R^m$ if we know that $H$ is generated by a regular sequence in $\hat R^m$.  This is equivalent to saying that $H/{H}^2$ is a locally free $ \hat
{R}^{m'}/H$-module. It is enough to produce a toric modification of
$R$ such that the strict transform of $H/H^2$ is equal to $H'/{H'}^2$ and is locally free on $\hbox{\rm Spec}\hat {R'}^{m'}/H'$. 	
Note that since we may assume that $H'$ is the strict transform of $H$ by the map $\hat R^m\to \hat {R'}^{m'}$,  the $ \hat {R'}^{m'}/H'$-module $H'/{H'}^2$
is equal to $H/H^2\otimes  \hat {R'}^{m'}/H'$ divided by its torsion. Now the module $H'/{H'}^2$ will be locally free if the map
$\hat R^m\to \hat{R'}^{m'}$  makes the Fitting ideal of the
$\hat R^{(\nu)}$-module
$H/H^2$ principal at the point picked by the valuation in the toric modification $R\to R'$. Using the Corollary to Proposition
\ref{principalization}, we see that after a toric modification of $R$ we can obtain that the
Fitting ideal is principal. Thus, we reach the situation where $R'$ is regular
and we are done. \par
 There will remain to translate the trace of the toric 
modification on the
completion $\hat R$ (or $\hat R^{(\nu)}$) into a finite sequence of blowing ups with non singular
centers coming from $R$. But one can prove a torus-equivariant version of Gruson-Raynaud's
([G-R]) and Hironaka's ([H2]) proper flattening theorem and use it to show that a toric
modification is dominated by a sequence of toric blowing ups, that is, blowing ups with a torus-invariant non singular center. In fact a more general
result is proved in [DC-P]; see also [O], Chap. 1,
\S 1.7, p.39. We will also have to show that the centers of these blowing-ups
can be chosen contained in the singular locus of $\hbox{\rm Spec} R$ and its transforms, whereas all
we know so far is that these centers are contained in the union of the coordinate
hyperplanes (and their transforms), and although we know by Proposition \ref{truereso} that our toric map can be chosen to be an isomorphism outside of
the singular locus of
$\hbox{\rm Spec}
\hbox{\rm gr}_{\hat\nu}\hat R^{(\nu)}$, we do not know yet the position of these hyperplanes with respect to the singular locus of
$\hbox{\rm Spec} \hat R^{(\nu)}$ .\par\medskip\noindent
 \begin{lemma}\label{monn}Let $\nu$ be a rational valuation of a complete equicharacteristic n\oe therian local ring $R$ with algebraically closed residue
field. Assume that for any
$j\notin F$, in the equality
$$\eta_{j+1} =d_{j+1}\big( \eta_j^{n_j}\prod_{k<j}\eta_k^{n^{(j)}_k}-\lambda_j
\prod_{k<j}\eta_k^{\ell^{(j)}_k}+\sum_sd_s^{(j)}\eta^s\big)$$ of Corollary \ref{semir}, there is no term $\prod_{k<j}\eta_k^{n^{(j)}_k}$. Then we may
assume also that there appear no powers of
$\eta_j$ greater than $n_j$. If we assign value $\nu_0(\eta_j)=\nu(\eta_j)$ for $j\in F$ and the value defined inductively by
$\nu_0(\eta_{j+1})=\nu_0(\eta_j^{n_j})$ for
$j\notin F$, there is in the series representing $\eta_{j+1}$ a unique monomial with the minimal $\nu_0$-value, namely $\eta_j^{n_j}$, and therefore
$\nu_0$ defines a monomial valuation on $R$.\end{lemma} 
\begin{proof}: for each monomial $\eta^s=\tilde\eta^p\eta_j^k$ in the series, we have
$$\nu(\tilde\eta^p)+k\nu(\eta_j)>n_j\nu(\eta_j)$$ and since $k\leq n_j$ this gives 
$(n_j-k)\nu(\eta_j)<\nu(\tilde\eta^p)$. The fact that $n_k\nu(\eta_k)<\nu(\eta_{k+1})$ for $k\notin F$ provides us with inequalities on the components
$p_1,\ldots ,p_\ell$ of $p$, which give us the inequality $\nu_0(\eta^s)>\nu_0(\eta_j^{n_j})$. Since we must have
$\nu_0(\eta_j^{n_j})<\nu_0(\prod_{k<j}\eta_k^{\ell^{(j)}_k})$, this gives the result.\end{proof}
\begin{corollary}
{\rm *}Given a rational valuation $\nu$ of a complete equicharacteristic n\oe therian local ring $R$ with algebraically closed residue field, assuming that
the conditions of Lemma \ref{monn} are satisfied, there exists a rational valuation $\nu_0$ of $R$ such that
$\hbox{\rm gr}_{\nu_0}R$ is finitely generated, $\nu_0(x)\leq\nu (x)$ for $x\in R$, and toric resolutions of the binomial variety $\hbox{\rm
Specgr}_{\nu_0}R$ provide (embedded) uniformizations of
$\nu_0$ which also uniformize
$\nu$ on
$R$.{\rm *}\end{corollary}
\begin{proof} Choose for $\nu_0$ the monomial valuation of $R$, with respect to $(\xi_i)_{i\in
F'}$, determined by
$\nu_0(\xi_i)=\nu (\xi_i)=\gamma_i$ for $i\in F'$ as in Lemma \ref{monn} and apply Corollary
\ref{finit}.\end{proof} It remains to see whether the assumption of lemma \ref{monn} is necessary and whether a similar result holds for any excellent
equicharacteristic local ring with an algebraically closed residue field.

{\textbf {A last example.}}
An example where $\hbox{\rm gr}_\nu R$ is finitely generated is obtained by taking the
first example after Proposition \ref{oneequicompletion}. The completion $\hat R^{(\nu)}$ is the
ring of a branch, the associated graded ring is that of the corresponding monomial
curve which may be desingularized as in subsections \ref{planebranches}
and \ref{branchrev} by a toric map. The transform $R'$ of $R$ by this map has the
property that it is a two-dimensional local ring with maximal ideal $m'$ such that the
maximal ideal of $\hat R'$ contains an element $f'$ such that $\hat R'/(f')$ is a regular
one-dimensional local ring and $(\hat R')_{(f')}$ is regular. So $\hat R'$ is regular and therefore so is $R'$. Essentially the same
argument will work whenever $\hbox{\rm dim}\hat R^{(\nu)} =1$, since once the
singularity of this ring is resolved we deal with a discrete valuation ring, a case where the
argument to bring the number of generators of $H$ down to its height presents no
difficulty.\par\medskip\noindent


\end{document}